\documentclass[12pt,twoside,leqno]{amsart}
\usepackage{amssymb,amsbsy,amsmath,amsfonts,amssymb,amscd,epsfig,
times,graphics,color,xypic}

\definecolor{blue}{cmyk}{1.,1.,0.,1}
\definecolor{red}{cmyk}{0.,1.,1.,1}
\definecolor{green}{cmyk}{1.,0.,1.,1}

\usepackage[latin1]{inputenc}
\newcommand{\C}{\mathbb{C}}\newcommand{\D}{\mathbb{D}}

\newcommand{\N}{\mathbb{N}}\renewcommand{\P}{\mathbb{P}}

\newcommand{\Z}{\mathbb{Z}}

\newcommand{\blue}{\textcolor{blue}}
\newcommand{\green}{\textcolor{green}}
\newcommand{\red}{\textcolor{red}}

\def\boiteepaisseavecuntitre#1{%
  \def\thickhrulefill{\leavevmode \leaders \hrule height 1pt \hfill \kern \z@}%
  \def\bkvz@before@breakbox{\ifhmode\par\fi\vskip\breakboxskip\relax}%
  \fboxrule=1pt
  \def\bkvz@set@linewidth{\advance\linewidth -2\fboxrule
                          \advance\linewidth -2\fboxsep}%
  \def\bkvz@left{\vrule \@width\fboxrule\hskip\fboxsep}%
  \def\bkvz@right{\hskip\fboxsep\vrule \@width\fboxrule}%
  \def\bkvz@top{\hbox to \hsize{%
      \vrule\@width\fboxrule\@height 1.2pt 
      \thickhrulefill{#1}\thickhrulefill
      \vrule\@width\fboxrule\@height 1.2pt}}%
  \def\bkvz@bottom{\hrule\@height\fboxrule}%
  \breakbox}

\newcommand{\THEOREM}{\smallskip\noindent}
\newcommand{\stopTHEOREM}{\smallskip}

\begin{document}

\title[
An algorithm to generate all Demailly-Semple invariants
]{
Application of computational invariant theory 
\\
to Kobayashi hyperbolicity
\\
and to Green-Griffiths 
algebraic degeneracy
}

\author{Jo\"el Merker}

\address{
D\'epartement de Math\'ematiques et Applications, UMR 8553
du CNRS, \'Ecole Normale
Sup\'erieure, 45 rue d'Ulm, F-75230 Paris Cedex 05, 
France. \ \
{\tt www.dma.univ-mrs.fr/$\sim$merker/}}

\email{merker@dma.ens.fr}



\begin{abstract}


A major unsolved problem (according to Demailly 1997) towards the
Kobayashi hyperbolicity conjecture in optimal degree is to understand
jet differentials of germs of holomorphic discs that are invariant
under any reparametrization of the source. The underlying group action
is not reductive, but we provide a complete algorithm to generate all
invariants, in arbitrary dimension $n$ and for jets of arbitrary order
$k$.

Two main new situations are studied in great details. For jets of
order 4 in dimension 4, we establish that the algebra of
Demailly-Semple invariants is generated by 2835 polynomials, while the
algebra of bi-invariants is generated by 16 mutually independent
polynomials sharing 41 gröbnerized syzygies. Nonconstant entire
holomorphic curves valued in an algebraic 3-fold (resp. 4-fold) $X^3
\subset \P^4 (\C)$ (resp. $X^4 \subset \P^5 (\C)$) of degree $d$
satisfy global differential equations as soon as $d \geqslant 72$
(resp. $d \geqslant 259$). A useful asymptotic formula for the
Euler-Poincaré characteristic of Schur bundles in terms
of Giambelli's determinants is derived.

For jets of order 5 in dimension 2, we establish that the algebra of
Demailly-Semple invariants is generated by 56 polynomials, while the
algebra of bi-invariants is generated by 17 mutually independent
polynomials sharing 105 gröbnerized syzygies.


\end{abstract}

\maketitle

\vspace{-0.5cm}

\begin{center}
\begin{minipage}[t]{11cm}
\baselineskip =0.35cm
{\scriptsize

\centerline{\bf Table of contents}

\smallskip

{\bf 1.~Introduction
\dotfill \pageref{Section-1}.}

{\bf 2.~Invariant polynomials and composite differentiation 
\dotfill \pageref{Section-2}.}

{\bf 3.~Bracketing process and syzygies: Jacobi, Plücker~1 and
Plücker~2 \dotfill \pageref{Section-3}.}

{\bf 4.~Survey of known descriptions of ${\sf E }_\kappa^n$ in low
dimensions for small jet levels
\dotfill \pageref{Section-4}.}

{\bf 5.~Initial invariants in dimension $n$ for arbitrary jet
level $\kappa \geqslant 1$ \dotfill \pageref{Section-5}.}

{\bf 6.~Description of the algorithm in dimension
$n=2$ for jet level $\kappa = 4$ \dotfill \pageref{Section-6}.}

{\bf 7.~Action of ${\sf GL}_n ( \C)$ and unipotent reduction 
\dotfill \pageref{Section-7}.}

{\bf 8.~Counter-expectation: insufficiency of bracket invariants 
\dotfill \pageref{Section-8}.}

{\bf 9.~Principle of the general algorithm \dotfill \pageref{Section-9}.}

{\bf 10.~Seventeen bi-invariant generators in dimension $n = 2$ 
for jet level $\kappa = 5$ \dotfill \pageref{Section-10}.}

{\bf 11.~Sixteen (fifteen) bi-invariant in dimension $n = 4$ 
($n=3$) for jet level $\kappa = 4$ \dotfill \pageref{Section-11}.}

{\bf 12.~Approximate Schur bundle decomposition of ${\sf E}_{4, 
m}^4 T_X^*$ \dotfill \pageref{Section-12}.}

{\bf 13.~Asymptotic expansion of the Euler characteristic $\chi \big(
X, \Gamma^{ (\ell_1, \dots, \ell_n)} T_X^* \big)$
\dotfill \pageref{Section-13}.}

{\bf 14.~Euler characteristic calculations
\dotfill \pageref{Section-14}.}

\smallskip

}\end{minipage}
\end{center}

\section*{ \S1.~Introduction}
\label{Section-1}

The Kobayashi hyperbolicity conjecture (1970), in optimal degree and
taking account of Brody's theorem (1978), expects that all entire
holomorphic curves $f : \C \to X$ into a complex projective
(algebraic, smooth) hypersurface $X = X^n \subset \P^{ n+1} (\C)$ must be
constant if $\deg X \geqslant 2n+1$, provided $X$ is generic. In
1980, Green and Griffiths conjectured that if $X \subset \P^{ n+1 } (\C)$
is of general type, which holds in degree $d \geqslant n + 3$, then
there is a proper algebraic subvariety $Y \subsetneq X$ which absorbs
the image of all nonconstant entire holomorphic maps $f : \C \to X$,
namely $f ( \C) \subset Y$ necessarily. Correspondingly, an entire
holomorphic $f : \C \to X$ will be called {\sl algebraically
degenerate} if its image is contained in some proper algebraic
subvariety (which might depend on $f$).

Publications up to 2008 are still quite far from approaching the two
optimal degrees $2n+1$ and $n+3$. For $X^2 \subset \P^3 (\C)$ very generic,
such entire $f$'s are known to be algebraically degenerate and even
constant, in degree $d \geqslant 21$ (resp. $d \geqslant 18$)
according to~\cite{ delg2000} (resp.~\cite{ pau2008}). For $X^3
\subset \P^4 (\C)$ very generic, algebraic degeneracy of such $f$'s holds
true in degree $d \geqslant 593$ according to~\cite{ rou2007a}. For
$X^4 \subset \P^5 (\C)$, a forthcoming work~\cite{ dmr2008} applying the
results of the present paper will obtain an effective degree lower
bound for algebraic degeneracy; other applications to the logarithmic
case also are imminent.

Quite unexpectedly, the two conjectures above and other similar
problems as well in complex algebraic geometry happened in the last
few years to pertain to purely algebraic problems, and not only to
rely upon the scope of some soft techniques (pluripotential theory,
currents, plurisubharmonic functions, {\em etc.}). Computational
invariant theory should be expressly invoked here, as the present
paper will show that what is at stake really is to find a complete
description of the algebra of polynomials that are invariant under a
certain Lie group action, which is {\em not}\, reductive.

\subsection*{ Green-Griffiths Jet differentials}
How can one figure out that a given nonconstant entire holomorphic map
$f = \C \to X^n \subset \P^{ n+1} (\C)$ is constrained to be somehow
degenerate by just being valued in $X$? Looking at its derivatives
$f', f'', \dots, f^{ (k)}$ (in some jet-chart), one may expect at
first to derive, by means of some suitable elimination process,
sufficiently many {\em differential equations} which might presumably
be due to the virtual guidance by some hidden $Y \subsetneq X$
absorbing $f ( \C)$.

For instance, for $X^2 \subset \P^3 (\C)$, the entire $f$'s do satisfy
(invariant) algebraic differential equations of order $k = 2$,
resp. $k = 3$, resp. $k = 4$ when $X$ has degree $d \geqslant 15$,
resp. $d \geqslant 11$, resp. $d \geqslant 9$ according to~\cite{
dem1997, delg2000}, resp.~\cite{ rou2006a}, resp.~\cite{
mer2008a}. For $X^3 \subset \P^4 (\C)$, differential equations of order $k
= 3$ enjoyed by any entire $f$ exist when $X$ has degree $d \geqslant
97$ (\cite{ rou2006b}).

Intrinsically speaking, consider the bundle $J_k$ of $k$-jets of
holomorphic curves $f : ( \D, 0) \longrightarrow (X, x)$ centered
at various points $x = f ( 0) \in X$. In the seminal article~\cite{
gg1980} (1980), Green and Griffiths introduced the fiber bundle ${\sf
E}_{ k, m}^{ GG} T_X^* \to X$ of jet polynomials of order $k$ and of
weighted degree $m$ whose fibers in some jet-chart are complex-valued
polynomials $Q \big( f', f'', \dots, f^{ ( k)} \big)$ satisfying the
weighted homogeneity:
\[
Q\big(\lambda f',\,\lambda^2 f'',\,\dots,\,\lambda^kf^{(k)}\big)
=
\lambda^m\,
Q\big(f',\,f'',\,\dots,f^{(k)}\big),
\]
for every $\lambda \in \C^*$. Global sections of ${\sf E}_{ k, m}^{
GG} T_X^*$ over $X$ are differential operators of order $k$.
Elementary reasonings show (\cite{ gg1980, dem1997, rou2007b,
div2007}) that ${\sf E}_{ k, m}^{ GG} T_X^*$ is in fact a graded
{\em vector}\, bundle isomorphic to the direct sum:
$$
\bigoplus_{\ell_1+2\ell_2+\cdots+k\ell_k=m}\,
{\rm Sym}^{\ell_1}T_X^*\otimes{\rm Sym}^{\ell_2}T_X^*
\otimes\cdots\otimes
{\rm Sym}^{\ell_k}T_X^*.
$$
Such a grading of ${\sf E}_{ k, m}^{ GG} T_X^*$
enables one (\cite{ gg1980}) to derive from Hirzebruch's Riemann-Roch
formula (\cite{hirz1966}) a sharp asymptotic estimate for its
Euler-Poincaré characteristic, namely:
\[
\small
\aligned
\chi\left(X,\,{\sf E}_{k,m}^{GG}T_X^*\right)
&
=
\frac{m^{(k+1)n-1}}{(k!)^n\,\big((k+1)n-1\big)!}\,
\bigg(
\frac{(-1)^n}{n!}\,
{\sf c}_1(X)^n\,({\rm log}\,k)^n
+
{\rm O}\big(({\rm log}\,k)^{n-1}\big)
\bigg)
+
\\
&\ \ \ \ \ \ \ \ \ \ \ \ \ \ \ \ \ \ \ \ \ \ \ \ \ \ \ \ \ \ \ \ 
\ \ \ \ \ \ \ \ \ \ \ \ \ \ \ \ \ \ \ \ \ \ \ \ \ \ \ \ \ \ \ \
\ \ \ \ \ \ \ \ \ \ \ \ \ \ \ \ \ \ \ \
+
{\rm O}\big(m^{(k+1)n-2}\big).
\endaligned
\]
This formula and the knowledge of the expression of the $n$-th power
of the first Chern class (implicitly integrated over $X$) in terms of
the degree:
\[
(-1)^n\,{\sf c}_1(X)^n
=
(d-n-2)^n\,d
\]
entails that, as the jet order $k$ tends to $\infty$, the
characteristic $\chi \big( X,\, {\sf E}_{k,m }^{GG } T_X^* \big)$
becomes eventually positive for $m$ large enough, as soon as $\deg X
\geqslant n + 3$. Thus, up to a constant factor, {\em ${\sf c}_1^n$
becomes the dominant term of $\chi \left( X, \,{\sf E }_{ k,m }^{ GG}
T_X^* \right)$ when $X^n \subset \P^{ n+1} (\C)$ is of general type}.

\subsection*{ Demailly-Semple invariant jet differentials}
In 1997, inspired also by an older paper of Semple, Demailly
introduced a subbundle of ${\sf E}_{ k,m}^{ GG} T_X^*$ having better
positivity properties and exhibiting a nice, stepwise
compactification process.

With $\D \subset \C$ denoting any nonempty open disc centered at
$0$ (possibly $\D = \C$), consider a nonconstant holomorphic curve
$f : \D \to X^n \subset \P^{ n+1} (\C)$. 
Of course, $f ' ( \zeta)$ then
belongs to the tangent space $T_{ X, f ( \zeta)}$ for every $\zeta \in
\D$. The projectivization $[ f' ( \zeta) ] \in PT_{ X, f(\zeta)}$
therefore belongs to the projectivized bundle of tangent lines to $X$, so
that one gratuitously obtains a lifting $f_{ [1]} := \big( f, \, [ f']
\big) : \D \longrightarrow P(T_X)$, at least for all $\zeta$ with
$f' ( \zeta) \neq 0$. Here, $P ( T_X)$ is $( 2n - 1)$-dimensional, but
the so lifted holomorphic curve $f_{ [ 1]}$ happens to be guided by a
certain $n$-dimensional subbundle of $P(T_X)$, better seen as follows.

Abstractly and generally speaking, let $Y$ be a complex manifold, let
$V \subset T_Y$ be any vector subbundle and call $(Y, V)$ a {\sl
directed manifold}. Define $Y' := P ( V)$ the projectivized bundle of
lines contained in the vector subbundle $V \subset TY$ with of course
$\dim Y' = \dim Y + {\rm rk} \, V - 1$. It is equipped with a natural
projection $\pi : Y' \to Y$ which enables one to introduce the {\sl
lifted subbundle} $V ' \subset T_{ Y'}$, the fiber of which, at an 
arbitrary point $(x, [ v]) \in Y'$, is precisely defined by:
\[
V_{(x,[v])}'
:=
\big\{
v'\in T_{X'}:\,d\pi(v')\in\C v
\big\},
\]
and the rank of which is clearly untouched: ${\rm rk} \, V' = {\rm rk} \,
V$. Most importantly, any nonconstant holomorphic $f : \D \to Y$
constrained to be $V$-tangent, namely to satisfy $f' ( \zeta) \in V_{
f ( \zeta)}$ for all $\zeta \in \D$, may be shown (\cite{
dem1997}) to lift automatically, even at points $\zeta$ where $f'
(\zeta)$ vanishes, as a map $f_{ [1]} : \D \to Y'$ which is also
constrained to be $V'$-tangent, namely which necessarily satisfies
$f_{ [1]} ' ( \zeta) \in V_{ f_{ [1]} ( \zeta)}'$ for all $\zeta \in
\D$. So lifting to a higher dimensional manifold keeps 
memory of the original guidance on the base.

Starting therefore with $Y = X^n \subset \P^{ n+1} (\C)$ and with $V =
T_X$, setting $X_0 := X$, $V_0 := T_X$, one defines first (\cite{
dem1997}) $X_1 := P ( T_X)$, $V_1 = V'$ and then inductively $(X_l,
V_l) := (X_{ l-1}', V_{ l-1}')$ with natural projections $\pi_{ l,
l-1} : X_l \to X_{ l-1}$. One then assembles everything for $l = 0$ to
$l = k$ as a tower of projectivized bundles with total projection
$\pi_{ 0, k} : X_k \to X$ and with intermediate projections $\pi_{
j,l} : X_l \to X_j$, for any $0 \leqslant j \leqslant l \leqslant
k$. By applying inductively the above lifting operator $f_{[ l]} :=
\big( f_{ [ l-1]} \big)_{[1]}$, every nonconstant holomorphic curve $f
: \D \to X$ gives rise to lifts $f_{ [l]} : \D \to X_l$ for
all $l = 0, 1, \dots, k$. Each of these lifts is guided by $V_l$,
namely $f_{ [ l]}' ( \zeta) \in V_{ l, f_{ [ l]}' ( \zeta)}$ for all
$\zeta \in \D$.

At each level $l \geqslant 1$, we have a tautological line bundle
$\mathcal{ O}_{ X_l} ( -1)$ over $X_l = P ( V_{ l-1})$ whose fiber at
a point $\big( x_{ l-1 }, [ v_{ l-1 }] \big) \in P ( V_{ l-1})$ just
consists of the line $\C\cdot v_{ l-1}$ directed by (a representative
of) $[v_{ l-1}]$, and similarly as in the projective spaces, one
may build the basic bundles $\mathcal{ O}_{ X_l} ( q)$ for every $q
\in \Z$.

Now, the bundle of invariant jet differentials of order $k$ and of
weighted degree $m$ is the subbundle\footnote{\, Because the dimension
$n$ will vary often in our study, it must be indicated as an exponent
in the notation of the Demailly-Semple bundle. } ${\sf E}_{ k, m}^n
T_X^*$ of ${\sf E}_{ k,m}^{ GG} T_X^*$ whose fibers at a point $x \in
X$ consist of polynomial differential operators $Q \big( f', f'',
\dots, f^{ ( k)} \big)$ which, under arbitrary local reparametrization
$\phi : (\C, 0) \longrightarrow ( \C, 0)$ of the source with $\phi (
0) = 0$, satisfy the general invariancy condition:
\[
Q\big((f\circ\phi)',\,(f\circ\phi)'',\dots,(f\circ\phi)^{(k)}\big)
=
\phi'(0)^m\,Q\big(f',f'',\dots,f^{(k)}\big),
\]
not only under rescaling-like changes of coordinates $\zeta \mapsto
\lambda \, \zeta$ with $\lambda \in \C^*$. This apparently neat
definition hides several algebraic objects which will be inspected and
explored in length throughout the present article. Comparing the two
bundles:
\[
\diagram 
X_k \rto^{\,} 
\dto_{\pi_{0,k}} 
& 
\mathcal{O}_{X_k}(m)
\dto^{\,} 
\\
X \rto^{\,} 
& 
{\sf E}_{k,m}^nT_X^*\,,
\enddiagram
\]
over $X$ and over $X_k$, one establishes (\cite{ dem1997})
the direct image formula 
\[
(\pi_{0,k})_*\mathcal{O}_{X_k}(m)
=
\mathcal{O}\big({\sf E}_{k,m}^nT_X^*).
\]

\subsection*{ Existence of global algebraic differential equations}
What then are the global algebraic differential equations that
nonconstant entire maps $f : \C \to X$ could satisfy? As the
hypersurface $X$ lives in $\P^{ n+1} (\C)$, it carries many ample line
bundles, {\em e.g.} any $\mathcal{ O}_{ X_k} ( q)$ with $q \geqslant
1$.

\smallskip
(\cite{ gg1980, dem1997})
{\em 
Fix an ample line bundle $A \to X$ and assume that ${\sf E}_{ k,m}^n
T_X^* \otimes A^{ -1}$ has nonzero sections, namely:
\[
h^0\big(X,\,{\sf E}_{k,m}^nT_X^*\otimes A^{-1}\big)
=
\dim H^0({\sf same})
\geqslant 1.
\]
Then for every global invariant operator $P \in \Gamma \big( X, \,
{\sf E}_{ k,m}^n T_X^* \otimes A^{ - 1} \big)$ valued in $A^{ -1}$, any
entire holomorphic curve $f$ must satisfy the algebraic differential
equation $P \big( f_{ [ k]} ) \equiv 0$. A similar result also holds true
for the larger bundle ${\sf E}_{ k, m}^{ GG} T_X^*$.
}\medskip

How then one can guarantee the existence of such sections $P$?
Because $X$ is elementarily seen to be of general type when $d
\geqslant n + 3$, it is expected (\cite{ dem1997, rou2006a}) in a
first moment that the Euler-Poincaré characteristic of the
Demailly-Semple subbundle ${\sf E}_{ k,m}^n T_X^*$ should behave in a
way quite similar to the Green-Griffiths bundle, with ${\sf
c}_1^n$ becoming the dominant term (up to a constant factor) as $k$
and $m$ both tend to $\infty$, so that $\chi \big(X,\, {\sf E}_{
k,m}^n T_X^* \big)$ should be eventually large, and furthermore in a
second moment, it is also expected that the dimension of the space of
global sections $H^0 \big( X, \, {\sf E}_{ k, m}^n T_X^* \big)$ should
be eventually large, due to some vanishing or to some control of the
higher order cohomology groups. The truth of such conjectural
expectations would presumably open new routes towards a solution in
optimal degree to the two above-mentioned conjectures.

\subsection*{ Seeking Schur bundle decomposition of ${\sf E}_{ k,m}^n 
T_X^*$} However, as is written in \cite{ dem1997}, it is a major
unsolved problem to find the decomposition of ${\sf E}_{ k, m}^n
T_X^*$ into direct sums of the irreducible Schur bundles $\Gamma^{ (
\ell_1, \ell_2, \dots, \ell_n )} T_X^*$ with $\ell_1 \geqslant \ell_2
\geqslant \cdots \geqslant \ell_n$ that are the basic bricks and whose
cohomology is somehow currently available. According to a possible
strategy developed for $k = n = 3$ mainly by Rousseau in~\cite{
rou2006a, rou2006b}, such a decomposition would yield access to the
Euler characteristic $\chi \big( X, \, {\sf E}_{ k, m}^n T_X^* \big)$,
and then afterwards, one would attain an effective estimate of $h^0
\big( X, \, {\sf E}_{ k, m}^n T_X^* \big)$, provided one controls the
other cohomology groups. In fact, the only decompositions known up to
now are the following; the second one (\cite{ rou2006a}) already
required a nontrivial argument based on a theorem of Popov about
polarization of multilinear invariants.

\smallskip$\bullet$\,\,
For $n = k = 2$ (\cite{dem1997}):
\[
\footnotesize
\aligned
{\sf E}_{2,m}^2T_X^*
=
\bigoplus_{a+3b=m}\,
\Gamma^{(a+b,\,b)}\,T_X^*.
\endaligned
\]

\smallskip$\bullet$\,\,
For $n=k=3$ and also for $n = 2$, $k = 3$ (\cite{rou2006a}):
\[
\footnotesize
\aligned
{\sf E}_{3,m}^3T_X^*
&
=
\bigoplus_{a+3b+5c+6d=m}\,
\Gamma^{(a+b+2c+d,\,b+c+d,\,d)}\,T_X^*,
\\
&
\text{\rm and}
\ \ \ \ \ \
{\sf E}_{3,m}^2T_X^*
=
\bigoplus_{a+3b+5c=m}\,
\Gamma^{(a+b+2c,\,b+c)}\,T_X^*.
\endaligned
\]

\smallskip$\bullet$\,\,
For $n=2$, $k = 4$ (\cite{ mer2008a}):
\[
\footnotesize
\aligned
{\sf E}_{4,m}^2T_X^*
&
=
\bigoplus_{a+3b+5c+8e=m}\,
\Gamma^{(a+b+2c+2e,\,b+c+2e)}\,T_X^*
\\
&\ \ \ \ \
\bigoplus_{7+a+5c+7d+8e=m}\,
\Gamma^{(3+a+2c+3d+2e,\,1+c+d+2e)}\,T_X^*.
\endaligned
\]
In this paper, we mainly attack the case $n = k = 4$. The complexity
increases suddenly and we seem to be still quite far from being able
to push the jet order $k$ to $\infty$.

\smallskip\noindent
{\bf Theorem.}
{\em 
On a smooth complex algebraic hypersurface $X^4 \subset \P^5 (\C)$, the
graduate $m$-th part ${\sf E}_{ 4, m}^4 T_X^*$ of the
Demailly-Semple bundle ${\sf E}_4^4 T_X^* = \oplus_m \, {\sf E}_{ 4,
m}^4 T_X^*$ has the following decomposition in direct sums of Schur
bundles:
\[
\aligned
&
{\sf E}_{4,m}^4T_X^*
=
\bigoplus_{(a,b,\dots,n)\in\N^{14}\backslash
(\square_1\cup\cdots\cup\square_{41}) 
\atop
o+3a+\cdots+21n+10p=m}\,
\\
&
{\scriptsize
\Gamma
\left(
\aligned
o+a+2b+3c+d+2e+3f+2g+2h+3i+4j+3k+3l+4m'+5n+p&
\\
a+b+c+d+e+f+2g+2h+2i+2j+2k+3l+3m'+3n+p&
\\
d+e+f+h+i+j+2k+2l+2m'+2n+p&
\\
p&
\endaligned
\right)}\,T_X^*,
\endaligned
\]
where the 41 subsets $\square_i$, $i = 1, 2, \dots, 41$, of $\N^{ 14}
\ni (a, b, \dots, l, m', n)$ are explicitly defined in the complete
statement on p.~\pageref{Schur-decomposition-4-4}.
}\medskip

It is known (\cite{ rou2006b}) that ${\sf E}_{ k, m}^3 T_X^*$ has no
nonzero sections for jet order $k = 1$ or $k = 2$. More generally
(\cite{ div2007}), for jet order $k \leqslant n - 1$ strictly smaller
than the dimension, sections are never available: $H^0 \big( X, \,
{\sf E}_{ k, m}^n T_X^* \big) = 0$. Consequently, even if one may
easily deduce from the above theorem a Schur decomposition of ${\sf
E}_{ 4, m}^3 T_X^*$, for applications to hyperbolicity in dimension
higher than 3, one should always start with jet order $k$ at least
equal to the dimension\footnote{\, Nonetheless, we ignore whether the
case $n = k = 5$ is accessible to us. }. The case $n = k = 4$ was the
first unknown one before.

\subsection*{ Asymptotic expansion of Euler-Poincaré characteristic}
Because the characteristic is just additive on direct sums of vector
bundles, knowing a representation of ${\sf E}_{ k, m}^n T_X^*$ (for
certain values of $n$, $k$, {\em e.g.} for $n = k = 4$) as a direct
sum of certain Schur bundle is very convenient, provided of course
that one already possesses an asymptotic for the Euler-Poincaré
characteristic of the $\Gamma^{ ( \ell_1, \dots, \ell_n ) } T_X^*$ as
$\ell_1 + \cdots + \ell_n \to \infty$. Section~13 will derive an
explicit asymptotic for which there seems to be no reference with a
precise enunciation (compare~\cite{bru1997, rou2004}). Because of the
relations ${\sf c}_k \big( T_X^* \big) = (-1)^k {\sf c}_k \big( T_X
\big)$, there is no loss of generality to express everything in terms
of the Chern classes of the tangent bundle $T_X$.

\smallskip\noindent{\bf Theorem.}
{\em 
The terms of highest order with respect to $\vert \ell \vert = \max_{
1 \leqslant i \leqslant n } \, \ell_i$ in the Euler-Poincaré
characteristic of the Schur bundle $\Gamma^{ ( \ell_1, \ell_2, \dots,
\ell_n )} \, T_X$ are homogeneous of order ${\rm O} \big( \vert \ell
\vert^{ \frac{ n ( n+1)}{ 2}} \big)$ and they are given by a sum of
$\ell_i'$-determinants indexed by all the partitions $(\lambda_1,
\dots, \lambda_n)$ of $n$:
\[
\small
\aligned
&
\chi\Big(X,\,\,
\Gamma^{(\ell_1,\ell_2,\dots,\ell_n)}\,T_X\Big)
=
\\
&
=
\sum_{\lambda\,\text{\rm partition of}\,\,n}\,
\frac{
{\sf C}_{\lambda^c}}{(\lambda_1+n-1)!\,\cdots\,\lambda_n!}\,
\left\vert
\begin{array}{cccc}
{\ell_1'}^{\lambda_1+n-1} & {\ell_2'}^{\lambda_1+n-1} &
\cdots & {\ell_n'}^{\lambda_1+n-1}
\\
{\ell_1'}^{\lambda_2+n-2} & {\ell_2'}^{\lambda_2+n-2} &
\cdots & {\ell_n'}^{\lambda_2+n-2}
\\
\vdots & \vdots & \ddots & \vdots
\\
{\ell_1'}^{\lambda_n} & {\ell_2'}^{\lambda_n} & 
\cdots & {\ell_n'}^{\lambda_n}
\end{array}
\right\vert
+
\\
&\ \ \ \ \ \ \ \ \ \
+
{\rm O}\Big(
\vert\ell\vert^{\frac{n(n+1)}{2}-1}
\Big),
\endaligned
\]
where $\ell_i ' := \ell_i + n - i$ for notational brevity, with
coefficients ${\sf C}_{\lambda^c}$ being expressed in terms of the
Chern classes ${\sf c}_k \big( T_X \big) = {\sf c}_k$ of $T_X$
by means of {\em Giambelli's determinantal expression} depending
upon the {\em conjugate} partition $\lambda^c$:
\[
{\sf C}_{\lambda^c}
=
{\sf C}_{(\lambda_1^c,\dots,\lambda_n^c)}
=
\left\vert
\begin{array}{cccccc}
{\sf c}_{\lambda_1^c} & {\sf c}_{\lambda_1^c+1} &
{\sf c}_{\lambda_1^c+2} & \cdots & {\sf c}_{\lambda_1^c+n-1}
\\
{\sf c}_{\lambda_2^c-1} & {\sf c}_{\lambda_2^c} &
{\sf c}_{\lambda_2^c+1} & \cdots & {\sf c}_{\lambda_2^c+n-2}
\\
\vdots & \vdots & \vdots & \ddots & \vdots
\\
{\sf c}_{\lambda_n^c-n+1} & {\sf c}_{\lambda_n^c-n+2} &
{\sf c}_{\lambda_n^c-n+3} & \cdots & 
{\sf c}_{\lambda_n^c}
\end{array}
\right\vert,
\]
on understanding by convention that ${\sf c}_k := 0$ 
for $k< 0$ or $k > n$, and that ${\sf c}_0 := 1$. 
}\medskip

\subsection*{ Effective calculations of characteristics in 
dimensions 3 and 4}
We then perform electronically assisted computations 
to obtain the desired, quite complicated value of the
characteristic of ${\sf E}_{ 4, m}^4 T_X^*$. 

\smallskip\noindent{\bf Theorem.}
{\em 
If $X^4 \subset \P^5 (\C)$ is a degree $d$ smooth algebraic 4-fold,
then as $m \to \infty$, one has the asymptotic:
\[
\footnotesize
\aligned
\chi\big(X,\,{\sf E}_{4,m}^4T_X^*\big)
&
=
\frac{m^{16}}{1313317832303894333210335641600000000000000}\,
\cdot\,d\,\cdot
\\
&\ \ \ \ \
\cdot
\big(
50048511135797034256235\,d^4
-
\\
&\ \ \ \ \ \ \ \ \ \
-
6170606622505955255988786\,d^3
-
\\
&\ \ \ \ \ \ \ \ \ \
-
928886901354141153880624704\,d
+
\\
&\ \ \ \ \ \ \ \ \ \
+
141170475250247662147363941\,d^2
+
\\
&\ \ \ \ \ \ \ \ \ \
+
1624908955061039283976041114
\big)
\\
&\ \ \ \ \ \
+
{\rm O}\big(m^{15}\big).
\endaligned
\]
Furthermore, the coefficient of $m^{ 16}$ here, a factorized
polynomial of degree
5 with respect to $d$, is positive in all degrees $d \geqslant 96$.

}\medskip

For $n = k = 3$, based on his above-mentioned Schur decomposition of
${\sf E}_{ 3, m}^3 T_X^*$,
Rousseau (\cite{ rou2006a}) showed that:
\[
\footnotesize
\aligned
\chi\big(X,\,{\sf E}_{3,m}^3T_X^*\big)
=
\frac{m^9}{81648000000}\,\cdot\,d\,\cdot\,
\big(
389d^3-20739d^2+185559d-358873
\big)
+
{\rm O}\big(m^8\big),
\endaligned
\]
and that the coefficient of $m^9$ is positive for all degrees $d \geqslant
43$. Furthermore, in~\cite{ rou2004}, Rousseau showed that $h^2 \big(
X, \, {\rm Sym }^m \, T_X^* \big) = \big( - \frac{ 7}{ 24} d + \frac{
1}{ 8} d^2 \big) \, m^5 + {\rm O} \big( m^4 \big)$ in any degree $d
\geqslant 6$, so that one cannot expect second
cohomology groups to vanish. Afterwards, as the main objective of the
paper~\cite{ rou2006b}, he first established the general majoration:
\[
\footnotesize
\aligned
h^2
\big(
X,\,\Gamma^{(\ell_1,\ell_2,\ell_3)}T_X^*
\big)
\leqslant
d(d+13)\,
\frac{3(\ell_1+\ell_2+\ell_3)^3}{2}\,
(\ell_1-\ell_2)(\ell_1-\ell_3)(\ell_2-\ell_3)
+
{\rm O}\big(\vert\ell\vert^5\big).
\endaligned
\]
he then deduced by summation from the cited decomposition
${\sf E}_{ 3, m}^3 T_X^* = \bigoplus_{ a + 3b +
5c + 6d = m} \, \Gamma^{ ( a+ b+ 2c+ d, \, b+ c + d, \, d)} T_X^*$
that:
\[
\footnotesize
\aligned
h^2\big(X,\,{\sf E}_{3,m}^3T_X^*\big)
\leqslant
\frac{49403}{252\cdot 10^7}\,
d(d+13)\,m^9
+
{\rm O}\big(m^8\big), 
\endaligned
\]
and finally, by applying the trivial minoration:
\[
\aligned
h^0\big(X,\,{\sf E}_{4,m}^3T_X^*\big)
\geqslant
\chi\big(X,\,{\sf E}_{4,m}^3T_X^*\big)
-
h^2\big(X,\,{\sf E}_{4,m}^3T_X^*\big),
\endaligned
\] 
stemming from the definition $\chi = h^0 - h^1 + h^2 - h^3$,
he immediately deduced the minoration:
\[
\footnotesize
\aligned
h^0\big(X,\,{\sf E}_{3,m}^3T_X^*\big)
\geqslant
\frac{m^9}{408240000000}\,\cdot\,d\,\cdot
\big(1945\,d^3-103695\,d^2-7075491\,d-105837083\big)
+
{\rm O}\big(m^8\big),
\endaligned
\]
in which the coefficient of $m^9$ is checked (again electronically) to
be positive in all degrees $d \geqslant 97$. As a result, nontrivial
sections of ${\sf E}_{ 3, m}^3 T_X^*$ exist when $\deg X \geqslant
97$.

For jets of order $4$ in dimension 3, when applying in dimension 3 our
decomposition of ${\sf E}_{ 4, m}^3 T_X^*$ into Schur bundles which
appears in the theorem on p.~\pageref{Schur-decomposition-4-4}, a
Maple computation using the cited majoration formula for 
$h^2 \big( X, \, \Gamma^{ (\ell_1, \ell_2, \ell_3)} T_X^* \big)$
then provides:
\[
h^2
\big(
X,\,{\sf E}_{4,m}^3T_X^*
\big)
\leqslant
d(d+13)\,
\frac{342988705758851}{29822568148961280000000}\,
m^{11}
+
{\rm O}\big(m^{10}\big).
\]

\smallskip\noindent{\bf Theorem.}
{\em
The asymptotic, as $m \to \infty$, of the Euler-Poincaré
characteristic of the Demailly-Semple bundle ${\sf E}_{ 4, m}^4 T_X^*$
on a degree $d$ smooth projective 
algebraic 3-fold $X^3 \subset \P^4 (\C)$ is given by:
\[
\footnotesize
\aligned
\chi\big(X,\,{\sf E}_{4,m}^3T_X^*\big)
&
=
\frac{m^{11}}{206133591045620367360000000}\,\cdot\,d\,\cdot\,
\big(
1029286103034112\,d^3
-
\\
&\ \ \ \ \
-
38980726828290305\,d^2
+
299551055917162501\,d
-
561169562618151944
\big)
+
\\
&\ \ \ \ \
+
{\rm O}\big(m^{10}\big),
\endaligned
\]
and the coefficient of $m^{ 11}$ here is positive in all degrees $d
\geqslant 29$. Furthermore, subtracting to this asymptotic
the above majorant of $h^2 \big( X, \, {\sf E}_{ 4, m}^3 
T_X^* \big)$:
\[
\footnotesize
\aligned
h^0\big(X,\,{\sf E}_{4,m}^3T_X^*\big)
&
\geqslant
\frac{m^{11}}{206133591045620367360000000}\,\cdot\,d\,\cdot\,
\big(
1029286103034112\,d^3
-
\\
&\ \ \ \ \
-
38980726828290305\,d^2
+
2071186878288015611\,d
-
31380762707285467400
\big)
+
\\
&\ \ \ \ \
+
{\rm O}\big(m^{10}\big),
\endaligned
\]
and the modified coefficient here of $m^{ 11}$
is now positive in all degrees $d\geqslant 72$.
}\medskip

This last condition $d \geqslant 72$ on the degree insuring the
existence of global invariant jet differentials of order $\kappa = 4$
on $X^3 \subset \P^4 ( \C)$ improves the condition $d \geqslant 97$
obtained in~\cite{ rou2006b} and appears to be slightly better than
the condition $d \geqslant 74$ obtained more recently in~\cite{
div2007} with another approach. A number of further numerical
applications shall appear soon
(\cite{ dmr2008}); as will be seen in a near future, in dimension 4,
the lower bound on the degree $d \geqslant 259$ for the existence of
sections which will based on the present approach will also improve
the bound $d \geqslant 298$ obtained in~\cite{ div2007}. Nonetheless,
we must stop at this point in order to describe the main contribution
of the present article. Last but not least, we cannot go beyond
without mentioning that Siu's strategy for establishing algebraic
degeneracy (\cite{ siu2004, pau2008, rou2007a}) will also bring
further fruits thanks to the recent construction of a global
meromorphic framing on the space of vertical $n$-jets tangent to the
universal hypersurface in arbitrary dimension $n$ (\cite{ mer2008b}).

\subsection*{ A problem in invariant theory}
Now, how does one obtain Schur decompositions of Demailly-Semple
bundles? To begin with, we show how one can understand the condition
of being invariant under reparametrization in terms of classical
invariant theory.

Let us from now on denote by $\kappa$ (instead
of $k$) the jet order and let us abbreviate
$j^\kappa f = \big( f', f'', \dots, f^{ (\kappa)} \big)$.

The group ${\sf G}_\kappa$ of $\kappa$-jets at the origin of local
reparametrizations $\phi ( \zeta) = \zeta + \phi'' ( 0) \, \frac{
\zeta^2}{ 2!} + \cdots + \phi^{ ( \kappa) } ( 0) \, \frac{
\zeta^\kappa }{ \kappa !} + \cdots$ that are tangent to the identity,
namely which satisfy $\phi' ( 0) = 1$, 
may be seen to act linearly on the $n
\kappa$-tuples $\big( f_{ j_1}', f_{ j_2}'', \dots, f_{ j_\kappa}^{ (
\kappa)} \big)$ by plain matrix multiplication, {\em i.e.} when we set
$g_i^{ ( \lambda)} := \big( f_i \circ \phi \big)^{ (\lambda)}$, a
computation applying the chain rule gives for each index $i$:
\[
\footnotesize
\aligned
\left(
\begin{array}{c}
g_i'
\\
g_i''
\\
g_i'''
\\
g_i''''
\\
\vdots
\\
g_i^{(\kappa)}
\end{array}
\right)
=
\left(
\begin{array}{ccccccc}
1 & 0 & 0 & 0 & \cdots & 0
\\
\phi'' & 1 & 0 & 0 & \cdots & 0
\\
\phi''' & 3\phi'' & 1 & 0 & \cdots & 0
\\
\phi'''' & 4\phi'''+3{\phi''}^2 & 6\phi'' & 1 & \cdots & 0
\\
\vdots & \vdots & \vdots & \vdots & \ddots & \vdots 
\\
\phi^{(\kappa)} & \cdots & \cdots & \cdots & \cdots & 1
\end{array}
\right)
\left(
\begin{array}{c}
f_i'\circ\phi
\\
f_i''\circ\phi
\\
f_i'''\circ\phi
\\
f_i''''\circ\phi
\\
\vdots
\\
f_i^{(\kappa)}\circ\phi
\end{array}
\right)
\ \ \ \ \ \ \ \ 
{\scriptstyle{(i\,=\,1\,\cdots\,n)}}.
\endaligned
\]
Polynomials ${\sf P} \big( j^\kappa f \big)$ invariant by
reparametrization satisfy by definition for some integer $m$: 
\[
{\sf P}\big(j^\kappa g)
=
{\sf P}\big(j^\kappa(f\circ\phi)\big) 
=
\phi'(0)^m\cdot
{\sf P}\big((j^\kappa f)\circ\phi\big)
=
{\sf P}\big((j^\kappa f)\circ\phi\big),
\] 
for any $\phi$. If we denote by ${\sf E}_{ \kappa, m}^n$ the vector
space consisting of such polynomials, the direct sum ${\sf E}_\kappa^n
= \oplus_{ m\geqslant 1} \, {\sf E}_{ \kappa, m}^n$ forms an algebra
graded by constancy of weights: ${\sf E}_{ \kappa, m_1}^n \cdot {\sf
E}_{ \kappa, m_2}^n \subset {\sf E}_{ \kappa, m_1 + m_2}^n$.

Then obviously when $\phi' (0) = 1$, the algebra ${\sf E}_\kappa^n$
just coincides with the algebra of invariants for the linear group
action represented by the group of matrices just written:
\[
{\sf P}\big(j^\kappa g\big)
=
{\sf P}\big({\sf M}_{\phi'',\phi''',\dots,\phi^{(\kappa)}}\cdot
j^\kappa f\big)
=
{\sf P}\big(j^\kappa f\big),
\]
with $\phi'', \phi''', \dots, \phi^{ ( \kappa)}$ interpreted as
arbitrary complex constants. Such a group clearly 
has dimension $\kappa - 1$.

But unfortunately, this group of matrices is a subgroup of the full
unipotent group, hence it is {\em non-reductive}, and for this reason,
it is impossible to apply almost anything from the so well developed
invariant theory of reductive actions (\cite{dk2002}). Moreover,
though the invariants of the full unipotent group are well understood,
as soon as one looks at a {\em proper} subgroup of it, formal
harmonies happen to be rapidly destroyed.

We ignore whether the algebra of invariants is finitely generated, in
general. But in all previously known cases (carefully reminded below)
and in all further new cases studied in this paper, ${\sf E}_\kappa^n$
is finitely generated. We will establish that the (graded) algebra
${\sf E}_4^4 = \oplus_m \, {\sf E}_{ 4, m}^4$ is generated by 2835
invariant polynomials and that ${\sf E}_5^2 = \oplus_m \, {\sf E}_{ 5,
m}^2$ is generated by 56 invariant polynomials. We will also provide,
in the theorem stated in length on p.~\pageref{normal-syzygies}, a
{\em general algorithm} which, {\em in arbitrary dimension $n$ and for
arbitrary jet order $\kappa$, generates all invariants by adding a new
invariant only when it cannot be expressed as a polynomial with
respect to the already known invariants, and which stops after a
finite number of loops if and only if ${\sf E}_\kappa^n = \oplus_m \,
{\sf E}_{ \kappa, m}^n$ is finitely generated as an algebra}.

\subsection*{ Insufficiency of bracketing}
By definition, a polynomial ${\sf P} \big( j^\kappa f \big)$ in the
$\kappa$-th order jet space which is invariant by reparametrization
must satisfy ${\sf P}\big( j^\kappa ( (f \circ \phi) \big) = {\phi
'}^m \, {\sf P} \big( j^\kappa f) \circ \phi \big)$ for every
biholomorphism $\phi : (\D, 0) \longrightarrow (\D, 0)$, where the
integer $m$ is called the {\sl weight} of ${\sf P}$, and where it is
implicitly understood that the base point is the origin. Also,
suppose next that ${\sf Q}$ is another invariant of weight $n$ in the
$\tau$-th order jet space, {\it i.e.} satisfying ${\sf Q} \big( j^\tau
( f\circ \phi) \big) = {\phi'}^n \, {\sf Q} \big( (j^\tau f)\circ \phi
\big)$. If ${\sf D} := \sum_{ k = 1}^n\, \sum_{ \lambda \in \N}\,
\frac{ \partial ( \bullet) }{ \partial f_k^{ ( \lambda) } } \, \cdot
f_k^{ ( \lambda + 1)}$ denotes the {\sl total differentiation
operator}, which acts on any polynomial in $f', f'', \dots, f^{ (
\kappa)}$ as if it differentiated it with respect to the (virtual)
source variable $\zeta \in \D$, then the {\sl bracket expression}:
\[
\big[{\sf P},\,{\sf Q}\big]
:=
n\,{\sf D}{\sf P}\cdot{\sf Q}
-
m\,{\sf P}\cdot{\sf D}{\sf Q}
\]
will easily be checked (in \S3) to provide gratuitously another
invariant of weight $m + n + 1$ in the jet space of order $1 + \max (
\kappa, \tau)$.

For jet order $\kappa = 1$, the algebra of invariants is just $\C
\big[ f_1', f_2', \dots, f_n' \big]$. For $\kappa = 2$, the algebra
${\sf E}_2^n$ is generated by the $f_i'$ together with the
two-dimensional Wronskians $f_i' f_j'' - f_i'' f_j'$ which identify to
the brackets $\big[ f_j', \, f_i' \big]$, where $1 \leqslant i, j
\leqslant n$.

For $\kappa = 3$ in dimension $n = 2$, the Demailly-Semple algebra
${\sf E}_3^2$ is generated by 5 mutually independent invariants:
\[
f_1',\ \ \ \ \
f_2',\ \ \ \ \ 
\Lambda^3
:=
\big[f_2',\,f_1'\big],\ \ \ \ \
\Lambda_1^5
:=
\big[\Lambda^3,\,f_1'\big],\ \ \ \ \
\Lambda_2^5
:=
\big[\Lambda^3,\,f_2'\big],
\]
which all are furnished by just bracketing, according to~\cite{ rou2006a};
(but bracketing did not enter the scene there). 

In the next dimension $n = 3$ for jets of the same order $\kappa = 3$,
the Demailly-Semple algebra ${\sf E}_3^3$ is generated by 16 mutually
independent invariants (\cite{ rou2006a}), namely the $3 + 3 + 9 = 15$
following ones:
\[
f_i',\ \ \ \ \
\Lambda_{i,j}^3
:=
\big[f_j',\,f_i'\big],\ \ \ \ \
\Lambda_{i,j;\,k}^5
:=
\big[\Lambda_{i,j}^3,\,f_k'\big],\ \ \ \ \
\]
where $1 \leqslant i < j \leqslant 3$ and where $1 \leqslant k
\leqslant 3$, which are clearly all obtained by bracketing some invariants
from the preceding jet level, together with the three-dimensional
Wronskian:
\[
D_{1,2,3}^6
:=
\left\vert
\begin{array}{ccc}
f_1' & f_2' & f_3'
\\
f_1'' & f_2'' & f_3'' 
\\
f_1''' & f_2''' & f_3'''
\end{array}
\right\vert,
\]
which also appears, though after some division by $f_1'$,
to come from the brackets, for one checks by direct 
calculation the three relations:
\[
\aligned
\big[\Lambda_{1,2}^3,\,\Lambda_{1,3}^3\big]
&
=
-3\,f_1'\,D_{1,2,3}^6,
\ \ \ \ \ \ \ \ \ \
\big[\Lambda_{1,2}^3,\,\Lambda_{2,3}^3\big]
=
-3\,f_2'\,D_{1,2,3}^6,
\\
&
\ \ \ \ \ \ \ \ \
\big[\Lambda_{1,3}^3,\,\Lambda_{2,3}^3\big]
=
-3\,f_3'\,D_{1,2,3}^6.
\endaligned
\]
Here, as the reader may have observed already, we always put the
weight of every invariant at the upper index place.

Lastly, coming back to the dimension $n = 2$, for jet order $\kappa =
4$, the algebra ${\sf E}_4^2$ is generated by the 9 mutually
independent invariants (\cite{ dem2007, mer2008a}):
\[
\small
\aligned
&
f_1',\ \ \ \ \ \ \ \ \ \
f_2',\ \ \ \ \ \ \ \ \ \
\Lambda_{1,2}^3,\ \ \ \ \ \ \ \ \ \ 
\Lambda_{1,2;\,1}^5,\ \ \ \ \ \ \ \ \ \ 
\Lambda_{1,2;\,2}^5,
\\
&
\Lambda_{1,1}^7
:=
\big[\Lambda_{1,2;\,1}^5,\,f_1'\big],
\ \ \ \ \ \ \ \ \ \ \ \ \
\Lambda_{1,2}^7
:=
\big[\Lambda_{1,2;\,1}^5,\,f_2'\big]
=
\big[\Lambda_{1,2;\,2}^5,\,f_1'\big]
=
\Lambda_{2,1}^7,
\\
&
\Lambda_{2,2}^7
:=
\big[\Lambda_{1,2;\,2}^5,\,f_2'\big],
\ \ \ \ \ \ \ \ \ \ \ \ \ \ 
M^8
:=
\frac{1}{f_1'}\,
\big[\Lambda_{1,2;\,1}^5,\,\Lambda_{1,2}^3\big],
\endaligned
\]
coming again all from bracketing, possibly allowing a division 
by $f_1'$. 

In view of all these positive results, one could believe that
bracketing (with possible division) always generate all invariants
when passing from one jet level to the subsequent one. In fact, the
two so-called {\sl sigma-} and {\sl Omega-processes} are known to
generate all the invariants of binary forms in any degree (\cite{
ol1999, dk2002, pro2007}).

Unfortunately, in~\cite{ mer2008a}, we discovered that in dimension $n
= 2$ for jet order $\kappa = 5$, many invariants exist which are
totally independent from the ones obtained by bracketing the
invariants existing at the inferior jet levels $\kappa \leqslant
4$. Section~8 will provide more explanations, emphasizing that
{\em it is by no means possible to derive these further invariants
by dividing any incoming bracket-invariant by any other already 
known (bracket) invariant}.

Nonetheless, there could exist a second (and even a third)
algebraically uniform process which would generate gratuitously many
other invariants, and which, in cooperation with the bracketing
process, would be complete, but regarding such an 
idea, we must confess our ignorance.

\subsection*{ Initial rational expression for invariants}
Hopefully, the algorithm we already devised (and hid slightly?) 
in~\cite{ mer2008a} provides another route. How does it work?

To begin with, we define $\Lambda_{ 1, i}^3 := \big[ f_i', f_1'\big]$
and then by induction for any $\lambda$ with $3 \leqslant \lambda
\leqslant \kappa$:
\[
\Lambda_{1,i;\,1^{\lambda-2}}^{2\lambda-1}
:=
\big[\Lambda_{1,i;\,1^{\lambda-3}}^{2\lambda-3},\,f_1'\big].
\]
Being built by bracketing, these are invariants of weight $2 \lambda
- 1$ for any $i = 1, \dots, n$. The power $\lambda - 2$ of $1$ counts the
number of brackets with $f_1'$, starting from the Wronskian
$\Lambda_{ 1, i}^3$.

The preliminary step is to establish a {\em rational} representation
of any invariant polynomial as a sum of polynomials in terms of $f_1'$
and of the $\Lambda_{ 1,i;\, 1^{ \lambda-2 }}^{ 2 \lambda-1}$, $2
\leqslant i \leqslant n$, $1 \leqslant \lambda \leqslant \kappa$, a
representation in which $f_1'$ is allowed to have possibly negative
powers $(f_1' )^a$ with $- \frac{ \kappa - 1}{ \kappa} m \leqslant a
\leqslant m$. The following basic statement will appear in \S5.

\smallskip\noindent{\bf Lemma.}
{\em 
In dimension $n \geqslant 1$ and for jets of order $\kappa \geqslant
1$, every polynomial ${\sf P} = {\sf P} \big( j^\kappa f\big)$
invariant by reparametrization writes under the form:
\[
\small
\aligned
\boxed{
{\sf P}\big(j^\kappa f\big)
=
\sum_{-\frac{\kappa-1}{\kappa}m\leqslant a\leqslant m}\,
(f_1')^a\,{\sf P}_a
\left(
\begin{array}{ccccc} 
f_2', & f_3', & f_4', & \dots\,, & f_n',
\\
\Lambda_{1,2}^3, & \Lambda_{1,3}^3, & \Lambda_{1,4}^3, &
\dots\,, & \Lambda_{1,n}^3, 
\\
\cdots & \cdots & \cdots & \cdots & \cdots
\\
\Lambda_{1,2;\,1^{\kappa-2}}^{2\kappa-1}, & 
\Lambda_{1,3;\,1^{\kappa-2}}^{2\kappa-1}, & 
\Lambda_{1,4;\,1^{\kappa-2}}^{2\kappa-1}, &
\dots,\, & 
\Lambda_{1,n;\,1^{\kappa-2}}^{2\kappa-1}
\end{array}
\right)
}\,,
\endaligned
\]
where the integer $a$ takes all possibly negative values in the
interval $\big[ - \frac{ \kappa - 1}{ \kappa} \, m, m \big]$, for
certain weighted homogeneous polynomials:
\[
{\sf P}_a
=
\sum_{b_2+\cdots+b_n
+
3c_2+\cdots+3c_n+
\atop
+\cdots+
(2\kappa-1)q_2+\cdots+(2\kappa-1)q_n=m-a}\!\!\!\!
{\sf coeff}\,\cdot\,
\prod_{i=2}^n\,\big(F_i\big)^{b_i}\,
\prod_{i=2}^n\,\big(A_i^3\big)^{c_i}\,
\cdots\,\,
\prod_{i=2}^n\,\big(A_i^{2\kappa-1}\big)^{q_i}
\]
of weighted degree $m - a$, namely satisfying: 
\[
{\sf P}_a
\Big(
\delta\,F_i,\,\delta^3\,A_i^3,\,\dots,\,
\delta^{2\kappa-1}\,A_i^{2\kappa-1}
\Big)
=
\delta^{m-a}\cdot
{\sf P}_a
\Big(
F_i,\,A_i^3,\,\dots,\,A_i^{2\kappa-1}
\Big).
\]

Conversely, for every collection of such weighted homogeneous
polynomials ${\sf P}_a$ in $\C \big[ F_i, A_i^3, \dots, A_i^{ 2\kappa
- 1} \big]$ of weighted degree $m - a$ indexed by an integer $a$
running in $\big[ - \frac{ \kappa - 1}{ \kappa } \, m, \, m \big]$
such that the reduction to the same denominator and the simplification
of the finite sum:
\[
{\sf R}\big(j^\kappa f\big)
=
\sum_{-\frac{\kappa-1}{\kappa}m\leqslant a\leqslant m}\,
(f_1')^a\,{\sf P}_a
\left(
\begin{array}{ccccc} 
f_2', & f_3', & f_4', & \dots\,, & f_n',
\\
\Lambda_{1,2}^3, & \Lambda_{1,3}^3, & \Lambda_{1,4}^3, &
\dots\,, & \Lambda_{1,n}^3, 
\\
\cdots & \cdots & \cdots & \cdots & \cdots
\\
\Lambda_{1,2;\,1^{\kappa-2}}^{2\kappa-1}, & 
\Lambda_{1,3;\,1^{\kappa-2}}^{2\kappa-1}, & 
\Lambda_{1,4;\,1^{\kappa-2}}^{2\kappa-1}, &
\dots,\, & 
\Lambda_{1,n;\,1^{\kappa-2}}^{2\kappa-1}
\end{array}
\right)
\]
yields a {\em true} jet \underline{\em polynomial} in $\C \big[
j^\kappa f \big]$, then ${\sf R} \big( j^\kappa f \big)$ is
a polynomial invariant by reparametrization belonging to
${\sf E}_{ \kappa, m}^n$.
}\medskip

Next, we summarize briefly the way how our algorithm works;
mathematical causalities, motivations and ``reasons-why'' shall
be transparent to any reader who will study the example
${\sf E}_4^2$ detailed in Section~6. 

Suppose that, setting aside the special invariant $f_1'$, we already
know a certain number of invariants $L^{ l_1}, \dots, L^{ l_{ k_1}}$,
for instance the very initial ones above $f_2', \dots, f_n'$ together
with all the $\Lambda_{ 1,i;\, 1^{ \lambda-2 }}^{ 2 \lambda-1}$. The
recipe is to compute the ideal of relations between these invariants
after setting $f_1' = 0$ in them:
\[
\text{\sf Ideal-Rel}
\Big(
L^{l_1}(j^\kappa f)\big\vert_{f_1'=0},\,\dots,\,
L^{l_{k_1}}(j^\kappa f)\big\vert_{f_1'=0}
\Big).
\]
Using any symbolic package for computing Gröbner bases, suppose that,
for some monomial ordering, we may dispose of a Gröbner basis for the
ideal of relations between these restricted invariants which we shall
represent shortly by the following collection of algebraic equations:
\[
0
\equiv
{\sf S}_i\Big(L^{l_1}\big\vert_0,\,\dots,\,L^{l_{k_1}}\big\vert_0\Big)
\ \ \ \ \ \ \ \ \ \ \ \ \
{\scriptstyle{(i\,=\,1\,\cdots\,N_1)}}.
\]
One checks that each ${\sf S}_i$ may be supposed to be of constant
homogeneous weight $\mu_i$, namely:
\[
{\sf S}_i\big(\delta^{l_1}A_1,\dots,\delta^{l_{k_1}}A_{k_1}\big)
=
\delta^{\mu_i}{\sf S}_i\big(A_1,\dots,A_{k_1}\big)
\ \ \ \ \ \ \ \ \ \ \ \ \
{\scriptstyle{(i\,=\,1\,\cdots\,N_1)}}.
\]
Since ${\sf S}_i \big( j^\kappa f \big)$ 
vanishes identically after setting $f_1' = 0$, when
we do not set $f_1 ' = 0$, there must exist certain (possibly zero)
polynomial remainders ${\sf R}_i \big( j^\kappa f \big)$ such that we
may write in $\C \big[ j^\kappa f \big]$:
\[
{\sf S}_i\big(L^{l_1},\dots,L^{l_{k_1}}\big)
=
(f_1')^{\nu_i}\,{\sf R}_i\big(j^\kappa f\big)
\ \ \ \ \ \ \ \ \ \ \ \ \
{\scriptstyle{(i\,=\,1\,\cdots\,N_1)}},
\]
with ${\sf R}_i \not\equiv 0$ when $1 \leqslant \nu_i < \infty$ and
with ${\sf R}_i = 0$ by convention when $\nu_i = \infty$.

Then one easily convinces oneself that every remainder ${\sf R}_i
\big( j^\kappa f\big)$ also is a polynomial invariant by
reparametrization.

Afterwards, one then tests whether the first remainder ${\sf R}_1$
belongs to the algebra generated by $L^{ l_1}, \dots, L^{ l_{ k_1}}$.
If not, ${\sf R}_1$ must be added to the list of invariants. Next, one
tests whether ${\sf R}_2$ belongs to the algebra generated by ${\sf
L}^{ l_1}, \dots, L^{ l_{ k_1}}, {\sf R}_1$. If not, one adds ${\sf
R}_2$ to the list, and so on.

At the end, one gets a new list of invariants $L^{ l_1}, 
\dots, L^{ k_1}, M^{ m_1}, \dots, M^{ m_{ k_2 }}$ and
then one restarts a second loop by computing a
Gröbner basis for the ideal of relations:
\[
\text{\sf Ideal-Rel}
\Big(
L^{l_1}\big\vert_0,\,\dots,\,L^{l_{k_1}}\big\vert_0,\,\,
M^{m_1}\big\vert_0,\,\dots,\,M^{m_{k_2}}\big\vert_0
\Big).
\]

\smallskip\noindent{\bf Theorem.}
{\em
For a certain dimension $n$ and for a certain jet order $\kappa$,
suppose that, after performing a finite number of loops of the
algorithm, one possesses a finite number $1 + M$ of mutually
independent invariants $f_1'$, $\Lambda^{l_1}$, \dots,
$\Lambda^{l_M} \in \C \big[ j^\kappa f_1, \dots, j^\kappa f_n
\big]$ of weights $1, l_1, \dots, l_M$ belonging to ${\sf
E}_\kappa^n$, whose restrictions to $\{ f_1 ' = 0 \}$ share an ideal
of relations:
\[
\text{\sf Ideal-Rel}
\Big(\
\Lambda^{l_1}\big\vert_0,\ \
\dots\dots,\,
\Lambda^{l_M}\big\vert_0\
\Big)
\]
generated by a finite number $N$ (often large) of homogeneous syzygies:
\[
0
\equiv
{\sf S}_i
\big(\Lambda^{l_1}\big\vert_0,\,\dots,\,
\Lambda^{l_M}\big\vert_0\big),
\ \ \ \ \ \ \ \ \ \ \ \ \
{\scriptstyle{(i\,=\,1\,\cdots\,N)}}
\]
of weight $\mu_i$
assumed to be represented by a certain reduced Gröbner basis $\big<
{\sf S}_i \big>_{ 1 \leqslant i \leqslant N }$ for a certain monomial
order, with the crucial property that {\em no} new invariant
appears behind $f_1'$, namely with the property that, without setting
$f_1' = 0$, one has $N$ identically satisfied relations:
\[
0
\equiv
{\sf S}_i
\big(\Lambda^{l_1},\,\dots,\,\Lambda^{l_M}\big)
-
f_1'\,{\sf R}_i
\big(f_1',\,\Lambda^{l_1},\,\dots,\, \Lambda^{l_M}\big)
\ \ \ \ \ \ \ \ \ \ \ \ \
{\scriptstyle{(i\,=\,1\,\cdots\,N)}},
\]
for some remainders ${\sf R }_i$ {\em which all depend polynomially
upon the same collection of invariants $f_1', \Lambda^{ l_1 },
\dots, \Lambda^{ l_M }$}, so that no new invariant appears at
this stage.

Then the algorithm terminates and the algebra of invariants
coincides with:
\[
\boxed{
{\sf E}_\kappa^n
=
\C\big[
f_1',\,\Lambda^{l_1},\,\dots\dots,\,\Lambda^{l_M}
\big]
\ \ \ 
\text{\sf modulo syzygies}\
}\,.
\] 
}\medskip

As a standard byproduct of basic Gröbner bases theory, one deduces a
unique representation of any polynomial invariant under
reparametrization modulo the syzygies.

Indeed, for these values of $n$ and of $\kappa$, if one denotes the
leading terms (with respect to the monomial order in question) of the
above $N$ syzygies by:
\[
{\sf LT}\big({\sf S}_i(\Lambda)\big)
=
\big(\Lambda^{l_1}\big)^{\alpha_1^i}
\cdots
\big(\Lambda^{l_M}\big)^{\alpha_M^i}
\ \ \ \ \ \ \ \ \ \ \ \ \
{\scriptstyle{(i\,=\,1\,\cdots\,N)}},
\]
for certain specific multiindices $\big( \alpha_1^i, \dots, \alpha_M^i
\big) \in \N^M$, and if for $i = 1, \dots, N$ one denotes by:
\[
\square_i
:=
\alpha^i+\N^M
=
\big\{
\big(\alpha_1^i+b_1,\,\dots,\,\alpha_M^i+b_M\big):\,
b_1,\dots,b_M\in\N^M
\big\}
\]
the positive quadrant of $\N^M$ having vertex at $\alpha^i$, then a
general, arbitrary invariant in ${\sf E}_{\kappa, m}^n$ of weight
$m$ writes {\em uniquely} under the {\em normal form}:
\[
\sum_{0\leqslant a\leqslant m}\,
(f_1')^a\,\widetilde{\sf P}_a
\big(
\Lambda^{l_1},\,\dots,\,\Lambda^{l_M}
\big),
\]
with summation containing {\em only positive powers} of $f_1'$, where
each $\widetilde{ P}_a$ is of weight $m - a$ and is put under {\em
Gröbner-normalized form}:
\[
\boxed{
\widetilde{\sf P}_a
=
\sum_{(b_1,\dots,b_M)\in\N^M\backslash
(\square_1\cup\cdots\cup\square_N)\
\atop
l_1b_1+\cdots+l_Mb_M=m-a}\,
{\sf coeff}_{a;\,b_1,\dots,b_M}\cdot
\big(\Lambda^{l_1}\big)^{b_1}\cdots
\big(\Lambda^{l_M}\big)^{b_M}}\,,
\]
with complex coefficients ${\sf coeff }_{a; \,b_1, \dots, b_M}$
subjected to no restriction at all. 

\subsection*{ The kernel algorithm}
We would like to mention that, after the paper~\cite{ mer2008a} was
completed and submitted, on the occasion of a Workshop about
holomorphic extension of CR functions and their removable
singularities organized by Berit Stens{\o}nes and John-Erik
Forn{\ae}ss at the university of Michigan (Ann Arbor, December 2007),
Harm Derksen indicated to us the so-called {\sl Van den Essen's kernel
algorithm} for locally nilpotent derivations, the goal of which is to
generate all invariants for certain one-dimensional non-reductive
actions (\cite{ ess2000, dk2002, fre2007}). Although applied here to
actions of any dimension, our algorithm here is in substance the same,
though some features will be dealt with here more explicitly in the quite
nontrivial explorations to which the paper is devoted: homogeneity of
syzygies; stepwise generation of relations; skirting of Gröbner bases
when they fail (due to oversizeness) to compute of the remainders
${\sf R}_i$; systematic restriction to $\{ f_1' = 0 \}$ to shorten
time computation. 

In a near future, we hope to set up a refined algorithm which 
would almost completely tame the disturbing expression swelling.

\subsection*{ Action of ${\sf GL}_n (\C)$ and unipotent reduction}
Lastly, we come back to explaining how one obtains Schur
decompositions of Demailly-Semple bundles.

On an arbitrary fiber ${\sf E}_{ \kappa, m}^n$ of ${\sf E}_{ \kappa ,
m}^n T_X^*$ consisting of polynomials ${\sf P} \big( j^\kappa f \big)
= {\sf P} \big( f', f'', \dots, f^{ ( \kappa)} \big)$ invariant by
reparametrization, one looks at the action of matrices ${\sf w} = (
w_{ ij} ) \in {\sf GL}_n ( \C)$ which, for each jet level $\lambda$
with $1 \leqslant \lambda \leqslant \kappa$, multiplies by ${\sf w}$
the $n$ jet-components $f^\lambda := \big( f_1^{ (\lambda )}, \dots,
f_n^{ ( \lambda )} \big)$, namely which transforms them into ${\sf w}
\cdot f^\lambda := \big( \sum_{ j=1}^n\, w_{ 1 j} f_j^{ (\lambda )},
\dots, \sum_{ j=1}^n \, w_{ nj} f_j^{ (\lambda )} \big)$ with the {\em
same} matrix for each jet level $\lambda = 1, 2, \dots, \kappa$.

According to elementary representation theory, ${\sf E}_{ k,m}^n$ then
decomposes into a certain direct sum of irreducible ${\sf
GL}_n$-representations, which are nothing but the Schur
representations $\Gamma^{ (\ell_1, \ell_2, \dots, \ell_n)}$ indexed by
integers $\ell_1 \geqslant \ell_2 \geqslant \cdots \geqslant \ell_n$.
General reasons (\cite{ dem1997}) insure that such a decomposition on
fibers globalizes coherently as a decomposition between bundles over
$X \subset \P^{ n+1} (\C)$. How then does one determine the appearing
Schur components? It suffices to look at the so-called {\sl vectors of
highest weight}, which in our situation are just the polynomials
invariant by reparametrization ${\sf P} \in {\sf E}_{ \kappa, m}^n$
which are {\em unipotent-invariant}, namely which are left untouched
after multiplication by any unipotent matrix:
\[
{\sf u}\cdot{\sf P}\big(j^\kappa f\big)
=
{\sf P}\big(j^\kappa f)
\ \ \ \ \ \ 
\text{\rm for every}\ \ 
{\sf u}
=
\left(
\begin{array}{cccc}
1 & 0 & \cdots & 0
\\
u_{21} & 1 & \cdots & 0
\\
\vdots & \vdots & \ddots & \vdots
\\
u_{n1} & u_{n2} & \cdots & 1
\end{array}
\right).
\]
Then the full space ${\sf E}_{ \kappa, m}^n$ is obtained as just the
${\sf GL}_n ( \C)$-orbit of ${\sf UE}_{ \kappa, m}^n$, and this will
correspond to somehow {\sl polarizing} the lower indices of
bi-invariants, {\em see} below. We then call {\sl bi-invariants} the
polynomials which are both invariant under reparametrization and under
the unipotent action:
\[
\boxed{
{\sf P}\big(j^\kappa(f\circ\phi)\big)
=
(\phi')^m\cdot
{\sf P}\big((j^\kappa f)\circ\phi\big)
\ \ \ \ \ \ \ \ \ \
\text{\rm and}
\ \ \ \ \ \ \ \ \ \
{\sf P}\big(\text{\sf u}\cdot j^\kappa f)
=
{\sf P}\big(j^\kappa f)}\,.
\]
Thus, the bi-invariants are nothing but vectors of highest weight for
this representation of ${\sf GL}_n ( \C)$. According to the general
theory, to each vector of highest weight corresponds one and only one
irreducible Schur representation $\Gamma^{ (\ell_1, \ell_2, \dots,
\ell_n)}$. How does one finds the integers $\ell_i$?

Suppose that, after executing the algorithm, one already knows that
${\sf UE}_\kappa^n$ is generated by a finite number $f_1',
\Lambda^{ l_1}, \dots, \Lambda^{ l_M}$ of bi-invariants of weights $1,
l_1, \dots, l_M$, and suppose that we have a {\em unique} writing:
\[
\sum_{(a,b_1,\dots,b_M)\in\mathcal{N}}\,
{\sf coeff}_{a,b_1,\dots,b_M}\,
(f_1')^a\,\big(\Lambda^{l_1}\big)^{b_1}\cdots
\big(\Lambda^{l_M}\big)^{b_M}
\]
of an arbitrary, general bi-invariant modulo the syzygies, for a
certain monomial order, where $\mathcal{ N} \subset \N^{ 1 + M}$
denotes the complement of the union of quadrants having vertex at
leading exponents. Then for every $(a, b_1, \dots, b_M)$, the single
monomial $(f_1')^a\, \big( \Lambda^{ l_1} \big)^{ b_1} \cdots \big(
\Lambda^{ l_M} \big)^{ b_M}$ is a vector of highest weight, and if one
lets a general diagonal matrix:
\[
{\sf x}
:=
\left(
\begin{array}{ccc}
x_1 & \cdots & 0
\\
\vdots & \ddots & \vdots
\\
0 & \cdots & x_n
\end{array}
\right)
\]
act on it, the theory says it necessarily is an eigenvector, 
and the eigenvalue:
\[
{\sf x}\cdot
(f_1')^a\,\big(\Lambda^{l_1}\big)^{b_1}\cdots
\big(\Lambda^{l_M}\big)^{b_M}
=
x_1^{\ell_1}\cdots x_n^{\ell_n}\,
(f_1')^a\,\big(\Lambda^{l_1}\big)^{b_1}\cdots
\big(\Lambda^{l_M}\big)^{b_M}, 
\]
exhibits the wanted $\ell_i$'s which necessarily satisfy $\ell_1
\geqslant \cdots \geqslant \ell_n$.

In conclusion, {\em both in order to understand invariants and in
order to make Euler-characteristic computations, the very main goal
is to explore algebras of bi-invariants}.

By requiring unipotent-invariance, the initial rational expression for
bi-invariants will depend upon certain determinants defined as follows
in terms of the initial invariants $\Lambda_{ 1, i: \, 1^{ \lambda -
2}}^{ 2 \lambda - 1}$.

\smallskip\noindent{\bf Theorem.}
{\em 
In dimension $n\geqslant 1$ and for jets of arbitrary order $\kappa
\geqslant 1$, every bi-invariant polynomial ${\sf BP} = {\sf BP} \big(
j^\kappa f \big)$ invariant by reparametrization and invariant under the
unipotent action writes under the form:
\[
\footnotesize
\aligned
{\sf BP}\big(j^\kappa f\big)
=\!\!
\sum_{-\frac{\kappa-1}{\kappa}m\leqslant a\leqslant m}\!\!
(f_1')^a\,
{\sf BP}_a
\bigg(
\left\vert
\begin{array}{cccc}
\Lambda_{1,2}^{2\lambda_2-1} & \Lambda_{1,3}^{2\lambda_2-1} &
\cdots & \Lambda_{1,n_1}^{2\lambda_2-1}
\\
\Lambda_{1,2}^{2\lambda_3-1} & \Lambda_{1,3}^{2\lambda_3-1} &
\cdots & \Lambda_{1,n_1}^{2\lambda_3-1}
\\
\vdots & \vdots & \ddots & \vdots
\\
\Lambda_{1,2}^{2\lambda_3-1} & \Lambda_{1,3}^{2\lambda_3-1} &
\cdots & \Lambda_{1,n_1}^{2\lambda_3-1}
\end{array}
\right\vert_{n_1=1,2\dots,n}^{2\leqslant
\lambda_2,\dots,\lambda_{n_1}\leqslant\kappa}
\bigg),
\endaligned
\]
for certain specific polynomials ${\sf BP}_a$ which depend upon ${\sf
BP} ( j^\kappa f)$.
}\medskip

\subsection*{ Algebras of bi-invariants}
As announced in the abstract, we finalized two main applications of
our algorithm. Only one bi-invariant, namely $Y^{ 27}$, was missed
in~\cite{ mer2008a}, an article which pointed out that bracketing was
insufficient.

\smallskip\noindent{\bf Theorem.}
{\em
In dimension $n = 2$ for jet order $\kappa = 5$, the algebra ${\sf
UE}_5^2$ of jet polynomials ${\sf P} \big( j^5 f_1, j^5 f_2 \big)$
invariant by reparametrization and invariant under the unipotent
action is generated by 17 mutually independent bi-invariants
explicitly defined in Section~10:
\[
\boxed{
\aligned
&
\ \ \ \ \ \ \ 
f_1',\ \ \ \ \
\Lambda^3,\ \ \ \ \
\Lambda^5,\ \ \ \ \
\Lambda^7,\ \ \ \ \
\Lambda^9,\ \ \ \ \
M^8,\ \ \ \ \
M^{10},\ \ \ \ \ 
K^{12},
\\
&
N^{12},\ \ \ \ 
H^{14},\ \ \ \ 
F^{16},\ \ \ \ 
X^{18},\ \ \ \ 
X^{19},\ \ \ \ 
X^{21},\ \ \ \ 
X^{23},\ \ \ \ 
X^{25},\ \ \ \ 
Y^{27}
\endaligned
}\,.
\]
As a consequence, the full algebra ${\sf E}_5^2$ of jet polynomials 
${\sf P} \big( j^5 f \big)$ invariant by reparametrization
is generated by the polarizations:
\[
\boxed{
\aligned
&\ \ \ \ \ 
f_i',\ \ \ \ \ 
\Lambda^3,\ \ \ \ \
\Lambda_i^5,\ \ \ \ \
\Lambda_{i,j}^7,\ \ \ \ \
\Lambda_{i,j,k}^9,\ \ \ \ \
M^8,\ \ \ \ \ 
M_i^{10},\ \ \ \ \
K_{i,j}^{12},
\\
&
N^{12},\ \ \ \ 
H_i^{14},\ \ \ \ 
F_{i,j}^{16},\ \ \ \ 
X_{i,j,k}^{18},\ \ \ \ 
X_i^{19},\ \ \ \ 
X^{21},\ \ \ \ 
X_i^{23},\ \ \ \ 
X_{i,j}^{25},\ \ \ \ 
Y_{i,j,k}^{27}
\endaligned
}
\]
of these 17 bi-invariants, where the indices $i, j, k$ vary in $\{ 1,
2 \}$, whence the total number of these invariants equals:
\[
2+1+2+4+8+1+2+4+1+2+4+8+2+1+2+4+8
=
\fbox{\bf 56}\,.
\]}

Secondly, we obtain the following new result in dimension 4.
We must confess that we were unable to discover some harmonious 
algebraic structures which could probably (in)exist?

\smallskip\noindent{\bf Theorem.}
{\em
In dimension $n = 4$ for jets of order $\kappa = 4$, the algebra ${\sf
UE}_4^4$ of jet polynomials ${\sf P} \big( j^4 f_1, j^4 f_2, j^4 f_3,
j^4f_4\big)$ invariant by reparametrization and invariant under the
unipotent action is generated by 16 mutually independent
bi-invariants explicitly defined\footnote{\, The bi-invariant
$X^{ 21}$ here is different from the $X^{ 21}$ of the preceding
theorem. } in Section~11:
\[
\boxed{
\aligned
&\ \ \ \
W^{10},\ \ \ \ \
f_1',\ \ \ \ \
\Lambda^3,\ \ \ \ \
\Lambda^5,\ \ \ \ \
\Lambda^7,\ \ \ \ \
D^6,\ \ \ \ \
D^8,\ \ \ \ \
N^{10},\ \ \ \ \ 
\\
&
M^8,\ \ \ \ \
E^{10},\ \ \ \ \
L^{12},\ \ \ \ \
Q^{14},\ \ \ \ \
R^{15},\ \ \ \ \
U^{17},\ \ \ \ \
V^{19},\ \ \ \ \
X^{21},
\endaligned}
\]
whose restriction to $\{ f_1' = 0 \}$ has a reduced gröbnerized ideal
of relations, for the Lexicographic ordering, which consists of the 41
syzygies written on p.~\pageref{41-syzygies}.

Furthermore, any bi-invariant of weight $m$ writes uniquely in the
finite polynomial form:
\[
\aligned
{\sf P}\big(j^\kappa f\big)
=
&
\sum_{o,\,p}\,
(f_1')^o\,\big(W^{10}\big)^p\,
\sum_{(a,\dots,n)\in\N^{14}\backslash
(\square_1\cup\cdots\cup\square_{41})
\atop
3a+\cdots+21n=m-o-10p}\,
{\sf coeff}_{a,\dots,n,o,p}\,\cdot
\\
&\ \ \ \ \
\cdot
\big(\Lambda^3\big)^a\,
\big(\Lambda^5\big)^b\,
\big(\Lambda^7\big)^c\,
\big(D^6\big)^d\,
\big(D^8\big)^e\,
\big(N^{10}\big)^f
\big(M^8\big)^g\,
\big(E^{10}\big)^h\,
\\
&\ \ \ \ \ \ \
\big(L^{12}\big)^i\,
\big(Q^{14}\big)^j\,
\big(R^{15}\big)^k\,
\big(U^{17}\big)^l\,
\big(V^{19}\big)^m\,
\big(X^{21}\big)^n,
\endaligned
\]
with coefficients ${\sf coeff}_{a,\dots,n,o,p}$ subjected to no
restriction, where $\square_1$, \dots, $\square_{ 41}$ denote the
quadrants in $\N^{ 14}$ having vertex at the leading terms of the 41
syzygies in question.

As a consequence, the full algebra ${\sf E}_4^4$ of jet polynomials
${\sf P} \big( j^4 f \big)$ invariant by reparametrization is
generated by the polarizations of the 16 bi-invariants:
\[
\boxed{
\aligned
&\ \ \ \ \ \ \ \
W^{10},\ \ \ \ \
f_i',\ \ \ \ \
\Lambda_{[i,j]}^3,\ \ \ \ \
\Lambda_{[i,j];\,\alpha}^5,\ \ \ \ \
\Lambda_{[i,j];\,\alpha,\beta}^7,\ \ \ \ \
D_{[i,j,k]}^6,
\\
&
D_{[i,j,k];\,\alpha}^8,\ \ \ \ \
N_{[i,j,k];\,\alpha,\beta}^{10},\ \ \ \ \ 
M_{[i,j],[k,l]}^8,\ \ \ \ \
E_{[i,j,k],[p,q]}^{10},\ \ \ \ \
L_{[i,j,k],[p,q];\,\alpha}^{12},
\\
&\ \ \ \ \ \ \ \ \ \ \
Q_{[i,j,k],[p,q];\,\alpha,\beta}^{14},\ \ \ \ \
R_{[i,j,k],[p,q,r];\,\alpha}^{15},\ \ \ \ \
U_{[i,j,k],[p,q,r],[s,t]}^{17},
\\
&\ \ \ \ \ \ \ \ \ \ \ \ \ \ \ \ \ \ \ 
V_{[i,j,k],[p,q,r],[s,t];\,\alpha}^{19},\ \ \ \ \
X_{[i,j,k],[p,q,r],[s,t];\,\alpha,\beta}^{21},
\endaligned}
\] 
These polarized invariants are skew-symmetric with respect to each
collection of bracketed indices $[i,j,k]$, $[p,q,r]$, $[s,t]$, where
the roman indices satisfy $1 \leqslant i < j < k \leqslant 4$, where
$1 \leqslant p < q < r \leqslant 4$, where $1 \leqslant s < r
\leqslant 4$ and where the two greek indices $\alpha, \beta$ satisfy
$1\leqslant \alpha , \beta\leqslant 4$ without restriction and finally
the total number of these invariants generating the Demailly-Semple
algebra ${\sf E}_4^4$ equals:
\[
\aligned
&
1
+
4+6+24+96+4+16+64
+
\\
&
+
36+24+96+384+64+96+384+1536
=
\fbox{\bf 2835}\,.
\endaligned
\]
}\medskip

\subsection*{ Acknowledgments} 
The author heartfully thanks Jacques Beigbeder (SPI, \'Ecole Normale
Supérieure) for having installed the package {\sf FGb} (Spiral Team,
LIP6) containing efficient Gröbner basis algorithms that were used to
capture the quite huge ideals of relations exhibited in the present
{\sf arxiv.org} electronic (pre)publication. He is also grateful to
Jean-Pierre Demailly, to Jahwer El Goul, to Erwan Rousseau, and to
Simone Diverio for friend(ful)ly sharing thoughts about the puzzling
complexity of invariant jet differentials. In March 2007, Tien-Cuong
Dinh and Nessim Sibony suggested the question to our expertise.

Finally, the theorem on p.~\pageref{bi-invariant-4-4} which describes
the structure of the algebra of bi-invariants for $n = k = 4$ was
firmly gained during the author's stay at the Mittag-Leffler
Institute in April 2008.

\section*{ \S2.~Invariant polynomials and composite differentation}
\label{Section-2}

\subsection*{ Fixing basic notations}
Let $X$ be a smooth $n$-dimensional complex algebraic hypersurface of
$\P^{ n+1} ( \C)$, let $\D$ be the unit disc in $\C$ and consider an
arbitrary holomorphic disc $f: \D \to X$ valued in $X$, for instance
the restriction to $\D$ of some entire holomorphic curve $\C \to X$.
In some local chart on $X \simeq \D^n$ centered at $f ( 0)$, the
$\kappa$-jet $j_0^\kappa f$ of $f$ at $0 \in \D$ is represented by the
collection of all the derivatives, with respect to the variable $\zeta
\in \D$, of the $n$ components $f_1, \dots, f_n$ of $f$, up to order
$\kappa$, that is to say:
\[
j^\kappa f
=
\big(
f_1',\dots,f_n',\,
f_1'',\dots,f_n'',\,
\dots\dots,\,
f_1^{(\kappa)},\dots,f_n^{(\kappa)}
\big);
\]
from the beginning and throughout this study, we shall in fact
constantly omit to denote the base point $0\in \D$.

\subsection*{ Polynomials invariant by reparametrization}
For $\kappa \geqslant 1$, we consider polynomials in all the
jet variables:
\[
{\sf P}
= 
{\sf P} 
\big(j^\kappa f\big)
=
{\sf P}
\big(
f_{j_1}',\,f_{j_2}'',\dots,\,
f_{j_\kappa}^{(\kappa)}
\big),
\] 
where the indices $j_1, j_2, \dots, j_\kappa$ run in $\{ 1, \dots,
n\}$. An open problem in Demailly's strategy towards the Kobayashi
hyperbolicity conjecture (\cite{ dem1997, delg2000}) was to describe
those polynomials ${\sf P} \big( j^\kappa f \big)$ enjoying the
property that a change of variable $\D \ni \zeta \longmapsto \phi (
\zeta) \in \C$ in the source affects the polynomial only through
multiplication by some power of the first derivative of $\phi$:
\[
{\sf P}
\big(j^\kappa(f\circ \phi)\big)
=
(\phi')^m\cdot
{\sc P}\big(
(j^\kappa f)\circ\phi
\big),
\] 
where $m \geqslant 1$ is an integer which shall be called here the
{\em weight}\, of ${\sf P}$.

Choosing in particular $\phi$ to be simply a dilation $\zeta \mapsto
\delta \cdot z$ by a constant nonzero complex factor $\delta$, one
sees that such polynomials must at least ({\em cf.} \cite{ gg1980}) be
{\sl weighted homogeneous of order $m$} with respect to the weighted
anisotropic dilations:
\[
{\sf P}
\big(
\delta\cdot f_{j_1'},\,\delta^2\cdot f_{j_2}'',\,\dots,\,
\delta^{\kappa}\cdot f_{j_\kappa}^{(\kappa)}
\big)
\equiv
\delta^m\cdot
{\sf P}\big(f_{j_1}',\,f_{j_2}'',\,\dots,\,
f_{j_\kappa}^{(\kappa)}
\big).
\]
As a useful mnemonic, weight therefore always counts the total number of
primes.

By ${\sf E}_{\kappa, m}^n$, we will thus denote the vector space
consisting of all such polynomials. The direct sum ${\sf E
}_\kappa^n := \bigoplus_{ m \geqslant 1} \, {\sf E}_{ \kappa, m
}^n$ forms an algebra which is graded by constancy of weights, for the
definition yields:
\[
{\sf E}_{\kappa,m_1}^n
\cdot 
{\sf E}_{\kappa,m_2}^n
\subset
{\sf E}_{\kappa,m_1+m_2}^n.
\]
Following a nowadays established terminology, a polynomial ${\sf P}
\big( j^\kappa f \big)$ in this algebra will be said to be {\sl
invariant by reparametrization}. The present article
aims to describe a complete algorithm generating all such polynomials, 
sometimes briefly called {\sl invariants}.

\subsection*{ Example}
For $\kappa = 1$, the components $f_i'$ for $i = 1, \dots, n$ of the
jet satisfy:
\[
\big(
f_i\circ\phi
\big)'
=
\phi'\cdot f_i',
\]
hence every polynomial ${\sf P} = {\sf P} \big( f_1', \dots, f_n'
\big)$ which depends only upon the first order jet $j^1 f$ is
invariant by reparametrization. So ${\sf E}_1^n$ coincides with the
plain polynomial algebra $\C \big[ f_1', \dots, f_n' \big]$.

\subsection*{ Example}
For $\kappa = 2$, aside from the monomials $f_1', \dots, f_n'$ coming
from the preceding jet level $\kappa = 1$, 
there are yet the $2 \times 2$
determinants (clearly of weight $3$):
\[
\Delta_{i,j}^{',\,''}
:=
\left\vert
\begin{array}{cc}
f_i' & f_j'
\\
f_i'' & f_j''
\end{array}
\right\vert,
\]
for one easily checks, thanks to row linear dependence, that:
\[
\left\vert
\begin{array}{cc}
(f_i\circ\phi)' & (f_j\circ\phi)'
\\
(f_i\circ\phi)'' & (f_j\circ\phi)''
\end{array}
\right\vert
=
\left\vert
\begin{array}{cc}
\phi'f_i' & \phi'f_j'
\\
\phi''f_i'+{\phi'}^2f_i'' & \phi''f_j'+{\phi'}^2f_j''
\end{array}
\right\vert
=
{\phi'}^3\cdot\left\vert
\begin{array}{cc}
f_i' & f_j'
\\
f_i'' & f_j''
\end{array}
\right\vert.
\]
It is a theorem, to be stated below, that the $f_i'$ and
the $\Delta_{ j,k}^{ ', \, ''}$ generate the
algebra ${\sf E}_n^2$.

\subsection*{ Composite differentiation up to order 
$\kappa = 5$} Setting $g_i := f_i \circ \phi$ for $i = 1, \dots, n$,
the elementary chain rule provides derivatives of $g_i$ with respect
to the source variable $\zeta \in \D$:
\[
\aligned
g_i'
&
=
\phi'f_i',
\\
g_i''
&
=
\phi''f_i'
+
{\phi'}^2f_i'',
\\
g_i'''
&
=
\phi'''f_i'
+
3\,\phi''\phi'f_i''
+
{\phi'}^3f_i''',
\\
g_i''''
&
=
\phi''''f_i'
+
4\,\phi'''\phi'f_i''
+
3\,{\phi''}^2f_i''
+
6\,\phi''{\phi'}^2f_i'''
+
{\phi'}^4f_i'''',
\\
g_i'''''
&
:=
\phi'''''f_i'
+
5\,\phi''''\phi'f_i''
+
10\,\phi'''\phi''f_i''
+
15\,{\phi''}^2\phi'f_i'''
+
\\
&
\ \ \ \ \ \ \ \ \ \ \ \ \ \ \ \ \ \ \ \ \ \ \ \
+
10\,\phi'''{\phi'}^2f_i'''
+
10\,\phi''{\phi'}^3f_i''''
+
{\phi'}^5f_i'''''.
\endaligned
\]
Thus with $\kappa = 5$ for instance,
the goal is to 
find all polynomials ${\sf P} = {\sf P} \big( j^5 g \big)$
which, after replacing $g_i'$, $g_i''$, $g_i'''$, $g_i''''$ and
$g_i'''''$ by these expressions, have the property of {\em
cancelling}\, the derivatives $\phi''$, $\phi'''$, $\phi''''$ and
$\phi '''''$ of $\phi$ whose order is $\geqslant 2$, so that ${\sf P}
\big( j^5 g \big) = {\phi '}^m {\sf P} \big( j^5 f)$ for a certain $m
\in \N$.

\smallskip

For the sake of completeness, let us present the classical {\sl Fa\`a
di Bruno}, well known in the case of one variable $\zeta \in \C$.

\smallskip\noindent{\bf Theorem.}
{\em
For every integer $\kappa \geqslant 1$, the derivative of order
$\kappa$ of each composite function $g_i (z) := f_i \circ \phi (z)$ $(
1 \leqslant i \leqslant n)$ with respect to the variable $\zeta \in
\C$ is a polynomial with integer coefficients in the derivatives of
$f_i$ (same index $i$) and in the
derivatives of $\phi${\rm :}}
\[
\boxed{
\aligned
g_i^{(\kappa)}
& 
=
\sum_{e=1}^\kappa
\
\sum_{1\leqslant\lambda_1<\cdots<\lambda_e\leqslant\kappa}
\
\sum_{\mu_1\geqslant 1,\dots,\mu_e\geqslant 1}
\
\sum_{\mu_1\lambda_1+\cdots+\mu_e\lambda_e=\kappa}
\\
& \
\ \ \ \ \ 
\frac{\kappa !}{(\lambda_1!)^{\mu_1}\ \mu_1! 
\cdots
(\lambda_e!)^{\mu_e}\ \mu_e!}
\
\big(
\phi^{(\lambda_1)}
\big)^{\mu_1}
\cdots\cdots
\big(
\phi^{(\lambda_e)}
\big)^{\mu_e}
\
f_i^{(\mu_1+\cdots+\mu_e)}
\endaligned
}\,.
\]

To read this general formula with the help of the formulas specialized
above, let us observe that the general monomial $\big(
\phi^{(\lambda_1)} \big)^{\mu_1} \cdots\cdots \big( \phi^{(\lambda_e)}
\big)^{\mu_e}$ in the reparametrization jet gathers derivatives of
increasing orders $\lambda_1 < \lambda_2 < \cdots < \lambda_e$, with
$\mu_1, \mu_2, \dots, \mu_e$ counting their respective numbers. Then
the function $f_i$ is subjected to a partial differentiation of order
$\mu_1 + \mu_2 + \cdots + \mu_e$, the total number of derivatives
$\phi^{ (\lambda_k )}$ in the monomial in question. Finally, in the
permutation group $\mathfrak{ S}_\kappa$ of $\{ 1, 2, \dots, \kappa
\}$ whose cardinality clearly equals $\kappa !$, the quantity
$(\lambda_1!)^{\mu_1} \mu_1! \cdots (\lambda_e!)^{\mu_e} \mu_e!$
counts the number of permutations which possess $\mu_1$ cycles of
length $\lambda_1$, $\mu_2$ cycles of length $\lambda_2$, {\em etc.},
$\mu_e$ cycles of length $\lambda_e$, so that the fractional
coefficient $\frac{ \kappa !}{( \lambda_1! )^{ \mu_1}\, \mu_1! \cdots
(\lambda_e! )^{ \mu_e }\, \mu_e! }$ with $\kappa = \mu_1 \lambda_1 +
\mu_2 \lambda_2 + \cdots + \mu_e \lambda_e$ is an integer which
provides the cardinality of the (left or right) coset of $\mathfrak{
S}_\kappa$ modulo such a subgroup permutations. Notice that all these
observations are confirmed by the formulas developed above up to
$\kappa = 5$.

\smallskip

With such a formula, the problem of finding all polynomials invariant
by reparametrization can be interpreted in terms of invariant theory
(\cite{ dem1997, rou2006a}).

Indeed, the group ${\sf G}_\kappa$ of $\kappa$-jets at the origin of
local reparametrizations:
\[
\phi
(\zeta)
=
\zeta
+
\phi''(0)\,
\frac{\zeta^2}{2!}
+\cdots+
\phi^{(\kappa)}(0)\,
\frac{\zeta^\kappa}{\kappa!}
+
\cdots
\] 
that are tangent to the identity, namely $\phi' ( 0 ) = 1$, may be
seen, thanks to the above formulas, to act linearly on the $n
\kappa$-tuples $\big( f_{ j_1}', f_{ j_2}'', \dots, f_{j_\kappa}^{ (
\kappa)} \big)$ just by matrix multiplication:
\[
\left(
\begin{array}{c}
g_i'
\\
g_i''
\\
g_i'''
\\
g_i''''
\\
\vdots
\\
g_i^{(\kappa)}
\end{array}
\right)
=
\left(
\begin{array}{ccccccc}
1 & 0 & 0 & 0 & \cdots & 0
\\
\phi'' & 1 & 0 & 0 & \cdots & 0
\\
\phi''' & 3\phi'' & 1 & 0 & \cdots & 0
\\
\phi'''' & 4\phi'''+3{\phi''}^2 & 6\phi'' & 1 & \cdots & 0
\\
\vdots & \vdots & \vdots & \vdots & \ddots & \vdots 
\\
\phi^{(\kappa)} & \cdots & \cdots & \cdots & \cdots & 1
\end{array}
\right)
\left(
\begin{array}{c}
f_i'
\\
f_i''
\\
f_i'''
\\
f_i''''
\\
\vdots
\\
f_i^{(\kappa)}
\end{array}
\right)
\ \ \ \ \ \ \ \ 
{\scriptstyle{(i\,=\,1\,\cdots\,n)}}.
\]
Polynomials ${\sf P} \big( j^\kappa f \big)$ invariant by
reparametrization coincide with 
the invariants for this linear group action, an
action which is clearly unipotent, hence non-reductive. In such a
context, no general theory or algorithm exists to decide whether the
algebra of invariants is finitely generated ({\em cf.} Problem~2 p.~2
in~\cite{ dk2007}). In fact, we will attack the problem from
another point of view.

\section*{ \S3.~Bracketing process and syzygies:
\\
Jacobi, Plücker~1 and Plücker~2
}
\label{Section-3}

\subsection*{ Cross product between two invariants polynomials}
A natural process known to Demailly and to El Goul ({\em cf.}~\cite{
dem2007} and~\cite{ mer2008a}) is as follows.
Suppose that we know two reparametrization-invariant polynomials ${\sf
P} = {\sf P} \big( j^\kappa g \big)$ of weight $m$ and ${\sf Q} = {\sf
Q} \big( j^\tau f)$ of weight $n$, namely we have:
\[
\aligned
{\sf P}
\big(
j^\kappa 
g
\big)
&
=
{\phi'}^m\,
{\sf P}
\big(
(j^\kappa f)
\circ
\phi
\big),
\\
{\sf Q}
\big(
j^\tau 
g
\big)
&
=
{\phi'}^n\,
{\sf Q}
\big(
(j^\tau f)
\circ
\phi
\big),
\endaligned
\]
where we have again set $g := f \circ \phi$. To differentiate a
polynomial with respect to the source variable $\zeta \in \C$ amounts
to apply to it the {\sl total differentiation operator}:
\[
{\sf D}
:=
\sum_{k=1}^n\,
\sum_{\lambda\in\N}\,
\frac{\partial (\bullet)}{\partial f_k^{(\lambda)}}\,
\cdot\,
f_k^{(\lambda+1)},
\]
which gives here:
\[
\aligned
\big[
{\sf D}{\sf P}
\big]
\big(
j^{\kappa+1}g
\big)
&
=
m\,\phi''\,
{\phi'}^{m-1}\,
{\sf P}
\big(
(j^\kappa f)\circ\phi
\big)
+
{\phi'}^m\,\phi'\,
\big[
{\sf D}{\sf P}
\big]
\big(
(j^{\kappa+1}f)\circ\phi
\big),
\\
\big[
{\sf D}{\sf Q}
\big]
\big(
j^{\tau+1}g
\big)
&
=
n\,\phi''\,
{\phi'}^{n-1}\,
{\sf Q}
\big(
(j^\kappa f)\circ\phi
\big)
+
{\phi'}^n\,\phi'\,
\big[
{\sf D}{\sf Q}
\big]
\big(
(j^{\tau+1}f)\circ\phi
\big),
\endaligned
\]
and in order to remove the second order derivative $\phi''$, it
suffices to perform a {\sl cross-product}, namely to form the $2
\times 2$ determinant:
\[
\aligned
&
\left\vert
\begin{array}{cc}
\big[{\sf D}{\sf P}\big]\big(j^{\kappa+1}g\big)\
&\
m\,{\sf P}\big(j^\kappa g\big)
\\
\big[{\sf D}{\sf Q}\big]\big(j^{\tau+1}g\big)\
&\
n\,{\sf Q}\big(j^\tau g\big)
\end{array}
\right\vert
=
\\
&
=
\left\vert
\begin{array}{cc}
m\,\phi''\,
{\phi'}^{m-1}\,
{\sf P}
\big(
(j^\kappa f)\circ\phi
\big)
+
{\phi'}^{m+1}\,
\big[
{\sf D}{\sf P}
\big]
\big(
(j^{\kappa+1}f)\circ\phi
\big)\
&\
m\,{\phi'}^m\,{\sf P}\big((j^\kappa f)\circ\phi\big)
\\
n\,\phi''\,
{\phi'}^{n-1}\,
{\sf Q}
\big(
(j^\kappa f)\circ\phi
\big)
+
{\phi'}^{n+1}\,
\big[
{\sf D}{\sf Q}
\big]
\big(
(j^{\tau+1}f)\circ\phi
\big)
&\
n\,{\phi'}^n\,{\sf Q}\big((j^\kappa f)\circ\phi\big)
\end{array}
\right\vert
\\
&
\ \ \ \ \ \ \ \ \ \ \ \ \ \ \ \ \ \ \ \ \ \
=
\left\vert
\begin{array}{cc}
{\phi'}^{m+1}\,
\big[
{\sf D}{\sf P}
\big]
\big(
(j^{\kappa+1}f)\circ\phi
\big)\
&\
m\,{\phi'}^m\,{\sf P}\big((j^\kappa f)\circ\phi\big)
\\
{\phi'}^{n+1}\,
\big[
{\sf D}{\sf Q}
\big]
\big(
(j^{\tau+1}f)\circ\phi
\big)
&\
n\,{\phi'}^n\,{\sf Q}\big((j^\kappa f)\circ\phi\big)
\end{array}
\right\vert
\\
&
\ \ \ \ \ \ \ \ \ \ \ \ \ \ \ \ \ \ \ \ \ \
\ \ \ \ \ \ \ \ \ \ \ \ \ \ \ \ \ \ \ \ \ \
=
{\phi'}^{m+n+1}\,
\left\vert
\begin{array}{cc}
\big[{\sf D}{\sf P}\big]\big(j^{\kappa+1}f\big)\
&\
m\,{\sf P}\big(j^\kappa f\big)
\\
\big[{\sf D}{\sf Q}\big]\big(j^{\tau+1}f\big)\
&\
n\,{\sf Q}\big(j^\tau f\big)
\end{array}
\right\vert,
\endaligned
\]
which therefore happens to constitute a new invariant of weight $m + n
+ 1$ in the jet space of order $1 + \max ( \kappa, \tau)$ {\em 
increased by
one unit}.

\subsection*{ Bracket operator $\big[ \cdot, \, \cdot\big]$ and
its accompanying syzygies} Thus, {\em every pair of invariants
automatically produces a new invariant}:
\[
\boxed{
\big[
{\sf P},\,{\sf Q}
\big]
:=
n\,{\sf D}{\sf P}\cdot{\sf Q}
-
m\,{\sf P}\cdot{\sf D}{\sf Q}}\,,
\]
which is obviously skew-symmetric with respect to the pair $\big( {\sf
P}, {\sf Q} \big)$. For instance, we recover in this way {\em all} the
invariants of (jet) order 2 mentioned above:
\[
\big[f_i',\,f_j'\big]
=
{\sf D}f_i'\cdot f_j'
-
f_i'{\sf D}f_j'
=
f_i''f_j'
-
f_i'f_j''
=
-\Delta_{i,j}^{',\,''},
\]
and again, we notice that bracketing increases jet order by one unit.

Certainly, as soon as at least 3 pairwise distinct invariants ${\sf
P}$, ${\sf Q}$ and ${\sf R}$ are known, a {\em Jacobi-type} identity
(checked on pp.~867--868 of \cite{ mer2008a}) must hold:
\[
\boxed{
\aligned
(\mathcal{J}ac):
\ \ \ \ \ \ \ \ \ \ \ \
0
&
\equiv
\big[\big[{\sf P},\,{\sf Q}\big],\,{\sf R}\big]
+
\big[\big[{\sf R},\,{\sf P}\big],\,{\sf Q}\big]
+
\big[\big[{\sf Q},\,{\sf R}\big],\,{\sf P}\big]
\endaligned}\,.
\]
Although such relations give nothing for jet order $\kappa = 2$,
because the jet order of an iterated bracket $\big[ [ \cdot, \, \cdot
], \, \cdot \big]$ is in any case $\geqslant 3$, if we introduce the
following new bracket-type invariants: 
\[
\aligned
\big[\Delta_{i,j}^{',\,''},\,f_k'\big]
&
=
{\sf D}\big(f_i'f_j''-f_i''f_j')\cdot f_k'
-
3\,(f_i'f_j''-f_i''f_j')\cdot f_k''
\\
&
=
(f_i'f_j'''-f_i'''f_j')\cdot f_k'
-
3\,(f_i'f_j''-f_i''f_j')\cdot f_k'',
\endaligned
\]
then we gratuitously have Jacobi-type relations 
which will hold true at the next jet level
$\kappa = 3$:
\[
0
\equiv
\big[\Delta_{i,j}^{',\,''},\,f_k'\big]
+
\big[\Delta_{k,i}^{',\,''},\,f_j'\big]
+
\big[\Delta_{j,k}^{',\,''},\,f_i'\big].
\]

On the other hand, we remind that 
the $2 \times 2$ minors ${\sf a}_{ 1,2}^{
j_1, j_2} := \det \big( {\sf a }_i^j \big)_{ i = 1, 2}^{ j= j_1, j_2}$
of an arbitrary $2 \times N$ complex-valued matrix $\big( {\sf a}_i^j
\big)_{ i=1, 2}^{ 1 \leqslant j \leqslant N}$ are known to enjoy
(\cite{ mer2008a}, p.~883) the so-called {\em quadratic Plücker
relations} which are usually organized in two families\footnote{\, In
the first line, the sum bears upon circular permutations of $(j_1,
j_2, j_3)$; in the second line, $j_3$ is fixed and the sum bears upon
circular permutations of $(j_1, j_2, j_4)$. Equivalently, 
one could have fixed $j_4$ and considered circular permutations
of $(j_1, j_2, j_3)$.}:
\[
\aligned
(\mathcal{P}lck_1):
\ \ \ \ \ \ \ \ \ \ \ \
0
&
\equiv
{\sf a}_1^{j_1}\cdot
{\sf a}_{1,2}^{j_2,j_3}
+
{\sf a}_1^{j_3}\cdot
{\sf a}_{1,2}^{j_1,j_2}
+
{\sf a}_1^{j_2}\cdot
{\sf a}_{1,2}^{j_3,j_1},
\\
(\mathcal{P}lck_2):
\ \ \ \ \ \ \ \ \ \ \ \
0
&
\equiv
{\sf a}_{1,2}^{j_1,j_2}
\cdot
{\sf a}_{1,2}^{j_3,j_4}
+
{\sf a}_{1,2}^{j_1,j_2}
\cdot
{\sf a}_{1,2}^{j_3,j_4}
+
{\sf a}_{1,2}^{j_1,j_2}
\cdot
{\sf a}_{1,2}^{j_3,j_4},
\endaligned
\]
and which may be checked by expanding the minors, just observing
cancellations\footnote{\, In fact, only these relations appear in the
ideal of syzygies between the ${\sf a}_1^j$ and the ${\sf a}_{ 1, 2}^{
j_1, j_2}$, for an
appropriate monomial order (\cite{ mist2005}, p.~277).}. We
then deduce that our bracketing process, when interpreted as computing
the minors of an auxiliary matrix:
\[
\left(
\begin{array}{ccccc}
m\,{\sf P} & n\,{\sf Q} & o\,{\sf R} & p\,{\sf S}
& \cdots
\\
{\sf D}{\sf P} & {\sf D}{\sf Q} & {\sf D}{\sf R} & {\sf D}{\sf S}
& \cdots
\end{array}
\right),
\] 
whose first line lists known invariants multiplied by their own
weight, and whose second line lists their total derivatives, we
immediately deduce that our bracketing process introduces the
following two supplementary families of identically satisfied {\sl
Plückerian-like} relations:
\[
\boxed{
\aligned
(\mathcal{P}lck_1):
\ \ \ \ \ \ \ \ \ \ \ \
0
&
\equiv
m\,{\sf P}\,\big[{\sf Q},\,{\sf R}\big]
+
o\,{\sf R}\,\big[{\sf P},\,{\sf Q}\big]
+
n\,{\sf Q}\,\big[{\sf R},\,{\sf P}\big],
\\
(\mathcal{P}lck_2):
\ \ \ \ \ \ \ \ \ \ \ \
0
&
\equiv
\big[{\sf P},\,{\sf Q}\big]
\cdot
\big[{\sf R},\,{\sf S}\big]
+
\big[{\sf S},\,{\sf P}\big]
\cdot
\big[{\sf R},\,{\sf Q}\big]
+
\big[{\sf Q},\,{\sf S}\big]
\cdot
\big[{\sf R},\,{\sf P}\big]
\endaligned
}\,.
\]
Throughout the text, identically satisfied relations between
polynomials will often be called {\sl syzygies}, following the
terminology of classical invariant theory
(\cite{ ol1999}).

\smallskip

For instance, at the jet level $\kappa = 2$, we plainly have:
\[
\aligned
0
&
\equiv
\Delta_{i,j}^{',\,''}\cdot f_k'
+
\Delta_{k,i}^{',\,''}\cdot f_j'
+
\Delta_{j,k}^{',\,''}\cdot f_i',
\\
0
&
\equiv
\Delta_{i,j}^{',\,''}\cdot\Delta_{k,l}^{',\,''}
+
\Delta_{l,i}^{',\,''}\cdot\Delta_{k,j}^{',\,''}
+
\Delta_{j,l}^{',\,''}\cdot\Delta_{k,i}^{',\,''},
\endaligned
\]
for all indices $i, j, k, l = 1, \dots, n$.

\subsection*{ A general notation for Wronskian-like determinants}
It will be quite useful to abbreviate the explicit denotation of the
further, rather complicated invariants that we shall have to deal with
in the sequel by introducing the minors:
\[
\Delta_{i,\,\,j}^{(\alpha),(\beta)}
:=
\left\vert
\begin{array}{cc}
f_i^{(\alpha)} & f_j^{(\alpha)}
\\
f_i^{(\beta)} & f_j^{(\beta)}
\end{array}
\right\vert,
\ \ \ \ \ \ \ \ \ \ \
\Delta_{i,\,\,j,\,\,k}^{(\alpha),(\beta),(\gamma)}
:=
\left\vert
\begin{array}{ccc}
f_i^{(\alpha)} & f_j^{(\alpha)} & f_k^{(\alpha)}
\\
f_i^{(\beta)} & f_j^{(\beta)} & f_k^{(\beta)}
\\
f_i^{(\gamma)} & f_j^{(\gamma)} & f_k^{(\gamma)}
\end{array}
\right\vert,
\ \ \ \ \ \ \ \ \ \
\text{\em etc.}
\]
extracted from the jet matrix $\big( f_i^{ (\lambda)}\big)$.
Top indices list derivative orders, appearing in rows.

Thanks to skew-symmetry, after some row or column permutations, one
can always write these determinants in such a way that the lower,
dimensional indices satisfy $1 \leqslant i < j < k \leqslant n$ and
similarly, the upper, derivative indices also satisfy $1 \leqslant
\alpha < \beta < \gamma \leqslant \kappa$ at the same time.

In fact, the already mentioned observation that $\Delta_{ i, j}^{ ',
\, ''}$ always provides an invariant easily generalizes, for if we set:
\[
g_i^{(\lambda)}
:=
\big(f_i\circ\phi\big)^{(\lambda)}, 
\]
then by either manipulating the Faà di Bruno formula written above, or
by using a less explicit intermediate inductive assertion 
in order to pass from one jet level to the next jet level, one may
subject the determinants to row linear combinations in order to
establish the following:

\smallskip\noindent{\bf Lemma.} {\em 
For every $\lambda$ with $1\leqslant \lambda \leqslant \kappa$ and
for all indices $i_1, i_2, \dots, i_\lambda = 1, \dots, n$, one
has:
\[
\footnotesize
\aligned
\left\vert
\begin{array}{cccc}
g_{i_1}' & g_{i_2}' & \cdots & g_{i_\lambda}'
\\
g_{i_1}'' & g_{i_2}'' & \cdots & g_{i_\lambda}''
\\
\vdots & \vdots & \ddots & \vdots
\\
g_{i_1}^{(\lambda)} & g_{i_2}^{(\lambda)} & \cdots &
g_{i_\lambda}^{(\lambda)}
\end{array}
\right\vert
=
(\phi')^{2\lambda-1}\cdot
\left\vert
\begin{array}{cccc}
f_{i_1}' & f_{i_2}' & \cdots & f_{i_\lambda}'
\\
f_{i_1}'' & f_{i_2}'' & \cdots & f_{i_\lambda}''
\\
\vdots & \vdots & \ddots & \vdots
\\
f_{i_1}^{(\lambda)} & f_{i_2}^{(\lambda)} & \cdots &
f_{i_\lambda}^{(\lambda)}
\end{array}
\right\vert
=
(\phi')^{2\lambda-1}\cdot
\Delta_{i_1,i_2,\dots,i_\lambda}^{',\,'',\,\dots,(\lambda)},
\endaligned
\]
hence all the Wronskian-like determinants $\Delta_{ i_1,i_2, \dots,
i_\lambda}^{ ', \, '', \dots, (\lambda)}$ always are invariant by
reparametrization.
}\medskip

Here, 
it is crucial that the derivative order starts from $1$ at the first
row and increases by one unit exactly while descending stepwise along
the rows; otherwise, we would {\em not in any case} get a true invariant;
for instance in the expression:
\[
\small
\aligned
\left\vert
\begin{array}{cc}
g_i'' & g_j''
\\
g_i''' & g_j'''
\end{array}
\right\vert
&
=
\left\vert
\begin{array}{cc}
\phi''f_i'+{\phi'}^2f_i''
&
\phi''f_j'+{\phi'}^2f_j''
\\
\phi'''f_i'+3\,\phi''\phi'f_i''+{\phi'}^3f_i'''\ \
&
\phi'''f_j'+3\,\phi''\phi'f_j''+{\phi'}^3f_j'''
\end{array}
\right\vert
\\
&
=
\left\vert
\begin{array}{cc}
\phi''f_i'+{\phi'}^2f_i''
&
\phi''f_j'+{\phi'}^2f_j''
\\
\phi'''f_i'-2{\phi'}^3f_i'''
&
\phi'''f_j'-2{\phi'}^3f_j'''
\end{array}
\right\vert,
\endaligned
\] 
no further simplification enables to get rid of $\phi''$, 
$\phi'''$ and such an obstruction happens to hold generally.

\subsection*{ Combinatorics of the subalgebra generated 
by the Wronskians} Thus at least, we know a large family of
invariants. The following statement goes back to 
the nineteenth century.

\smallskip\noindent{\bf Proposition.}
(\cite{ krpr1996, dem1997, mist2005})
{\em 
For jets of order $\kappa = 2$ in arbitrary dimension $n\geqslant 2$,
the algebra ${\sf E}_2^n$ consists of the algebra
generated by the $n + \frac{ n ( n-1)}{ 2}$
fundamental invariants: 
\[
f_k'
\ \ \ \ \ \ \ \ \ \
\text{\rm and}
\ \ \ \ \ \ \ \ \ \
\Delta_{i,j}^{',\,''}
=
\left\vert
\begin{array}{cc}
f_i' & f_j'
\\
f_i'' & f_j''
\end{array}
\right\vert
\]
and their syzygy ideal is generated by the two families of 
Plückerian relations written above:
\[
\left\{\,
\aligned
0
&
\equiv
\Delta_{i,j}^{',\,''}\cdot f_k'
+
\Delta_{k,i}^{',\,''}\cdot f_j'
+
\Delta_{j,k}^{',\,''}\cdot f_i',
\\
0
&
\equiv
\Delta_{i,j}^{',\,''}\cdot\Delta_{k,l}^{',\,''}
+
\Delta_{l,i}^{',\,''}\cdot\Delta_{k,j}^{',\,''}
+
\Delta_{j,l}^{',\,''}\cdot\Delta_{k,i}^{',\,''}.
\endaligned\right.
\]
}\medskip

\section*{\S4.~Survey of known descriptions of ${\sf E}_\kappa^n$ 
\\
in low dimensions for small jet levels} 
\label{Section-4}

The above-defined algebra ${\sf E}_\kappa^n$ of jet polynomials ${\sf P}
\big( j^\kappa f \big)$ invariant by reparametrization is understood
only in certain specific situations.

\subsection*{Demailly 1997} 
At first, in dimension $n \geqslant 2$ for jet level $\kappa = 2$, the
$n + \frac{ n ( n-1) }{ 2}$ generators of the proposition
just above appear on
p.~341 of~\cite{ dem1997}, namely every polynomial ${\sf P}$ in ${\sf
E}_2^n$ writes:
\[
{\sf P}
\big(j^2f)
\equiv
\mathcal{P}_{\sf P}
\big(
f_1',\dots,f_n',\,\Delta_{1,2}^{',\,''},\dots,
\Delta_{n-1,n}^{',\,\,''}\big)
\]
having as arguments the basic invariants in question.

In the particular case of surfaces, namely for $n = 2$, no syzygy
exists between $f_1'$, $f_2'$ and $\Delta_{ 1, 2}^{ ', \, ''}$, hence
${\sf E}_2^2$ coincides with a plain polynomial algebra:
\[
{\sf E}_2^2
=
\C\big[f_1',f_2',\Delta_{1,2}^{',\,''}\big].
\]

\subsection*{ Basic notions of invariant theory}
For higher $n$'s and $\kappa$'s, unpredictable syzygies will
obscure the picture, but before pursuing, we must fix a suitable
terminology. We formulate these concepts for ${\sf E}_\kappa^n$, 
but they hold quite more generally.

\smallskip\noindent{\bf Definition.}
If, for certain values of $n$ and $\kappa$, there are finitely many
invariants $\Lambda_1, \dots, \Lambda_{\sf last}$ in ${\sf
E}_\kappa^n$ with the property that every polynomial ${\sf P} \big(
j^\kappa f \big) \in {\sf E}_\kappa^n$ invariant by reparametrization
can be written as a polynomial:
\[
{\sf P}\big(j^\kappa f\big)
\equiv
\mathcal{P}_{\sf P}
\big(
\Lambda_1,\dots,\Lambda_{\sf last}
\big)
\]
having $\Lambda_1, \dots, \Lambda_{\sf last }$ as arguments,
we shall say that ${\sf E}_\kappa^n$ is {\sl generated}
(as an algebra) {\sl by $\Lambda_1, \dots, \Lambda_{ \sf last}$}. 

\smallskip\noindent{\bf Definition.}
Further, we shall say that $\Lambda_1, \dots, \Lambda_{\sf last}$ are
{\sl mutually independent} if, for every middle index with $1
\leqslant {\sf middle} \leqslant {\sf last}$, there does not exist any
polynomial $\mathcal{ P}$ such that $\Lambda_{\sf middle}$ identifies
to a polynomial:
\[
\Lambda_{\sf middle}
=
\mathcal{P}\big(
\Lambda_1,\dots,\widehat{\Lambda_{\sf middle}},\dots,\Lambda_{\sf last}
\big)
\]
in the other remaining invariants. Then $\Lambda_1, \dots,
\Lambda_{\sf last}$ will be called {\sl fundamental invariants
generating} ${\sf E}_\kappa^n$ (for such values of $n$, $\kappa$) and
an indivivual $\Lambda_{ \sf middle}$ will be called a {\sl basic
invariant}.

\smallskip
For a fixed ${\sf E}_\kappa^n$, 
all sets of fundamental invariants, either finite
or infinite, have the same cardinality.

\subsection*{Weights always appear as upper indices}
Also, we want for later use to introduce the new notation:
\[
\Lambda_{1,2}^3
:=
\Delta_{1,2}^{',\,''},
\]
where we specify the row indices $1, 2$ and where we specially
emphasize the weight $3$, counting the total number of primes. In
fact, throughout the whole paper, {\em we shall systematically write
the weight of every basic invariant as its upper index}. We thus can
continue the survey.

\subsection*{Demailly 2004; Rousseau 2006} 
Next, in dimension $n = 2$ for jet level $\kappa = 3$, it is
shown\footnote{\, The result was known to Demailly (unpublished). }
in~\cite{ rou2006a} that the algebra ${\sf E}_3^2$ is generated by the
three invariants $f_1'$, $f_2'$ and $\Delta_{ 1, 2}^{ ', \, ''}$
(already known from the preceding jet level) to which one adds the two
further invariants of weight $5$:
\[
\aligned
\Lambda_{1,2;\,1}^5
:=
&\
\big[\Lambda_{1,2}^3,\,f_1'\big]
\ \ \ \ \ \ \ \ \ \ \ \ \ \ \
\text{\rm and}
\ \ \ \ \ \ \ \ \ \ \ \ \ \ \
\Lambda_{1,2;\,2}^5
:=
\big[\Lambda_{1,2}^3,\,f_2'\big]
\\
=
&\
\Delta_{1,3}^{',\,''}\,f_1'
-
3\,\Delta_{1,2}^{',\,''}\,f_1''
\ \ \ \ \ \ \ \ \ \ \ \ \ \ \ \ \,
\ \ \ \ \ \ \ \ \ \ \ \ \ \ \
=
\Delta_{1,3}^{',\,''}\,f_2'
-
3\,\Delta_{1,2}^{',\,''}\,f_2'',
\endaligned
\]
the only possible brackets, as one checks. Moreover, these five
invariants $f_1'$, $f_2'$, $\Lambda_{ 1, 2}^3$, $\Lambda_{ 1, 2; \, 1}^5$ and $\Lambda_{
1, 2; \, 2}^5$ are mutually
independent and their syzygy ideal is
{\em principal}, generated by the single quadratic relation:
\[
0
\equiv
3\,\Lambda_{1,2}^3\Lambda_{1,2}^3
-
f_2'\Lambda_{1,2;\,1}^5
+
f_1'\Lambda_{1,2;\,2}^5.
\] 
One sees that this syzygy just comes 
($\mathcal{ P}lck_2$). In fact, 
($\mathcal{ J}ac$) and
($\mathcal{ P}lck_1$) give nothing.

\subsection*{Rousseau 2006} 
Now, in dimension $n = 3$ and for jet level
$\kappa = 3$, applying a theorem of Popov, Rousseau (\cite{ rou2006a},
p.~403) deduced that the algebra ${\sf E}_3^3$ is generated by all the
invariants known in dimension $n - 1 = 2$ whose lower indices are {\em
polarized} in all possible ways, namely the 15 invariants:
\[
\small
\aligned
&
f_1',\ \ \ \ \
f_2',\ \ \ \ \
f_3',\ \ \ \ \
\Lambda_{1,2}^3,\ \ \ \ \
\Lambda_{1,3}^3,\ \ \ \ \
\Lambda_{2,3}^3,
\\
\Lambda_{1,2;\,1}^5,\ \ \ \ \
\Lambda_{1,3;\,1}^5,\ \ \ \ \
&
\Lambda_{2,3;\,1}^5,\ \ \ \ \
\Lambda_{1,2;\,2}^5,\ \ \ \ \
\Lambda_{1,3;\,2}^5,\ \ \ \ \
\Lambda_{2,3;\,2}^5,\ \ \ \ \
\Lambda_{1,2;\,3}^5,\ \ \ \ \
\Lambda_{1,3;\,3}^5,\ \ \ \ \
\Lambda_{2,3;\,3}^5,
\endaligned
\]
together with a single further invariant, the Wronskian:
\label{theorem-3-3}
\[
D_{1,2,3}^6
:=
\left\vert
\begin{array}{ccc}
f_1' & f_2' & f_3' 
\\
f_1'' & f_2'' & f_3''
\\
f_1''' & f_2''' & f_3'''
\end{array}
\right\vert.
\]
This makes 16 invariants in sum. An alternative, direct proof of this
result may be found in~\cite{ mer2008a}, and will also pop up again in
the present paper.

We must mention that the Wronskian $D_{ 1, 2, 3}^6$ also appears in
fact in terms of brackets, for one checks the following three relations:
\[
\aligned
\big[\Lambda_{1,2}^3,\,\Lambda_{1,3}^3\big]
&
=
-3\,f_1'\,D_{1,2,3}^6,
\ \ \ \ \ \ \ \ \ \
\big[\Lambda_{1,2}^3,\,\Lambda_{2,3}^3\big]
=
-3\,f_2'\,D_{1,2,3}^6,
\\
&
\ \ \ \ \ \ \ \ \
\big[\Lambda_{1,3}^3,\,\Lambda_{2,3}^3\big]
=
-3\,f_3'\,D_{1,2,3}^6.
\endaligned
\]
 
A Maple computation (\cite{ rou2004}) also provided the ideal of
relations between these 16 invariants. Among the 62 generators of the
reduced
Gröbner basis supplied by Maple after $\sim 15$ hours of symbolic
computations, 30 appear to be minimal generators of the ideal of
relations, the 32 remaining ones being further automatically generated
S-polynomials which are required to complete the basis. Remarkably,
it may be checked (\cite{ mer2008a}) that
{\em each one of the 30 minimal syzygies} in question
is included among the collection of
syzygies deduced from our three fundamental families by inserting
$f_i'$ and $\Lambda_{ j,k}^3$ in all possible ways in place of ${\sf P}$,
${\sf Q}$, ${\sf R}$ and ${\sf T}$:
\[
\label{syzygies-3-3}
(\mathcal{J}ac):
\ \ \ \ \ 
\Big\{
\aligned
0
&
\equiv
\big[[f_i',\,f_j'],\,f_k'\big]
+
\big[[f_k',\,f_i'],\,f_j'\big]
+
\big[[f_j',\,f_k'],\,f_i'\big],
\endaligned
\]
\[
\small
(\mathcal{P}lck_1):
\ \ \ \ \ 
\left\{
\aligned
0
&
\equiv
f_i'\big[f_j',\,f_k']
+
f_k'\big[f_i',\,f_j']
+
f_j'\big[f_k',\,f_i'],
\\
0
&
\equiv
f_i'\big[f_j',\,\Lambda_{k,l}^3\big]
+
3\,\Lambda_{k,l}^3\big[f_i',\,f_j'\big]
+
f_j'\big[\Lambda_{k,l}^3,\,f_i'\big],
\\
0
&
\equiv
f_i'\big[\Lambda_{j,k}^3,\,\Lambda_{l,m}^3\big]
+
\Lambda_{l,m}^3\big[f_i',\,\Lambda_{j,k}^3\big]
+
\Lambda_{j,k}^3\big[\Lambda_{l,m}^3,\,f_i'\big],
\\
0
&
\equiv
\Lambda_{i,j}^3\big[\Lambda_{k,l}^3,\,\Lambda_{m,n}^3\big]
+
\Lambda_{m,n}^3\big[\Lambda_{i,j}^3,\,\Lambda_{k,l}^3\big]
+
\Lambda_{k,l}^3\big[\Lambda_{m,n}^3,\,\Lambda_{i,j}^3\big],
\endaligned\right.
\]
\[
\footnotesize
(\mathcal{P}lck_2):
\ \ \ \ \ 
\left\{
\aligned
0
&
\equiv
\big[f_i',\,f_j'\big]\cdot\big[f_k',\Lambda_{l,m}^3\big]
+
\big[\Lambda_{l,m}^3,\,f_i'\big]\cdot\big[f_k',f_j'\big]
+
\big[f_j',\,\Lambda_{l,m}^3\big]\cdot\big[f_k',f_i'\big],
\\
0
&
\equiv
\big[f_i',\,f_j'\big]\cdot\big[\Lambda_{k,l}^3,\,\Lambda_{m,n}^3\big]
+
\big[\Lambda_{m,n}^3,\,f_i'\big]\cdot\big[\Lambda_{k,l}^3,\,f_j'\big]
+
\big[f_j',\,\Lambda_{m,n}^3\big]\cdot\big[\Lambda_{k,l}^3,\,f_i'\big],
\\
0
&
\equiv
\big[f_i',\,\Lambda_{j,k}^3\big]\cdot
\big[\Lambda_{l,m}^3,\,\Lambda_{n,p}^3\big]
+
\big[\Lambda_{n,p}^3,\,f_i'\big]\cdot
\big[\Lambda_{l,m}^3,\,\Lambda_{j,k}^3\big]
+
\big[\Lambda_{j,k}^3,\,\Lambda_{n,p}^3\big]
\cdot\big[\Lambda_{l,m}^3,\,f_i'\big],
\endaligned\right.
\]
where the indices $i$, $j$, $k$, $l$, $m$, $n$, and $p$ take all the
values $1, 2, 3$.

\subsection*{Demailly-El Goul 2004; Rousseau
2007; M. 2007} Finally\footnote{\, The result was known (unpublished)
to experts; a proof appears in \cite{ mer2008a}. }, for jets of order
$\kappa = 4$ in dimension $n = 2$, the algebra ${\sf E}_4^2$ is
generated by the five invariants:
\[
f_1',\ \ \ \ \
f_2',\ \ \ \ \
\Lambda_{1,2}^3,\ \ \ \ \ 
\Lambda_{1,2;\,1}^5,\ \ \ \ \
\Lambda_{1,2;\,2}^5
\]
already known from the preceding jet level, 
to which one adds the four further invariants
gently provided by bracketing:
\[
\small
\aligned
\Lambda_{1,1}^7
:=
\big[\Lambda_{1,2;\,1}^5,\,f_1'\big],
\ \ \ \ \ \ \ \ 
\Lambda_{1,2}^7
&
:=
\big[\Lambda_{1,2;\,1}^5,\,f_2'\big]
=
\big[\Lambda_{1,2;\,2}^5,\,f_1'\big]
=
\Lambda_{2,1}^7,
\ \ \ \ \ \ \ \ 
\Lambda_{2,2}^7
:=
\big[\Lambda_{1,2;\,2}^5,\,f_2'\big],
\\
M^8
&
:=
\frac{1}{f_1'}\,
\big[\Lambda_{1,2;\,1}^5,\,\Lambda_{1,2}^3\big].
\endaligned
\]
This in sum makes 9 fundamental invariants. Notice the (necessary)
division by $f_1'$ to get $M^8$. The two missing brackets\footnote{\,
Details of computations may be found in~\cite{ mer2008a}, pp.~870--871
and also pp.~882--886. }:
\[
\big[\Lambda_{1,2;\,2}^5,\,\Lambda_{1,2}^3\big]
=
f_2'\,M^8
\ \ \ \ \ \ \ \ 
\text{\rm and}
\ \ \ \ \ \ \ \ 
\big[\Lambda_{1,2;\,1}^5,\,\Lambda_{1,2;\,2}^5\big]
=
\Lambda_{1,2}^3\,M^8
\]
appear to in fact belong already to the algebra generated by these 9
invariants.

Now, we lighten a little the notation by dropping some of the lower
indices, especially in the $\Delta_{ 1, 2}^{ ( \alpha), (\beta)}
\equiv : \Delta^{ (\alpha), (\beta)}$, because in dimension $n = 2$,
by skew-symmetry of determinants, only $(1,2)$ can appear at 
the bottom.

\smallskip\noindent{\bf Theorem.} 
\label{theorem-2-4}
(\cite{ mer2008a})
{\em 
For jets of order $\kappa = 4$ in dimension $n = 2$, the algebra ${\sf
E}_4^2$ is generated by 9 mutually 
independent fundamental invariants
explicitly defined by:
\[
\small
\aligned
&
f_1',\ \ \ \ \ \ \ \ \ \ \ \ \
f_2', \ \ \ \ \ \ \ \ \ \ \ \ \ \ \
\Lambda^3
:=
\Delta^{',\,''},
\endaligned
\]
\[
\small
\aligned
\Lambda_1^5
&
:=
\Delta^{',\,'''}\,f_1'
-
3\,\Delta^{',\,''}\,f_1'',
\\
\Lambda_2^5
&
:=
\Delta^{',\,'''}\,f_2'
-
3\,\Delta^{',\,''}\,f_2'',
\endaligned
\]
\[
\small
\aligned
\Lambda_{1,1}^7
&
:=
\big(
\Delta^{',\,''''}
+
4\,
\Delta^{'',\,'''}\big)\,f_1'f_1'
-
10\,\Delta^{',\,'''}\,f_1'f_1''
+
15\,\Delta^{',\,''}\,f_1''f_1'',
\\
\Lambda_{1,2}^7
&
:=
\big(
\Delta^{',\,''''}
+
4\,
\Delta^{'',\,'''}\big)\,f_1'f_2'
-
5\,\Delta^{',\,'''}\big(f_1''f_2'
+
f_2''f_1'\big)
+
15\,\Delta^{',\,''}\,f_1''f_2'',
\\
\Lambda_{2,2}^7
&
:=
\big(\Delta^{',\,''''}
+
4\,\Delta^{'',\,'''}\big)\,f_2'f_2'
-
10\,\Delta^{',\,'''}\,f_2'f_2''
+
15\,\Delta^{',\,''}\,f_2''f_2'',
\endaligned
\]
\[
\small
M^8
:=
3\,\Delta^{',\,''''}\,\Delta^{',\,''}
+
12\,\Delta^{'',\,'''}\,\Delta^{',\,''}
-
5\,\Delta^{',\,'''}\,\Delta^{',\,'''}
\]
whose ideal of relations is generated by 9 
fundamental syzygies:
\[
\aligned
\left[
0
\overset{1}{\equiv}
f_2'\,\Lambda_1^5
-
f_1'\,\Lambda_2^5
-
3\,\Lambda^3\,\Lambda^3,
\right.
\endaligned
\]
\[
\left[
\aligned
0
&
\overset{2}{\equiv}
f_2'\,\Lambda_{1,1}^7
-
f_1'\,\Lambda_{1,2}^7
-
5\,\Lambda^3\,\Lambda_1^5,
\\
0
&
\overset{3}{\equiv}
f_2'\,\Lambda_{1,2}^7
-
f_1'\,\Lambda_{2,2}^7
-
5\,\Lambda^3\,\Lambda_2^5,
\endaligned\right.
\]
\[
\left[
\aligned
0
&
\overset{4}{\equiv}
f_1'\,f_1'\,M^8
-
3\,\Lambda^3\,\Lambda_{1,1}^7
+
5\,\Lambda_1^5\,\Lambda_1^5,
\\
0
&
\overset{5}{\equiv}
f_1'\,f_2'\,M^8
-
3\,\Lambda^3\,\Lambda_{1,2}^7
+
5\,\Lambda_1^5\,\Lambda_2^5,
\\
0
&
\overset{6}{\equiv}
f_2'\,f_2'\,M^8
-
3\,\Lambda^3\,\Lambda_{2,2}^7
+
5\,\Lambda_2^5\,\Lambda_2^5,
\endaligned\right.
\]
\[
\left[
\aligned
0
&
\overset{7}{\equiv}
f_1'\,\Lambda^3\,M^8
-
\Lambda_1^5\,\Lambda_{1,2}^7
+
\Lambda_2^5\,\Lambda_{1,1}^7,
\\
0
&
\overset{8}{\equiv}
f_2'\,\Lambda^3\,M^8
-
\Lambda_1^5\,\Lambda_{2,2}^7
+
\Lambda_2^5\,\Lambda_{1,2}^7,
\endaligned\right.
\]
\[
\left[
\aligned
0
&
\overset{9}{\equiv}
5\,\Lambda^3\,\Lambda^3\,M^8
-
\Lambda_{2,2}^7\,\Lambda_{1,1}^7
+
\Lambda_{1,2}^7\,\Lambda_{1,2}^7,
\endaligned\right.
\]
which are all obtained by means of the three families of automatic
relations $(\mathcal{ J} ac )$, $(\mathcal{ P}lck_1)$ and $(\mathcal{
P}lck_2 )$.
}

\subsection*{ Summary and induction}
Thus, all known descriptions of algebras of jet polynomials
invariant by reparametrization were obtained by 
starting with the trivial
list:
\[
f_1',\,f_2',\,\dots,\,f_n'
\]
of invariants of order 1, and bracketing them again and again in order
to lift oneself to higher jet levels. The principle of induction
then dictates to continue such a process. 

\subsection*{ Jets of order $\kappa = 5$ in dimension $n = 2$}
Bracketing all invariants from the preceding jet level $\kappa = 4$
amounts to compute all the $2 \times 2$ minors of the following $2
\times 9$ matrix:
\[
\small
\left\vert\!\left\vert
\begin{array}{ccccccccc}
f_1' \ \ & \ \ f_2' \ \ & \ \ 3\,\Lambda^3 \ \ & 
\ \ 5\,\Lambda_1^5 \ \ & \ \ 5\,\Lambda_2^5 \ \ & \ \ 
7\,\Lambda_{1,1}^7 \ \ & \ \ 7\,\Lambda_{1,2}^7 \ \ & \ \ 
7\,\Lambda_{2,2}^7 \ \ & \ \ 8\,M^8
\\
{\sf D}f_1' \ \ & \ \ {\sf D}f_2' \ \ & \ \
{\sf D}\Lambda^3 \ \ & \ \ {\sf D}\Lambda_1^5 \ \ & \ \
{\sf D}\Lambda_2^5 \ \ & \ \ {\sf D}\Lambda_{1,1}^7 \ \ & \ \ 
{\sf D}\Lambda_{1,2}^7 \ \ & \ \ {\sf D}\Lambda_{2,2}^7 \ \ & \ \
{\sf D}M^8
\end{array}
\right\vert\!\right\vert,
\]
which in sum makes a total of $\frac{ 9!}{ 2! \, 7!} = 36$ brackets.
But taking account of the fact that the $\frac{ 5!}{ 2! \, 3!} = 10$
minors of the first 5 columns correspond to the already known passage
from $\kappa = 3$ to $\kappa = 4$, just a few less brackets, namely
$36 - 10 = 26$ have to be computed, namely the
eight families: 
\[
\begin{array}{cc}
\big[\Lambda_{i,j}^7,\,f_k'\big],
\ \ \ \ \ \ \ \ \ \ \ \ \ \ \ 
&
\ \ \ \ \ \ \ \ \ \ \ \ \ \ \ 
\big[M^8,\,f_i'\big],
\\
\big[\Lambda_{i,j}^7,\,\Lambda^3\big],
\ \ \ \ \ \ \ \ \ \ \ \ \ \ \ 
&
\ \ \ \ \ \ \ \ \ \ \ \ \ \ \ 
\big[M^8,\,\Lambda^3\big],
\\
\big[\Lambda_{i,j}^7,\,\Lambda_k^5\big],
\ \ \ \ \ \ \ \ \ \ \ \ \ \ \ 
&
\ \ \ \ \ \ \ \ \ \ \ \ \ \ \ 
\big[M^8,\,\Lambda_i^5\big],
\\
\big[\Lambda_{i,j}^7,\,\Lambda_{k,l}^7\big],
\ \ \ \ \ \ \ \ \ \ \ \ \ \ \ 
&
\ \ \ \ \ \ \ \ \ \ \ \ \ \ \ 
\big[M^8,\,\Lambda_{i,j}^7\big].
\end{array}
\]
In~\cite{ mer2008a}, this task was achieved, thoroughly and in great
details, the obtained brackets being all written in terms of the
$\Delta^{ ( \alpha), (\beta)}$. Furthermore, by inspecting
systematically the first fundamental family\footnote{\, The other two
families of syzygies $(\mathcal{ P} lck_1)$ and $(\mathcal{ P}lck_2)$
having all their terms quadratic, no resolved relation for any bracket
invariant $\Pi$ of the form $\Pi = {\sf polynomial} \, ( \Lambda^1,
\dots, \Lambda_{ \sf last})$ can arise from them. } of syzygies
$(\mathcal{ J}ac)$, some superfluous brackets that are certain
polynomials in terms of previously known invariants were left out.

\smallskip\noindent{\bf Theorem.}
\label{bracket-5-2}
(\cite{ mer2008a})
{\em For jets of order $\kappa = 5$ in dimension $n = 2$, the algebra
of bracket invariants in ${\sf E}_5^2$
is generated by exactly {\bf 24} mutually independent
fundamental invariants:
\[
\boxed{
\aligned
f_1',\ \
f_2',\ \
\Lambda^3,\ \
\Lambda_1^5,\ \
\Lambda_2^5,\ \
\Lambda_{1,1}^7,\ \
\Lambda_{1,2}^7,\ \
\Lambda_{2,2}^7,\ \
M^8,\
\\
\Lambda_{1,1,1}^9,\ \
\Lambda_{1,2,1}^9,\ \
\Lambda_{2,1,2}^9,\ \
\Lambda_{2,2,2}^9,\ \
M_1^{10},\ \
M_2^{10},\
\\
N^{12}, \ \
K_{1,1}^{12},\ \
K_{1,2}^{12}=
K_{2,1}^{12},\ \
K_{2,2}^{12},\
\\
H_1^{14},\ \
H_2^{14},\ \
F_{1,1}^{16},\ \
F_{1,2}^{16},\ \
F_{2,2}^{16},\
\endaligned}
\]
among which the pure order 5 brackets are
defined by:
\[
\blue{
\aligned
\green{\Lambda_{i,j,k}^9}
&
:=
\big[\Lambda_{i,j}^7,\,f_k'\big]
\\
\green{M_i^{10}}
&
:=
\big[M^8,\,f_i'\big]
\\
\green{N^{12}}
&
:=
\big[M^8,\,\Lambda^3\big]
\\
\green{K_{i,j}^{12}}
&
:=
\big[\Lambda_{i,j}^7,\,\Lambda_1^5\big]\red{\big/f_i'}
\\
\green{H_i^{14}}
&
:=
\big[M^8,\,\Lambda_i^5\big]
\\
\green{F_{i,j}^{16}}
&
:=
\big[M^8,\,\Lambda_{i,j}^7\big]
\endaligned}
\]
\[
\]
and are explicitly given by the following normalized formulas:
\[
\small
\aligned
\Lambda_{i,j,k}^9
&
:=
\Delta^{',\,'''''}\,f_i'f_j'f_k'
+
5\,\Delta^{'',\,''''}\,f_i'f_j'f_k'
-
\\
&
\ \ \ \ \
-
4\,\Delta^{',\,''''}\big(f_i''f_j'+f_i'f_j'')\,f_k'
-
7\,\Delta^{',\,''''}\,f_i'f_j'f_k''
-
\\
&
\ \ \ \ \
-
16\,\Delta^{'',\,'''}\big(f_i''f_j'+f_i'f_j''\big)f_k'
-
28\,\Delta^{'',\,'''}\,f_i'f_j'f_k''
-
\\
&
\ \ \ \ \
-
5\,\Delta^{',\,'''}\big(f_i'''f_j'+f_i'f_j''')\,f_k'
+
35\,\Delta^{',\,'''}\big(f_i''f_j''f_k'+f_i''f_j'f_k''
+f_i'f_j''f_k''\big)
-
\\
&
\ \ \ \ \
-
105\,\Delta^{',\,''}\,f_i''f_j''f_k'',
\endaligned
\]
\[
\small
\aligned
M_i^{10}
&
:=
\big[
3\,\Delta^{',\,'''''}\,\Delta^{',\,''}
+
15\,\Delta^{'',\,''''}\,\Delta^{',\,''}
-
7\,\Delta^{',\,''''}\,\Delta^{',\,'''}
+
2\,\Delta^{'',\,'''}\,\Delta^{',\,'''}
\big]\,f_i'
-
\\
&
\ \ \ \ \
-
\big[
24\,\Delta^{',\,''''}\,\Delta^{',\,''}
+
96\,\Delta^{'',\,'''}\,\Delta^{',\,''}
-
40\,\Delta^{',\,'''}\,\Delta^{',\,'''}
\big]f_i'',
\endaligned
\]
\[
\small
\aligned
N^{12}
&
:=
9\,\Delta^{',\,'''''}\,\Delta^{',\,''}\,\Delta^{',\,''}
+
45\,\Delta^{'',\,''''}\,\Delta^{',\,''}\,\Delta^{',\,''}
-
45\,\Delta^{',\,''''}\,\Delta^{',\,'''}\,\Delta^{',\,''}
-
\\
&
\ \ \ \
-
90\,\Delta^{'',\,'''}\,\Delta^{',\,'''}\,\Delta^{',\,''}
+
40\,\Delta^{',\,'''}\,\Delta^{',\,'''}\,\Delta^{',\,'''},
\endaligned
\]
\[
\footnotesize
\aligned
K_{i,j}^{12}
&
:=
f_i'f_j'
\Big(
5\,\Delta^{',\,'''''}\,\Delta^{',\,'''}
+
25\,\Delta^{'',\,''''}\,\Delta^{',\,'''}
-
7\,\Delta^{',\,''''}\,\Delta^{',\,''''}
-
56\,\Delta^{'',\,'''}\,\Delta^{',\,''''}
-
112\,\Delta^{'',\,'''}\,\Delta^{'',\,'''}
\Big)
+
\\
&
\ \ \ \ \
+
\frac{(f_i'f_j''+f_i''f_j')}{2}
\Big(
-15\,\Delta^{',\,'''''}\,\Delta^{',\,''}
-
75\,\Delta^{'',\,''''}\,\Delta^{',\,''}
+
65\,\Delta^{',\,''''}\,\Delta^{',\,'''}
+
110\,\Delta^{'',\,'''}\,\Delta^{',\,'''}
\Big)
+
\\
&
\ \ \ \ \
+
\frac{(f_i'f_j'''+f_i'''f_j')}{2}
\Big(
-50\,\Delta^{',\,'''}\,\Delta^{',\,'''}
\Big)
+
\\
&
\ \ \ \ \
+
f_i''f_j''
\Big(
-25\,\Delta^{',\,'''}\,\Delta^{',\,'''}
+
15\,\Delta^{',\,''''}\,\Delta^{',\,''}
+
60\,\Delta^{'',\,'''}\,\Delta^{',\,''}
\Big),
\endaligned
\]
\[
\footnotesize
\aligned
H_i^{14}
&
:=
\Big(
15\,\Delta^{',\,'''''}\,\Delta^{',\,'''}\,\Delta^{',\,''}
+
75\,\Delta^{'',\,''''}\,\Delta^{',\,'''}\,\Delta^{',\,''}
+
5\,\Delta^{',\,''''}\,\Delta^{',\,'''}\,\Delta^{',\,'''}
+
\\
&
\ \ \ \ \
+
170\,\Delta^{'',\,'''}\,\Delta^{',\,'''}\,\Delta^{',\,'''}
-
24\,\Delta^{',\,''''}\,\Delta^{',\,''''}\,\Delta^{',\,''}
-
192\,\Delta^{',\,''''}\,\Delta^{'',\,'''}\,\Delta^{',\,''}
-
\\
&
\ \ \ \ \
-
384\,\Delta^{'',\,'''}\,\Delta^{'',\,'''}\,\Delta^{',\,''}
\Big)\,f_i'
+
\Big(
-
45\,\Delta^{',\,'''''}\,\Delta^{',\,''}\,\Delta^{',\,''}
-
225\,\Delta^{'',\,''''}\,\Delta^{',\,''}\,\Delta^{',\,''}
+
\\
&
\ \ \ \ \
+
225\,\Delta^{',\,''''}\,\Delta^{',\,'''}\,\Delta^{',\,''}
+
450\,\Delta^{'',\,'''}\,\Delta^{',\,'''}\,\Delta^{',\,''}
-
200\,\Delta^{',\,'''}\,\Delta^{',\,'''}\,\Delta^{',\,'''}
\Big)\,f_i'',
\endaligned
\]
\[
\footnotesize
\aligned
F_{i,j}^{16}
&
:=
\Big(
-3\,\Delta^{',\,'''''}\,\Delta^{',\,''''}\,\Delta^{',\,''}
-
15\,\Delta^{'',\,''''}\,\Delta^{',\,''''}\,\Delta^{',\,''}
-
12\,\Delta^{',\,'''''}\,\Delta^{'',\,'''}\,\Delta^{',\,''}
+
\\
&
\ \ \ \ \ 
+
40\,\Delta^{',\,'''''}\,\Delta^{',\,'''}\,\Delta^{',\,'''}
-
60\,\Delta^{'',\,''''}\,\Delta^{'',\,'''}\,\Delta^{',\,''}
+
200\,\Delta^{'',\,''''}\,\Delta^{',\,'''}\,\Delta^{',\,'''}
-
\\
&
\ \ \ \ \
-
49\,\Delta^{',\,''''}\,\Delta^{',\,''''}\,\Delta^{',\,'''}
-
422\,\Delta^{',\,''''}\,\Delta^{'',\,'''}\,\Delta^{',\,'''}
-
904\,\Delta^{'',\,'''}\,\Delta^{'',\,'''}\,\Delta^{',\,'''}
\Big)f_i'f_j'
+
\endaligned
\]
\[
\footnotesize
\aligned
\\
&
\ \ \ \ \ 
+
\Big(
-105\,\Delta^{',\,'''''}\,\Delta^{',\,'''}\,\Delta^{',\,''}
-
525\,\Delta^{'',\,''''}\,\Delta^{',\,'''}\,\Delta^{',\,''}
+
205\,\Delta^{',\,''''}\,\Delta^{',\,'''}\,\Delta^{',\,'''}
-
\\
&
\ \ \ \ \
-
230\,\Delta^{'',\,'''}\,\Delta^{',\,'''}\,\Delta^{',\,'''}
+
96\,\Delta^{',\,''''}\,\Delta^{',\,''''}\,\Delta^{',\,''}
+
768\,\Delta^{',\,''''}\,\Delta^{'',\,'''}\,\Delta^{',\,''}
+
\\
&
\ \ \ \ \
+
1536\,\Delta^{'',\,'''}\,\Delta^{'',\,'''}\,\Delta^{',\,''}
\Big)
\big(
f_i''f_j'+f_i'f_j''
\big)
+
\endaligned
\]
\[
\footnotesize
\aligned
&
\ \ \ \ \
+
\Big(
-200\,\Delta^{',\,'''}\,\Delta^{',\,'''}\,\Delta^{',\,'''}
\Big)\big(
f_i'''f_j'+f_i'f_j'''
\big)
+
\\
&
\ \ \ \ \ 
+
\Big(
315\,\Delta^{',\,'''''}\,\Delta^{',\,''}\,\Delta^{',\,''}
+
1575\,\Delta^{'',\,''''}\,\Delta^{',\,''}\,\Delta^{',\,''}
-
1575\,\Delta^{',\,''''}\,\Delta^{',\,'''}\,\Delta^{',\,''}
-
\\
&
\ \ \ \ \
-
3150\,\Delta^{'',\,'''}\,\Delta^{',\,'''}\,\Delta^{',\,''}
+
1400\,\Delta^{',\,'''}\,\Delta^{',\,'''}\,\Delta^{',\,'''}
\Big)f_i''f_j'',
\endaligned
\]
where the indices $i$, $j$ and $k$ run in $\{ 1, 2 \}$. Furthermore,
the ideal of relations between these 24 fundamental bracket invariants
consists of all the syzygies that one obtains\footnote{\, The data of our
manuscript are not reproduced here. } by substituting in $( \mathcal{
P } lck_1 )$ or in $( \mathcal{ P }lck_2 )$ for ${\sf P}$, ${\sf Q}$,
${\sf R}$, ${\sf T}$ three or four among the nine invariants $f_1'$,
$f_2'$, $\Lambda^3$, $\Lambda_1^5$, $\Lambda_2^5$, $\Lambda_{1, 1}^7$,
$\Lambda_{ 1, 2}^7$, $\Lambda_{ 2, 2}^7$, $M^8$, in all possible ways,
which makes in sum:
\[
{\textstyle{\frac{9!}{3!\,6!}+\frac{9!}{4!\,5!}}}
=
84+126
=
210
\]
generating syzygies.
}\medskip

It is now great time to offer ideas, arguments, principles of
computations, and also proofs.

\section*{\S5.~Initial invariants in dimension $n$ 
\\
for arbitrary jet level $\kappa \geqslant 1$}
\label{Section-5}

\subsection*{ Reparametrizing by $f_1^{ -1}$}
To fix ideas and to better offer the intuition of our computations, we
shall firstly work in dimension $n = 2$ until everything about the
first basic step becomes clear, so that afterwards, the description of
the birth of the initial invariants in the higher dimensions $n = 3$
and $n = 4$ shall present no real difficulty.

Thus, let ${\sf P} \big( j^\kappa f_1, \, j^\kappa f_2)$ be a
polynomial of weight $m$ that is invariant by reparametrization. By
definition,
\def\theequation{$*$}\begin{equation}
{\sf P}\big(
j^\kappa(f\circ\phi)\big) 
= 
{\phi'}^m\,{\sf P}\big((j^\kappa f)
\circ\phi\big),
\end{equation}
for every local biholomorphism of $\C$ fixing $0$. Following a trick
of Rousseau (\cite{ rou2006a}), we will apply this formula to the
inverse mapping $\phi := f_1^{ -1}$ of the first coordinate map $f_1 :
\C \to \C$, assuming that $f_1' ( 0 ) \neq 0$, whence $\phi' = \frac{
1}{ f_1'} \circ f_1^{ -1}$. We will explain in a moment that the
assumption $f_1 '( 0) \neq 0$ is harmless for the result.

At first, we trivially have $f_1 \circ f_1^{ -1} = {\rm Id}$, whence
$\big( f_1 \circ f_1^{ -1} \big) ' = 1$ and $\big( f_1 \circ f_1^{ -1}
\big)^{ (\lambda)} = 0$ for all $\lambda \geqslant 2$. Next, by some
direct computations, the derivatives of the reparametrization of $f_2$
happen to be:
\[
\footnotesize
\aligned
\big(f_2\circ f_1^{-1}\big)'
&
=
\bigg(
\frac{f_2'}{f_1'}
\bigg)
\circ f_1^{-1},
\\
\big(f_2\circ f_1^{-1}\big)''
&
=
\bigg[
\frac{f_2''}{(f_1')^2}
-
\frac{f_1''f_2'}{(f_1')^3}
\bigg]
\circ f_1^{-1}
=
\bigg[
\frac{f_1'f_2''
-
f_1''f_2'}{(f_1')^3}
\bigg]
\circ f_1^{-1}
\\
&
=
\frac{\Lambda^3}{(f_1')^3}\circ f_1^{-1},
\endaligned
\]
where we recognize here our favorite Wronskian $\Lambda^3 = \Delta_{ 1,
2}^{ ', \, ''}$. Furthermore, by pursuing as we should the
computations with the help of our beloved total differentiation
operator, we next get:
\[
\footnotesize
\aligned
\big(f_2\circ f_1^{-1}\big)'''
&
=
\bigg(
\frac{{\sf D}\Lambda^3}{(f_1')^4}
-
3\,\frac{\Lambda^3\,f_1''}{(f_1')^5}
\bigg)
\circ f_1^{-1}
=
\frac{\big[\Lambda^3,\,f_1'\big]}{(f_1')^5}
\circ f_1^{-1}
\\
&
=
\frac{\Lambda_1^5}{(f_1')^5}
\circ f_1^{-1},
\\
\big(f_2\circ f_1^{-1}\big)''''
\circ f_1^{-1}
&
=
\frac{\big[\Lambda_1^5,\,f_1'\big]}{(f_1')^7}
\circ f_1^{-1}
\\
&
=
\frac{\Lambda_{1,1}^7}{(f_1')^7}
\circ f_1^{-1},
\endaligned
\]
and so on, with the now clear formal facts that numerators should be
constructed by successively bracketing with $f_1'$, their weight
being visible as just the power of $f_1'$ in the denominator.

With indices, we may therefore define inductively the collection of
{\sl initial invariants} (including $f_1'$ and $\Lambda^3$):
\[
\aligned
\Lambda_{1^{\lambda-2}}^{2\lambda-1}
:=
&\
\big[
\Lambda_{1^{\lambda-3}}^{2\lambda-3},\,f_1'
\big]
\\
=
&\
{\sf D}\Lambda_{1^{\lambda-3}}^{2\lambda-3}\cdot f_1'
-
(2\lambda-3)\Lambda_{1^{1-\lambda_3}}^{2\lambda-3}\cdot f_1'',
\endaligned
\]
for all $\lambda$ with $3 \leqslant \lambda \leqslant \kappa$, where
at the bottom of $\Lambda_{ 1^\ell }^\bullet$, the
notation $1^\ell$ stands for $\ell$ copies of $1$. We then 
get by induction:
\[
\boxed{
\big(f_2\circ f_1^{-1}\big)^{(\lambda)}
=
\frac{\Lambda_{1^{\lambda-2}}^{2\lambda-1}}{
(f_1')^{2\lambda-1}}
\circ f_1^{-1}
}\,.
\]

It would not be a so straightforward task to find a general explicit
expression of these invariants $\Lambda_{ 1^{ \kappa - 2 }}^{2 \kappa
-1}$ in terms of $j^\kappa f$ for arbitrary jet order. For instance,
the invariant $\Lambda_{ 1, 1, 1}^9$, obtained by specializing $i = j
= k = 1$ in the expression given in the theorem stated above (and by
simplifying) reads:
\[
\aligned
\Lambda_{1,1,1}^9
&
=
\big(\Delta^{',\,'''''}
+
5\,\Delta^{',\,''''}\big)\,f_1'f_1'f_1'
-
\big(15\,\Delta^{',\,''''}
+
60\,\Delta^{'',\,'''}\big)\,f_1'f_1'f_1''
-
\\
&
\ \ \ \ \
-
10\,\Delta^{',\,'''}\,f_1'f_1'f_1'''
+
105\,\Delta^{',\,'''}\,f_1'f_1''f_1''
-
105\,\Delta^{',\,''}\,f_1''f_1''f_1''.
\endaligned
\]
Nonetheless, we will in fact not really need to expand the expressions
of these initial invariants.

\smallskip\noindent{\bf Fact.}
{\em The invariants $f_1'$, $\Lambda^3$, $\Lambda_1^5$, $\Lambda_{
1,1}^7$,
\dots, $\Lambda_{1^{\kappa -2}}^{
2 \kappa - 1}$ are mutually algebraically
independent}.

\medskip
This is just because $\Lambda_{ 1^{ \lambda-2}}^{ 2
\lambda - 1}$ is a polynomial in
$\big( j^\lambda f_1, j^\lambda f_2 \big)$ while $\Lambda_{ 
1^{ \lambda-1}}^{ 2 \lambda
+ 1}$ contains the higher jet monomial $f_2^{ (\lambda + 1)}
[ f_1' ]^\lambda$.

\subsection*{ Initial rational expression for invariant
polynomials}
A general polynomial ${\sf P} \big( j^\kappa f \big)$ of weight $m$ in
${\sf E}_{ \kappa, m}^2$ writes in expanded form:
\[
\small
\aligned
{\sf P}\big(j^\kappa f_1,\,j^\kappa f_2\big)
=
\sum_{a_1^1+a_2^1+2a_1^2+2a_2^2
+\cdots+
\kappa a_1^\kappa+\kappa a_2^\kappa=m}\,
{\sf coeff}\cdot
&
\big(f_1'\big)^{a_1^1}\,
\big(f_2'\big)^{a_2^1}\,
\big(f_1''\big)^{a_1^2}\,
\big(f_2''\big)^{a_2^2}\,
\\
&\ \
\cdots\,
\big(f_1^{(\kappa)}\big)^{a_1^\kappa}\,
\big(f_2^{(\kappa)}\big)^{a_2^\kappa},
\endaligned
\]
where by ``${\sf coeff}$'' we mean varying, but notationally
unspecified complex numbers. Reparametrizing by $\phi := f_1^{ -1}$
by an application 
of the definition $(*)$, we should have the relation:
\[
\frac{1}{(f_1'\circ f_1^{-1})^m}\cdot
{\sf P}
\big(
j^\kappa f_1,\,j^\kappa f_2
\big)\circ f_1^{-1}
=
{\sf P}\big(j^\kappa(f_1\circ f_1^{-1}),\,
j^\kappa(f_2\circ f_1^{-1})\big),
\]
in the open subset $\{ f_1' \neq 0 \}$ of the jet space
$J^\kappa ( \C, \, \C^n)$. Thanks to
the preparatory computations above, we 
may replace each monomial in the right hand side, and
this gives us a quite interesting representation:
\[
\small
\aligned
\frac{1}{(f_1'\circ f_1^{-1})^m}\cdot
{\sf P}
\big(
j^\kappa f_1,\,j^\kappa f_2
\big)\circ f_1^{-1}
&
=
\\
=
\bigg[
\sum_{a_1^1+a_2^1+2a_1^2+2a_2^2
+\cdots+
\kappa a_1^\kappa+\kappa a_2^\kappa=m}\,
{\sf coeff}\cdot
&
\big(1\big)^{a_1^1}\,
\Big(\frac{f_2'}{f_1'}\Big)^{a_2^1}\,
\big(0\big)^{a_1^2}\,
\bigg(\frac{\Lambda^3}{(f_1')^3}
\bigg)^{a_2^2}\,
\\
&\ \ \ \ \ \ \ \
\cdots\,
\big(0\big)^{a_1^\kappa}\,
\bigg(
\frac{\Lambda_{1^{\kappa-2}}^{2\kappa-1}}{(f_1')^{2\kappa-1}}
\bigg)^{a_2^\kappa}
\bigg]
\circ f_1^{-1}.
\endaligned
\]
Immediately, we reparametrize this identity by $f_1$, which then
simply erases all the appearing $f^{-1}$, we see that monomials with
positive exponent $a_1^\lambda \geqslant 1$ for some 
$\lambda$ with $2 \leqslant
\lambda \leqslant \kappa$ automatically vanish, and we reduce
monomials to the same denominator:
\[
\aligned
&
{\sf P}\big(j^\kappa f_1,\,j^\kappa f_2\big)
=
\sum_{a_1^1+a_2^1+2a_2^2+\cdots+\kappa a_2^\kappa=m}\,
{\sf coeff}\cdot
\frac{\big(f_2'\big)^{a_2^1}\,
\big(\Lambda^3\big)^{a_2^2}\,
\cdots
\big(\Lambda_{1^{\kappa-2}}^{2\kappa-1}\big)^{a_2^\kappa}}{
(f_1')^{-m+a_2^1+3a_2^2+\cdots+(2\kappa-1)a_2^\kappa}}.
\endaligned
\]
What is the largest power of $f_1'$ as a denominator in the monomials
of the right hand side? Supposing for a while that the quantities 
$a_i^j$ are
nonnegative real numbers, instead of integers,
we may simplify step by step the definition
of this maximum:
\[
\small
\aligned
&
\max_{a_1^1+a_2^1+2a_2^2+\cdots+\kappa a_2^\kappa=m}\,
\Big(-m+a_2^1+3a_2^2+\cdots+(2\kappa-1)a_2^\kappa\Big)
=
\\
&
=
\max\Big\vert_{a_1^1=0}
\Big(\text{\rm substitute}\,\,a_2^1=m-2a_2^2-\cdots-
\kappa a_2^\kappa\,\,\text{\rm in the same quantity}\Big)
\ \ \ \ \ \ \ \ \ \
\ \ \ \ \ \ \ \ \ \
\\
&
=
\max_{a_2^1+2a_2^2+3a_2^3+\cdots+\kappa a_2^\kappa=m}\,
\Big(a_2^2+2a_2^3+\cdots+(\kappa-1)a_2^\kappa\Big)
\ \ \ \ \ \ \ \ \ \
\text{\rm [divide and multiply by 2]}
\\
&
=
\frac{1}{2}\cdot
\max_{2a_2^2+3a_2^3+\cdots+\kappa a_2^\kappa=m}\,
\Big(2a_2^2+4a_2^3+\cdots+(2\kappa-2)a_2^\kappa\Big)
\ \ \ \ \ \ \ \ \ \
\text{\rm [substitute}\,\,2a_2^2\text{\rm ]}
\\
&
=
\frac{m}{2}
+
\frac{1}{2}\cdot
\max_{3a_2^3+4a_2^4+\cdots+\kappa a_2^\kappa=m}
\Big(a_2^3+2a_2^4+\cdots+(\kappa-2)a_2^\kappa\Big)
\ \ \ \ \ 
\text{\rm [divide and multiply by 3]}
\\
&
=
\frac{m}{2}
+
\frac{1}{2\cdot 3}\cdot\,
\max_{3a_2^3+4a_2^4+\cdots+\kappa a_2^\kappa=m}
\Big(3a_2^3+6a_2^4+\cdots+3(\kappa-2)a_2^\kappa\Big)
\ \ \ \ \ \ \ \ \ \
\text{\rm [substitute}\,\,3a_2^3\text{\rm ]}
\\
&
=
\frac{m}{2}
+
\frac{m}{6}
+
\frac{1}{3\cdot 4}\cdot\,
\max_{4a_2^4+5a_2^5+\cdots+\kappa a_2^\kappa=m}
\Big(4a_2^4+8a_2^5+\cdots+4(\kappa-3)a_2^\kappa\Big)
\\
&
=
\frac{2}{3}\,m
+
\frac{1}{12}\,m
+
\frac{1}{4\cdot 5}\cdot\,
\max_{5a_2^5+6a_2^6+\cdots+\kappa a_2^\kappa=m}
\Big(5a_2^6+10a_2^6+\cdots+5(\kappa-4)a_2^\kappa\Big)
\\
&
=
\frac{3}{4}\,m
+
\frac{1}{20}\,m
+
\frac{1}{5\cdot 6}\cdot\,
\max_{6a_2^6+7a_2^7+\cdots+\kappa a_2^\kappa=m}
\Big(6a_2^6+12a_2^7+\cdots+6(\kappa-5)a_2^\kappa\Big)
\\
&
=
\frac{4}{5}\,m
+\cdots
\ \ \ \ \ \ \ \ \ \ \ \ \ \ \ 
\text{\rm [observe the induction]}
\\
&
=
\frac{\kappa-1}{\kappa}\,m.
\endaligned
\]
Thus, when the $a_2^i$ are restricted to be integers, we in any case
deduce that the maximally negative power of $f_1'$ is $\geqslant -
\frac{ ( \kappa - 1)}{ \kappa}\, m$. Reorganizing the result, we then
obtain a representation of ${\sf P} \big( j^\kappa f\big)$, valid by
construction in the subset $\{ f_1' \neq 0 \}$ of the jet space
$J^\kappa \big( C, \, \C^n)$, as a sum of powers of $f_1'$ :
\[
\aligned
&
{\sf P}\big(j^\kappa f\big)
=
\sum_{-\frac{\kappa-1}{\kappa}\,m\leqslant a\leqslant m}\,\,\,\,
(f_1')^a\cdot
{\sf P}_a\Big(
f_2',\,\Lambda^3,\,
\Lambda_1^5,\,\Lambda_{1,1}^7,\,\dots,\,
\Lambda_{1^{\kappa-2}}^{2\kappa-1}\Big),
\endaligned
\]
multiplied by certain polynomials ${\sf P}_a$ which depend upon ${\sf
P}$ and are {\em not} arbitrary. In fact, by reduction to the same
denominator, we may write:
\[
{\sf P}\big(j^\kappa f\big)
=
\frac{{\sf Q}\big(f_1',\,f_2',\,\Lambda^3,\,\dots,\,
\Lambda_{1^{\kappa-2}}^{2\kappa-1}\big)}{(f_1')^{-a_0}},
\]
where $a_0$ is the smallest exponent $a$ of $f_1'$ above. Chasing the
denominator in the case where $a_0$ is negative (this would be
unnecessary if $a_0 \geqslant 0$), we get an identity $(f_1')^{ a_0 }
\cdot {\sf P} \equiv {\sf Q}$ between two polynomials valid in $\{
f_1' \neq 0\}$, hence everywhere by the principle of analytic
continuation. Thus, the restriction $f_1' \neq 0$ is removed.

\subsection*{ Weighted homogeneities}
Let $\mu \in \Z$ be an integer, possibly negative. A rational
expression ${\sf R} \big( j^\kappa f \big) \in \text{\rm Frac} \big(
\C [ j^\kappa f ] \big)$ will be said to be of {\sl weighted
homogeneous degree} $\mu$ when for every complex weighted
$\delta$-dilation which acts in accordance with the number of primes:
\[
\delta\cdot j^\kappa f
:=
\big(
\delta\,f_1',\,\delta f_2',\,\delta^2\,f_1'',\,\delta^2\,f_2'',\cdots,
\delta^\kappa\,f_1^{(\kappa)},\,\delta^\kappa\,f_2^{(\kappa)}
\big),
\]
the dilation factor escapes the parentheses to exactly the $\mu$-th
power:
\[
{\sf R}\big(\delta\cdot j^\kappa f\big)
=
\delta^\mu\cdot
{\sf R}\big(j^\kappa f\big).
\]
When ${\sf R}$ is a polynomial, $\mu$ is then the total, constant
number of primes of each monomial.
 
By choosing the reparametrization $\phi$ to just be a
$\delta$-dilation in the source, with nonzero $\delta \in \C$, we
immediataly see that our original jet polynomial ${\sf P} \in {\sf
E}_{ \kappa, m}^2$\,\,---\,\,hence 
also its rational expression obtained above\,\,---\,\,must 
in particular be weighted homogeneous of degree $m$:
\[
{\sf P}\big(\delta\cdot j^\kappa f\big)
=
\delta^m\cdot{\sf P}\big(j^\kappa f\big).
\]
In addition and in particular, using the definition $\Lambda_{ 1^{
\lambda - 2}}^{ 2 \lambda - 1} = \big[ \Lambda_{ 1^{ \lambda - 3}}^{ 2
\lambda - 3}, \, f_1' \big]$, one easily verifies by induction that
the invariant $\Lambda_{1^{ \lambda - 2}}^{ 2 \lambda - 1}$ is
homogeneous of degree equal to its weight $2 \lambda - 1$, 
an integer which we
had already specified as the upper index:
\[
\Lambda_{1^{\lambda-2}}^{2\lambda-1}
\big(\delta\cdot j^\kappa f\big)
=
\delta^{2\lambda-1}\cdot\Lambda_{1^{\lambda-2}}^{2\lambda-1}
\big(j^\lambda f\big).
\]

In an analogous fashion, introducing some new extra independent
variables $F_1, \, F_2, \, A^3, \, \dots, \, A^{ 2 \kappa - 1}$
corresponding to $f_1', f_2', \Lambda^3, \, \dots, \, \Lambda_{ 1^{
\kappa - 2}}^{ 2 \kappa - 1}$, a rational expression ${\sf T} \in
\text{\rm Frac} \big( \C [ F_1, F_2, A^3, \dots, A^{ 2 \kappa - 1} ]
\big)$ will be said to be of {\sl weighted homogeneous degree}\, $\mu$
when it enjoys:
\[
{\sf T}
\big(\delta F_1,\,\delta F_2,\,\delta^3A^3,\,\delta^5A^5,\,\dots,\,
\delta^{2\kappa-1}A^{2\kappa-1}\big)
=
\delta^\mu\cdot
{\sf T}
\big(F_1,\,F_2,\,A^3,\,A^5,\,\dots,\,A^{2\kappa-1}\big),
\] 
for every $\delta \in \C$.

\smallskip\noindent{\bf Lemma.}
\label{initial-2-k}
{\em 
In dimension $n = 2$ for jets of order $\kappa \geqslant 2$, every jet
polynomial ${\sf P} = {\sf P} \big( j^\kappa f)$ invariant by
reparametrization writes under the form:
\[
\boxed{
{\sf P}\big(j^\kappa f\big)
=
\sum_{-\frac{\kappa-1}{\kappa}m\leqslant a\leqslant m}\,
(f_1')^a\,{\sf P}_a\Big(
f_2',\,\Lambda^3,\,\Lambda_1^5,\,\Lambda_{1,1}^7,\,
\cdots,\,
\Lambda_{1^{\kappa-2}}^{2\kappa-1}
\Big)}\,,
\] 
where the integer $a$ takes possibly negative values in the interval
$\big[ - \frac{ \kappa-1 }{ \kappa } \, m, \, m \big]$, for certain
weighted homogeneous polynomials:
\[
{\sf P}_a
=
\sum_{b_2+3c_3+\cdots+(2\kappa-1)
c_{2\kappa-1}=m-a}\,\,
\text{\sf coeff}\cdot\,
\big(F_2\big)^{b_2}\,
\big(A^3\big)^{c_3}\,\cdots\,
\big(A^{2\kappa-1}\big)^{c_{2\kappa-1}}
\]
of weighted degree $m-a$.

Conversely, for every collection of such weighted homogeneous
polynomials ${\sf P}_a$ in $\C \big[ F_2, \, A^3, \, \dots, \,
A^{2\kappa -1 } \big]$ of weighted degree $m - a$ indexed by an
integer $a$ running in $\big[- \frac{ \kappa-1 }{ \kappa }\, m, \, m
\big]$ such that the reduction to the same denominator and the
simplification of the finite sum:
\[
{\sf R}\big(j^\kappa f\big)
:=
\sum_{-\frac{\kappa-1}{\kappa}m\leqslant a\leqslant m}\,
(f_1')^a\,{\sf P}_a\Big(
f_2',\,\Lambda^3,\,\Lambda_1^5,\,\Lambda_{1,1}^7,\,
\cdots,\,
\Lambda_{1^{\kappa-2}}^{2\kappa-1}
\Big)
\]
yields a {\em true} jet \underline{\em polynomial} in $\C \big[
j^\kappa f \big]$, then ${\sf R} ( j^\kappa f)$ is a polynomial
invariant by reparametrization belonging to ${\sf E}_{\kappa, m}^2$. }

\proof
We saw that ${\sf P}$ is homogeneous of degree $m$, namely:
\[
\aligned
\sum_{-\frac{\kappa-1}{\kappa}m\leqslant a\leqslant m}\,
&
\delta^a\big(f_1'\big)^a\,
{\sf P}_a\Big(\delta f_2',\,\delta^3\Lambda^3,\,\dots,\,
\delta^{2\kappa-1}\Lambda_{1^{\kappa-2}}^{2\kappa-1}\Big)
=
\\
&
=
\delta^m\cdot\sum_{-\frac{\kappa-1}{\kappa}m\leqslant a\leqslant m}\,
\big(f_1'\big)^a\,{\sf P}_a
\Big(f_2',\,\Lambda^3,\,\dots,\,\Lambda_{1^{\kappa-2}}^{2\kappa-1}\Big).
\endaligned
\]
By algebraic independency of $f_1'$ with respect to $\C \big[ f_2',
\Lambda^3, \dots, \Lambda_{1^{\kappa-2}}^{ 
2 \kappa - 1} \big]$, we then may identify
powers of $f_1'$, getting for each $a$:
\[
{\sf P}_a\Big(\delta f_2',\,\delta^3\Lambda^3, \dots, 
\delta^{2\kappa-1}\Lambda_{1^{\kappa-2}}^{2\kappa-1}\Big)
=
\delta^{m-a}\,{\sf P}_a
\Big(
f_2',\Lambda^3,\dots,\Lambda_{1^{\kappa-2}}^{2\kappa-1}
\Big).
\]
Further, the algebraic independency of $f_2', \Lambda^3, \dots,
\Lambda_{ 1^{\kappa-2}}^{ 2 \kappa - 1}$ then entails that the
homogeneities:
\[
{\sf P}_a\Big(\delta F_2,\delta^3A^3,\dots,
\delta^{2\kappa-1}A^{2\kappa-1}\Big)
=
\delta^{m-a}\cdot{\sf P}_a\Big(F_2,A^3,\dots,A^{2\kappa-1}\Big)
\]
hold in the polynomial algebra $\C \big[ F_2, A^3, \dots, A^{ 2 \kappa
- 1} \big]$. This gives the claimed representation of any ${\sf P}
\in {\sf E}_{ \kappa, m}^2$.

Conversely, assuming that the ${\sf P}_a$ are homogeneous in this way,
then for any reparametrization $\phi$, setting $g_i = f_i \circ \phi$
for $i = 1, 2$ and recalling:
\[
\Lambda_{1^{\lambda-2}}^{2\lambda-1}\big(j^\lambda g\big)
=
(\phi')^{2\lambda-1}\,
\Lambda_{1^{\lambda-2}}^{2\lambda-1}\big(j^\lambda f),
\]
we immediately deduce that:
\[
\aligned
{\sf P}_a\Big(g_2',\,
&
\Lambda^3\big(j^3g\big),\,\dots,\,
\Lambda_{1^{\kappa-2}}^{2\kappa-1}\big(j^\kappa g\big)\Big)
=
\\
&
=
{\sf P}_a\Big(\phi'\cdot f_2'\circ\phi,\,
(\phi')^3\cdot\Lambda^3\circ\phi,\,\dots,\,
(\phi')^{2\kappa-1}\cdot\Lambda_{1^{\kappa-2}}^{2\kappa-1}\circ\phi
\Big)
\\
&
=
(\phi')^{m-a}\cdot{\sf P}_a\Big(f_2',\,\Lambda^3,\,\dots,\,
\Lambda_{1^{\kappa-2}}^{2\kappa-1}
\Big),
\endaligned
\]
whence multiplication by $(g_1')^a = (\phi')^a \, (f_1')^a$ and
summation gives:
\[
\aligned
\sum_{-\frac{\kappa-1}{\kappa}m\leqslant a\leqslant m}\,
{g_1'}^a\,
&
{\sf P}_a\Big(g_2',\,\Lambda^3\big(j^3g\big),\,\dots,\,
\Lambda_{1^{\kappa-2}}^{2\kappa-1}\big(j^\kappa g\big)\Big)
=
\\
&
(\phi')^m\cdot
\sum_{-\frac{\kappa-1}{\kappa}m\leqslant a\leqslant m}\,
(f_1')^a\,{\sf P}_a\Big(
f_2',\,\Lambda^3,\,\dots,\,\Lambda_{1^{\kappa-2}}^{2\kappa-1}
\Big)\circ\phi,
\endaligned
\]
which exactly means, as soon as such a rational sum represents a true
polynomial, that it belongs to ${\sf E}_{ \kappa, m}^2$, {\em 
quod erat demonstrandum}.
\endproof

\subsection*{ Arbitrary dimension}
To generalize the preceding proposition, suppose now that $n \geqslant
2$ is arbitrary. The same trick of reparametrizing each $f_i$ by $\phi
= f_1^{ -1}$ gives birth to a collection of {\sl initial invariants}
appearing as numerators of:
\[
\big(f_i\circ f_1^{-1}\big)^{(\lambda)}
=
\frac{\Lambda_{1,i;\,1^{\lambda-2}}^{2\lambda-1}\big(j^\lambda f\big)}{
(f_1')^{2\lambda-1}},
\]
for $i = 2, 3, \dots, n$, where the $\Lambda$-invariants 
depending on $i$ and on $\lambda$ are
defined inductively by successively bracketing with $f_1'$:
\[
\aligned
&
\Lambda_{1,i}^3
:=
\big[f_i',\,f_1'\big],
\ \ \ \ \ \ \ \ \ \ \ \ \ \
\ \ \ \ \ \ \ \ \ \ \ \ \ \
\Lambda_{1,i;\,1}^5
:=
\big[\Lambda_{1,i}^3,\,f_1'\big],
\\
&
\text{\rm and generally:} \ \ \ \ \ \ \
\Lambda_{1,i;\,1^{\lambda-2}}^{2\lambda-1}
:=
\big[\Lambda_{1,i;\,1^{\lambda-3}}^{2\lambda-3},\,f_1'\big],
\ \ \ \ \ \ \
\text{\rm for}\ \ 3 \leqslant \lambda \leqslant \kappa.
\endaligned
\]
Our considerations about brackets show that these polynomials
are effectively invariant by reparametrization. Furthermore:

\smallskip\noindent{\bf Fact.}
{\em 
The $n + (n-1) (\kappa-1)$ invariants:
\[
\label{fact-n-k}
\begin{array}{cccccc}
f_1', & f_2', & f_3', & f_4', & \dots\,, & f_n',
\\
& \Lambda_{1,2}^3, & \Lambda_{1,3}^3, & \Lambda_{1,4}^3, &
\dots\,, & \Lambda_{1,n}^3, 
\\
& \Lambda_{1,2;\,1}^5, & \Lambda_{1,3;\,1}^5, & \Lambda_{1,4;\,1}^5, &
\dots\,, & \Lambda_{1,n;\,1}^5,
\\
& \cdots & \cdots & \cdots & \cdots & \cdots 
\\
&
\Lambda_{1,2;\,1^{\kappa-2}}^{2\kappa-1}, & 
\Lambda_{1,3;\,1^{\kappa-2}}^{2\kappa-1}, & 
\Lambda_{1,4;\,1^{\kappa-2}}^{2\kappa-1}, &
\dots,\, & 
\Lambda_{1,n;\,1^{\kappa-2}}^{2\kappa-1}
\end{array}
\]
are mutually algebraically independent.
}\medskip

Indeed, $\Lambda_{ 1, i; \, 1^{ \lambda - 2}}^{ 2 \lambda - 1}$
contains the monomial $f_i^{ ( \lambda)} { f_1'}^{ \lambda - 1}$,
while the invariants $\Lambda_{ 1, j; \, 1^{ \lambda - 3}}^{ 2 \lambda
- 3}$ only depend upon $j^{ \lambda - 1}f$.

\smallskip

Reasonings similar to the ones developed above yield the following
lemma, valuable for any $n \geqslant 1$ and any $\kappa
\geqslant 1$.

\smallskip\noindent{\bf Lemma.}
\label{initial-n-k}
{\em 
In dimension $n \geqslant 1$ and for jets of order $\kappa \geqslant
1$, every polynomial ${\sf P} = {\sf P} \big( j^\kappa f\big)$
invariant by reparametrization writes under the form:
\[
\small
\aligned
\boxed{
{\sf P}\big(j^\kappa f\big)
=
\sum_{-\frac{\kappa-1}{\kappa}m\leqslant a\leqslant m}\,
(f_1')^a\,{\sf P}_a
\left(
\begin{array}{ccccc} 
f_2', & f_3', & f_4', & \dots\,, & f_n',
\\
\Lambda_{1,2}^3, & \Lambda_{1,3}^3, & \Lambda_{1,4}^3, &
\dots\,, & \Lambda_{1,n}^3, 
\\
\cdots & \cdots & \cdots & \cdots & \cdots
\\
\Lambda_{1,2;\,1^{\kappa-2}}^{2\kappa-1}, & 
\Lambda_{1,3;\,1^{\kappa-2}}^{2\kappa-1}, & 
\Lambda_{1,4;\,1^{\kappa-2}}^{2\kappa-1}, &
\dots,\, & 
\Lambda_{1,n;\,1^{\kappa-2}}^{2\kappa-1}
\end{array}
\right)
}\,,
\endaligned
\]
where the integer $a$ takes all possibly negative values in the
interval $\big[ - \frac{ \kappa - 1}{ \kappa} \, m, m \big]$, for
certain weighted homogeneous polynomials:
\[
{\sf P}_a
=
\sum_{b_2+\cdots+b_n
+
3c_2+\cdots+3c_n+
\atop
+\cdots+
(2\kappa-1)q_2+\cdots+(2\kappa-1)q_n=m-a}\!\!\!\!
{\sf coeff}\,\cdot\,
\prod_{i=2}^n\,\big(F_i\big)^{b_i}\,
\prod_{i=2}^n\,\big(A_i^3\big)^{c_i}\,
\cdots\,\,
\prod_{i=2}^n\,\big(A_i^{2\kappa-1}\big)^{q_i}
\]
of weighted degree $m - a$, namely satisfying: 
\[
{\sf P}_a
\Big(
\delta\,F_i,\,\delta^3\,A_i^3,\,\dots,\,
\delta^{2\kappa-1}\,A_i^{2\kappa-1}
\Big)
=
\delta^{m-a}\cdot
{\sf P}_a
\Big(
F_i,\,A_i^3,\,\dots,\,A_i^{2\kappa-1}
\Big).
\]

Conversely, for every collection of such weighted homogeneous
polynomials ${\sf P}_a$ in $\C \big[ F_i, A_i^3, \dots, A_i^{ 2\kappa
- 1} \big]$ of weighted degree $m - a$ indexed by an integer $a$
running in $\big[ - \frac{ \kappa - 1}{ \kappa } \, m, \, m \big]$
such that the reduction to the same denominator and the simplification
of the finite sum:
\[
{\sf R}\big(j^\kappa f\big)
=
\sum_{-\frac{\kappa-1}{\kappa}m\leqslant a\leqslant m}\,
(f_1')^a\,{\sf P}_a
\left(
\begin{array}{ccccc} 
f_2', & f_3', & f_4', & \dots\,, & f_n',
\\
\Lambda_{1,2}^3, & \Lambda_{1,3}^3, & \Lambda_{1,4}^3, &
\dots\,, & \Lambda_{1,n}^3, 
\\
\cdots & \cdots & \cdots & \cdots & \cdots
\\
\Lambda_{1,2;\,1^{\kappa-2}}^{2\kappa-1}, & 
\Lambda_{1,3;\,1^{\kappa-2}}^{2\kappa-1}, & 
\Lambda_{1,4;\,1^{\kappa-2}}^{2\kappa-1}, &
\dots,\, & 
\Lambda_{1,n;\,1^{\kappa-2}}^{2\kappa-1}
\end{array}
\right)
\]
yields a {\em true} jet \underline{\em polynomial} in $\C \big[
j^\kappa f \big]$, then ${\sf R} \big( j^\kappa f \big)$ is
a polynomial invariant by reparametrization belonging to
${\sf E}_{ \kappa, m}^n$.
}\medskip

\section*{\S6.~Description of the algorithm 
\\
in dimension $n=2$ for jet level $\kappa = 4$}
\label{Section-6}

\subsection*{Necessity of negative powers of $f_1'$}
Our aim now is to prove\footnote{\, An alternative proof was provided
in~\cite{ mer2008a}. } the theorem which describes the algebraic
structure of ${\sf E}_4^2$, {\em see} p.~\pageref{theorem-2-4}. We
will thus illustrate in a concrete case the general algorithm which
will be presented in Section~9 below. We hope this will make
the general considerations intuitively clearer.

\proof
Compared to the initial rational representation:
\[
{\sf P}\big(j^4f\big)
=
\sum_{-\frac{3}{4}m\leqslant a\leqslant m}\,
(f_1')^a\,{\sf P}_a
\big(f_2',\,\Lambda^3,\,\Lambda_1^5,\,\Lambda_{1,1}^7\big)
\]
of an arbitrary polynomial ${\sf P} \in {\sf E}_4^2$ 
that was furnished by the lemma on
p.~\pageref{initial-2-k}, the theorem on p.~\pageref{theorem-2-4}
states that 4 further invariants, namely $\Lambda_2^5$, $\Lambda_{ 1,
2}^7$, $\Lambda_{ 2, 2}^7$ and $M^8$, are necessary to generate the
full algebra ${\sf E}_4^2$. In fact, by looking at the 9 syzygies
listed in the theorem in question, one may easily obtain the
expression of these 4 further invariants in $\C \big[ f_2', \,
\Lambda^3, \, \Lambda_1^5, \, \Lambda_{ 1, 1}^7 \big] \big[ \frac{ 1}{
f_1'} \big]$, namely:
\[
\aligned
\Lambda_2^5
&
=
\frac{f_2'\Lambda_1^5-3\,\Lambda^3\Lambda^3}{f_1'},
\\
\Lambda_{1,2}^7
&
=
\frac{f_2'\Lambda_{1,1}^7-5\,\Lambda^3\Lambda_1^5}{f_1'},
\\
\Lambda_{2,2}^7
&
=
\frac{f_2'f_2'\Lambda_{1,1}^7-10\,f_2'\Lambda^3\Lambda_1^5
+15\,\Lambda^3\Lambda^3\Lambda^3}{f_1'f_1'f_1'},
\\
M^8
&
=
\frac{3\,\Lambda^3\Lambda_{1,1}^7-5\,\Lambda_1^5\Lambda_1^5}{f_1'f_1'}.
\endaligned
\]

\smallskip\noindent{\bf Crucial observation.}
{\em So when one substitutes, in an arbitrary polynomial:
\[
{\sf P}\big(
f_1',\,f_2',\,\Lambda^3,\,\Lambda_1^5,\,\Lambda_{1,1}^7,\,
\Lambda_{1,2}^7,\,\Lambda_{2,2}^7,\,M^8\big)
\]
these rational representations of $\Lambda_2^5$, $\Lambda_{ 1, 2}^7$,
$\Lambda_{ 2, 2}^7$, $M^8$, one indeed obtains a rational expression
as the one above which {\em necessarily and unavoidably} incorporates
negative powers of $f_1'$}.
\medskip

Well, how then should we interpret our initial rational expression?
Why are the 4 further invariants $\Lambda_2^5$, $\Lambda_{ 1, 2}^7$,
$\Lambda_{ 2, 2}^7$ and $M^8$ invisible in it?

First of all, as a preliminary, we must at least show that the 9
fundamental
invariants $f_1'$, $f_2'$, $\Lambda^3$, $\Lambda_1^5$, $\Lambda_2^5$,
$\Lambda_{1,1}^7$, $\Lambda_{ 1, 2}^7$, $\Lambda_{ 2, 2}^7$ and $M^8$
are mutually independent.

On this purpose, we set $f_1' = 0$ in these 9 fundamental invariants,
and this then leaves us with the 8 invariants:
\[
f_2',\ \ \ 
\Lambda^3\big\vert_0,\ \ \ \ \
\Lambda_1^5\big\vert_0,\ \ \ \ \
\Lambda_2^5\big\vert_0,\ \ \ \ \
\Lambda_{1,1}^7\big\vert_0,\ \ \ \ \
\Lambda_{1,2}^7\big\vert_0,\ \ \ \ \
\Lambda_{2,2}^7\big\vert_0,\ \ \ \ \
M^8\big\vert_0,
\]
which we shall shortly call {\sl restricted invariants}, our notation
being self-evident. The following assertion is simply checked by
inspecting the explicit expressions.

\smallskip\noindent{\bf Fact.}
{\em 
The four restricted invariants:
\[
\aligned
f_2',\ \ \ \ \ \ \ \ \ \ \ \ \ \ \
\Lambda^3\big\vert_0
&
=
-f_1''f_2',\ \ \ \ \ \ \ \ \ \ \ \ \ \ \ 
\Lambda_2^5\big\vert_0
=
3\,f_1''f_2'f_1''\ \ \ \ \ \ \ \ \
\text{\rm and}
\\
\Lambda_{2,2}^7\big\vert_0
&
=
\big(-f_1''''f_2'+\Delta^{'',\,'''}\big)\,f_2'f_2'
+
10\,f_1'''f_2'f_2'f_2''
-
15\,f_1''f_2'f_2''f_2''
\endaligned
\]
are mutually algebraically independent.
}\medskip

It follows that $f_1'$, $f_2'$, $\Lambda^3$, $\Lambda_2^5$ and
$\Lambda_{ 2, 2}^7$ are algebraically independent. Switching lower
indices, $f_1'$, $f_2'$, $\Lambda^3$, $\Lambda_1^5$ and $\Lambda_{ 1,
1}^7$ are also algebraically independent.

\smallskip

Next, by looking at the 9 syzygies listed in the theorem,
we may express each one of the 8 restricted invariants
by means of the above four algebraically independent
restricted invariants, provided that one allows a
division by $f_2'$:
\[
\aligned
&
f_2',
\\
&
\Lambda^3\big\vert_0,
\\
&
\Lambda_1^5\big\vert_0
=
3\,
\frac{\Lambda^3\vert_0\,\Lambda^3\vert_0}{f_2'},
\\
&
\Lambda_2^5\big\vert_0,
\endaligned
\]
\[
\aligned
&
\Lambda_{1,1}^7\big\vert_0
=
15\,\frac{\Lambda^3\vert_0\,\Lambda^3\vert_0\,\Lambda^3\,\vert_0}{
f_2'f_2'},
\\
&
\Lambda_{1,2}^7\big\vert_0
=
5\,\frac{\Lambda^3\vert_0\,\Lambda_2^5\vert_0}{f_2'},
\\
&
\Lambda_{2,2}^7\big\vert_0,
\endaligned
\]
\[
\aligned
&
M^8\big\vert_0
=
\frac{3\,\Lambda^3\vert_0\,\Lambda_{2,2}^7-
5\,\Lambda_2^5\vert_0\,\Lambda_2^5\vert_0}{f_2'f_2'}.
\endaligned
\]
In fact, all divisions by $f_2'$ or
by $f_2' f_2'$ do cancel out after simplification.

\smallskip\noindent{\bf Lemma.}
\label{lemma-independence}
{\em 
The 9 fundamental invariants $f_1'$, $f_2'$, $\Lambda^3$,
$\Lambda_1^5$, $\Lambda_2^5$, $\Lambda_{1,1}^7$, $\Lambda_{ 1, 2}^7$,
$\Lambda_{ 2, 2}^7$ and $M^8$ are mutually independent. More
precisely, $f_1'$, $f_2'$, $\Lambda^3$, $\Lambda_1^5$ and $\Lambda_{
1, 1}^7$ are algebraically independent and there exist no polynomial
representation of either one of the following four forms:
\[
\aligned
\Lambda_2^5
&
=
{\sf polynomial}
\big(f_1',\,f_2',\,\Lambda^3,\,\Lambda_1^5,\,\Lambda_{1,1}^7\big),
\\
\Lambda_{1,2}^7
&
=
{\sf polynomial}
\big(f_1',\,f_2',\,\Lambda^3,\,\Lambda_1^5,\,
\Lambda_2^5,\,\Lambda_{1,1}^7\big),
\\
\Lambda_{2,2}^7
&
=
{\sf polynomial}
\big(f_1',\,f_2',\,\Lambda^3,\,\Lambda_1^5,\,
\Lambda_2^5,\,\Lambda_{1,1}^7,\,\Lambda_{1,2}^7\big),
\\
M^8
&
=
{\sf polynomial}
\big(f_1',\,f_2',\,\Lambda^3,\,\Lambda_1^5,\,
\Lambda_2^5,\,\Lambda_{1,1}^7,\,\Lambda_{1,2}^7,\,\Lambda_{2,2}^7\big).
\endaligned
\]
}\medskip

\proof
By setting $f_1' = 0$ in a polynomial representation such as the
first one and by replacing the values of some of the restricted
invariants, we get:
\[
\aligned
\Lambda_2^5\big\vert_0
&
=
\sum\,{\sf coeff}\cdot
\big(f_2'\big)^b\,\big(\Lambda^3\big)^c\,\big(\Lambda_1^5\big)^d\,
\big(\Lambda_{1,1}^7\big)^e\Big\vert_0
\\
&
=
\sum\,{\sf coeff}\cdot
\big(f_2'\big)^b\,\big(\Lambda^3\big)^c\,
\bigg(3\,\frac{\Lambda^3\,\Lambda^3}{f_2'}\bigg)^d\,
\bigg(15\,\frac{\Lambda^3\,\Lambda^3\,\Lambda^3}{f_2'f_2'}
\bigg)^e\Big\vert_0
\\
&
=
\sum\,{\sf coeff}\cdot
\big(f_2'\big)^{b-d-2e}\,\big(\Lambda^3\big\vert_0\big)^{c+2d+3e},
\endaligned
\]
where the exponents
$b, c, d, e \geqslant 0$ are nonnegative integers.
But this is impossible, because $\Lambda_2^5 \big\vert_0$ is
transcendental over $\C \big[ f_2', \, \Lambda^3 \big\vert_0 \big]$.

Next, for the second hypothetical representation, the same kind of
replacement yields:
\[
5\,
\frac{\Lambda^3\,\Lambda_2^5}{f_2'}
\Big\vert_0
=
\sum\,
{\sf coeff}\cdot
\big(f_2'\big)^{b-d-2f}\,
\big(\Lambda^3\big)^{c+2d+3f}\,
\big(\Lambda_2^5\big)^e
\Big\vert_0.
\]
So, by identifying the powers of the restricted algebraically
independent invariants $f_2'$, $\Lambda^3 \big\vert_0$, $\Lambda_2^5
\big\vert_0$, we get 
three equations between integers:
\[
-1=b-d-2f,\ \ \ \ \ \ \ \ \ \ \ \
1=c+2d+3f,\ \ \ \ \ \ \ \ \ \ \ \
1=e,
\]
which are seen to be impossible, since $b, c, d, e, f \geqslant 0$,
the second one yielding $d = f = 0$, while the first one then reads $-
1 = b$.

Similarly as for the first one, the third hypothetical representation
is {\em a priori}\, excluded, because the right hand side does not
depend upon $\Lambda_{ 2, 2}^7 \big\vert_0$ at all.

Finally, the fourth hypothetical representation amounts to:
\[
\frac{3\,\Lambda^3\,\Lambda_{2,2}^7
-
5\,\Lambda_2^5\,\Lambda_2^5}{f_2'f_2'}
\Big\vert_0
=
\sum\,{\sf coeff}\cdot
\big(f_2')^{b-d-2f-g}\,
\big(\Lambda^3\big)^{c+2d+3f+g}\,
\big(\Lambda_2^5\big)^{e+g}\,
\big(\Lambda_{2,2}^7\big)^h
\Big\vert_0,
\]
hence looking at the representation of the first term $\frac{ 3\,
\Lambda^3 \, \Lambda_{ 2, 2}^7}{ f_2' f_2'}$ of the left hand side,
and identifying powers, we get three equations:
\[
-2=b-d-2f-g,\ \ \ \ \ \ \ \ \ \ \ \
1=c+2d+3f+g,\ \ \ \ \ \ \ \ \ \ \ \
1=h.
\]
The second one implies $d = f = 0$ and $g = 0$ or $1$, whence
the first one then becomes impossible.
\endproof

\subsection*{ First loop of the algorithm}
The initial expression:
\[
{\sf P}\big(j^4f\big)
=
\sum_{-\frac{3}{4}m\leqslant a\leqslant m}\,
(f_1')^a\,{\sf P}_a
\big(f_2',\,\Lambda^3,\,\Lambda_1^5,\,\Lambda_{1,1}^7\big)
\]
shows five invariants $f_1'$, $f_2'$, $\Lambda^3$, $\Lambda_1^5$,
$\Lambda_{ 1, 1}^7$, and four restricted invariants $f_2'$, $\Lambda^3
\big\vert_0$, $\Lambda_1^5 \big\vert_0$, $\Lambda_{ 1, 1}^7
\big\vert_0$. To determine the structure of ${\sf E}_4^2$, 
here is the first loop of our algorithm.

\smallskip$\bullet$\
Compute the ideal of relations\footnote{\, We will discuss in a
while two ways of computing ideal of relations. 
The data reproduced here are obtained by means of
Gröbner bases computations. } between the 4 known
restricted invariants:
\[
{\sf Ideal-Rel}\Big(
f_2'\big\vert_0,\,
\Lambda^3\big\vert_0,\,
\Lambda_1^5\big\vert_0,\,
\Lambda_{1,1}^7\big\vert_0
\Big),
\]
namely a generating set of the ideal of all polynomials
$\mathcal{ Q} \big( F_2, \, A^3, \, A^5, \, A^7 \big)$ in four
variables that give zero, identically, after substituting these four
restricted invariants.

\smallskip$\bullet$\
Get as generators of this ideal of relations the three relations, 
valuable for $f_1' = 0$:
\[
\aligned
0
&
\equiv
3\,\Lambda^3\Lambda^3
-
f_2'\Lambda_1^5\Big\vert_0,
\\
0
&
\equiv
5\,\Lambda^3\Lambda_1^5
-
f_2'\Lambda_{1,1}^7\Big\vert_0,
\\
0
&
\equiv
5\,\Lambda_1^5\Lambda_1^5
-
3\,\Lambda^3\Lambda_{1,1}^7\Big\vert_0.
\endaligned
\]

\smallskip$\bullet$\
Consequently, without setting $f_1' = 0$, there should exist
remainders that are a multiple of $f_1'$:
\[
\aligned
\\
0
&
\equiv
3\,\Lambda^3\Lambda^3
-
f_2'\Lambda_1^5
+
\green{f_1'}
\times\blue{\sf something},
\\
0
&
\equiv
5\,\Lambda^3\Lambda_1^5
-
f_2'\Lambda_{1,1}^7
+
\green{f_1'}
\times\blue{\sf something},
\\
0
&
\equiv
5\,\Lambda_1^5\Lambda_1^5
-
3\,\Lambda^3\Lambda_{1,1}^7
+
\green{f_1'}
\times\blue{\sf something}.
\endaligned
\]

\smallskip$\bullet$\
Each ``\blue{something}'' necessarily also is an invariant belonging
to ${\sf E}_4^2$, because it is a polynomial and we can write it as
$\frac{ 1}{ f_1'}$ times a corresponding quadratic expression in the
already known invariants $f_2'$, $\Lambda^3$, $\Lambda_1^5$,
$\Lambda_{ 1, 1}^7$.

\smallskip$\bullet$\
Find the maximal power by which $f_1'$ factors each remaining
``\blue{something}''.

\smallskip$\bullet$\
Get the three identically satisfied relations:
\[
\aligned
0
&
\equiv
3\,\Lambda^3\Lambda^3
-
f_2'\Lambda_1^5
\blue{+f_1'\Lambda_2^5},
\\
0
&
\equiv
5\,\Lambda^3\Lambda_1^5
-
f_2'\Lambda_{1,1}^7
\blue{+f_1'\Lambda_{1,2}^7},
\\
0
&
\equiv
5\,\Lambda_1^5\Lambda_1^5
-
3\,\Lambda^3\Lambda_{1,1}^7
\blue{+f_1'f_1'M^8},
\endaligned
\]
where the appearing new invariants are already 
known from the statement of the theorem.

\smallskip$\bullet$\
Test whether or not the so obtained three invariants:
\[
\blue{
\Lambda_2^5,\ \ \ \ \ \ \ \ \ \ \ \ \
\Lambda_{1,2}^7,\ \ \ \ \ \ \ \ \ \ \ \ \
M^8,
}
\]
belong or do not belong to the algebra generated by the
previously known invariants. Here in fact, neither $\Lambda_2^5$ nor
$\Lambda_{ 1, 2}^7$, nor $M^8$ belongs to $\C\big[ f_1', f_2',
\Lambda^3, \Lambda_1^5, \Lambda_{ 1, 1}^7 \big]$, as
we have already verified.

\subsection*{ Second loop of the algorithm}
We now restart the process with our new, extended list of 7 invariants
$f_1'$, $f_2'$, $\Lambda^3$, $\Lambda_1^5$, $\Lambda_2^5$,
$\Lambda_{1, 1}^7$, $\Lambda_{ 1, 2}^7$ and $M^8$.

\smallskip$\bullet$\
Compute the ideal of relations between the
6 restricted invariants known at this stage:
\[
{\sf Ideal-Rel}\Big(
f_2'\big\vert_0,\,
\Lambda^3\big\vert_0,\,
\Lambda_1^5\big\vert_0,\,
\Lambda_2^5\big\vert_0,\,
\Lambda_{1,1}^7\big\vert_0,\,
\Lambda_{1,2}^7\big\vert_0,\,
M^8\big\vert_0
\Big).
\]

\smallskip$\bullet$\
Get the 6 equations:
\[
\small
\aligned
0
&
\equiv
3\,\Lambda^3\Lambda^3
-
f_2'\Lambda_1^5\Big\vert_0,
\\
0
&
\equiv
5\,\Lambda^3\Lambda_1^5
-
f_2'\Lambda_{1,1}^7\Big\vert_0,
\\
0
&
\equiv
5\,\Lambda_1^5\Lambda_1^5
-
3\,\Lambda^3\Lambda_{1,1}^7\Big\vert_0,
\endaligned
\]
\[
\small
\aligned
0
&
\equiv
5\,\Lambda^3\Lambda_2^5
-
f_2'\Lambda_{1,2}^7
\Big\vert_0,
\\
0
&
\equiv
5\,\Lambda_1^5\Lambda_2^5
-
3\,\Lambda^3\Lambda_{1,2}^7
\Big\vert_0,
\\
0
&
\equiv
\Lambda_{1,1}^7\Lambda_2^5
-
\Lambda_1^5\Lambda_{1,2}^7
\Big\vert_0.
\endaligned
\]

\smallskip$\bullet$\
Compute the remainders behind a power of $f_1'$:
\[
\aligned
0
&
\equiv
3\,\Lambda^3\Lambda^3
-
f_2'\Lambda_1^5
\green{+f_1'\Lambda_2^5},
\\
0
&
\equiv
5\,\Lambda^3\Lambda_1^5
-
f_2'\Lambda_{1,1}^7
\green{+f_1'\Lambda_{1,2}^7},
\\
0
&
\equiv
5\,\Lambda_1^5\Lambda_1^5
-
3\,\Lambda^3\Lambda_{1,1}^7
\green{+f_1'f_1'M^8},
\endaligned
\]
\[
\small
\aligned
0
&
\equiv
5\,\Lambda^3\Lambda_2^5
-
f_2'\Lambda_{1,2}^7
\blue{+f_1'\Lambda_{2,2}^7},
\\
0
&
\equiv
5\,\Lambda_1^5\Lambda_2^5
-
3\,\Lambda^3\Lambda_{1,2}^7
\green{+f_1'f_2'M^8},
\\
0
&
\equiv
\Lambda_{1,1}^7\Lambda_2^5
-
\Lambda_1^5\Lambda_{1,2}^7
\green{+f_1'\Lambda^3M^8}.
\endaligned
\]

\smallskip$\bullet$\
Get only one new invariant $\Lambda_{ 2, 2}^7$ not belonging to
the algebra generated by already known invariants $\C \big[ f_1', f_2',
\Lambda^3, \Lambda_1^5, \Lambda_2^5, \Lambda_{ 1, 1}^7, \Lambda_{ 1,
2}^7, M^8 \big]$.

\subsection*{ Third loop of the algorithm}
The final list of syzygies, after filling in the remainders and
testing whether new invariants appear, reads:
\[
\small
\aligned
0
&
\equiv
3\,\Lambda^3\Lambda^3
-
f_2'\Lambda_1^5
\green{+f_1'\Lambda_2^5},
\\
0
&
\equiv
5\,\Lambda^3\Lambda_1^5
-
f_2'\Lambda_{1,1}^7
\green{+f_1'\Lambda_{1,2}^7},
\\
0
&
\equiv
5\,\Lambda_1^5\Lambda_1^5
-
3\,\Lambda^3\Lambda_{1,1}^7
\green{+f_1'f_1'M^8},
\endaligned
\]
\[
\small
\aligned
0
&
\equiv
5\,\Lambda^3\Lambda_2^5
-
f_2'\Lambda_{1,2}^7
\green{+f_1'\Lambda_{2,2}^7},
\\
0
&
\equiv
5\,\Lambda_1^5\Lambda_2^5
-
3\,\Lambda^3\Lambda_{1,2}^7
\green{+f_1'f_2'M^8},
\\
0
&
\equiv
\Lambda_{1,1}^7\Lambda_2^5
-
\Lambda_1^5\Lambda_{1,2}^7
\green{+f_1'\Lambda^3M^8},
\endaligned
\]
\[
\small
\aligned
0
&
\equiv
5\,f_2'\Lambda_1^5M^8
+
3\,\Lambda_{1,2}^7\Lambda_{1,2}^7
-
3\,\Lambda_{1,1}^7\Lambda_{2,2}^7
\red{+0},
\\
0
&
\equiv
f_2'\Lambda^3M^8
+
\Lambda_2^5\Lambda_{1,2}^7
-
\Lambda_1^5\Lambda_{2,2}^7
\red{+0},
\\
0
&
\equiv
f_2'f_2'M^8
+
5\,\Lambda_2^5\Lambda_2^5
-
3\,\Lambda^3\Lambda_{2,2}^7
\red{+0}.
\endaligned
\]
Three new syzygies only appear, namely the last three ones above, and
for each of them, the remainders that are a multiple of $f_1'$ are
identically zero, which we specify explicitly by writing ``$\red{ +0
}$''. Importantly, {\em no new invariant appears at this stage}. 

We then claim that the algorithm stops ({\em cf.} also Section~9), and
that the following proposition holds true. In fact, 
the arguments of proof will follow from the general 
theorem of \S9.

\smallskip\noindent{\bf Proposition.}
{\em 
An arbitrary polynomial ${\sf P} = {\sf P} \big( j^4 j \big)$ in ${\sf
E}_4^2$ invariant by reparametrization writes uniquely under the form:
\[
\small
\aligned
{\sf P}\big(j^4j\big)
&
=
\mathcal{Q}
\big(f_1',\,f_2',\,\Lambda_{1,1}^7,\,\Lambda_{2,2}^7,\,M^8\big)
+
\Lambda^3\,
\mathcal{R}
\big(f_1',\,f_2',\,\Lambda_{1,1}^7,\,\Lambda_{2,2}^7,\,M^8\big)
+
\\
&
+
\Lambda_1^5\,
\mathcal{S}
\big(f_1',\,f_2',\,\Lambda_{1,1}^7,\,\Lambda_{2,2}^7,\,M^8\big)
+
\Lambda_2^5\,
\mathcal{T}
\big(f_1',\,f_2',\,\Lambda_{1,1}^7,\,\Lambda_{2,2}^7,\,M^8\big)
+
\\
&
+
\Lambda_{1,2}^7\,
\mathcal{U}
\big(f_1',\,f_2',\,\Lambda_{1,1}^7,\,\Lambda_{2,2}^7,\,M^8\big)
+
\Lambda^3\Lambda_{1,2}^7\,
\mathcal{V}
\big(f_1',\,f_2',\,\Lambda_{1,1}^7,\,\Lambda_{2,2}^7,\,M^8\big),
\endaligned
\]
where $\mathcal{ Q}$, $\mathcal{ R}$, $\mathcal{ S}$, $\mathcal{ T}$,
$\mathcal{ U}$ and $\mathcal{ V}$ are complex polynomials in five
variables subjected to no restriction.
}\medskip

\section*{ \S7.~Action of ${\sf GL}_n(\C)$ and unipotent
reduction}
\label{Section-7}

\subsection*{Sums of irreducible Schur representations} 
The cohomology of Schur bundles $\Gamma^{( \ell_1, \ell_2, \dots,
\ell_n) }\, T_X^*$ on a complex algebraic projective hypersurface $X^n
\subset \P^{ n+1} ( \C)$ being available through Hirzebruch's
Riemann-Roch formula (\S13 below), we should look for a
decomposition of the Demailly-Semple bundle ${\sf E}_{ \kappa, m}^n
T_X^*$ as a direct sum of Schur bundles, at least in the cases where
we understand the algebraic structure of the fiber algebras ${\sf E}_{
\kappa, m}^n$. We recall that according to a fundamental theorem of
representation theory (\cite{ fuha1991}), any group action of ${\sf
GL}_n ( \C)$ on a space of polynomials is 
isomorphic to a certain direct sum of
irreducible Schur representations.

\subsection*{ Action of ${\sf GL}_n ( \C)$ on the jet space}
On this purpose, similarly as in~\cite{ rou2006a}, we therefore define an
appropriate linear action of ${\sf GL}_n ( \C)$ on the $\kappa$-th jet
space $J^\kappa ( \C, \, \C^n)$. By definition, an arbitrary element
{\sf w} of ${\sf GL}_n ( \C)$ written in matrix form:
\[
\text{\sf w}
=
\left(
\begin{array}{ccc}
w_{11} & \cdots & w_{1n}
\\
\vdots & \ddots & \vdots
\\
w_{n1} & \cdots & w_{nn}
\end{array}
\right)
\]
shall transform the collection $\big( f_1^{ ( \lambda)},
\dots, f_n^{ ( \lambda)} \big)$ of the $n$ components of
a $\kappa$-jet $j^\kappa f$
at each $\lambda$-th jet level, just by matrix multiplication:
\[
\left\{
\aligned
\text{\sf w}\cdot f_1^{(\lambda)}
&
=
w_{11}\,f_1^{(\lambda)}
+\cdots+
w_{1n}\,f_n^{(\lambda)}
\\
\cdots\ \ \ \ 
& \ \ \ \ \ \ \ \
\cdots \ \ \ \ \ \ \ \ \ \ \ \ \ \ \ \ \ \ \ \ \
\cdots
\\
\text{\sf w}\cdot f_n^{(\lambda)}
&
=
w_{n1}\,f_1^{(\lambda)}
+\cdots+
w_{nn}\,f_n^{(\lambda)},
\endaligned\right.
\]
with the same matrix {\sf w} at each jet level $\lambda$
with $1 \leqslant \lambda \leqslant \kappa$.

\smallskip\noindent{\bf Definition.}
A polynomial ${\sf P} \big( j^\kappa f \big)$ invariant by
reparametrization will be called a \underline{\sl bi-invariant} if it
is a {\sl vector of highest weight} for this representation of ${\sf
GL}_n ( \C)$, namely if it is invariant by the {\sl unipotent
subgroup} ${\sf U}_n ( \C) \subset {\sf GL}_n ( \C)$ constituted by
(unipotent) matrices of the form:
\[
\text{\sf u}
=
\left(
\begin{array}{ccccc}
1 & 0 & 0 & \cdots & 0
\\
u_{21} & 1 & 0 & \cdots & 0
\\
u_{31} & u_{32} & 1 & \cdots & 0
\\
\vdots & \vdots & \vdots & \ddots & \vdots
\\
u_{n1} & u_{n2} & u_{n3} & \cdots & 1
\end{array}
\right).
\] 
The vector space of bi-invariant polynomials ${\sf P}$ thus
satisfies:
\[
\boxed{
{\sf P}\big(j^\kappa(f\circ\phi)\big)
=
(\phi')^m\cdot
{\sf P}\big((j^\kappa f)\circ\phi\big)
\ \ \ \ \ \ \ \ \ \
\text{\rm and}
\ \ \ \ \ \ \ \ \ \
{\sf P}\big(\text{\sf u}\cdot j^\kappa f)
=
{\sf P}\big(j^\kappa f)}\,.
\]
In the sequel, the vector space of bi-invariants of weight $m$ will be
denoted by ${\sf UE}_{ \kappa, m}^n$. Also, one defines the graded
algebra of bi-invariants ${\sf UE}_\kappa^n := \bigoplus_{ m \geqslant
1} \, {\sf UE}_{ \kappa, m}^n$ with of course ${\sf UE}_{ \kappa,
m_1}^n \cdot {\sf UE}_{ \kappa, m_2}^n \subset {\sf UE}_{ \kappa, m_1
+ m_2}^n$.

\medskip

Without delay, we emphasize four fundamental observations.

\medskip$\bullet$
The full space ${\sf E}_{ \kappa, m}^n$ is obtained as just the ${\sf
GL}_n ( \C)$-orbit of ${\sf UE}_{ \kappa, m}^n$.

\medskip$\bullet$
The algebraic structure of ${\sf UE}_\kappa^n$
is {\em always} much simpler than that of ${\sf E}_\kappa^n$. 
For instance:

\medskip\ \ \ \ \ ---\ \ \ 
{\small ${\sf UE}_3^3$ is generated by only 4 bi-invariant
polynomials\footnote{\, {\it See} the proposition on
p.~\pageref{bi-invariant-3-3} below, or the considerations on
pp.~931--932 in~\cite{ mer2008a}. } $f_1'$, $\Lambda_{ 1,2}^3$,
$\Lambda_{1,2; \, 1}^5$ and $D_{1,2,3}^6$ which are {\em algebraically
independent} (no syzygy!), whereas, according to~\cite{
rou2004, rou2006a} or to the description given
on p.~\pageref{theorem-3-3} here, the full algebra ${\sf E}_3^3$ is
generated by 16 invariants, submitted to the three complicated
families of syzygies developed on p.~\pageref{syzygies-3-3}.}

\medskip\ \ \ \ \ ---\ \ \ 
{\small ${\sf UE}_4^2$ is generated by the 5 bi-invariant
polynomials $f_1'$, $\Lambda^3$, $\Lambda_1^5$, $\Lambda_{ 1,
1}^7$ and $M^8$, whose ideal of relations is principal, generated by
the single syzygy:
\[
0
\overset{4}{\equiv}
f_1'f_1'M^8
-
3\,\Lambda^3\Lambda_{1,1}^7
+
5\,\Lambda_1^5\Lambda_1^5,
\]
while, according to the theorem on p.~\pageref{theorem-2-4}, the
full algebra ${\sf E}_4^2$ is generated by 9 invariants submitted to 9
fundamental syzygies.}

\medskip\ \ \ \ \ ---\ \ \ 
{\small
We will establish that ${\sf UE}_4^4$ is generated by 16 mutually
independent bi-invariant polynomials, while ${\sf E}_4^4$ is generated
by 2835 polynomials invariant by reparametrization. Also, we will
show 41 syzygies generate the ideal of relations between (the
restriction to $\{ f_1' = 0\}$ of) these 16 generators of ${\sf
UE}_4^4$, while we ignore the structure of the (presumably out of
human scale) ideal of relations between the 2835 
generators of ${\sf E}_4^4$.
}

\medskip\ \ \ \ \ ---\ \ \ 
{\small
We will establish that ${\sf UE}_5^2$ is generated by 17 mutually
independent bi-invariant polynomials, while ${\sf E}_5^2$ is generated
by 56 polynomials invariant by reparametrization. We will show 66
syzygies generating the ideal of relations between (the restriction to
$\{ f_1' = 0\}$ of) these 17 generators of ${\sf UE}_5^2$.}

\medskip$\bullet$
In any case, if we can show that ${\sf UE}_\kappa^n$ is, for a certain
$n$ and for a certain $\kappa$, generated as an algebra by a finite
number of bi-invariants, we may easily deduce as a corollary finite
generation of the full algebra ${\sf E}_\kappa^n$. For instance:

\medskip\ \ \ \ \ ---\ \ \
{\small For $n = \kappa = 3$, computing the ${\sf GL}_3( \C)$-orbit of
the 4 bi-invariants $f_1'$, $\Lambda_{ 1, 2}^3$, $\Lambda_{ 1, 2; \,
1}^5$ and $D_{ 1, 2, 3}^6$ amounts to polarize their lower indices,
which yields the invariants $f_i'$, $\Lambda_{ i,j}^3$, $\Lambda_{
i,j;\, k}^5$ and $D_{i,j,k}^6$ generating ${\sf E}_3^3$. }

\medskip\ \ \ \ \ ---\ \ \
{\small For $n = 2$ and $\kappa = 4$, computing the ${\sf GL}_2(
\C)$-orbit of the 5 bi-invariants $f_1'$, $\Lambda^3$, $\Lambda_1^5$,
$\Lambda_{ 1, 1}^7$ and $M^8$ again amounts to polarize their lower
indices, which yields the invariants $f_i'$, $\Lambda^3$,
$\Lambda_i^5$, $\Lambda_{ i,j}^7$ and $M^8$ generating
${\sf E}_4^2$.}

\medskip$\bullet$
Finally, for applications to Kobayashi hyperbolicity (which involves
estimating the Euler-Poincaré characteristic of ${\sf E}_{\kappa, m}^n
T_X^*$), it is useless to look for a complete understanding of the
algebraic structure of ${\sf E}_\kappa^n$, and it only suffices to
possess a complete description of the algebra of bi-invariants ${\sf
UE}_\kappa^n$. In fact, as will be (re)explained in \S12, each
bi-invariant will correspond to one and to only one Schur bundle.

\medskip
\centerline{\fbox{\sf So from now on, 
we focus our attention on bi-invariants}}

\subsection*{ Initial representation of bi-invariants}
We now restart with the initial, rational expression of any polynomial
invariant by reparametrization provided by the lemma on
p.~\pageref{initial-n-k} and we want to determine when such a polynomial
is, in addition, invariant by the unipotent action.

To begin with, we consider the subgroup ${\sf U}_n^* (\C)$ of ${\sf
U}_n ( \C)$ generated by matrices of the form:
\[
\text{\sf u}^*
=
\left(
\begin{array}{ccccc}
1 & 0 & 0 & \cdots & 0
\\
u_{21} & 1 & 0 & \cdots & 0
\\
u_{31} & 0 & 1 & \cdots & 0
\\
\vdots & \vdots & \vdots & \ddots & \vdots
\\
u_{n1} & 0 & 0 & \cdots & 1
\end{array}
\right).
\]
Clearly, the components of the first order jet $j^1 f$ are modified by
the action of $\text{\sf u}^*$:
\[
\left\{
\aligned
\text{\sf u}^*\cdot f_1'
&
=
f_1',
\\
\text{\sf u}^*\cdot f_2'
&
=
f_2'+u_{21}f_1',
\\
\text{\sf u}^*\cdot f_3'
&
=
f_3'+u_{31}f_1',
\\
\cdots\ \ \,
&\ \ \ \ \ \ \ \ \ \
\cdots
\\
\text{\sf u}^*\cdot f_n'
&
=
f_n'+u_{n1}f_1'.
\endaligned\right.
\]
On the other hand, all the $\Lambda_{ 1, i}^3$ are 
left invariant:
\[
\text{\sf u}^*\cdot\Lambda_{1,i}^3
=
\text{\sf u}^*\cdot\big[f_i',\,f_1'\big]
=
\big[f_i'+u_{i1}f_1',\,f_1'\big]
=
\big[f_i',\,f_1'\big]
=
\Lambda_{1,i}^3,
\]
and in fact, more generally, one may verify that the same is
true of higher $\Lambda$'s:
\[
\text{\sf u}^*\cdot\Lambda_{1,i;\,1}^5
=
\Lambda_{1,i;\,1}^5,
\ \ \ \ \ \ \ \
\text{\sf u}^*\cdot\Lambda_{1,i;\,1,1}^7
=
\Lambda_{1,i;\,1,1}^7,\,\,
\dots\dots,\,\,
\text{\sf u}^*\cdot\Lambda_{1,i;\,1^{\kappa-2}}^{2\kappa-1}
=
\Lambda_{1,i;\,1^{\kappa-2}}^{2\kappa-1},
\]
for any $i = 2, 3, \dots, n$. Consequently, the requirement that a
polynomial invariant by reparametrization 
${\sf P} \big( j^\kappa f\big) \in {\sf E}_{ \kappa, m}^n$
be in addition also invariant by the
unipotent subgroup ${\sf U}_n^* ( \C) \subset {\sf U}_n ( \C)$, namely
${\sf u}^* \cdot {\sf P} \big( j^\kappa f \big) = {\sf P} \big(
j^\kappa f \big)$, shall be written in length as follows, when
employing the mentioned representation given on
p.~\pageref{initial-n-k}:
\[
\small
\aligned
&
\sum_a\,
(f_1')^a\,{\sf P}_a
\Big(
f_2'+u_{21}f_1',\,f_3'+u_{31}f_1',\,\dots,\,f_n'+u_{n1}f_1',
\\
&
\ \ \ \ \ \ \ \ \ \ \ \ \ \ \ \ \ \ \ \ \ \ \ \ 
\ \ \ \ \ \ \ \ \ \ \ \ \ 
\Lambda_{1,2}^3,\,\dots,\,\Lambda_{1,n}^3,\,
\dots\dots,\,
\Lambda_{1,2;\,1^{\kappa-2}}^{2\kappa-1},\,\dots,\,
\Lambda_{1,n;\,1^{\kappa-2}}^{2\kappa-1}
\Big)
=
\endaligned
\]
\[
\small
\aligned
&
=
\sum_a\,
(f_1')^a\,{\sf P}_a
\Big(f_2',\,f_3',\,\dots,\,f_n',
\\
&
\ \ \ \ \ \ \ \ \ \ \ \ \ \ \ \ \ \ \ \ \ \ \ \ 
\ \ \ \ \ \ \ \ \ \ \ \ \ \ \ \
\Lambda_{1,2}^3,\,\dots,\,\Lambda_{1,n}^3,\,
\dots\dots,\,
\Lambda_{1,2;\,1^{\kappa-2}}^{2\kappa-1},\,\dots,\,
\Lambda_{1,n;\,1^{\kappa-2}}^{2\kappa-1}
\Big).
\endaligned
\]
Because the $n + (n-1)(\kappa -1)$ invariants $f_1', \dots, f_n'$,
$\Lambda_{ 1, i; \, 1^{ \lambda-2}}^{ 2\lambda - 1}$, $2 \leqslant i
\leqslant n$, $2 \leqslant \lambda \leqslant \kappa$, are algebraically
independent, we deduce
that each ${\sf P}_a$ must be independent of $f_2', f_3', \dots,
f_n'$, so that we come to the simpler rational expression:
\[
{\sf R}
=
\sum_a\,
(f_1')^a\,{\sf P}_a
\Big(
\Lambda_{1,2}^3,\,\dots,\,\Lambda_{1,n}^3,\,
\dots\dots,\,
\Lambda_{1,2;\,1^{\kappa-2}}^{2\kappa-1},\,\dots,\,
\Lambda_{1,n;\,1^{\kappa-2}}^{2\kappa-1}
\Big),
\]
which is however not yet invariant under the full unipotent action.

\subsection*{ Second unipotent subgroup}
Next, we consider the subgroup ${\sf U}_n^\sharp ( \C) \subset {\sf
U}_n (\C)$ constituted by matrices of the form:
\[
\text{\sf u}^\sharp
=
\left(
\begin{array}{ccccc}
1 & 0 & 0 & \cdots & 0
\\
0 & 1 & 0 & \cdots & 0
\\
0 & u_{32} & 1 & \cdots & 0
\\
\vdots & \vdots & \vdots & \ddots & \vdots
\\
0 & u_{n2} & u_{n3} & \cdots & 1
\end{array}
\right).
\]
Since ${\sf U}_n^* ( \C)$ and ${\sf U}_n^\sharp ( \C)$ clearly
generate the full unipotent group ${\sf U}_n (\C)$, it now only
remains to require the ${\sf U}_n^\sharp ( \C)$-invariance for the
rational expression ${\sf R}$ obtained just above.

The requirement $\text{\sf u}^\sharp \cdot {\sf R} = {\sf R}$ can
in turn be written in length as follows:
:
\[
\aligned
\text{\sf u}^\sharp
\bigg(
\sum_a\,
(f_1')^a\,{\sf P}_a
\Big(
\Lambda_{1,2}^3,\,
\dots,\,
\Lambda_{1,n;\,1^{\kappa-2}}^{2\kappa-1}
\Big)
\bigg)
&
=
\sum_a\,(f_1')^a\,{\sf P}_a
\Big(
\text{\sf u}^\sharp\cdot\Lambda_{1,2}^3,\,\dots,\,
\text{\sf u}^\sharp\cdot\Lambda_{1,n;\,1^{\kappa-2}}^{2\kappa-1}
\Big)
\\
&
=
\sum_a\,(f_1')^a\,{\sf P}_a
\Big(
\Lambda_{1,2}^3,\,\dots,\,
\Lambda_{1,n;\,1^{\kappa-2}}^{2\kappa-1}
\Big).
\endaligned
\]
But on the other hand, for any $\lambda$ with $2 \leqslant \lambda
\leqslant \kappa$, one may verify that the action of $\text{\sf
u}^\sharp$ on the initial $\Lambda$-invariants appearing as arguments
of ${\sf R}$ is given by the triangular formulas:
\[
\aligned
\text{\sf u}^\sharp\cdot
\Lambda_{1,2;\,1^{\lambda-2}}^{2\lambda-1}
&
=
\Lambda_{1,2;\,1^{\lambda-2}}^{2\lambda-1},
\\
\text{\sf u}^\sharp\cdot
\Lambda_{1,3;\,1^{\lambda-2}}^{2\lambda-1}
&
=
\Lambda_{1,3;\,1^{\lambda-2}}^{2\lambda-1}
+
u_{32}\,
\Lambda_{1,2;\,1^{\lambda-2}}^{2\lambda-1},
\\
\text{\sf u}^\sharp\cdot
\Lambda_{1,4;\,1^{\lambda-2}}^{2\lambda-1}
&
=
\Lambda_{1,4;\,1^{\lambda-2}}^{2\lambda-1}
+
u_{43}\,
\Lambda_{1,3;\,1^{\lambda-2}}^{2\lambda-1}
+
u_{42}\,
\Lambda_{1,2;\,1^{\lambda-2}}^{2\lambda-1},
\\
\cdots \ \ \ \ \ \ \
& 
\ \ \ \ \ \ \ \ \ 
\cdots 
\ \ \ \ \ \ \ \ \ \ \ \ \ \ \ \ \ \ 
\cdots 
\ \ \ \ \ \ \ \ \ \ \ \ \ \ \ \ \ \ \ \
\cdots
\\
\text{\sf u}^\sharp\cdot
\Lambda_{1,n;\,1^{\lambda-2}}^{2\lambda-1}
&
=
\Lambda_{1,n;\,1^{\lambda-2}}^{2\lambda-1}
+
u_{n,n-1}\,
\Lambda_{1,n-1;\,1^{\lambda-2}}^{2\lambda-1}
+\cdots+
u_{n2}\,
\Lambda_{1,2;\,1^{\lambda-2}}^{2\lambda-1}.
\endaligned
\]
The algebraic independency of $f_1'$, $\Lambda_{ 1, i; \, 1^{\lambda -
2}}^{ 2 \lambda - 1}$ then implies that such an ${\sf R}$ is ${\sf
U}_n^\sharp (\C)$-invariant if and only if every ${\sf P}_a$ is so,
namely if and only if the following identity holds:
\[
\small
\aligned
{\sf P}_a
&
\Big(
A_{1,2}^3,\,A_{1,3}^3+u_{32}\,A_{1,2}^3,\,\dots,\,
A_{1,n}^3+u_{n,n-1}\,A_{1,n-1}^3+\cdots+u_{n2}\,A_{1,2}^3,
\\
&\ \ \
A_{1,2}^5,\,A_{1,3}^5+u_{32}\,A_{1,2}^5,\,\dots,\,
A_{1,n}^5+u_{n,n-1}\,A_{1,n-1}^5+\cdots+u_{n2}\,A_{1,2}^5,
\\
& \ \ \ 
\cdots 
\ \ \ \ \ \ \ \ \ \
\cdots
\ \ \ \ \ \ \ \ \ \ \ \ \ \ \ \ \ \ \ \ \ \ \ \ \ \ \ \ \ \
\cdots
\\
&\ \ \
A_{1,2}^{2\kappa-1},\,
A_{1,3}^{2\kappa-1}+u_{32}\,A_{1,2}^{2\kappa-1},\,\dots,\,
A_{1,n}^{2\kappa-1}+u_{n,n-1}\,A_{1,n-1}^{2\kappa-1}
+\cdots+
u_{n2}\,A_{1,2}^{2\kappa-1}
\Big)
=
\endaligned
\]
\[
\small
\aligned
&
=
{\sf P}_a
\Big(
A_{1,2}^3,\,A_{1,3}^3,\,\dots,\,A_{1,n}^3,
\\
&\ \ \ \ \ \ \ \ \ \ \ \
A_{1,2}^5,\,A_{1,3}^5,\,\dots,\,A_{1,n}^5,
\\
&\ \ \ \ \ \ \ \ \ \ \ \
\cdots \ \ \ \
\cdots \ \ \ \ \ \ \ \ \ \ \ \
\cdots
\\
&\ \ \ \ \ \ \ \ \ \ \ \
A_{1,2}^{2\kappa-1},\,A_{1,3}^{2\kappa-1},\,\dots,\,
A_{1,n}^{2\kappa-1}
\Big),
\endaligned
\]
as polynomials in $\C \big[ A_{ 1, 2}^3, \dots, A_{ 1, n}^3, \dots
\dots, A_{ 1, 2}^{ 2\kappa - 1}, \dots, A_{ 1, n}^{ 2 \kappa - 1}
\big]$, for every $\text{\sf u}^\sharp$, and for every $a$ with $-
\frac{ m- 1}{ m} \kappa \leqslant a \leqslant m$.

Here, we recognize a full unipotent action, acted by means of a 
general $(n-1) \times (n-1)$ unipotent matrice of the form:
\[
\left(
\begin{array}{cccc}
1 & 0 & \cdots & 0
\\
u_{32} & 1 & \cdots & 0
\\
\vdots & \vdots & \ddots & \vdots
\\
u_{n2} & u_{n3} & \cdots & 1
\end{array}
\right)
\in{\sf U}_{n-1}(\C),
\]
on the set of the $\kappa - 1$ vectors of $\C^{ n-1}$ defined by:
\[
\big(
A_{1,2}^{2\lambda-1},\,
A_{1,3}^{2\lambda-1},\,\dots,\,
A_{1,n}^{2\lambda-1}
\big)
\ \ \ \ \ \ \ \ \ \ \ \ \ \
{\scriptstyle{(2\,\leqslant\,\lambda\,\leqslant\,\kappa)}}.
\]
It is known (\cite{ krpr1996, pro2007}) that the invariants for such
an action are constituted by all the minors:
\[
\small
\aligned
\Pi_2^{\lambda_2}
:=
A_{1,2}^{2\lambda_2-1},
\ \ \ \ \ \ \ \ \
\Pi_{2,3}^{\lambda_2,\lambda_3}
&
:=
\left\vert
\begin{array}{cc}
A_{1,2}^{2\lambda_2-1} & A_{1,3}^{2\lambda_2-1}
\\
A_{1,2}^{2\lambda_3-1} & A_{1,3}^{2\lambda_3-1}
\end{array}
\right\vert,
\\
\Pi_{2,3,4}^{\lambda_2,\lambda_3,\lambda_4}
&
:=
\left\vert
\begin{array}{ccc}
A_{1,2}^{2\lambda_2-1} & A_{1,3}^{2\lambda_2-1} & A_{1,4}^{2\lambda_2-1}
\\
A_{1,2}^{2\lambda_3-1} & A_{1,3}^{2\lambda_3-1} & A_{1,4}^{2\lambda_3-1}
\\
A_{1,2}^{2\lambda_4-1} & A_{1,3}^{2\lambda_4-1} & A_{1,4}^{2\lambda_4-1}
\end{array}
\right\vert,
\endaligned
\]
and generally:
\[
\Pi_{2,3,4,\dots,n_1}^{\lambda_2,\lambda_3,\lambda_4,\dots,\lambda_{n_1}}
:=
\left\vert
\begin{array}{ccccc}
A_{1,2}^{2\lambda_2-1} & A_{1,3}^{2\lambda_2-1} & A_{1,4}^{2\lambda_2-1}
& \cdots & A_{1,n_1}^{2\lambda_2-1}
\\
A_{1,2}^{2\lambda_3-1} & A_{1,3}^{2\lambda_3-1} & A_{1,4}^{2\lambda_3-1}
& \cdots & A_{1,n_1}^{2\lambda_3-1}
\\
A_{1,2}^{2\lambda_4-1} & A_{1,3}^{2\lambda_4-1} & A_{1,4}^{2\lambda_4-1}
& \cdots & A_{1,n_1}^{2\lambda_4-1}
\\
\vdots & \vdots & \vdots & \ddots & \vdots
\\
A_{1,2}^{2\lambda_{n_1}-1} & A_{1,3}^{2\lambda_{n_1}-1} &
A_{1,4}^{2\lambda_{n_1}-1} & \cdots &
A_{1,n_1}^{2\lambda_{n_1}-1}
\end{array}
\right\vert,
\]
for all $n_1$ from $n_1 = 1$ up to $n_1 = n$, and for arbitrary
$\lambda_j$ with $2 \leqslant \lambda_j \leqslant \kappa$. In fact,
one immediately sees that these minors are obviously invariant by the
unipotent action of ${\sf U}_{ n-1} ( \C)$, thanks to the fact that
column linear dependence leaves untouched any determinant.

\THEOREM

\smallskip\noindent\fbox{\bf THEOREM}
\label{bi-invariant-n-k}
{\sf\em 
In dimension $n\geqslant 1$ and for jets of arbitrary order $\kappa
\geqslant 1$, every bi-invariant polynomial ${\sf BP} = {\sf BP} \big(
j^\kappa f \big)$ invariant by reparametrization and invariant under the
unipotent action writes under the form:
\[
\footnotesize
\aligned
{\sf BP}\big(j^\kappa f\big)
=\!\!
\sum_{-\frac{\kappa-1}{\kappa}m\leqslant a\leqslant m}\!\!
(f_1')^a\,
{\sf BP}_a
\bigg(
\left\vert
\begin{array}{cccc}
\Lambda_{1,2}^{2\lambda_2-1} & \Lambda_{1,3}^{2\lambda_2-1} &
\cdots & \Lambda_{1,n_1}^{2\lambda_2-1}
\\
\Lambda_{1,2}^{2\lambda_3-1} & \Lambda_{1,3}^{2\lambda_3-1} &
\cdots & \Lambda_{1,n_1}^{2\lambda_3-1}
\\
\vdots & \vdots & \ddots & \vdots
\\
\Lambda_{1,2}^{2\lambda_3-1} & \Lambda_{1,3}^{2\lambda_3-1} &
\cdots & \Lambda_{1,n_1}^{2\lambda_3-1}
\end{array}
\right\vert_{n_1=1,2\dots,n}^{2\leqslant
\lambda_2,\dots,\lambda_{n_1}\leqslant\kappa}
\bigg),
\endaligned
\]
for certain specific polynomials ${\sf BP}_a$ which 
depend upon ${\sf BP} ( j^\kappa f)$.
}

\stopTHEOREM

\subsection*{ The case $n = \kappa = 3$}
After ${\sf U}_3^* (\C)$-reduction, an arbitrary 
element of ${\sf UE}_{3, m}^3$ writes:
\[
{\sf R}
=
\sum_{-\frac{2}{3}m\leqslant a\leqslant m}\,(f_1')^a\,{\sf P}_a
\Big(
\Lambda_{1,2}^3,\,\Lambda_{1,3}^3,\,
\Lambda_{1,2;\,1}^5,\,\Lambda_{1,3;\,1}^5
\Big).
\]
Then the ${\sf U}_3^\sharp ( \C)$-reduction presented above shows that
there are four initial bi-invariants, namely the three obvious ones
$f_1'$, $\Lambda_{ 1, 2}^3$, $\Lambda_{ 1, 2; \, 1}^5$ together with:
\[
\left\vert
\begin{array}{cc}
\Lambda_{1,2}^3 & \Lambda_{1,3}^3
\\
\Lambda_{1,2;\,1}^5 & \Lambda_{1,3;\,1}^5
\end{array}
\right\vert
=
f_1'f_1'\cdot
\left\vert
\begin{array}{ccc}
f_1' & f_2' & f_3' 
\\
f_1'' & f_2'' & f_3''
\\
f_1''' & f_2''' & f_3'''
\end{array}
\right\vert
=:
f_1'f_1'\cdot
D_{1,2,3}^6,
\]
where the first equality, which follows from a direct calculation,
gives birth to the three-dimensional Wronskian. By pluging this minor
in the above rational expression of ${\sf R}$, we obtain that any
bi-invariant polynomial in ${\sf UE}_{3, m}^3$ writes under the form:
\[
{\sf BP}\big(j^3f\big)
=
\sum_{-\frac{2}{3}m\leqslant a\leqslant m}\,
(f_1')^a\,\widetilde{\sf P}_a
\Big(
\Lambda_{1,2}^3,\,\Lambda_{1,2;\,1}^5,\,D_{1,2,3}^6
\Big),
\]
for certain (new) polynomials $\widetilde{\sf P}_a$. More is true, for
we claim that there are no negative powers of $f_1'$ anymore in 
such a rational representation.

\smallskip\noindent{\bf Proposition.}
\label{bi-invariant-3-3}
{\em 
Any bi-invariant polynomial ${\sf BP} \in {\sf UE}_{3, m}^3$ writes
uniquely under the form:
\[
{\sf BP}\big(j^3f\big)
=
\sum_{0\leqslant a\leqslant m}\,
(f_1')^a\,{\sf BP}_a
\Big(
\Lambda_{1,2}^3,\,\Lambda_{1,2;\,1}^5,\,D_{1,2,3}^6
\Big),
\] 
where the ${\sf BP}_a$ are arbitrary polynomials. In fact:
\[
\boxed{
{\sf UE}\big(j^3f\big)
=
\C\big[f_1',\,\Lambda_{1,2}^3,\,\Lambda_{1,2;\,1}^5,\,
D_{1,2,3}^6\big]}\,.
\]
}\medskip

\proof
One verifies at first sight that, after setting $f_1' = 0$, the 3
restricted invariants:
\[
\Lambda_{1,2}^3\big\vert_0
=
-f_1''f_2',
\ \ \ \ \ \ \ \ \ \ \ \ \ \
\Lambda_{1,2;\,1}^5\big\vert_0
=
3\,f_1''f_2'f_1''
\ \ \ \ \ \ \ \
\text{\rm and}
\ \ \ \ \ \ \ \
D_{1,2,3}^6\big\vert_0
=
\left\vert
\begin{array}{ccc}
0 & f_2' & f_3' 
\\
f_1'' & f_2'' & f_3'' 
\\
f_1''' & f_2''' & f_3'''
\end{array}
\right\vert
\]
are mutually algebraically independent. Suppose then by contradiction
that the expression:
\[
{\sf BP}\big(j^3f\big)
=
\sum_{-a_0\leqslant a\leqslant m}\,
(f_1')^a\,
\widetilde{\sf P}_a\Big(\Lambda_{1,2}^3,\,
\Lambda_{1,2;\,1}^5,\,D_{1,2,3}^6\Big),
\]
starts with a not identically zero $\widetilde{ \sf P}_{- a_0} \big(
A^3, A^5, \Delta^6) \not \equiv 0$ for some smallest negative power
$-a_0< 0$ of $f_1'$. Multiplying both sides by $(f_1')^{ a_0}$ and
setting $f_1' = 0$ afterwards, the left term $(f_1')^{ a_0} \, {\sf
BP} \big( j^3 f \big)$ then vanishes, hence one would derive an
identity:
\[
0
\equiv
\widetilde{\sf P}_{-a_0}
\Big(\Lambda_{1,2}^3\big\vert_0,\,\Lambda_{1,2;\,1}^5\big\vert_0,\,
D_{1,2,3}^6\big\vert_0\Big)
\]
between restricted bi-invariants which would then entail $\widetilde{
\sf P}_{- a_0} \equiv 0$ because the arguments are algebraically
independent, a contradiction.

Consequently, the rational expression for ${\sf BP} \big( j^3 f\big)$
was already polynomial and inversely, every arbitrary polynomial in
$\C \big[ f_1', \, \Lambda_{ 1, 2}^3, \, \Lambda_{ 1, 2; \, 1}^5, \,
D_{ 1, 2, 3}^6 \big]$ obviously is a bi-invariant.
\endproof

\subsection*{ The case $n = \kappa = 4$}
After ${\sf U}_4^* (\C)$-reduction, an arbitrary 
element of ${\sf UE}_{4, m}^4$ writes under the form:
\[
\aligned
{\sf R}
=
\sum_{-\frac{3}{4}m\leqslant a\leqslant m}\,(f_1')^a\,{\sf P}_a
\Big(
\Lambda_{1,2}^3,\,\Lambda_{1,3}^3,\,\Lambda_{1,4}^3,\,
\Lambda_{1,2;\,1}^5,\,\Lambda_{1,3;\,1}^5,\,\Lambda_{1,4;\,1}^5,\,
\Lambda_{1,2;\,1,1}^7,\,\Lambda_{1,3;\,1,1}^7,\,\Lambda_{1,4;\,1,1}^7
\Big).
\endaligned
\]
Then the ${\sf U}_4^\sharp ( \C)$-reduction presented above shows that
there are the 4 obvious initial bi-invariants: 
\[
f_1',
\ \ \ \ \ \ \ \ \ 
\Lambda_{1,2}^3,
\ \ \ \ \ \ \ \ \ 
\Lambda_{1,2;\,1}^5
\ \ \ \ \ \ \ \
\text{\rm and}
\ \ \ \ \ \ \ \ \
\Lambda_{1,2;\,1,1}^7,
\]
together with the 4 further ones: 
\[
D^6,
\ \ \ \ \ \ \ \ \ 
D^8
=
\big[D^6,\,f_1'\big],
\ \ \ \ \ \ \ \ \ 
N^{10}
\ \ \ \ \ \ \ \
\text{\rm and}
\ \ \ \ \ \ \ \
W^{10},
\]
that are obtained by dividing the 4 minors involving the $\Lambda$'s
by the maximal power of $f_1'$ which appears in factor, namely:
\[
\aligned
\left\vert
\begin{array}{cc}
\Lambda_{1,2}^3 & \Lambda_{1,3}^3
\\
\Lambda_{1,2;1}^5 & \Lambda_{1,3;1}^5 
\end{array}
\right\vert
&
\equiv
f_1'f_1'D^6,
\\
\left\vert
\begin{array}{cc}
\Lambda_{1,2}^3 & \Lambda_{1,3}^3
\\
\Lambda_{1,2;1,1}^7 & \Lambda_{1,3;1,1}^7 
\end{array}
\right\vert
&
\equiv
f_1'f_1'D^8,
\endaligned
\]
\[
\aligned
\left\vert
\begin{array}{cc}
\Lambda_{1,2;1}^5 & \Lambda_{1,3;1}^5
\\
\Lambda_{1,2;1,1}^7 & \Lambda_{1,3;1,1}^7 
\end{array}
\right\vert
&
\equiv
f_1'f_1'N^{10},
\\
\left\vert
\begin{array}{ccc}
\Lambda_{1,2}^3 & \Lambda_{1,3}^3 & \Lambda_{1,4}^3
\\
\Lambda_{1,2;1}^5 & \Lambda_{1,3;1}^5 & \Lambda_{1,4;1}^5
\\
\Lambda_{1,2;1,1}^7 & \Lambda_{1,3;1,1}^7 & \Lambda_{1,4;1,1}^7 
\end{array}
\right\vert
&
\equiv
f_1'f_1'f_1'f_1'f_1'\,W^{10},
\endaligned
\]
where the last one behind $(f_1')^5$ appears to be equal to the 
four-dimensional {\sl Wronskian}:
\[
W^{10}
:=
\left\vert
\begin{array}{cccc}
f_1' & f_2' & f_3' & f_4 '
\\
f_1'' & f_2'' & f_3'' & f_4''
\\
f_1''' & f_2''' & f_3''' & f_4'''
\\
f_1'''' & f_2'''' & f_3'''' & f_4''''
\end{array}
\right\vert,
\]
and where the first three ones are explicitly defined by:
\[
D^6
:=
\left\vert
\begin{array}{ccc}
f_1' & f_2' & f_3'
\\
f_1'' & f_2'' & f_3''
\\
f_1''' & f_2''' & f_3'''
\end{array}
\right\vert,
\]
\[
D^8
:=
f_1'
\left\vert
\begin{array}{ccc}
f_1' & f_2' & f_3'
\\
f_1'' & f_2'' & f_3''
\\
f_1'''' & f_2'''' & f_3''''
\end{array}
\right\vert
-
6\,f_1''
\left\vert
\begin{array}{ccc}
f_1' & f_2' & f_3'
\\
f_1'' & f_2'' & f_3''
\\
f_1''' & f_2''' & f_3'''
\end{array}
\right\vert,
\]
\[
\aligned
N^{10}
&
:=
f_1'f_1'\,
\left\vert
\begin{array}{ccc}
f_1' & f_2' & f_3'
\\
f_1''' & f_2''' & f_3'''
\\
f_1'''' & f_2'''' & f_3''''
\end{array}
\right\vert
-
3\,f_1'f_1''\,
\left\vert
\begin{array}{ccc}
f_1' & f_2' & f_3'
\\
f_1'' & f_2'' & f_3''
\\
f_1'''' & f_2'''' & f_3''''
\end{array}
\right\vert
+
\\
&\ \ \ \ \
+
4\,f_1'f_1'''\,
\left\vert
\begin{array}{ccc}
f_1' & f_2' & f_3'
\\
f_1'' & f_2'' & f_3''
\\
f_1''' & f_2''' & f_3'''
\end{array}
\right\vert
+
3\,f_1''f_1''
\left\vert
\begin{array}{ccc}
f_1' & f_2' & f_3'
\\
f_1'' & f_2'' & f_3''
\\
f_1''' & f_2''' & f_3'''
\end{array}
\right\vert.
\endaligned
\]
By pluging these 8 bi-invariants in the rational expression
written on p.~\pageref{bi-invariant-n-k}, we obtain that any
bi-invariant polynomial in ${\sf UE}_{ 4, m}^4$ writes under the form:
\[
\boxed{
{\sf BP}\big(j^4j\big)
=
\sum_{-\frac{3}{4}m\leqslant a\leqslant m}\,
(f_1')^a\,\widetilde{\sf P}_a
\Big(\Lambda^3,\,\Lambda^5,\,\Lambda^7,\,D^6,\,D^8,\,
N^{10},\,
W^{10}\Big)}\,.
\]
This expression will be the very starting point for the application of
our general algorithm, to be presented in Section~9 below. In fact, as
in the case $n = 2$, $\kappa = 4$ of Section~6, there will appear
further independent {\sl ghost bi-invariants hidden behind powers of
$f_1'$}.

\section*{ \S8.~Counterexpectation: insufficiency of
bracket invariants}
\label{Section-8}

According to the unexpected, main outcome of~\cite{ mer2008a}, the
theorem for $n = 2$ and
$\kappa = 5$ on p.~\pageref{bracket-5-2} about bracket invariants does {\em
not} capture all Demailly-Semple (bi-)invariants. This was striking,
because brackets were sufficient to capture all invariants in all
previously known studies\footnote{\, On observes that ${\sf UE}_3^3$
is {\em not} obtained by bracketing bi-invariants in ${\sf UE}_2^3$
(think of $D^6$), but nevertheless ${\sf UE}_3^3$ is the
unipotent-invariant subalgebra of ${\sf E}_3^3$, and ${\sf E}_3^3$
itself is obtained by bracketing invariants from the preceding jet
level. }, namely for ${\sf E}_2^n$, for ${\sf E}_3^2$, for ${\sf
E}_3^3$ and for ${\sf E}_4^2$.

Aside from the 11 bi-invariants $f_1'$, $\Lambda^3$, $\Lambda_1^5$,
$\Lambda_{ 1, 1}^7$, $M^8$, $\Lambda_{ 1, 1, 1 }^9$, $M_1^{ 10}$, $N^{
12}$, $K_{ 1, 1}^{ 12}$, $H_1^{ 14}$ and $F_{ 1, 1}^{ 16}$, there are
yet the following 6 bi-invariants $X^{ 18}$, $X^{ 19}$, $X^{ 21}$,
$X^{ 23}$, $X^{ 25}$ and $Y^{ 27}$ that are defined by dividing by
$f_1'$ some appropriate quadratic combinations between already known
bi-invariants. We provide here the complete explicit expressions. It
is shown in~\cite{ mer2008a} that the 16 first bi-invariants are
mutually independent and it would be easy, by using the same method,
to verify that when one adds the last, 17-th bi-invariant $Y^{ 27}$, one
still gets a list of 17 mutually independent bi-invariants.

Importantly, we emphasize that by no means any of these 6 further
bi-invariants can come from inspecting the bracket invariants, by
dividing them either by $f_1'$, or by $\Lambda^3$ or by anything based
in brackets, because in~\cite{ mer2008a}, {\em all the possible
bracket invariants} were computed thoroughly, were simplified and were
analyzed at a piece. {\em The existence of $X^{ 18}$, $X^{ 19}$, $X^{
21}$, $X^{ 23}$, $X^{ 25}$, $Y^{ 27}$ really shows that bracketing
does not generate the algebra of bi-invariants ${\sf UE}_5^2$}. A
similar phenomenon will appear to take place in dimension $n = 3$ for
jet order $\kappa = 4$.

Before reading the formulas, we would like to mention that the
invariant $X^{ 21}$ of ${\sf UE}_5^2$ below is not the same as the
invariant $X^{ 21}$ of ${\sf UE}_3^4$ appearing in \S11. Our
manuscript sheets used the same notation, and we hope this should not
cause any confusion.
\[
\scriptsize
\aligned
X^{18}
:=
&\
\frac{-5\,\Lambda_{1,1,1}^9\,M_1^{10}
+
56\,\Lambda_{1,1}^7\,K_{1,1}^{12}}{
f_1'}
\\
=
&\
f_1'f_1'f_1'\Big(
-
18816\,\Delta^{',''''}\,\big[\Delta^{'','''}\big]^2 - 25088\,
\big[\Delta^{'','''}\big]^3 - 15\,\big[
\Delta^{','''''}\big]^2\,\Delta^{',''} - 150\,
\Delta^{','''''}\,\Delta^{'',''''}\,\Delta^{',''} 
\\
&
+ 315\,\Delta^{','''''}\,\Delta^{',''''}\,\Delta^{','''} + 960\,
\Delta^{','''''}\,\Delta^{'','''}\,\Delta^{','''} - 375\,\big[\Delta^{'',''''}\big]^2
\,\Delta^{',''} + 1575\,\Delta^{'',''''}\,\Delta^{',''''}\,\Delta^{','''}
\\
&
+ 4800\,\Delta^{'',''''}\,\Delta^{'','''}\,\Delta^{','''} - 392\,
\big[\Delta^{',''''}\big]^3 - 4704\,\big[\Delta^{',''''}\big]^2\,\Delta^{'','''}
\Big)
-
f_1'f_1'f_1''\,
\Big(
- 2475\,\Delta^{'',''''}\,\Delta^{',''''}\,\Delta^{',''}
\\
&
- 9900\,\Delta^{'',''''}\,\Delta^{'','''}\,\Delta^{',''} - 2850\,
\Delta^{','''''}\,\big[\Delta^{','''}\big]^2 + 51330\,\Delta^{',''''}\,
\Delta^{'','''}\,\Delta^{','''}
\\
&
+ 92760\,\big[\Delta^{'','''}\big]^2\,\Delta^{','''} - 14250\,\Delta^{'',''''}\,
\big[\Delta^{','''}\big]^2 + 7035\,\big[\Delta^{',''''}\big]^2\,\Delta^{','''} - 
495\,\Delta^{','''''}\,\Delta^{',''''}\,\Delta^{',''}
\endaligned
\]
\[
\scriptsize
\aligned
&
- 1980\,\Delta^{','''''}\,\Delta^{'','''}\,\Delta^{',''}
\Big) 
- 
f_1'f_1'f_1'''\,
\Big(
- 11100\,\Delta^{'','''}\,\big[\Delta^{','''}\big]^2 - 3150\,
\Delta^{',''''}\,\big[\Delta^{','''}\big]^2
\Big)
\\
&
+
f_1'f_1''f_1''\,
\Big(
- 109440\,\big[\Delta^{'','''}\big]^2\,\Delta^{',''} - 19050\,\Delta^{'','''}\,
\big[\Delta^{','''}\big]^2 - 32325\,\Delta^{',''''}\,\big[\Delta^{','''}\big]^2 
\\
&
+ 11025\,\Delta^{','''''}\,\Delta^{','''}\,\Delta^{',''} + 55125
\,\Delta^{'',''''}\,\Delta^{','''}\,\Delta^{',''} - 6840\,\big[\Delta^{',''''}\big]^2
\,\Delta^{',''}
\\
&
- 54720\,\Delta^{',''''}\,\Delta^{'','''}\,\Delta^{',''}
\Big) 
-
f_1'f_1''f_1'''\,\Big(
+ 30000\,
\big[\Delta^{','''}\big]^3
\Big) 
- 
f_1''f_1''f_1''\,
\Big(11025\,\Delta^{','''''}\,
\big[\Delta^{',''}\big]^2
\\
&
- 55125\,\Delta^{'',''''}\,
\big[\Delta^{',''}\big]^2 + 55125\,\Delta^{',''''}\,
\Delta^{','''}\,\Delta^{',''} + 110250\,\Delta^{'','''}\,
\Delta^{','''}\,\Delta^{',''} 
\\
&
- 49000\,\big[\Delta^{','''}\big]^3
\Big).
\endaligned
\]
\[
\scriptsize
\aligned
X^{19}
:=
&\
\frac{-5\,M_1^{10}\,M_1^{10}+64\,M^8\,K_{1,1}^{12}}{f_1'}
\\
=
&\,
f_1'\, 
\Big(
1170\,\Delta^{','''''}\,\Delta^{',''''}\,\Delta^{','''}\,
\Delta^{',''} - 45\,\big[\Delta^{','''''}\big]^2\,\big[
\Delta^{',''}\big]^2 - 450\,
\Delta^{','''''}\,\Delta^{'',''''}\,\big[\Delta^{',''}\big]^2 
\\
&
+ 74220\,\big[\Delta^{'','''}\big]^2\,\big[\Delta^{','''}\big]^2 + 3780\,
\Delta^{','''''}\,\Delta^{'','''}\,\Delta^{','''}\,\Delta^{',''} - 1600\,
\Delta^{','''''}\,\big[\Delta^{','''}\big]^3 
\\
&
- 1125\,\big[\Delta^{'',''''}\big]^2\,\big[\Delta^{',''}\big]^2 + 5850\,
\Delta^{'',''''}\,\Delta^{',''''}\,\Delta^{','''}\,\Delta^{',''} + 18900\,
\Delta^{'',''''}\,\Delta^{'','''}\,\Delta^{','''}\,\Delta^{',''} 
\\
&
- 8000\,\Delta^{'',''''}\,\big[\Delta^{','''}\big]^3 - 1344\,
\big[\Delta^{',''''}\big]^3\,
\Delta^{',''} - 16128\,\big[\Delta^{',''''}\big]^2\,\Delta^{'','''}\,
\Delta^{',''} + 1995\,\big[\Delta^{',''''}\big]^2\,\big[\Delta^{','''}\big]^2 
\endaligned
\]
\[
\scriptsize
\aligned
&
- 64512\,\Delta^{',''''}\,\big[\Delta^{'','''}\big]^2\,\Delta^{',''} + 
27660\,\Delta^{',''''}\,\Delta^{'','''}\,\big[\Delta^{','''}\big]^2 - 86016\,
\big[\Delta^{'','''}\big]^3\,\Delta^{',''}
\Big) 
\\
&
+
f_1''\,
\Big(
- 74400\,\Delta^{'','''}\,\big[\Delta^{','''}\big]^3 - 10800\,\Delta^{'',''''}\,
\Delta^{',''''}\,\big[\Delta^{',''}\big]^2 - 2160\,\Delta^{','''''}\,\Delta^{',''''}
\,\big[\Delta^{',''}\big]^2 
\\
&
- 8640\,\Delta^{','''''}\,\Delta^{'','''}\,\big[\Delta^{',''}\big]^2 + 
3600\,\Delta^{','''''}\,\big[\Delta^{','''}\big]^2\,\Delta^{',''} + 64800\,
\Delta^{',''''}\,\Delta^{'','''}\,\Delta^{','''}\,\Delta^{',''} 
\\
&
- 43200\,\Delta^{'',''''}\,\Delta^{'','''}\,\big[\Delta^{',''}\big]^2 + 
18000\,\Delta^{'',''''}\,\big[\Delta^{','''}\big]^2\,\Delta^{',''} + 10800\,
\big[\Delta^{',''''}\big]^2\,\Delta^{','''}\,\Delta^{',''} 
\\
&
- 27600\,\Delta^{',''''}\,\big[\Delta^{','''}\big]^3 + 86400\,
\big[\Delta^{'','''}\big]^2\,\Delta^{','''}\,\Delta^{',''}
\Big)
+ 
f_1'''\,\Big(
16000\,\big[\Delta^{','''}\big]^4\Big).
\endaligned
\]
\[
\scriptsize
\aligned
X^{21}
:=
&\
\frac{-5\,M_1^{10}\,N^{12}+8\,M^8\,H_1^{14}}{f_1'}
\\
=
&\, 
- 135\,\big[\Delta^{','''''}\big]^2\,\big[\Delta^{',''}\big]^3 - 1350
\,\Delta^{','''''}\,\Delta^{'',''''}\,\big[\Delta^{',''}\big]^3 + 1350\,
\Delta^{','''''}\,\Delta^{',''''}\,\Delta^{','''}\,\big[\Delta^{',''}\big]^2 
\\
&
+ 2700\,\Delta^{','''''}\,\Delta^{'','''}\,\Delta^{','''}\,
\big[\Delta^{',''}\big]^2 - 1200\,\Delta^{','''''}\,\big[\Delta^{','''}\big]^3\,
\Delta^{',''} - 3375\,\big[\Delta^{'',''''}\big]^2\,\big[\Delta^{',''}\big]^3 
\\
&
+ 6750\,\Delta^{'',''''}\,\Delta^{',''''}\,\Delta^{','''}\,
\big[\Delta^{',''}\big]^2 + 13500\,\Delta^{'',''''}\,\Delta^{'','''}\,
\Delta^{','''}\,\big[\Delta^{',''}\big]^2 - 6000\,\Delta^{'',''''}\,
\big[\Delta^{','''}\big]^3\,
\Delta^{',''} 
\endaligned
\]
\[
\scriptsize
\aligned
&
- 576\,\big[\Delta^{',''''}\big]^3\,\big[\Delta^{',''}\big]^2 - 6912\,
\big[\Delta^{',''''}\big]^2\,\Delta^{'','''}\,\big[\Delta^{',''}\big]^2 - 495\,
\big[\Delta^{',''''}\big]^2\,\big[\Delta^{','''}\big]^2\,\Delta^{',''} 
\\
&
- 27648\,\Delta^{',''''}\,\big[\Delta^{'','''}\big]^2\,\big[\Delta^{',''}\big]^2
+ 9540\,\Delta^{',''''}\,\Delta^{'','''}\,\big[\Delta^{','''}\big]^2\,
\Delta^{',''} + 1200\,\Delta^{',''''}\,\big[\Delta^{','''}\big]^4 
\\
&
- 36864\,\big[\Delta^{'','''}\big]^3\,\big[\Delta^{',''}\big]^2 + 32580\,
\big[\Delta^{'','''}\big]^2\,\big[\Delta^{','''}\big]^2\,\Delta^{',''} - 7200\,
\Delta^{'','''}\,\big[\Delta^{','''}\big]^4.
\endaligned
\]
\[
\scriptsize
\aligned
X^{23}
:=
&\
\frac{-7\,N^{12}\,K_{1,1}^{12}+M^8\,F_{1,1}^{16}}{f_1'}
\\
=
&\, 
f_1'\,
\Big(432\,\Delta^{','''''}\,\big[\Delta^{',''''}\big]^2\,\big[\Delta^{',''}\,
\big]^2 + 3456\,\Delta^{','''''}\,\Delta^{',''''}\,\Delta^{'','''}\,
\big[\Delta^{',''}\big]^2 + 1710\,\Delta^{','''''}\,\Delta^{',''''}\,
\big[\Delta^{','''}\big]^2\,
\Delta^{',''} 
\\
&
- 3150\,\Delta^{','''''}\,\Delta^{'',''''}\,\Delta^{','''}\,
\big[\Delta^{',''}\big]^2 + 540\,\Delta^{','''''}\,\Delta^{'','''}\,\big[\Delta^{','''}
\big]^2\,\Delta^{',''} - 1600\,\Delta^{','''''}\,\big[\Delta^{','''}\big]^4 
\\
&
- 7875\,\big[\Delta^{'',''''}\big]^2\,\Delta^{','''}\,\big[\Delta^{',''}\big]^2
+ 6912\,\Delta^{','''''}\,\big[\Delta^{'','''}\big]^2\,\big[
\Delta^{',''}\big]^2 - 8000
\,\Delta^{'',''''}\,\big[\Delta^{','''}\big]^4 
\\
&
- 2352\,\big[\Delta^{',''''}\big]^3\,\Delta^{','''}\,\Delta^{',''} - 
23904\,\big[\Delta^{',''''}\big]^2\,\Delta^{'','''}\,\Delta^{','''}\,\Delta^{',''}
+ 2205\,\big[\Delta^{',''''}\big]^2\,\big[\Delta^{','''}\big]^3 
\endaligned
\]
\[
\scriptsize
\aligned
&
- 78336\,\Delta^{',''''}\,\big[\Delta^{'','''}\big]^2\,\Delta^{','''}\,
\Delta^{',''} + 34740\,\Delta^{',''''}\,\Delta^{'','''}\,\big[\Delta^{','''}\big]^3
- 81408\,\big[\Delta^{'','''}\big]^3\,\Delta^{','''}\,\Delta^{',''} 
\\
&
+ 72180\,\big[\Delta^{'','''}\big]^2\,\big[\Delta^{','''}\big]^3 + 2160\,
\Delta^{'',''''}\,\big[\Delta^{',''''}\big]^2\,\big[\Delta^{',''}\big]^2 + 17280\,
\Delta^{'',''''}\,\Delta^{',''''}\,\Delta^{'','''}\,\big[\Delta^{',''}\big]^2 
\\
&
+ 8550\,\Delta^{'',''''}\,\Delta^{',''''}\,\big[\Delta^{','''}\big]^2\,
\Delta^{',''} + 34560\,\Delta^{'',''''}\,\big[\Delta^{'','''}\big]^2\,
\big[\Delta^{',''}\big]^2 + 2700\,\Delta^{'',''''}\,\Delta^{'','''}\,
\big[\Delta^{','''}\big]^2\,
\Delta^{',''} 
\\
&
- 315\,\big[\Delta^{','''''}\big]^2\,\Delta^{','''}\,\big[\Delta^{',''}\big]^2
\Big)
+
f_1''\,
\Big(
23625\,\big[\Delta^{'',''''}\big]^2\,\big[\Delta^{',''}\big]^3
 - 47250\,\Delta^{'',''''}\,\Delta^{',''''}\,\Delta^{','''}\,\big[\Delta^{',''}
\big]^2 
\\
&
- 94500\,\Delta^{'',''''}\,\Delta^{'','''}\,\Delta^{','''}\,
\big[\Delta^{',''}\big]^2 + 42000\,\Delta^{'',''''}\,\big[\Delta^{','''}\big]^3\,
\Delta^{',''} + 576\,\big[\Delta^{',''''}\big]^3\,\big[\Delta^{',''}\big]^2 
\endaligned
\]
\[
\scriptsize
\aligned
&
+ 6912\,\big[\Delta^{',''''}\big]^2\,\Delta^{'','''}\,\big[\Delta^{',''}\big]^2
+ 20745\,\big[\Delta^{',''''}\big]^2\,\big[\Delta^{','''}\big]^2\,\Delta^{',''} + 
27648\,\Delta^{',''''}\,\big[\Delta^{'','''}\big]^2\,\big[\Delta^{',''}\big]^2 
\\
&
+ 945\,\big[\Delta^{','''''}\big]^2\,\big[\Delta^{',''}\big]^3 + 9450\,
\Delta^{','''''}\,\Delta^{'',''''}\,\big[\Delta^{',''}\big]^3 - 9450\,\Delta^{','''''}
\,\Delta^{',''''}\,\Delta^{','''}\,\big[\Delta^{',''}\big]^2 
\\
&
- 18900\,\Delta^{','''''}\,\Delta^{'','''}\,\Delta^{','''}\,
\big[\Delta^{',''}\big]^2 + 8400\,\Delta^{','''''}\,\big[\Delta^{','''}\big]^3\,
\Delta^{',''} + 71460\,\Delta^{',''''}\,\Delta^{'','''}\,\big[\Delta^{','''}\big]^2\,
\Delta^{',''} 
\\
&
- 37200\,\Delta^{',''''}\,\big[\Delta^{','''}\big]^4 + 36864\,\big[\Delta^{'','''}\,
\big]^3\,\big[\Delta^{',''}\big]^2 + 48420\,\big[\Delta^{'','''}\big]^2\,
\big[\Delta^{','''}\big]^2\,\Delta^{',''} 
\\
&
- 64800\,\Delta^{'','''}\,\big[\Delta^{','''}\big]^4
\Big)
+ 
f_1'''\,\Big(
16000\,\big[\Delta^{','''}\big]^5\Big).
\endaligned
\]
\[
\scriptsize
\aligned
X^{25}
:=
&\
\frac{-56\,K_{1,1}^{12}\,H_1^{14}+5\,M_1^{10}\,F_{1,1}^{16}}{f_1'}
\\
=
&\, 
f_1'\,f_1'\,
\Big(
-45\,\big[\Delta^{','''''}\big]^2\,\Delta^{',''''}\,
\big[\Delta^{',''}\big]^2 - 180\,\big[\Delta^{','''''}\big]^2\,\Delta^{'','''}\,
\big[\Delta^{',''}\big]^2
- 3600\,\big[\Delta^{','''''}\big]^2\,\big[\Delta^{','''}\big]^2\,\Delta^{',''} 
\\
&
 - 2800\,\Delta^{','''''}\,\Delta^{',''''}\,\big[\Delta^{','''}\big]^3 - 
83200\,\Delta^{','''''}\,\Delta^{'','''}\,\big[\Delta^{','''}\big]^3 - 1125\,
\big[\Delta^{'',''''}\big]^2\,\Delta^{',''''}\,\big[\Delta^{',''}\big]^2 
\\
&
- 4500\,\big[\Delta^{'',''''}\big]^2\,\Delta^{'','''}\,\big[\Delta^{',''}\big]^2
- 90000\,\big[\Delta^{'',''''}\big]^2\,\big[\Delta^{','''}\big]^2\,\Delta^{',''} - 
14000\,\Delta^{'',''''}\,\Delta^{',''''}\,\big[\Delta^{','''}\big]^3 
\\
&
- 416000\,\Delta^{'',''''}\,\Delta^{'','''}\,\big[\Delta^{','''}\big]^3 - 
150528\,\big[\Delta^{',''''}\big]^3\,\Delta^{'','''}\,\Delta^{',''} - 903168\,
\big[\Delta^{',''''}\big]^2\,\big[\Delta^{'','''}\big]^2\,\Delta^{',''} 
\\
&
+ 163800\,\big[\Delta^{',''''}\big]^2\,\Delta^{'','''}\,\big[\Delta^{','''}\big]^2 
-2408448\,\Delta^{',''''}\,\big[\Delta^{'','''}\big]^3\,\Delta^{',''} + 
1129500\,\Delta^{',''''}\,\big[\Delta^{'','''}\big]^2\,\big[\Delta^{','''}\big]^2 
\\
&
- 9408\,\big[\Delta^{',''''}\big]^4\,\Delta^{',''} + 3675\,\big[\Delta^{',''''}\,
\big]^3\,\big[\Delta^{','''}\big]^2 - 2408448\,\big[\Delta^{'','''}\big]^4\,
\Delta^{',''} + 2132400\,\big[\Delta^{'','''}\big]^3\,\big[\Delta^{','''}\big]^2 
\\
&
- 450\,\Delta^{','''''}\,\Delta^{'',''''}\,\Delta^{',''''}\,
\big[\Delta^{',''}\big]^2 - 1800\,\Delta^{','''''}\,\Delta^{'',''''}\,\Delta^{'','''}\,
\big[\Delta^{',''}\big]^2 
\endaligned
\]
\[
\scriptsize
\aligned
&
- 36000\,\Delta^{','''''}\,\Delta^{'',''''}\,\big[\Delta^{','''}\big]^2\,
\Delta^{',''} + 11970\,\Delta^{','''''}\,\big[\Delta^{',''''}\big]^2\,
\Delta^{','''}\,\Delta^{',''} 
\\
&
+ 187920\,\Delta^{','''''}\,\big[\Delta^{'','''}\big]^2\,\Delta^{','''}\,
\Delta^{',''} + 59850\,\Delta^{'',''''}\,\big[\Delta^{',''''}\big]^2\,
\Delta^{','''}\,\Delta^{',''} 
\\
&
+ 939600\,\Delta^{'',''''}\,\big[\Delta^{'','''}\big]^2\,\Delta^{','''}\,
\Delta^{',''} + 474300\,\Delta^{'',''''}\,\Delta^{',''''}\,\Delta^{'','''}\,
\Delta^{','''}\,\Delta^{',''} 
\\
&
+ 94860\,\Delta^{','''''}\,\Delta^{',''''}\,\Delta^{'','''}\,
\Delta^{','''}\,\Delta^{',''}
\Big)
+ 
f_1'f_1''\,
\Big(
- 2556600
\,\Delta^{',''''}\,\Delta^{'','''}\,\big[\Delta^{','''}\big]^3 
\endaligned
\]
\[
\scriptsize
\aligned
&
- 5014200\,\big[\Delta^{'','''}\big]^2\,\big[\Delta^{','''}\big]^3 - 187950\,
\big[\Delta^{',''''}\big]^2\,\big[\Delta^{','''}\big]^3 + 5621760\,\Delta^{',''''}\,
\big[\Delta^{'','''}\big]^2\,\Delta^{','''}\,\Delta^{',''} 
\\
&
+ 5652480\,\big[\Delta^{'','''}\big]^3\,\Delta^{','''}\,\Delta^{',''}
- 2764800\,\Delta^{'',''''}\,\big[\Delta^{'','''}\big]^2\,\big[\Delta^{',''}\big]^2 
\\
&
+ 99000\,\Delta^{'',''''}\,\Delta^{'','''}\,\big[\Delta^{','''}\big]^2\,
\Delta^{',''} + 500000\,\Delta^{'',''''}\,\big[\Delta^{','''}\big]^4 + 174720\,
\big[\Delta^{',''''}\big]^3\,\Delta^{','''}\,\Delta^{',''} 
\\
&
+ 1751040\,\big[\Delta^{',''''}\big]^2\,\Delta^{'','''}\,\Delta^{','''}\,
\Delta^{',''} - 276480\,\Delta^{','''''}\,\Delta^{',''''}\,\Delta^{'','''}\,
\big[\Delta^{',''}\big]^2 
\\
&
- 105300\,\Delta^{','''''}\,\Delta^{',''''}\,\big[\Delta^{','''}\big]^2\,
\Delta^{',''} - 552960\,\Delta^{','''''}\,\big[\Delta^{'','''}\big]^2\,
\big[\Delta^{',''}\big]^2 
\endaligned
\]
\[
\scriptsize
\aligned
&
+ 19800\,\Delta^{','''''}\,\Delta^{'','''}\,\big[\Delta^{','''}\big]^2\,
\Delta^{',''} + 100000\,\Delta^{','''''}\,\big[\Delta^{','''}\big]^4 + 551250\,
\big[\Delta^{'',''''}\big]^2\,\Delta^{','''}\,\big[\Delta^{',''}\big]^2 
\\
&
- 172800\,\Delta^{'',''''}\,\big[\Delta^{',''''}\big]^2\,\big[\Delta^{',''}\big]^2
- 1382400\,\Delta^{'',''''}\,\Delta^{',''''}\,\Delta^{'','''}\,
\big[\Delta^{',''}\big]^2 
\\
&
- 526500\,\Delta^{'',''''}\,\Delta^{',''''}\,\big[\Delta^{','''}\big]^2\,
\Delta^{',''} - 34560\,\Delta^{','''''}\,\big[\Delta^{',''''}\big]^2\,
\big[\Delta^{',''}\big]^2 + 22050\,\big[\Delta^{','''''}\big]^2\,\Delta^{','''}\,
\big[\Delta^{',''}\big]^2
\\
&
+ 220500\,\Delta^{','''''}\,\Delta^{'',''''}\,\Delta^{','''}\,
\big[\Delta^{',''}\big]^2
\Big)
\\
&
+ 
f_1'f_1'''\,
\Big(
28000\,\Delta^{',''''}\,\big[\Delta^{','''}\big]^4 + 472000\,
\Delta^{'','''}\,\big[\Delta^{','''}\big]^4
\Big)
\endaligned
\]
\[
\scriptsize
\aligned
&
+
f_1''f_1''\,
\Big(
330750\,\Delta^{','''''}\,\Delta^{',''''}\,\Delta^{','''}\,\big[\Delta^{',''}\big]^2
+ 661500\,\Delta^{','''''}\,\Delta^{'','''}\,\Delta^{','''}\,
\big[\Delta^{',''}\big]^2 
\\
&
- 294000\,\Delta^{','''''}\,\big[\Delta^{','''}\big]^3\,\Delta^{',''} - 
330750\,\Delta^{','''''}\,\Delta^{'',''''}\,\big[\Delta^{',''}\big]^3 
\\
&
+ 1653750\,\Delta^{'',''''}\,\Delta^{',''''}\,\Delta^{','''}\,
\big[\Delta^{',''}\big]^2 + 3307500\,\Delta^{'',''''}\,\Delta^{'','''}\,
\Delta^{','''}\,\big[\Delta^{',''}\big]^2 
\\
&
- 1470000\,\Delta^{'',''''}\,\big[\Delta^{','''}\big]^3\,\Delta^{',''}
- 2880\,\big[\Delta^{',''''}\big]^3\,\big[
\Delta^{',''}\big]^2 - 34560\,\big[\Delta^{',''''}
\big]^2\,\Delta^{'','''}\,\big[\Delta^{',''}\big]^2 
\endaligned
\]
\[
\scriptsize
\aligned
&
- 812475\,\big[\Delta^{',''''}\big]^2\,\big[\Delta^{','''}\big]^2\,\Delta^{',''}\,
- 138240\,\Delta^{',''''}\,\big[\Delta^{'','''}\big]^2\,\big[\Delta^{',''}\big]^2 - 
33075\,\big[\Delta^{','''''}\big]^2\,\big[\Delta^{',''}\big]^3 
\\
&
+ 1446000\,\Delta^{',''''}\,\big[\Delta^{','''}\big]^4 - 184320\,
\big[\Delta^{'','''}\big]^3\,\Delta^{',''}\big]^2 
- 3077100\,\big[\Delta^{'','''}\big]^2\,
\big[\Delta^{','''}\big]^2\,\Delta^{',''} 
\\
&
+ 2844000\,\Delta^{'','''}\,\big[\Delta^{','''}\big]^4 - 826875\,
\big[\Delta^{'',''''}\big]^2\,\big[\Delta^{',''}\big]^3 - 3192300\,\Delta^{',''''}\,
\Delta^{'','''}\,\big[\Delta^{','''}\big]^2\,\Delta^{',''}
\Big)
\\
&
+
f_1''f_1'''\,
\Big( - 640000\,\big[\Delta^{','''}\big]^5
\Big).
\endaligned
\]
\[
\scriptsize
\aligned
Y^{27}
:=
&\
\frac{-56\,K_{1,1}^{12}F_{1,1}^{16}+M_1^{10}X^{18}}{f_1'}
\\
=
&\
f_1'f_1'f_1'\Big(
572820\,\Delta^{','''''}\,\big[\Delta^{'','''}\big]^2\,\big[\Delta^{','''}\big]^2
-
5343744\,\Delta^{',''''}\,\big[\Delta^{'','''}\big]^3\,\Delta^{','''}
-
752640\,\Delta^{'',''''}\,\big[\Delta^{'','''}\big]^3\,\Delta^{',''}
-
\\
&\
-
286944\,\big[\Delta^{',''''}\big]^3\,\Delta^{'','''}\,\Delta^{','''}
+
2864100\,\Delta^{'',''''}\,\big[\Delta^{'','''}\big]^2\,\big[\Delta^{','''}\big]^2
-
1862784\,\big[\Delta^{',''''}\big]^2\,\big[\Delta^{'','''}\big]^2\,\Delta^{','''}
+
\\
&\
+
27195\,\Delta^{','''''}\,\big[\Delta^{',''''}\big]^2\,\big[\Delta^{','''}\big]^2
-
112000\,\Delta^{','''''}\,\Delta^{'',''''}\,\big[\Delta^{','''}\big]^3
-
150528\,\Delta^{','''''}\,\big[\Delta^{'','''}\big]^3\,\Delta^{',''}
-
\\
&
-2352\,\Delta^{','''''}\,\big[\Delta^{',''''}\big]^3\,\Delta^{',''}
+
135975\,\Delta^{'',''''}\,\big[\Delta^{',''''}\big]^2\,\big[\Delta^{','''}\big]^2
-
11760\,\Delta^{'',''''}\,\big[\Delta^{',''''}\big]^3\,\Delta^{',''}
-
\\
&
-3375\,\Delta^{','''''}\,\big[\Delta^{'',''''}\big]^2\,\big[\Delta^{',''}\big]^2
-
675\,\big[\Delta^{','''''}\big]^2\,\Delta^{'',''''}\,\big[\Delta^{',''}\big]^2
-
45\,\big[\Delta^{','''''}\big]^3\,\big[\Delta^{',''}\big]^2
-
\endaligned
\]
\[
\scriptsize
\aligned
&
-11200\,\big[\Delta^{','''''}\big]^2\,\big[\Delta^{','''}\big]^3
-
5625\,\big[\Delta^{'',''''}\big]^3\,\big[\Delta^{',''}\big]^2
-
280000\,\big[\Delta^{'',''''}\big]^2\,\big[\Delta^{','''}\big]^3
-
\\
&
-
16464\,\big[\Delta^{',''''}\big]^4\,\Delta^{','''}
-
5720064\,\big[\Delta^{'','''}\big]^4\,\Delta^{','''}
+
1890\,\big[\Delta^{','''''}\big]^2\,\Delta^{',''''}\,\Delta^{','''}\,\Delta^{',''}
+
\\
&
+6210\,\big[\Delta^{','''''}\big]^2\,\Delta^{'','''}\,\Delta^{','''}\,\Delta^{',''}
-
28224\,\Delta^{','''''}\,\big[\Delta^{',''''}\big]^2\,\Delta^{'','''}\,\Delta^{',''}
-
112896\,\Delta^{','''''}\,\Delta^{',''''}\,\big[\Delta^{'','''}\big]^2\,\Delta^{',''}
+
\\
&
+
255360\,\Delta^{','''''}\,\Delta^{',''''}\,\Delta^{'','''}\,\big[\Delta^{','''}\big]^2
+
47250\,\big[\Delta^{'',''''}\big]^2\,\Delta^{',''''}\,\Delta^{','''}\,\Delta^{',''}
+
155250\,\big[\Delta^{'',''''}\big]^2\,\Delta^{'','''}\,\Delta^{','''}\,\Delta^{',''}
-
\\
&
-141120\,\Delta^{'',''''}\,\big[\Delta^{',''''}\big]^2\,\Delta^{'','''}\,\Delta^{',''}
-
564480\,\Delta^{'',''''}\,\Delta^{',''''}\,\big[\Delta^{'','''}\big]^2\,\Delta^{',''}
+
1276800\,\Delta^{'',''''}\,\Delta^{',''''}\,\Delta^{'','''}\,\big[\Delta^{','''}\big]^2
+
\\
&
+18900\,\Delta^{','''''}\,\Delta^{'',''''}\,\Delta^{',''''}\,\Delta^{','''}\,\Delta^{',''}
+
62100\,\Delta^{','''''}
\Big)
+
\endaligned
\]
\[
\scriptsize
\aligned
+
&\
\Big(\Big(
-36450\,\Delta^{','''''}\,\Delta^{'',''''}\,\Delta^{',''''}\,\big[\Delta^{',''}\big]^2
-
145800\,\Delta^{','''''}\,\Delta^{'',''''}\,\Delta^{'','''}\,\big[\Delta^{',''}\big]^2
+
832500\,\Delta^{','''''}\,\big[\Delta^{','''}\big]^2\,\Delta^{',''}
-
\\
&
-149310\,\Delta^{','''''}\,\big[\Delta^{',''''}\big]^2\,
-
2680560\,\Delta^{','''''}\,\big[\Delta^{'','''}\big]^2\,\Delta^{','''}\,\Delta^{',''}
-
746550\,\Delta^{'',''''}\,\big[\Delta^{',''''}\big]^2\,\Delta^{','''}\,\Delta^{',''}
-
\\
&
-13402800\,\Delta^{'',''''}\,\big[\Delta^{'','''}\big]^2\,\Delta^{','''}\,\Delta^{',''}
-
3645\,\big[\Delta^{','''''}\big]^2\,\Delta^{',''''}\,\big[\Delta^{',''}\big]^2
-
14580\,\big[\Delta^{','''''}\big]^2\,\Delta^{'','''}\,\big[\Delta^{',''}\big]^2
+
\\
&
+83250\,\big[\Delta^{','''''}\big]^2\,\big[\Delta^{','''}\big]^2\,\Delta^{',''}
-
245700\,\Delta^{','''''}\,\Delta^{',''''}\,\big[\Delta^{','''}\big]^3
+
682200\,\Delta^{','''''}\,\Delta^{'','''}\,\big[\Delta^{','''}\big]^3
-
\endaligned
\]
\[
\scriptsize
\aligned
&
-91125\,\big[\Delta^{'',''''}\big]^2\,\Delta^{',''''}\,\big[\Delta^{',''}\big]^2
-
364500\,\big[\Delta^{'',''''}\big]^2\,\Delta^{'','''}\,\big[\Delta^{',''}\big]^2
+
2081250\,\big[\Delta^{'',''''}\big]^2\,\big[\Delta^{','''}\big]^2\,\Delta^{',''}
-
\\
&
-1228500\,\Delta^{'',''''}\,\Delta^{',''''}\,\big[\Delta^{','''}\big]^3
+
3411000\,\Delta^{'',''''}\,\Delta^{'','''}\,\big[\Delta^{','''}\big]^3
+
1354752\,\big[\Delta^{',''''}\big]^3\,\Delta^{'','''}\,\Delta^{',''}
+
\\
&
+8128512\,\big[\Delta^{',''''}\big]^2\,\big[\Delta^{'','''}\big]^2\,\Delta^{',''}
+
1796760\,\big[\Delta^{',''''}\big]^2\,\Delta^{'','''}\,\big[\Delta^{','''}\big]^2
+
21676032\,\Delta^{',''''}\,\big[\Delta^{'','''}\big]^3\,\Delta^{',''}
+
\\
&
+850140\,\Delta^{',''''}\,\big[\Delta^{'','''}\big]^2\,\big[\Delta^{','''}\big]^2
+
84672\,\big[\Delta^{',''''}\big]^4\,\Delta^{',''}
+
274155\,\big[\Delta^{',''''}\big]^3\,\big[\Delta^{','''}\big]^2
+
\\
&
+21676032\,\big[\Delta^{'','''}\big]^4\,\Delta^{',''}
-
7801680\,\big[\Delta^{'','''}\big]^3\,\big[\Delta^{','''}\big]^2
-
6336900\,\Delta^{'',''''}\,\Delta^{',''''}\,\Delta^{'','''}\,\Delta^{','''}\,\Delta^{',''}
-
\\
&
-1267380\,\Delta^{','''''}\,\Delta^{',''''}\,\Delta^{'','''}\,\Delta^{','''}\,\Delta^{',''}
\Big)\,f_1''
+
\endaligned
\]
\[
\scriptsize
\aligned
&
+
\Big(1120000\,\Delta^{'',''''}\,\big[\Delta^{','''}\big]^4
-
5062200\,\big[\Delta^{'','''}\big]^2\,\big[\Delta^{','''}\big]^3
+
224000\,\Delta^{','''''}\,\big[\Delta^{','''}\big]^4
-
\\
&
-271950\,\big[\Delta^{',''''}\big]^2\,\big[\Delta^{','''}\big]^3
-
2364600\,\Delta^{',''''}\,\Delta^{'','''}\,\big[\Delta^{','''}\big]^3
\Big)\,f_1'''
\Big)\,f_1'f_1'
+
\\
&
+
\Big(
\Big(
34044300\,\big[\Delta^{'','''}\big]^2\,\big[\Delta^{','''}\big]^3
+
108675\,\big[\Delta^{',''''}\big]^2\,\big[\Delta^{','''}\big]^3
-
231525\,\big[\Delta^{','''''}\big]^2\,\Delta^{','''}\,\big[\Delta^{',''}\big]^2
+
\\
&
+278640\,\Delta^{','''''}\,\big[\Delta^{',''''}\big]^2\,\big[\Delta^{',''}\big]^2
+
4458240\,\Delta^{','''''}\,\big[\Delta^{'','''}\big]^2\,\big[\Delta^{',''}\big]^2
-
5788125\,\big[\Delta^{'',''''}\big]^2\,\Delta^{','''}\,\big[\Delta^{',''}\big]^2
+
\endaligned
\]
\[
\scriptsize
\aligned
&
+
1393200\,\Delta^{'',''''}\,\big[\Delta^{',''''}\big]^2\,\big[\Delta^{',''}\big]^2
+
22291200\,\Delta^{'',''''}\,\big[\Delta^{'','''}\big]^2\,\big[\Delta^{',''}\big]^2
-
1396080\,\big[\Delta^{',''''}\big]^3\,\Delta^{','''}\,\Delta^{',''}
+
\\
&
+14733900\,\Delta^{',''''}\,\Delta^{'','''}\,\big[\Delta^{','''}\big]^3
-
44766720\,\big[\Delta^{'','''}\big]^3\,\Delta^{','''}\,\Delta^{',''}
-
1284000\,\Delta^{','''''}\,\big[\Delta^{','''}\big]^4
-
\\
&
-6420000\,\Delta^{'',''''}\,\big[\Delta^{','''}\big]^4
-
2315250\,\Delta^{','''''}\,\Delta^{'',''''}\,\Delta^{','''}\,\big[\Delta^{',''}\big]^2
+
2229120\,\Delta^{','''''}\,\Delta^{',''''}\,\Delta^{'','''}\,\big[\Delta^{',''}\big]^2
+
\endaligned
\]
\[
\scriptsize
\aligned
&
+1386450\,\Delta^{','''''}\,\Delta^{',''''}\,\big[\Delta^{','''}\big]^2\,\Delta^{',''}
+
915300\,\Delta^{','''''}\,\Delta^{'','''}\,\big[\Delta^{','''}\big]^2\,\Delta^{',''}
+
11145600\,\Delta^{'',''''}\,\Delta^{',''''}\,\Delta^{'','''}\,\big[\Delta^{',''}\big]^2
+
\\
&
+6932250\,\Delta^{'',''''}\,\Delta^{',''''}\,\big[\Delta^{','''}\big]^2\,\Delta^{',''}
+
4576500\,\Delta^{'',''''}\,\Delta^{'','''}\,\big[\Delta^{','''}\big]^2\,\Delta^{',''}
-
13966560\,\big[\Delta^{',''''}\big]^2\,\Delta^{'','''}\,\Delta^{','''}\,\Delta^{',''}
-
\\
&
-44720640\,\Delta^{',''''}\,\big[\Delta^{'','''}\big]^2\,\Delta^{','''}\,\Delta^{',''}
\Big)\,f_1''f_1''
+
\Big(
2268000\,\Delta^{',''''}\,\big[\Delta^{','''}\big]^4
+
792000\,\Delta^{'','''}\,\big[\Delta^{','''}\big]^4
\Big)\,f_1''f_1'''
-
\\
&
-1120000\,\big[\Delta^{','''}\big]^5\,f_1'''f_1'''
\Big)\,f_1'
+
\endaligned
\]
\[
\scriptsize
\aligned
&
+
\Big(
-4630500\,\Delta^{','''''}\,\Delta^{'','''}\,\Delta^{','''}\,\big[\Delta^{',''}\big]^2
+
2058000\,\Delta^{','''''}\,\big[\Delta^{','''}\big]^3\,\Delta^{',''}
-
11576250\,\Delta^{'',''''}\,\Delta^{',''''}\,\Delta^{','''}\,\big[\Delta^{',''}\big]^2
-
\\
&
-23152500\,\Delta^{'',''''}\,\Delta^{'','''}\,\Delta^{','''}\,\big[\Delta^{',''}\big]^2
+
10290000\,\Delta^{'',''''}\,\big[\Delta^{','''}\big]^3\,\Delta^{',''}
+
2880\,\big[\Delta^{',''''}\big]^3\,\big[\Delta^{',''}\big]^2
+
\\
&
+34560\,\big[\Delta^{',''''}\big]^2\,\Delta^{'','''}\,\big[\Delta^{',''}\big]^2
+
5773725\,\big[\Delta^{',''''}\big]^2\,\big[\Delta^{','''}\big]^2\,\Delta^{',''}
+
138240\,\Delta^{',''''}\,\big[\Delta^{'','''}\big]^2\,\big[\Delta^{',''}\big]^2
+
\\
&
+231525\,\big[\Delta^{','''''}\big]^2\,\big[\Delta^{',''}\big]^3
+
2315250\,\Delta^{','''''}\,\Delta^{'',''''}\,\big[\Delta^{',''}\big]^3
-
2315250\,\Delta^{','''''}\,\Delta^{',''''}\,\Delta^{','''}\,\big[\Delta^{',''}\big]^2
+
\\
&
+22922100\,\big[\Delta^{'','''}\big]^2\,\big[\Delta^{','''}\big]^2\,\Delta^{',''}
-
20484000\,\Delta^{'','''}\,\big[\Delta^{','''}\big]^4
+
5788125\,\big[\Delta^{'',''''}\big]^2\,\big[\Delta^{',''}\big]^3
-
\\
&
-10266000\,\Delta^{',''''}\,\big[\Delta^{','''}\big]^4
+
184320\,\big[\Delta^{'','''}\big]^3\,\big[\Delta^{',''}\big]^2
+
23037300\,\Delta^{',''''}\,\Delta^{'','''}\,\big[\Delta^{','''}\big]^2\,\Delta^{',''}
\Big)\,f_1''f_1''f_1''
+
\\
&
+
4560000\,\big[\Delta^{','''}\big]^5\,f_1''f_1''f_1'''.
\endaligned
\]
It will be a theorem, to be established in \S10 below, that the 17
mutually independent bi-invariants $f_1'$, $\Lambda^3$, $\Lambda_1^5$,
$\Lambda_{ 1, 1}^7$, $M^8$, $\Lambda_{ 1, 1, 1}^9$, $M_1^{ 10}$, $N^{
12}$, $K_{ 1, 1}^{ 12}$, $H_1^{ 14}$, $F_{ 1, 1}^{ 16}$, $X^{ 18}$,
$X^{ 19}$, $X^{ 21}$, $X^{ 23}$, $X^{ 25}$ and $Y^{ 27}$ generate the
algebra ${\sf UE}_5^2$.

\section*{\S9.~Principle of the general algorithm}
\label{Section-9}

\subsection*{ Initializing the algorithm}
We now explain a general algorithm which generates all bi-invariants,
which stops after a finite number of steps if and only if the algebra
of bi-invariants is finitely generated and which, in such a
circumstance, yields a complete generating family of mutually
independent bi-invariants together with a complete generating family
of syzygies between these bi-invariants. The same algorithm would
work equally well for Demailly-Semple invariants, but as we already
observed, in the desired applications, the complexity and the
cardinality of generators and of syzygies being much higher, only the
exploration of bi-invariants seems accessible.

Fix the dimension $n$ and the jet order $\kappa$, both arbitrary.
Start from the representation of an arbitrary bi-invariant of weight
$m$ gained previously thanks to the proposition on
p.~\pageref{bi-invariant-n-k}:
\[
{\sf P}
=
{\sf P}\big(j^\kappa f\big)
=
\sum_{-\frac{\kappa-1}{\kappa}m\leqslant a\leqslant m}\,
(f_1')^a\,
{\sf P}\big(L^{l_1},\dots,L^{l_{k_1}}\big),
\]
where the $L^{ l_i}$, $i = 1, \dots, k_1$, have weight $l_i$ and come
from the $\Lambda$-minors written there, after a division by an
appropriate maximal factoring power of $f_1'$, {\em cf.} the two
special cases analyzed after the general proposition. Call $f_1', L^{
l_1}, \dots, L^{ l_{ k_1}}$ the {\sl initial bi-invariants}.

\subsection*{ First loop of the algorithm}
The first step of the algorithm consists in computing a reduced
Gröbner basis (for a certain monomial order) of the ideal of relations
of the restrictions to $\{ f_1' = 0 \}$ of these initial
bi-invariants:
\[
\text{\sf Ideal-Rel}
\Big(
L^{l_1}\big\vert_0,\,\dots,\,
L^{l_{k_1}}\big\vert_0
\Big).
\]
In some favorables circumstances, this task may be done by symbolic
Gröbner bases packages, although it is well known that due to
exponentiality of time computation and to expression swelling, Gröbner
bases often appear to be
frustratingly unusable. Write as follows the so obtained
gröbnerized syzygies:
\[
0
\equiv
{\sf S}_i
\Big(
L^{l_1}(j^\kappa f)\big\vert_0,\,\dots,\,
L^{l_{k_1}}(j^\kappa f)\big\vert_0
\Big)
\ \ \ \ \ \ \ \ \ \ \ \ \
{\scriptstyle{(i\,=\,1\,\cdots\,N_1)}}.
\]

At first, we claim that, without loss of generality, one may assume
that each syzygy polynomial ${\sf S}_i$ is {\sl weighted homogeneous},
say of weight $\mu_i$, namely satisfies:
\[
{\sf S}_i
\big(
\delta^{l_1}A_1,\,\dots,\,\delta^{l_{k_1}}A_{k_1}
\big)
=
\delta^{\mu_i}\,
{\sf S}_i
\big(
A_1,\,\dots,\,A_{k_1}
\big),
\]
in $\C \big[ A_1, \dots, A_{ k_1} \big]$ for every weighted dilation
factor $\delta \in \C$. Indeed, dilating $j^\kappa f$ as usual:
\[
\delta\cdot j^\kappa f
:=
\big(\delta^\lambda f_i^{(\lambda)}\big)_{1\leqslant i\leqslant n}^{
1\leqslant\lambda\leqslant\kappa},
\]
since the syzygies hold for any collection of $n\kappa$ components
$\big( f_i^{ (\lambda )}\big)_{ 1 \leqslant i\leqslant n}^{
1\leqslant\lambda\leqslant\kappa}$ in the jet space, they must hold
too with $\big( \delta^\lambda f_i^{ (\lambda )}\big)_{ 1 \leqslant
i\leqslant n}^{ 1\leqslant\lambda\leqslant\kappa}$, namely:
\[
\aligned
0
&
\equiv
{\sf S}_i
\Big(
L^{l_1}(\delta\cdot j^\kappa f)\big\vert_0,\,\dots,\,
L^{l_{k_1}}(\delta\cdot j^\kappa f)\big\vert_0
\Big)
\\
&
=
{\sf S}_i
\Big(
\delta^{l_1}\,
L^{l_1}(j^\kappa f)\big\vert_0,\,\dots,\,
\delta^{l_{k_1}}\,L^{l_{k_1}}(j^\kappa f)\big\vert_0
\Big)
\ \ \ \ \ \ \ \ \ \ \ \ \
{\scriptstyle{(i\,=\,1\,\cdots\,N_1)}},
\endaligned
\]
and we may use the fact that the $L^{l_i}$ are invariant under
reparametrization. Therefore, if we gather together, in each syzygy
polynomial ${\sf S}_i$, all terms which have equal, constant 
weight $\mu$:
\[
{\sf S}_i
=
\sum_\mu\,{\sf S}_i^\mu,
\ \ \ \ \ \ \ \ \
\text{\rm with}
\ \ \ \ \ \ \ \ \
{\sf S}_i^\mu
\big(\delta^{l_1}A_1,\,\dots,\,\delta^{l_{k_1}}A_{k_1}\big)
=
\delta^\mu\,
{\sf S}_i^\mu
\big(
A_1,\,\dots,\,A_{k_1}
\big),
\]
we may expand according to weight the obtained relations under the
specific form:
\[
0
\equiv
\sum_\mu\,
\delta^\mu\,
{\sf S}_i^\mu
\Big(
L^{l_1}(j^\kappa f)\big\vert_0,\,\dots,\,
L^{l_{k_1}}(j^\kappa f)\big\vert_0
\Big)
\ \ \ \ \ \ \ \ \ \ \ \ \
{\scriptstyle{(i\,=\,1\,\cdots\,N_1)}}.
\]
Because these identities then hold in $\C \big[ \delta, \, j^\kappa
f\big]$, they are equivalent to the (possibly larger) collection of
{\em constantly weighted} syzygies:
\[
0
\equiv
{\sf S}_i^\mu
\Big(
L^{l_1}(j^\kappa f)\big\vert_0,\,\dots,\,
L^{l_{k_1}}(j^\kappa f)
\Big)
\ \ \ \ \ \ \ \ \ \ \ \ \
{\scriptstyle{(i\,=\,1\,\cdots\,N_1\,;\,\forall\,\mu)}},
\]
and this justifies the claim.

So let $\mu_i$ be the weight of the (homogeneous) syzygy ${\sf S}_i$,
for $i = 1, \dots, N_1$. Because by assumption each polynomial ${\sf
S}_i \big( L^{ l_1} ( j^\kappa f) , \dots, L^{ l_{ k_1}} ( j^\kappa f)
\big)$ vanishes identically in $\C \big[ j^\kappa f\big]$ after
setting $f_1' = 0$, there are maximal factoring powers $(f_1')^{
\nu_i}$ of $f_1'$, with $1 \leqslant \nu_i \leqslant \infty$, and
there are certain (possibly zero) polynomial remainders ${\sf R}_i
\big( j^\kappa f\big)$ such that we may write in $\C \big[ j^\kappa f
\big]$:
\[
{\sf S}_i\big(
L^{l_1},\,\dots,\,L^{l_{k_1}}
\big)
=
(f_1')^{\nu_i}\,
{\sf R}_i\big(j^\kappa f\big)
\ \ \ \ \ \ \ \ \ \ \ \ \
{\scriptstyle{(i\,=\,1\,\cdots\,N_1)}},
\]
with ${\sf R}_i \not \equiv 0$ when $1 
\leqslant \nu_i < \infty$ and with ${\sf
R}_i = 0$ by convention when $\nu_i = \infty$.

We claim that each such ${\sf R}_i \big( j^\kappa f\big)$ is then a
bi-invariant. In fact, it is a polynomial by definition,
and its representation as a quotient:
\[
{\sf R}_i\big(j^\kappa f\big)
=
\frac{{\sf S}_i\big(L^{l_1},\dots,L^{l_{k_1}}\big)}{
(f_1')^{\nu_i}}
\]
of two polynomials invariant by reparametrizations and invariant under
the unipotent action shows at once that ${\sf R}_i$ too enjoys
bi-invariancy.

The second step of the algorithm consists in testing, for each $i$,
whether or not ${\sf R}_i$ belongs to the algebra $\C \big[ f_1', L^{
l_1}, \dots, L^{ l_{ k_1}} \big]$ generated by the initial
bi-invariants. In the case where no new bi-invariant appears, the
algorithm will be shown to terminate, so let us assume that at least
one ${\sf R}_i$ provides a new bi-invariant, independent of $f_1', L^{
l_1}, \dots, L^{ l_{ k_1}}$. It is then clear that after renumbering
the ${\sf R}_i$ if necessary, one may assume that:
\[
\left\{
\aligned
&
{\sf R}_1\ \ \
\text{\rm is independent of}\ \
f_1',L^{l_1},\dots,L^{l_{k_1}},
\\
&
{\sf R}_2\ \ \
\text{\rm is independent of}\ \
f_1',L^{l_1},\dots,L^{l_{k_1}},{\sf R}_1,
\\
&\cdots\cdots\cdots
\\
&
{\sf R}_{k_2}\ \ \
\text{\rm is independent of}\ \
f_1',L^{l_1},\dots,L^{l_{k_1}},{\sf R}_1,\dots,{\sf R}_{k_2-1},
\endaligned\right.
\]
while for the next indices $i = k_2+1, \dots, N_1$:
\[
\Big\{
{\sf R}_i\ \ \ 
\text{\rm belongs to the algebra}\ \
\C\big[f_1',L^{l_1},\dots,L^{l_{k_1}},{\sf R}_1,\dots,{\sf R}_{k_2}\big].
\]
Denoting instead by $M^{ m_1}, \dots, M^{ m_{k_2 }}$ these ${\sf R}_i$
for $i = 1, \dots, k_2$ which provide new mutually independent
bi-invariants, where as usual the weights $m_i := \mu_i - \nu_i$, for
$i=1, \dots, k_2$ are put in exponent place, we can therefore
write down in more explicit form the filled syzygy polynomials
(without setting $f_1' = 0$):
\[
\left\{
\aligned
0
&
\equiv
{\sf S}_i\big(L^{l_1},\dots,L^{l_{k_1}}\big)
+
(f_1')^{\nu_i}\,M^{m_i}
\ \ \ \ \ \ \ \ \ \ \ \ \
{\scriptstyle{(i\,=\,1\,\cdots\,k_2)}},
\\
0
&
\equiv
{\sf S}_i\big(L^{l_1},\dots,L^{l_{k_1}}\big)
+
(f_1')^{\nu_i}\,
{\sf R}_i\big(L^{l_1},\dots,L^{l_{k_1}},M^{m_1},\dots,M^{m_{k_2}}\big)
\ \ \ \ \ \ \ \
{\scriptstyle{(i\,=\,k_2+1\,\cdots\,N_1)}},
\endaligned\right.
\]
from which we recover at once, by setting $f_1'$, the original
syzygies:
\[
0
\equiv
{\sf S}_i
\big(L^{ l_1}\big\vert_0,\,\dots,\, 
L^{l_{ k_1}}\big\vert_0\big)
\ \ \ \ \ \ \ \ \ \ \ \ \
{\scriptstyle{(i\,=\,1\,\cdots\,N_1)}}.
\]
So the equations above, when written explicitly in specific
applications below, shall show both the collection of new appearing
bi-invariants $M^{ m_1}, \dots, M^{ m_{ k_2}}$ (without setting $f_1'
= 0$) and (after setting $f_1'$) a reduced Gröbner basis for the ideal of
relations between the initial bi-invariants $L^{ l_1 } \big\vert_0,
\dots, L^{ L_{ k_1}} \big\vert_0$.

\subsection*{ Second and further loops of the algorithm}
Next, we restart the process with the new, larger collection of
bi-invariants, namely we compute a reduced Gröbner basis (for a
certain monomial order compatible with the preceding loop):
\[
\text{\sf Ideal-Rel}
\Big(
L^{l_1}\big\vert_0,\,\dots,\,L^{l_{k_1}}\big\vert_0,\,\,
M^{m_1}\big\vert_0,\,\dots,\,M^{m_{k_2}}\big\vert_0
\Big).
\]
Write as follows the so obtained gröbnerized syzygies, after filling
the remainders behind a power of $f_1'$ and after testing whether
these remainders provide new bi-invariants:
\[
\left\{
\aligned
0
&
\equiv
{\sf S}_i\big(L^{l_1},\dots,L^{l_{k_1}}\big)
+
(f_1')^{\nu_i}\,M^{m_i}
\ \ \ \ \ \ \ \ \ \ \ \ \
{\scriptstyle{(i\,=\,1\,\cdots\,k_2)}},
\\
0
&
\equiv
{\sf S}_i\big(L^{l_1},\dots,L^{l_{k_1}}\big)
+
(f_1')^{\nu_i}\,
{\sf R}_i\big(L^{l_1},\dots,L^{l_{k_1}},M^{m_1},\dots,M^{m_{k_2}}\big)
\ \ \ \ \ \ \ \
{\scriptstyle{(i\,=\,k_2+1\,\cdots\,N_1)}},
\\
0
&
\equiv
{\sf T}_j
\big(L^{l_1},\dots,L^{l_{k_1}},M^{m_1},\dots,M^{m_{k_2}}\big)
+
(f_1')^{\nu_j}\,N^{n_j}
\ \ \ \ \ \ \ \ \ \ \ \ \
{\scriptstyle{(j\,=\,1\,\cdots\,k_3)}},
\\
0
&
\equiv
{\sf T}_j
\big(L^{l_1},\dots,L^{l_{k_1}},M^{m_1},\dots,M^{m_{k_2}}\big)
+
\\
&\ \ \ \ \
+
(f_1')^{\nu_j}\,
R_j\big(L^{l_1},\dots,L^{l_{k_1}},M^{m_1},\dots,M^{m_{k_2}},
N^{n_1},\dots,N^{n_{k_3}}\big)
\ \ \ \ \ \ \ \ 
{\scriptstyle{(j\,=\,k_3+1\,\cdots\,N_2)}}.
\endaligned\right.
\]
with $N^{ n_1}, \dots, N^{ n_{ k_3}}$ denoting the
new appearing bi-invariants, of weight $n_1, \dots, n_{ k_3}$.

Successively, continue to perform further loops as long as new
bi-invariants appear which do not belong to the algebra generated by
already known bi-invariants.

\subsection*{ Termination of the algorithm}
Either there always appear new bi-invariants or, after a finite number
of loops, we come to a situation which falls under the scope of the
following important statement.

\THEOREM

\smallskip\noindent\fbox{\bf THEOREM}\ \ 
\label{normal-syzygies}
{\sf\em
For a certain dimension $n$ and for a certain jet order $\kappa$,
suppose that, after performing a finite number of loops of the
algorithm, one possesses a finite number $1 + M$ of mutually
independent bi-invariants $f_1'$, $\Lambda^{\ell_1}$, \dots,
$\Lambda^{\ell_M} \in \C \big[ j^\kappa f_1, \dots, j^\kappa f_n
\big]$ of weights $1, \ell_1, \dots, \ell_M$ belonging to ${\sf
UE}_\kappa^n$, whose restrictions to $\{ f_1 ' = 0 \}$ share an ideal
of relations:
\[
\text{\sf Ideal-Rel}
\Big(\
\Lambda^{\ell_1}\big\vert_0,\ \
\dots\dots,\,
\Lambda^{\ell_M}\big\vert_0\
\Big)
\]
generated by a finite number $N$ (often large) of homogeneous syzygies:
\[
0
\equiv
{\sf S}_i
\big(\Lambda^{\ell_1}\big\vert_0,\,\dots,\,
\Lambda^{\ell_M}\big\vert_0\big),
\ \ \ \ \ \ \ \ \ \ \ \ \
{\scriptstyle{(i\,=\,1\,\cdots\,N)}}
\]
of weight $\mu_i$
assumed to be represented by a certain reduced Gröbner basis $\big<
{\sf S}_i \big>_{ 1 \leqslant i \leqslant N }$ for a certain monomial
order, with the crucial property that {\em no} new bi-invariant
appears behind $f_1'$, namely with the property that, without setting
$f_1' = 0$, one has $N$ identically satisfied relations:
\[
0
\equiv
{\sf S}_i
\big(\Lambda^{\ell_1},\,\dots,\,\Lambda^{\ell_M}\big)
-
f_1'\,{\sf R}_i
\big(f_1',\,\Lambda^{\ell_1},\,\dots,\, \Lambda^{\ell_M}\big)
\ \ \ \ \ \ \ \ \ \ \ \ \
{\scriptstyle{(i\,=\,1\,\cdots\,N)}},
\]
for some remainders ${\sf R }_i$ {\em which all depend polynomially
upon the same collection of invariants $f_1', \Lambda^{ \ell_1 },
\dots, \Lambda^{ \ell_M }$}, so that no new bi-invariant appears at
this stage.

Then the algorithm terminates and the algebra of bi-invariants
coincides with:
\[
\boxed{
{\sf UE}_\kappa^n
=
\C\big[
f_1',\,\Lambda^{\ell_1},\,\dots\dots,\,\Lambda^{\ell_M}
\big]
\ \ \ 
\text{\sf modulo syzygies}\
}\,.
\]

In addition, for these values of $n$ and of $\kappa$, if one denotes
the leading terms of the syzygies by:
\[
{\sf LT}\big({\sf S}_i(\Lambda)\big)
=
\big(\Lambda^{\ell_1}\big)^{\alpha_1^i}
\cdots
\big(\Lambda^{\ell_M}\big)^{\alpha_M^i}
\ \ \ \ \ \ \ \ \ \ \ \ \
{\scriptstyle{(i\,=\,1\,\cdots\,N)}},
\]
for certain specific multiindices $\big( \alpha_1^i, \dots, \alpha_M^i
\big) \in \N^M$, and if for $i = 1, \dots, N$ one denotes by:
\[
\square_i
:=
\alpha^i+\N^M
=
\big\{
\big(\alpha_1^i+b_1,\,\dots,\,\alpha_M^i+b_M\big):\,
b_1,\dots,b_M\in\N^M
\big\}
\]
the positive quadrant of $\N^M$ having vertex at $\alpha^i$, then a
general, arbitrary bi-invariant in ${\sf UE}_{\kappa, m}^n$ of weight
$m$ writes {\em uniquely} under the {\em normal form}:
\[
\sum_{0\leqslant a\leqslant m}\,
(f_1')^a\,\widetilde{\sf P}_a
\big(
\Lambda^{\ell_1},\,\dots,\,\Lambda^{\ell_M}
\big),
\]
with summation containing {\em only positive powers} of $f_1'$, where
each $\widetilde{ P}_a$ is of weight $m - a$ and is put under {\em
Gröbner-normalized form}:
\[
\boxed{
\widetilde{\sf P}_a
=
\sum_{(b_1,\dots,b_M)\in\N^M\backslash
(\square_1\cup\cdots\cup\square_N)\
\atop
\ell_1b_1+\cdots+\ell_Mb_M=m-a}\,
{\sf coeff}_{a;\,b_1,\dots,b_M}\cdot
\big(\Lambda^{\ell_1}\big)^{b_1}\cdots
\big(\Lambda^{\ell_M}\big)^{b_M}}\,,
\]
with complex coefficients ${\sf coeff }_{a; \,b_1, \dots, b_M}$
subjected to no restriction at all. }

\stopTHEOREM

\proof
We start with the list of initial bi-invariants $f_1', L^{ l_1 },
\dots, L^{ l_{ k_1 }}$ and with the initial, rational representation of
an arbitrary bi-invariant ${\sf P } \big( j^\kappa f \big) 
\in {\sf UE}_{ \kappa, m }^n$ which was obtained previously:
\[
\aligned
{\sf P}\big(j^\kappa f\big)
&
=
\sum_{-\frac{\kappa-1}{\kappa}m\leqslant a\leqslant m}\,
(f_1')^a\,{\sf P}_a\big(L^{l_1},\dots,L^{l_{k_1}}\big)
\\
&
=
(f_1')^{a_0}\,{\sf P}_{a_0}
+
\sum_{a_0+1\leqslant a\leqslant m}\,(f_1')^a\,{\sf P}_a,
\endaligned
\]
and we denote by $a_0$ the smallest appearing exponent of $f_1'$.
Clearly, the final list of bi-invariants $\Lambda^{ \ell_1 }, \dots,
\Lambda^{ \ell_M}$ stabilized after a finite number of loops of the
algorithm contains $L^{ l_1}, \dots, L^{ l_{ k_1 }}$ as its first
$k_1$ terms. Working in the polynomial ring $\C \big[ A^1, \dots, A^{
k_1}, \dots, A^M \big]$, we may then divide ${\sf P}_{ a_0}$ by the
ideal of relations $\big< {\sf S}_i ( A) \big>_{ 1\leqslant i
\leqslant N}$:
\[
{\sf P}_{a_0}
\big(A^1,\dots,A^{k_1}\big)
=
\widetilde{\sf P}_{a_0}
\big(A^1,\dots,A^{k_1},\dots,A^M\big)
+
\sum_{i=1}^N\,q_i(A)\cdot{\sf S}_i(A),
\]
with multiplicands $q_i ( A)$ of weight $m - a_0 - \mu_i$, getting a
remainder $\widetilde{ \sf P }_{ a_0}$ of weight $m - a_0$ which in
general will depend upon all the variables $A^1, \dots, A^{ k_1 },
\dots, A^M$ and which is {\em unique} (while the multiplicands $q_i$
cannot be unique, as soon as $N \geqslant 2$), by virtue of a
classical feature of Gröbner bases. Consequently, replacing the
independent variables $A^l$ by the bi-invariants in the arguments and
then substituting each ${\sf S }_i ( \Lambda)$ by $f_1' \, {\sf R }_i
( f_1', \Lambda)$\,\,---\,\,thanks to the main assumption that in
filled syzygies, all the remainders behind $f_1'$ depend polynomially
upon the same bi-invariants $f_1', \Lambda^{ \ell_1 }, \dots,
\Lambda^{ \ell_M }$\,\,---, we then get a normalized representation of
$\widetilde{\sf P}_{ a_0}$:
\[
\aligned
{\sf P}_{a_0}
\big(L^{l_1},\dots,L^{l_{k_1}}\big)
&
=
\widetilde{\sf P}_{a_0}
\big(\Lambda^{\ell_1},\dots,\Lambda^{\ell_M}\big)
+
\sum_{i=1}^N\,q_i(\Lambda)\cdot{\sf S}_i(\Lambda)
\\
&
=
\widetilde{\sf P}_{a_0}
\big(\Lambda^{\ell_1},\dots,\Lambda^{\ell_M}\big)
+
\sum_{i=1}^N\,q_i(\Lambda)\cdot
f_1'\,{\sf R}_i\big(f_1',\Lambda\big)
\\
&
=
\widetilde{\sf P}_{a_0}
\big(\Lambda^{\ell_1},\dots,\Lambda^{\ell_M}\big)
+
f_1'\,\widetilde{\sf R}_{a_0}
\big(f_1',\,\Lambda^{\ell_1},\dots,\Lambda^{\ell_M}\big),
\endaligned
\]
(modulo an uncontrolled remainder $\widetilde{\sf R}_{ a_0}$ which
hopefully, lies behind $f_1'$) which we may therefore inject in our
rational representation:
\[
{\sf P}\big(j^\kappa f\big)
=
(f_1')^{a_0}\,\widetilde{\sf P}_{a_0}(\Lambda)
+
(f_1')^{a_0+1}\,\widetilde{\sf R}_{a_0}
\big(f_1',\,\Lambda\big)
+
\sum_{a_0+1\leqslant a\leqslant m}\,
(f_1')^a\,{\sf P}_a(L).
\]
But both $\widetilde{ \sf P}_{ a_0}$ and $f_1' \, \widetilde{ \sf R}_{
a_0}$ being of weight $m-a_0$ as was ${\sf P}_{ a_0}$, it follows
that, when developping the perturbing term $(f_1')^{ a_0+1
}\,\widetilde{\sf R}_{ a_0 } \big( f_1', \Lambda)$ in powers of
$f_1'$, the fact that this remainder is of weight $m$ guarantees that
the sum does not go beyond $(f_1')^m$, and thus, we come to an
expression:
\[
{\sf P}\big(j^\kappa f\big)
=
(f_1')^{a_0}\,\widetilde{\sf P}_{a_0}(\Lambda)
+
\sum_{a_0+1\leqslant a\leqslant m}\,
(f_1')^a\,
{\sf Q}_a(\Lambda)
\]
entirely similar to the one we started with, whose first term:
\[
\widetilde{\sf P}_{a_0}
=
\sum_{(b_1,\dots,b_M)\in\N^M\backslash
(\square_1\cup\cdots\cup\square_N)\
\atop
\ell_1b_1+\cdots+\ell_Mb_M=m-a_0}\,
{\sf coeff}_{a_0;\,b_1,\dots,b_M}\cdot
\big(\Lambda^{\ell_1}\big)^{b_1}\cdots
\big(\Lambda^{\ell_M}\big)^{b_M},
\]
is normalized modulo the syzygies. But we can then subject the next
term ${\sf Q}_{ a_0 + 1} ( \Lambda)$ to the same process, and
consequently by induction, after a finite number of steps, we come to
an expression in which {\em all} multiplicands of a power of $f_1'$
have been normalized:
\[
{\sf P}\big(j^\kappa f\big) 
=
\sum_{a_0'\leqslant a\leqslant m}\,
(f_1')^a\,
\sum_{(b_1,\dots,b_M)\in\N^M\backslash
(\square_1\cup\cdots\cup\square_N)\
\atop
\ell_1b_1+\cdots+\ell_Mb_M=m-a}\,
{\sf coeff}_{a;\,b_1,\dots,b_M}\cdot
\big(\Lambda^{\ell_1}\big)^{b_1}\cdots
\big(\Lambda^{\ell_M}\big)^{b_M},
\]
with a possibly larger $a_0'$ (in case $\widetilde{ P}_{ a_0}$
vanishes identically). However, the smallest $a_0$ in the initial
expression for ${\sf P} \big( j^\kappa f\big)$ was possibly negative
and hence our $a_0'$ here can still be negative too, and our gained
representation of ${\sf P} \big( j^\kappa f \big)$ can still be not
polynomial.

Hopefully, we may now claim that there are no negative powers of
$f_1'$ anymore in such a normalized expression, so that the right hand
side is a true polynomial.

Indeed, suppose that $a_0' < 0$ with $\widetilde{ P }_{ a_0' } \not
\equiv 0$. Multiply both sides by $(f_1' )^{ - a_0'}$, set afterwards
$f_1' = 0$ and then get in such a way a nontrivial identity:
\[
0
\equiv
\sum_{(b_1,\dots,b_M)\in\N^M\backslash
(\square_1\cup\cdots\cup\square_N)\
\atop
\ell_1b_1+\cdots+\ell_Mb_M=m-a_0'}\,
{\sf coeff}_{a_0';\,b_1,\dots,b_M}\cdot
\big(\Lambda^{\ell_1}\big\vert_0\big)^{b_1}\cdots
\big(\Lambda^{\ell_M}\big\vert_0\big)^{b_M}.
\]
This equation would then represent a syzygy between bi-invariants
restricted to $\{ f_1' = 0\}$ whose leading term is strictly smaller
than the leadings terms of the syzygies ${\sf S}_i$. This would
contradict the assumption that the collection $\big< {\sf S}_i \big>_{
1 \leqslant i \leqslant N}$ is a Gröbner basis for the ideal of
relations between $\Lambda^{ \ell_1} \big\vert_0, \dots, \Lambda^{
\ell_M} \big \vert_0$. So $a_0' \geqslant 0$, namely the normalized
representation is polynomial.

The same argument shows that the normalized representation is unique.

Finally, it suffices to say, if not remarked stealthily before, that
any polynomial in $f_1', \, \Lambda^{ \ell_1}, \dots, \Lambda^{
\ell_M}$ obviously is a bi-invariant. The proof is now complete.
\endproof

\section*{\S10.~Seventeen bi-invariant generators 
\\
in dimension $n = 2$ for jet level $\kappa = 5$}
\label{Section-10}

\subsection*{ First loop of the algorithm}
According to these general principles, in the case $n = 2$, $\kappa =
5$, we should therefore start with the initial rational
representation:
\[
{\sf P}\big(j^5f\big)
=
\sum_{-\frac{4}{5}m\leqslant a\leqslant m}\,
(f_1')^a\,{\sf P}_a
\big(
\Lambda^3,\,\Lambda^5,\,\Lambda^7,\,\Lambda^9
\big)
\]
of an arbitrary bi-invariant ${\sf P} \big( j^5 f \big) \in {\sf
UE}_5^2$. For simplicity reasons, we shall denote without any lower index
each one of the appearing bi-invariants. In fact, among all the
invariants explicitly defined in the theorem on
p.~\pageref{bracket-5-2}, bi-invariants correspond to lower indices
being constantly equal to 1, and one has also to consider the {\em
non-bracket} bi-invariants introduced in \S8.

So according to the general algorithm, we have to start by computing
the ideal of relations:
\[
\text{\sf Ideal-Rel}
\Big(
\Lambda^3\big\vert_0,\
\Lambda^5\big\vert_0,\
\Lambda^7\big\vert_0,\
\Lambda^9\big\vert_0
\Big).
\]
For this easy first step, we may use any
Gröbner bases package\footnote{\, {\em See}
dim-2-order-5-step-1-with-FGb.mw at~\cite{ mer2008c}. }. For the
Reverse Degree Lexicographic Ordering, the result provided is:
\[
\aligned
0
&
\equiv
-7\,\Lambda^7\big\vert_0\Lambda^7\big\vert_0
+
5\,\Lambda^5\big\vert_0\Lambda^9\big\vert_0,
\\
0
&
\equiv
-7\,\Lambda^5\big\vert_0\Lambda^7\big\vert_0
+
3\,\Lambda^3\big\vert_0\Lambda^9\big\vert_0,
\\
0
&
\equiv
-5\,\Lambda^5\big\vert_0\Lambda^5\big\vert_0
+
3\,\Lambda^3\big\vert_0\Lambda^7\big\vert_0.
\endaligned
\]
Then we compute the remainder bi-invariants appearing behind a power
of $f_1'$. Here, for the three syzygies, the maximal factoring power
of $f_1'$ is the same, equal to 2, and three new bi-invariants appear:
\[
\aligned
0
&
\equiv
-7\,\Lambda^7\Lambda^7
+
5\,\Lambda^5\Lambda^9
\blue{-f_1'f_1'K^{12}},
\\
0
&
\equiv
-7\,\Lambda^5\Lambda^7
+
3\,\Lambda^3\Lambda^9
\blue{-f_1'f_1'M^{10}},
\\
0
&
\equiv
-5\,\Lambda^5\Lambda^5
+
3\,\Lambda^3\Lambda^7
\blue{-f_1'f_1'M^8},
\endaligned
\]
namely: $M^8$, $M^{ 10}$ and $K^{ 12 }$. Either looking at the
syzygies of the second loop below, or computing directly by hand, or
playing a bit with Maple, we find the values of the restrictions to
$\{ f_1' = 0 \}$ of all the bi-invariants obtained so far, expressed
in (rational) terms of the three restricted bi-invariants, $\Lambda^3
\vert_0$, $\Lambda^5 \vert_0$ and $M^8 \vert_0$ which are 
easily checked to be algebraically independent:
\[
\aligned
&
\underline{\Lambda^3}\big\vert_0
\\
&
\underline{\Lambda^5}\big\vert_0
\\
&
\Lambda^7\big\vert_0
=
{\textstyle{\green{\frac{5}{3}}\,
\frac{\Lambda^5\vert_0\,\Lambda^5\vert_0}{
\Lambda^3\vert_0\ \ \ \ \ \ }}},
\\
&
\Lambda^9\big\vert_0
=
{\textstyle{\green{\frac{35}{9}}\,
\frac{\Lambda^5\vert_0\,\Lambda^5\vert_0\,\Lambda^5\vert_0}{
\Lambda^3\vert_0\,\Lambda^3\vert_0\ \ \ \ \ \ }}},
\endaligned
\]
\[
\aligned
&
M^8\big\vert_0
\\
&
M^{10}\big\vert_0
=
{\textstyle{\green{\frac{8}{3}}\,
\frac{\Lambda^5\vert_0\,M^8\vert_0}{
\Lambda^3\vert_0\ \ \ \ \ \ }}},
\\
&
K^{12}\big\vert_0
=
{\textstyle{\green{\frac{5}{9}}\,
\frac{\Lambda^5\vert_0\,\Lambda^5\vert_0\,M^8\vert_0}{
\Lambda^3\vert_0\,\Lambda^3\vert_0\ \ \ \ \ \ }}}.
\endaligned
\]
Proceeding then as in the lemma on p.~\pageref{lemma-independence} and
using these rational expressions, one may establish that the 8
bi-invariants known so far, namely $f_1'$, $\Lambda^3$, $\Lambda^5$,
$\Lambda^7$, $\Lambda^9$, $M^8$, $M^{10}$ and $K^{ 12}$, are mutually
independent.

\subsection*{ Second loop of the algorithm}
Afterwards, we must compute the ideal of relations between the
7 restricted bi-invariants in question:
\[
\text{\sf Ideal-Rel}
\Big(
\Lambda^3\big\vert_0,\
\Lambda^5\big\vert_0,\
\Lambda^7\big\vert_0,\
\Lambda^9\big\vert_0,\
M^8\big\vert_0,\
M^{10}\big\vert_0,\
K^{12}\big\vert_0
\Big).
\]
For the Degree Reverse Lexicographic Ordering, a Gröbner basis for
this ideal of relations consists of the following 10
polynomials\footnote{\, {\em See} dim-2-order-5-step-2-with-FGb.mw
at~\cite{ mer2008c}. } (in which the remainders behind a power of
$f_1'$ have already been filled):
\[
\aligned
0
&
\equiv
-5\,M^{10}M^{10}
+
64\,M^8K^{12}
\blue{-f_1'X^{19}},
\\
0
&
\equiv
-5\,\Lambda^9M^{10}
+
56\,\Lambda^7K^{12}
\blue{-f_1'X^{18}},
\\
0
&
\equiv
-8\,\Lambda^9M^8
+
7\,\Lambda^7M^{10}
\green{-f_1'F^{16}},
\\
0
&
\equiv
-\Lambda^9M^8
+
7\,\Lambda^5K^{12}
\blue{-f_1'F^{16}},
\\
0
&
\equiv
-8\,\Lambda^7M^8
+
5\,\Lambda^5M^{10}
\green{-f_1'H^{14}},
\endaligned
\]
\[
\aligned
0
&
\equiv
-\Lambda^7M^8
+
3\,\Lambda^3K^{12}
\blue{-f_1'H^{14}},
\\
0
&
\equiv
-8\,\Lambda^5M^8
+
3\,\Lambda^3M^{10}
\blue{-f_1'N^{12}},
\\
0
&
\equiv
-7\,\Lambda^7\Lambda^7
+
5\,\Lambda^5\Lambda^9
\green{-f_1'f_1'K^{12}},
\\
0
&
\equiv
-7\,\Lambda^5\Lambda^7
+
3\,\Lambda^3\Lambda^9
\green{-f_1'f_1'M^{10}},
\\
0
&
\equiv
-5\,\Lambda^5\Lambda^5
+
3\,\Lambda^3\Lambda^7
\green{-f_1'f_1'M^8}.
\endaligned
\]
How exactly do we manage to fill in what appears at the end of each
syzygy behind any power of $f_1'$?

\subsection*{ A standard obstacle: unavailability because of size
computations}
A natural idea would be to automatically apply the {\em Algebra
Membership Algorithm} based on Gröbner bases (\cite{ ess2000},
p.~289), but this would be (at least for us) impossible, because this
test would rely upon the (unavalaible to us) knowledge of a full
Gröbner basis for the ideal generated by the 8 equations:
\[
t_1-f_1',\
l_3-\Lambda^3,\
l_5-\Lambda^5,\
l_7-\Lambda^7,\
l_9-\Lambda^9,\
m_8-M^8,\
m_{10}-M^{10},\
k_{12}-K^{12},
\]
in the ring of 18 variables: 
\[
\C\big[ 
j^5f_1, j^5 f_2,\, 
t_1,l_3,l_5,l_7,l_9,m_8,m_{10},k_{12}
\big]
\]
with any monomial ordering having the only property that each jet
variable $f_i^{ ( \lambda)}$ is bigger than any monomial written with
only the 8 auxiliary variables $t_1, l_3, l_5, l_7, l_9, m_8, m_{ 10},
k_{ 12}$. Indeed, according to Proposition C.2.3 in the reference
cited, any remainder behind a power of $f_1'$, for instance the one
appearing in the sixth syzygy above: 
\[
{\sf rem}_6
:=
\frac{1}{f_1'}\,\big(8\,\Lambda^5M^8 
-
3\,\Lambda^3M^{10}\big),
\]
would then belong to the algebra generated by the 8 already known
bi-invariants: $f_1'$, $\Lambda^3$, $\Lambda^5$, $\Lambda^7$,
$\Lambda^9$, $M^8$, $M^{ 10}$, $K^{ 12}$, if and only if the {\em
normal form} of ${\sf rem}_6$ with respect to such a Gröbner basis
would belong to $\C \big[ t_1, \dots, k_{ 12} \big]$, and in such a
case, the (unique) normal form in question ${\sf rem}_6$ would provide
without any further effort the corresponding polynomial.

\smallskip

\centerline{\fbox{\sf However, Gröbner bases
here are blocked due to oversizeness 
}} 

\smallskip
Hence to bypass such a (usual, forseeable) 
drawback of Gröbner bases, we have to proceed differently.

What we do using Maple is a little bit tricky, and it works
well. After division by $f_1'$ (most often, and sometimes also by
$(f_1')^2$, but never by $(f_1')^3$), we start by computing each one
of the 10 remainder; in fact, since 3 of them were already treated in
the first loop, only 7 remainders have to be studied here. On the
other hand and as an independent preparation, we may check by
inspecting the explicit expressions given at the end of \S4, that
$\Lambda^3 \big\vert_0$, $\Lambda^5 \big\vert_0$, $M^8 \big\vert_0$
and $N^{ 12 } \big\vert_0$ (our ${\sf rem}_6$ itself!) are mutually
algebraically independent. Subsequently, we compute a Gröbner basis
for the four polynomial:
\[
l_{30}-\Lambda^3\big\vert_0,\
l_{50}-\Lambda^5\big\vert_0,\
m_{80}-M^8\big\vert_0,\
n_{120}-N^{12}\big\vert_0,
\]
in the ring $\C \big[ j^5f_1 \vert_0, j^5 f_2, l_{ 30}, l_{ 50}, m_{
80}, n_{ 120} \big]$, where $l_{ 30}$, $l_{ 50}$, $m_{ 80}$ and $n_{
120}$ denote auxiliary, supplementary variables, with any monomial
order having the property that each jet variable $f_i^{ ( \lambda)}$
is bigger than any monomial written with only the 4 auxiliary
variables $l_{ 30}$, $l_{ 50}$, $m_{ 80}$ and $n_{ 120}$. This then
is available to the computer: size is reasonable and it costs less
than 5 minutes on any computer. Then we set $f_1' = 0$ in each
remainder ${\sf rem}_k$, getting ${\sf rem}_k \big\vert_0$. We then
multiply each restricted remainder ${\sf rem}_k \big\vert_0$ for $k =
1, 2, \dots, 10$ by a suitable power of $\Lambda^3 \big\vert_0$
choosen by head, for instance if one looks at the third remainder:
\[
\Lambda^3\big\vert_0\Lambda^3\big\vert_0\cdot
{\sf rem}_3\big\vert_0
=
\Lambda^3\big\vert_0\Lambda^3\big\vert_0
\cdot
\bigg[
\frac{1}{f_1'}
\big(8\,\Lambda^9M^8-7\,\Lambda^7M^{10}\big)\bigg]_{f_1'=0}.
\]
Then we compute the normal form of this latter polynomial with respect
to the mentioned auxiliary Gröbner basis. For instance, our computer
yields for the third remainder the normal form:
\[
{\textstyle{\frac{35}{9}}}\,
l_{50}l_{50}n_{120}.
\] 
This result therefore means that the third unknown remainder ${\sf
rem}_3$ (appearing in the third syzygy) which we denoted in advance by
$F^{ 16}$, has the following value after setting $f_1' = 0$:
\[
F^{16}\big\vert_0
=
{\textstyle{\green{\frac{35}{9}}\,
\frac{\Lambda^5\vert_0\,\Lambda^5\vert_0\,N^{12}\vert_0}{
\Lambda^3\vert_0\,\Lambda^3\vert_0\ \ \ \ \ \ \ \ }}}.
\]
Then we test by hand and by head whether such a value for $f_1' = 0$
can be obtained as a polynomial in terms of the 7 previously known
restricted bi-invariants $\Lambda^3 \vert_0, \dots, K^{ 12}
\vert_0$. Here, it is easy to convince oneself that this cannot be the
case, so that $F^{ 16}$ really is a new bi-invariant.

On the other hand, we should do the same work for the fourth remainder
${\sf rem}_4$. It then happens that we find the {\em same} value at
$f_1' = 0$ in terms of $\Lambda^3 \vert_0, \Lambda^5 \vert_0, M^8
\vert_0, N^{ 12} \vert_0$. So we suspect that {\em without setting
$f_1' = 0$}, the two remainders ${\sf rem}_3$ and ${\sf rem}_4$ could
be identical and finally, a simple computation with Maple verifies
that this is indeed the case. Other remainders are computed similary,
and we thus have fully explained all our trick to bypass the
unavailability of full Gröbner bases due to oversizeness in this
problem. 

However, we would like to mention that achieving such a kind of task
took hours and days of patience. Hopefully, checking {\em a
posteriori} with Maple that a syzygy effectively holds is much, much 
more rapid and the reader will find in the Maple worksheets
referenced here the declaration of new bi-invariants at each 
step and the checking (at a piece) of all syzygies by means
of the basic ``{\small {\sf simplify}}'' command of Maple.

Finally, to finish with the second loop,
we give the values, restricted to $\{ f_1' = 0\}$, of the 5
appearing new bi-invariants at this stage:
\[
\aligned
&
N^{12}\big\vert_0
\\
&
H^{14}\big\vert_0
=
{\textstyle{\green{\frac{5}{3}}\,
\frac{\Lambda^5\vert_0\,N^{12}\vert_0}{
\Lambda^3\vert_0\,\ \ \ \ \ \ }}},
\\
&
F^{16}\big\vert_0
=
{\textstyle{\green{\frac{35}{9}}\,
\frac{\Lambda^5\vert_0\,\Lambda^5\vert_0\,N^{12}\vert_0}{
\Lambda^3\vert_0\,\Lambda^3\vert_0\ \ \ \ \ \ \ \ }}},
\\
&
X^{18}\big\vert_0
=
{\textstyle{\green{\frac{1225}{27}}\,
\frac{\Lambda^5\vert_0\,
\Lambda^5\vert_0\,\Lambda^5\vert_0\,N^{12}\vert_0}{
\Lambda^3\vert_0\,\Lambda^3\vert_0\,\Lambda^3\vert_0\ \ \ \ \ \ \ \ }}},
\\
&
X^{19}\big\vert_0
=
{\textstyle{\green{\frac{80}{3}}\,
\frac{
\Lambda^5\vert_0\,M^8\vert_0\,N^{12}\vert_0}{
\Lambda^3\vert_0\,\Lambda^3\vert_0\ \ \ \ \ \ \ \ \ }}}.
\endaligned
\]

\subsection*{ Third loop of the algorithm}
Now that we have explained how we proceed, we can offer directly the
32 filled syzygies appearing at the next step\footnote{\, {\em See}
dim-2-order-5-step-3-with-FGb.mw at~\cite{ mer2008c}. }, 
again for the Degree Reverse Lexicographic Ordering.
\[
\aligned
0
&
\equiv
-5\,F^{16}F^{16}
+
H^{14}X^{18}
\green{-f_1'K^{12}X^{19}},
\\
0
&
\equiv
-7\,H^{14}F^{16}
+
N^{12}X^{18}
\green{-f_1'M^{10}X^{19}},
\\
0
&
\equiv
-7\,H^{14}H^{14}
+
5\,N^{12}F^{16}
\green{-f_1'M^8X^{19}},
\\
0
&
\equiv
-56\,K^{12}F^{16}
+
M^{10}X^{18}
\blue{-f_1'Y^{27}},
\endaligned
\]
\[
\aligned
0
&
\equiv
-56K^{12}H^{14}
+
5\,M^{10}F^{16}
\blue{-f_1'X^{25}},
\\
0
&
\equiv
-8\,K^{12}N^{12}
+
M^{10}H^{14}
\blue{-f_1'X^{23}},
\\
0
&
\equiv
-49\,K^{12}H^{14}
+
M^8X^{18}
\green{-f_1'X^{25}},
\\
0
&
\equiv
-7\,K^{12}N^{12}
+
M^8F^{16}
\green{-f_1'X^{23}},
\\
0
&
\equiv
-5\,M^{10}N^{12}
+
8\,M^8H^{14}
\blue{-f_1'X^{21}},
\endaligned
\]
\[
\aligned
\\
0
&
\equiv
-48\,K^{12}F^{16}
+
\Lambda^9X^{19}
\green{-f_1'Y^{27}},
\\
0
&
\equiv
-48\,K^{12}H^{14}
+
\Lambda^7X^{19}
\green{-f_1'X^{25}}
\\
0
&
\equiv
-5\,\Lambda^9F^{16}
+
\Lambda^7X^{18}
\green{+8f_1'K^{12}K^{12}},
\\
0
&
\equiv
-\Lambda^9H^{14}
+
\Lambda^7F^{16}
\green{+f_1'M^{10}K^{12}},
\\
0
&
\equiv
-5\,\Lambda^9N^{12}
+
7\,\Lambda^7H^{14}
\green{+56f_1'M^8K^{12}-f_1'f_1'X^{19}},
\\
0
&
\equiv
-48K^{12}N^{12}
+
\Lambda^5X^{19}
\green{-7\,f_1'X^{23}},
\endaligned
\]
\[
\aligned
0
&
\equiv
-7\,\Lambda^9H^{14}
+
\Lambda^5X^{18}
\green{+8\,f_1'M^{10}K^{12}},
\\
0
&
\equiv
-\Lambda^9N^{12}
+
\Lambda^5F^{16}
\green{+f_1'M^{10}M^{10}},
\\
0
&
\equiv
-\Lambda^7N^{12}
+
\Lambda^5H^{14}
\green{+f_1'M^8M^{10}},
\\
0
&
\equiv
-10\,M^{10}N^{12}
+
\Lambda^3X^{19}
\green{-{\textstyle{\frac{7}{3}}}\,f_1'X^{21}},
\\
0
&
\equiv
-35\,\Lambda^9N^{12}
+
3\,\Lambda^3X^{18}
\green{+{\textstyle{\frac{285}{8}}}\,f_1'M^{10}M^{10}
-{\textstyle{\frac{7}{8}}}\,f_1'f_1'X^{19}},
\\
0
&
\equiv
-7\,\Lambda^7N^{12}
+
3\,\Lambda^3F^{16}
\green{+8f_1'M^8M^{10}},
\\
0
&
\equiv
-5\,\Lambda^5N^{12}
+
3\,\Lambda^3H^{14}
\green{+8\,f_1'M^8M^8},
\endaligned
\]
\[
\aligned
0
&
\equiv
-5\,M^{10}M^{10}
+
64\,M^8K^{12}
\green{-f_1'X^{19}},
\\
0
&
\equiv
-5\,\Lambda^9M^{10}
+
56\,\Lambda^7K^{12}
\green{-f_1'X^{18}},
\\
0
&
\equiv
-8\,\Lambda^9M^8
+
7\,\Lambda^7M^{10}
\green{-f_1'F^{16}},
\\
0
&
\equiv
-\Lambda^9M^8
+
7\,\Lambda^5K^{12}
\green{-f_1'F^{16}},
\\
0
&
\equiv
-8\,\Lambda^7M^8
+
5\,\Lambda^5M^{10}
\green{-f_1'H^{14}},
\endaligned
\]
\[
\aligned
0
&
\equiv
-\Lambda^7M^8
+
3\,\Lambda^3K^{12}
\green{-f_1'H^{14}},
\\
0
&
\equiv
-8\,\Lambda^5M^8
+
3\,\Lambda^3M^{10}
\green{-f_1'N^{12}},
\\
0
&
\equiv
-7\,\Lambda^7\Lambda^7
+
5\,\Lambda^5\Lambda^9
\green{-f_1'f_1'K^{12}},
\\
0
&
\equiv
-7\,\Lambda^5\Lambda^7
+
3\,\Lambda^3\Lambda^9
\green{-f_1'f_1'M^{10}},
\\
0
&
\equiv
-5\,\Lambda^5\Lambda^5
+
3\,\Lambda^3\Lambda^7
\green{-f_1'f_1'M^8}.
\endaligned
\]
Here, 4 new bi-invariants appear:
\[
X^{21},\ \ \ \ \
X^{23},\ \ \ \ \
X^{25},\ \ \ \ \
Y^{17}.
\]
Their values restricted to $\{ f_1 ' = 0\}$ are:
\[
\aligned
&
X^{21}\big\vert_0
=
{\textstyle{\green{-\frac{5}{3}}\,
\frac{N^{12}\vert_0\,N^{12}\vert_0}{
\Lambda^3\vert_0\,\ \ \ \ \ \ }}}
{\textstyle{\green{-\frac{64}{3}}\,
\frac{M^8\vert_0\,M^8\vert_0\,M^8\vert_0}{
\Lambda^3\vert_0\,\ \ \ \ \ \ \ \ \ \ \ \ \ }}},
\\
&
X^{23}\big\vert_0
=
{\textstyle{\green{-\frac{35}{3}}\,
\frac{\Lambda^5\vert_0\,N^{12}\vert_0\,N^{12}\vert_0}{
\Lambda^3\vert_0\,\Lambda^3\vert_0\, \ \ \ \ \ \ \ \ }}}
{\textstyle{\green{-\frac{64}{9}}\,
\frac{\Lambda^5\vert_0\,M^8\vert_0\,M^8\vert_0\,M^8\vert_0}{
\Lambda^3\vert_0\,\Lambda^3\vert_0\,\ \ \ \ \ \ \ \ \ \ \ \ \ \ }}},
\\
&
X^{25}\big\vert_0
=
{\textstyle{\green{-\frac{1225}{27}}\,
\frac{\Lambda^5\vert_0\,\Lambda^5\vert_0\,N^{12}\vert_0\,N^{12}\vert_0}{
\Lambda^3\vert_0\,\Lambda^3\vert_0\,\Lambda^3\vert_0\, 
\ \ \ \ \ \ \ \ }}}
{\textstyle{\green{-\frac{320}{27}}\,
\frac{\Lambda^5\vert_0\,\Lambda^5\vert_0\,
M^8\vert_0\,M^8\vert_0\,M^8\vert_0}{
\Lambda^3\vert_0\,\Lambda^3\vert_0\,\Lambda^3\vert_0\,
\ \ \ \ \ \ \ \ \ \ \ \ \ \ }}},
\\
&
Y^{27}\big\vert_0
=
{\textstyle{\green{-\frac{8575}{81}}\,
\frac{\Lambda^5\vert_0\,\Lambda^5\vert_0\,\Lambda^5\vert_0\,
N^{12}\vert_0\,N^{12}\vert_0}{
\Lambda^3\vert_0\,\Lambda^3\vert_0\,\Lambda^3\vert_0\,\Lambda^3\vert_0\, 
\ \ \ \ \ \ \ \ }}}
{\textstyle{\green{-\frac{320}{81}}\,
\frac{\Lambda^5\vert_0\,\Lambda^5\vert_0\,\Lambda^5\vert_0\,
M^8\vert_0\,M^8\vert_0\,M^8\vert_0}{
\Lambda^3\vert_0\,\Lambda^3\vert_0\,\Lambda^3\vert_0\,\Lambda^3\vert_0\,
\ \ \ \ \ \ \ \ \ \ \ \ \ \ }}}.
\endaligned
\]

\subsection*{ Fourth loop of the algorithm}
The Gröbner basis of syzygies between the restriction to $\{ f_1' = 0
\}$ of the 17 bi-invariants known so far consists here of 105
equations. By an independent calculation, we checked that 39 among
these 105 generators belong to the ideal of the 66 remaining ones. We
could fill in the remainders behind a power of $f_1'$. To test whether
there appear new bi-invariants, it is in fact useless to fill in the
39 left out remainders. Here are the 66 syzygies\footnote{\, {\em See}
dim-2-order-5-step-4-with-FGb.mw at~\cite{ mer2008c}. } in question:
\[
\aligned
0
&
\equiv
X^{18}X^{23}-8\,F^{16}X^{25}+7\,H^{14}Y^{27}
\red{+0},
\\
0
&
\equiv
5\,F^{16}X^{23}-8\,H^{14}X^{25}+5\,N^{12}Y^{27}
\green{+f_1'\,X^{19}X^{19}},
\\
0
&
\equiv
7\,K^{12}X^{23}
-
M^{10}X^{25}
+
M^8Y^{27}
\red{+0},
\\
0
&
\equiv
5\,\Lambda^9X^{23}
-
8\,\Lambda^7X^{25}
+
5\,\Lambda^5Y^{27}
\green{-8\,f_1'\,K^{12}X^{19}},
\\
0
&
\equiv
7\,\Lambda^7X^{23}
-
8\,\Lambda^5X^{25}
+
3\,\Lambda^3Y^{27}
\green{-f_1'\,M^{10}X^{19}},
\\
0
&
\equiv
X^{18}X^{21}-57\,H^{14}X^{25}+40\,N^{12}Y^{27}
\green{+7\,f_1'\,X^{19}X^{19}},
\endaligned
\]
\[
\aligned
0
&
\equiv
F^{16}X^{21}-8\,H^{14}X^{23}+N^{12}X^{25}
\red{+0},
\\
0
&
\equiv
7\,K^{12}X^{21}-5\,M^{10}X^{23}+M^8X^{25}
\red{+0},
\\
0
&
\equiv
7\,\Lambda^9X^{21}-57\,\Lambda^5X^{25}+24\,\Lambda^3Y^{27}
\green{-15\,f_1'\,M^{10}X^{19}},
\\
0
&
\equiv
7\,\Lambda^7X^{21}-40\,\Lambda^5X^{23}+3\,\Lambda^3X^{25}
\green{-8\,f_1'\,M^8X^{19}},
\\
0
&
\equiv
X^{18}X^{19}-8\,K^{12}X^{25}+5\,M^{10}Y^{27}
\red{+0},
\\
0
&
\equiv
7\,F^{16}X^{19}-M^{10}X^{25}+8\,M^8Y^{27}
\red{+0},
\endaligned
\]
\[
\aligned
0
&
\equiv
7\,H^{14}X^{19}-5\,M^{10}X^{23}+8\,M^8X^{25}
\red{+0},
\\
0
&
\equiv
N^{12}X^{19}-M^{10}X^{21}+8\,M^8X^{23}
\red{+0},
\\
0
&
\equiv
6\,F^{16}X^{18}-\Lambda^9X^{25}+7\,\Lambda^7Y^{27}
\red{+0},
\\
0
&
\equiv
6\,H^{14}X^{18}-\Lambda^7X^{25}+5\,\Lambda^5Y^{27}
\green{-7\,f_1'\,K^{12}X^{19}},
\\
0
&
\equiv
6\,N^{12}X^{18}-\Lambda^5X^{25}+3\,\Lambda^3Y^{27}
\green{-7\,f_1'\,M^{10}X^{19}},
\\
0
&
\equiv
6\,M^{10}X^{18}-7\,\Lambda^9X^{19}
\green{+f_1'Y^{27}},
\endaligned
\]
\[
\aligned
0
&
\equiv
48\,M^8X^{18}-49\,\Lambda^7X^{19}
\green{+f_1'X^{25}},
\\
0
&
\equiv
30\,F^{16}F^{16}-\Lambda^7X^{25}+5\,\Lambda^5Y^{27}
\green{-f_1'\,K^{12}X^{19}},
\\
0
&
\equiv
42\,H^{14}F^{16}-\Lambda^5X^{25}
+3\,\Lambda^3Y^{27}
\green{-f_1'\,M^{10}X^{19}},
\\
0
&
\equiv
30\,N^{12}F^{16}-5\,\Lambda^5X^{23}+3\,\Lambda^3X^{25}
\green{-f_1'\,M^8X^{19}},
\\
0
&
\equiv
48\,K^{12}F^{16}-\Lambda^9X^{19}
\green{+f_1'Y^{27}},
\\
0
&
\equiv
30\,M^{10}F^{16}-7\,\Lambda^7X^{19}
\green{+f_1'X^{25}},
\endaligned
\]
\[
\aligned
0
&
\equiv
48\,M^8F^{16}-7\,\Lambda^5X^{19}
\green{+f_1'X^{23}},
\\
0
&
\equiv
5\,\Lambda^9F^{16}-\Lambda^7X^{18}
\green{-8\,f_1'\,K^{12}K^{12}},
\\
0
&
\equiv
7\,\Lambda^7F^{16}-\Lambda^5X^{18}
\green{-f_1'\,M^{10}K^{12}},
\\
0
&
\equiv
35\,\Lambda^5F^{16}-3\,\Lambda^3X^{18}
\green{-8\,f_1'\,M^8K^{12}+f_1'f_1'\,X^{19}},
\\
0
&
\equiv
42\,H^{14}H^{14}-5\,\Lambda^5X^{23}+3\,\Lambda^3X^{25}
\green{-f_1'M^8X^{19}},
\\
0
&
\equiv
6\,N^{12}H^{14}-\Lambda^5X^{21}+3\,\Lambda^3X^{23}
\red{+0},
\endaligned
\]
\[
\aligned
0
&
\equiv
48\,K^{12}H^{14}-\Lambda^7X^{19}
\green{+f_1'X^{25}},
\\
0
&
\equiv
6\,M^{10}H^{14}-\Lambda^5X^{19}
\green{+f_1'X^{23}},
\\
0
&
\equiv
16\,M^8H^{14}-\Lambda^3X^{19}
\green{+{\textstyle{\frac{1}{3}}}\,f_1'\,X^{21}},
\\
0
&
\equiv
7\,\Lambda^9H^{14}-\Lambda^5X^{18}
\green{-8\,f_1'\,M^{10}K^{12}},
\\
0
&
\equiv
49\,\Lambda^7H^{14}-3\,\Lambda^3X^{18}
\green{-5\,f_1'\,M^{10}M^{10}},
\\
0
&
\equiv
7\,\Lambda^5H^{14}-3\,\Lambda^3F^{16}
\green{-f_1'M^8M^{10}},
\endaligned
\]
\[
\aligned
0
&
\equiv
48\,K^{12}N^{12}-\Lambda^5X^{19}
\green{+7\,f_1'X^{23}},
\\
0
&
\equiv
10\,M^{10}N^{12}-\Lambda^3X^{19}
\green{+{\textstyle{\frac{7}{3}}}\,f_1'\,X^{21}},
\\
0
&
\equiv
35\,\Lambda^9N^{12}-3\,\Lambda^3X^{18}
\green{-{\textstyle{\frac{285}{8}}}\,f_1'\,M^{10}M^{10}
+{\textstyle{\frac{7}{8}}}\,f_1'f_1'\,X^{19}},
\\
0
&
\equiv
7\,\Lambda^7N^{12}-3\,\Lambda^3F^{16}
\green{-8\,f_1'\,M^8M^{10}},
\\
0
&
\equiv
5\,\Lambda^5N^{12}-3\,\Lambda^3H^{14}
\green{-8\,f_1'\,M^8M^8},
\\
0
&
\equiv
5\,M^{10}M^{10}-64\,M^8K^{12}
\green{+f_1'X^{19}},
\endaligned
\]
\[
\aligned
0
&
\equiv
5\,\Lambda^9M^{10}-56\,\Lambda^7K^{12}
\green{+f_1'X^{18}},
\\
0
&
\equiv
\Lambda^7M^{10}-8\,\Lambda^5K^{12}
\green{+f_1'F^{16}},
\\
0
&
\equiv
5\,\Lambda^5M^{10}-24\,\Lambda^3K^{12}
\green{+f_1'H^{14}},
\\
0
&
\equiv
\Lambda^9M^8-7\,\Lambda^5K^{12}
\green{+f_1'F^{16}},
\\
0
&
\equiv
\Lambda^7M^8-3\,\Lambda^3K^{12}
\green{+f_1'H^{14}},
\\
0
&
\equiv
8\,\Lambda^5M^8-3\,\Lambda^3M^{10}
\green{+f_1'N^{12}},
\endaligned
\]
\[
\aligned
0
&
\equiv
7\,\Lambda^7\Lambda^7-5\,\Lambda^5\Lambda^9
\green{+f_1'f_1'\,K^{12}},
\\
0
&
\equiv
7\,\Lambda^5\Lambda^7-3\,\Lambda^3\Lambda^9
\green{+f_1'f_1'M^{10}},
\\
0
&
\equiv
5\,\Lambda^5\Lambda^5-3\,\Lambda^3\Lambda^7
\green{+f_1'f_1'\,M^8},
\\
0
&
\equiv
7\,K^{12}X^{19}X^{19}+X^{25}X^{25}-5\,X^{23}Y^{27}
\red{+0},
\\
0
&
\equiv
M^{10}X^{19}X^{19}+X^{23}X^{25}-X^{21}Y^{27}
\red{+0},
\\
0
&
\equiv
M^8X^{19}X^{19}+5\,X^{23}X^{23}-X^{21}X^{25}
\red{+0},
\endaligned
\]
\[
\aligned
0
&
\equiv
56\,K^{12}K^{12}X^{19}+X^{18}X^{25}-5\,F^{16}Y^{27}
\red{+0},
\\
0
&
\equiv
M^{10}K^{12}X^{19}+F^{16}X^{25}-H^{14}Y^{27}
\red{+0},
\\
0
&
\equiv
8\,M^8K^{12}X^{19}+7\,H^{14}X^{25}-5\,N^{12}Y^{27}
\green{-f_1'X^{19}X^{19}},
\\
0
&
\equiv
M^8M^{10}X^{19}+7\,H^{14}X^{23}-N^{12}X^{25}
\red{+0},
\\
0
&
\equiv
8\,M^8M^8X^{19}+7\,H^{14}X^{21}-5\,N^{12}X^{23}
\red{+0},
\\
0
&
\equiv
448\,K^{12}K^{12}K^{12}+X^{18}X^{18}+5\,\Lambda^9Y^{27}
\red{+0},
\endaligned
\]
\[
\aligned
0
&
\equiv
48\,M^{10}K^{12}K^{12}+\Lambda^9X^{25}-\Lambda^7Y^{27}
\red{+0},
\\
0
&
\equiv
384\,M^8K^{12}K^{12}+7\,\Lambda^7X^{25}-5\,\Lambda^5Y^{27}
\green{+f_1'K^{12}X^{19}},
\\
0
&
\equiv
48\,M^8M^{10}K^{12}+7\,\Lambda^5X^{25}-3\,\Lambda^3Y^{27}
\green{+f_1'\,M^{10}X^{19}},
\\
0
&
\equiv
384\,M^8M^8K^{12}+35\,\Lambda^5X^{23}-3\,\Lambda^3X^{25}
\green{+f_1'M^8X^{19}},
\\
0
&
\equiv
48\,M^8M^8M^{10}+7\,\Lambda^5X^{21}-3\,\Lambda^3X^{23}
\red{+0},
\\
0
&
\equiv
64\,M^8M^8M^8+5\,N^{12}N^{12}+3\,\Lambda^3X^{21}
\red{+0}.
\endaligned
\]
Remarkably, no new bi-invariant appears at this fourth stage. 
According to the general principle, we may therefore
conclude that the algorithm stops.

\THEOREM

\smallskip\noindent\fbox{\bf THEOREM}\ \
{\sf\em 
In dimension $n = 2$ for jet order $\kappa = 5$, the algebra ${\sf
UE}_5^2$ of jet polynomials ${\sf P} \big( j^5 f_1, j^5 f_2 \big)$
invariant by reparametrization and invariant under the unipotent
action is generated by the 17 mutually independent bi-invariants
explicitly defined above:
\[
\boxed{
\aligned
&
\ \ \ \ \ \ \ 
f_1',\ \ \ \ \
\Lambda^3,\ \ \ \ \
\Lambda^5,\ \ \ \ \
\Lambda^7,\ \ \ \ \
\Lambda^9,\ \ \ \ \
M^8,\ \ \ \ \
M^{10},\ \ \ \ \ 
K^{12},
\\
&
N^{12},\ \ \ \ 
H^{14},\ \ \ \ 
F^{16},\ \ \ \ 
X^{18},\ \ \ \ 
X^{19},\ \ \ \ 
X^{21},\ \ \ \ 
X^{23},\ \ \ \ 
X^{25},\ \ \ \ 
Y^{27}
\endaligned
}
\]
whose restriction to $\{ f_1' = 0 \}$ has a reduced gröbnerized ideal
of relations for the Degree Reverse Lexicographic ordering which
consists of 105 equations, 66 of which generate the ideal in question
and whose remainders behind a power of $f_1'$ have been filled just
above. 

As a consequence, the full algebra ${\sf E}_5^2$ of jet polynomials 
${\sf P} \big( j^5 f \big)$ invariant by reparametrization
is generated by the polarizations:
\[
\boxed{
\aligned
&\ \ \ \ \ 
f_i',\ \ \ \ \ 
\Lambda^3,\ \ \ \ \
\Lambda_i^5,\ \ \ \ \
\Lambda_{i,j}^7,\ \ \ \ \
\Lambda_{i,j,k}^9,\ \ \ \ \
M^8,\ \ \ \ \ 
M_i^{10},\ \ \ \ \
K_{i,j}^{12},
\\
&
N^{12},\ \ \ \ 
H_i^{14},\ \ \ \ 
F_{i,j}^{16},\ \ \ \ 
X_{i,j,k}^{18},\ \ \ \ 
X_i^{19},\ \ \ \ 
X^{21},\ \ \ \ 
X_i^{23},\ \ \ \ 
X_{i,j}^{25},\ \ \ \ 
Y_{i,j,k}^{27}
\endaligned
}
\]
of these 17 bi-invariants, where the indices $i, j, k$ vary in $\{ 1,
2 \}$, whence the {\bf total number} of these
invariants equals:
\[
2+1+2+4+8+1+2+4+1+2+4+8+2+1+2+4+8
=
\fbox{\bf 56}\,.
\]

}

\stopTHEOREM

\section*{\S11.~Sixteen (fifteen) bi-invariant
\\
in dimension $n = 4$ ($n=3$) for jet level $\kappa = 4$}
\label{Section-11}

\subsection*{ First loop of the algorithm}
Coming back to the end of \S7, we start with the seven initial
bi-invariants:
\[
\aligned
\Lambda^3
&
=
\Delta_{1,2}^{',\,''},
\\
\Lambda^5
&
=
\Delta_{1,2}^{',\,'''}\,f_1'
-
3\,\Delta_{1,2}^{',\,''}\,f_1'',
\\
\Lambda^7
&
=
\Delta_{1,2}^{',\,''''}\,f_1'f_1'
+
\Delta_{1,2}^{'',\,'''}\,f_1'f_1'
-
10\,\Delta_{1,2}^{',\,'''}\,f_1'f_1''
+
15\,\Delta_{1,2}^{',\,''}\,f_1''f_1'',
\endaligned
\]
\[
\aligned
D^6
&
=
\Delta_{1,2,3}^{',\,'',\,'''},
\\
D^8
&
=
\Delta_{1,2,3}^{',\,''',\,''''}\,f_1'
-
3\,\Delta_{1,2,3}^{',\,'',\,''''}\,f_1'',
\\
N^{10}
&
=
\Delta_{1,2,3}^{',\,''',\,''''}\,f_1'f_1'
-
3\,\Delta_{1,2,3}^{',\,'',\,''''}\,f_1'f_1''
+
4\,\Delta_{1,2,3}^{',\,'',\,'''}\,f_1'f_1'''
+
3\,\Delta_{1,2,3}^{',\,'',\,'''}\,f_1''f_1'',
\\
W^{10}
&
=
\Delta_{1,2,3,4}^{',\,'',\,''',\,''''}.
\endaligned
\]
Then we compute the ideal of relations between these
bi-invariants, after setting $f_1 ' = 0$ in them: 
\[
{\sf Ideal-Rel}\,
\Big(
\Lambda^3\big\vert_0,\ \
\Lambda^5\big\vert_0,\ \
\Lambda^7\big\vert_0,\ \
D^6\big\vert_0,\ \
D^8\big\vert_0,\ \
N^{10}\big\vert_0,\ \
W^{10}\big\vert_0
\Big).
\]
We should observe that the first six initial bi-invariants
$\Lambda^3$, $\Lambda^5$, $\Lambda^7$, $D^6$, $D^8$ and $N^{ 10}$
depend only upon the first three jet components $\big( j^4 f_1, \, j^4
f_2, \, j^4 f_3)$ of $j^4 f$, while $W^{ 10}$ and $W^{ 10} \big
\vert_0$\,\,---\,\,which both contain the monomial $-f_4'''' f_3'
f_2'' f_1'''$\,\,---\,\,really depend upon the fourth jet component
$j^4 f_4$. It follows that $W^{ 10} \big \vert_0$ is algebraically
independent of $\Lambda^3 \big \vert_0$, $\Lambda^5 \big \vert_0$,
$\Lambda^7 \big \vert_0$, $D^6 \big \vert_0$, $D^8\big \vert_0$, $N^{
10 }\big \vert_0$, so it cannot intervene in the ideal of
relations. Without loss of generality, we therefore have to consider:
\[
{\sf Ideal-Rel}\,
\Big(
\Lambda^3\big\vert_0,\ \
\Lambda^5\big\vert_0,\ \
\Lambda^7\big\vert_0,\ \
D^6\big\vert_0,\ \
D^8\big\vert_0,\ \
N^{10}\big\vert_0
\Big).
\]
A Maple computation with the Degree Reverse Lexicographic ordering
yields a reduced
Gröbner basis for this ideal consisting of the following 6
generators\footnote{\, {\it See} {\sf
dim-3-order-4-step-1-with-FGb.mw} at \cite{ mer2008c}. }:
\[
\aligned
0
&
\overset{a}{\equiv}
5\,\Lambda^5\Lambda^5
-
3\,\Lambda^3\Lambda^7
\blue{+f_1'f_1'M^8},
\\
0
&
\overset{b}{\equiv}
2\,\Lambda^5D^6
-
\Lambda^3D^8
\blue{
+{\textstyle{\frac{1}{3}}}\,f_1'E^{10}},
\\
0
&
\overset{c}{\equiv}
\Lambda^7D^6
-
5\,\Lambda^3N^{10}
\blue{+f_1'L^{12}},
\\
0
&
\overset{d}{\equiv}
\Lambda^5D^8-6\,\Lambda^3N^{10}
\green{+f_1'L^{12}},
\\
0
&
\overset{e}{\equiv}
\Lambda^7D^8
-
10\,\Lambda^5N^{10}
\blue{-f_1'Q^{14}},
\\
0
&
\overset{f}{\equiv}
D^8D^8
-
12\,D^6N^{10}
\blue{-f_1'R^{15}}.
\endaligned
\]
To read these equations ({\em cf.} \S9), one should at first set $f_1
' = 0$ virtually in one's head and then consider that further
computations show what are the remainders behind a power of $f_1'$.
Five new bi-invariants appear which are implicitly defined by five
among these six sizygies and we provide their explicit expression in
terms of $\Delta$ determinants, after mild simplifications: 
\[
\aligned
M^8
&
:=
\frac{-5\,\Lambda^5\Lambda^5+3\,\Lambda^3\Lambda^7}{f_1'f_1'}
\\
&
=
3\,\Delta_{1,2}^{',\,''''}\,\Delta_{1,2}^{',\,''}
+
12\,\Delta_{1,2}^{'',\,'''}\,\Delta_{1,2}^{',\,''}
-
5\,\Delta_{1,2}^{',\,'''}\,\Delta_{1,2}^{',\,'''},
\endaligned
\]
\[
\aligned
E^{10}
&
:=
\frac{-6\,\Lambda^5\,D^6+3\,\Lambda^3\,D^8}{f_1'}
\\
&
=
3\,\Delta_{1,2,3}^{',\,'',\,''''}\,\Delta_{1,2}^{',\,''}
-
6\,\Delta_{1,2,3}^{',\,'',\,'''}\,\Delta_{1,2}^{',\,'''},
\endaligned
\]
\[
\aligned
L^{12}
&
:=
\frac{-\Lambda^7D^6+5\,\Lambda^3N^{10}}{f_1'}
\\
&
=
-\Delta_{1,2,3}^{',\,'',\,'''}\,\Delta_{1,2}^{',\,''''}\,f_1'
-
4\,\Delta_{1,2,3}^{',\,'',\,'''}\,\Delta_{1,2}^{'',\,'''}\,f_1'
+
5\,\Delta_{1,2,3}^{',\,''',\,''''}\,\Delta_{1,2}^{',\,''}\,f_1'
+
10\,\Delta_{1,2,3}^{',\,'',\,'''}\,\Delta_{1,2}^{',\,'''}\,f_1''
-
\\
&
\ \ \ \ \ \ \ \ \ \ \ \ \ \ \ \ \ \ \ \ \ \ \ \ \ \ \ \ \ \ \ \ \ \ \ \
-
15\,\Delta_{1,2,3}^{',\,'',\,''''}\,\Delta_{1,2}^{',\,''}\,f_1''
+
20\,\Delta_{1,2,3}^{',\,'',\,'''}\,\Delta_{1,2}^{',\,''}\,f_1''',
\endaligned
\]
\[
\aligned
Q^{14}
&
:=
\frac{\Lambda^7D^8-10\,\Lambda^5N^{10}}{f_1'}
\\
&
=
-10\,\Delta_{1,2,3}^{',\,''',\,''''}\,\Delta_{1,2}^{',\,'''}\,f_1'f_1'
+
\Delta_{1,2,3}^{',\,'',\,''''}\,\Delta_{1,2}^{',\,''''}\,f_1'f_1'
+
4\,\Delta_{1,2,3}^{',\,''',\,''''}\,\Delta_{1,2}^{'',\,'''}\,f_1'f_1'
+
\\
&\ \ \ \ \
+20\,\Delta_{1,2,3}^{',\,'',\,''''}\,\Delta_{1,2}^{',\,'''}\,f_1'f_1''
+
30\,\Delta_{1,2,3}^{',\,''',\,''''}\,\Delta_{1,2}^{',\,''}\,f_1'f_1''
-
6\,\Delta_{1,2,3}^{',\,'',\,'''}\,\Delta_{1,2}^{',\,''''}\,f_1'f_1''
-
\\
&\ \ \ \ \
-24\,\Delta_{1,2,3}^{',\,'',\,'''}\,\Delta_{1,2}^{'',\,'''}\,f_1'f_1''
-
40\,\Delta_{1,2,3}^{',\,'',\,'''}\,\Delta_{1,2}^{',\,'''}\,f_1'f_1'''
-
75\,\Delta_{1,2,3}^{',\,'',\,''''}\,\Delta_{1,2}^{',\,''}\,f_1''f_1''
+
\\
&\ \ \ \ \ 
+30\,\Delta_{1,2,3}^{',\,'',\,'''}\,\Delta_{1,2}^{',\,'''}\,f_1''f_1''
+
120\,\Delta_{1,2,3}^{',\,'',\,'''}\,\Delta_{1,2}^{',\,''}\,f_1''f_1''',
\endaligned
\]
\[
\aligned
R^{15}
&
:=
\Delta_{1,2,3}^{',\,'',\,''''}\,\Delta_{1,2,3}^{',\,'',\,''''}\,f_1'
-
12\,\Delta_{1,2,3}^{',\,''',\,''''}\,\Delta_{1,2,3}^{',\,'',\,'''}\,f_1'
+
24\,\Delta_{1,2,3}^{',\,'',\,''''}\,\Delta_{1,2,3}^{',\,'',\,'''}\,f_1''
-
\\
&\ \ \ \ \
-48\,\Delta_{1,2,3}^{',\,'',\,'''}\,
\Delta_{1,2,3}^{',\,'',\,'''}\,f_1''',
\endaligned
\] 
and as usual, the weights are denoted by an upper index. Setting $W^{
10}$ apart, in order to verify that these 11 bi-invariants are
mutually independent, one computes at first which value they have
after setting $f_1 ' = 0$:
\[
\aligned
&
\underline{\Lambda^3}\big\vert_0
\\
&
\underline{\Lambda^5}\big\vert_0
\\
&
\Lambda^7\big\vert_0
=
{\textstyle{\green{\frac{5}{3}}\,
\frac{\Lambda^5\vert_0\,\Lambda^5\vert_0}{
\Lambda^3\vert_0\ \ \ \ \ \ }}}
\endaligned
\]
\[
\aligned
&
\underline{D^6}\big\vert_0
\\
&
D^8\big\vert_0
\ \,
=
{\textstyle{\green{2}\,
\frac{\Lambda^5\vert_0\,D^6\vert_0}{
\Lambda^3\vert_0\ \ \ \ \ \ }}}
\\
&
N^{10}\big\vert_0
=
{\textstyle{\green{\frac{1}{3}}\,
\frac{\Lambda^5\vert_0\,\Lambda^5\vert_0\,D^6\vert_0}{
\Lambda^3\vert_0\,\Lambda^3\vert_0\ \ \ \ \ \ }}}
\endaligned
\]
\[
\aligned
&
\underline{M^8}\big\vert_0
\\
&
\underline{E^{10}}\big\vert_0
\\
&
L^{12}\big\vert_0
=
{\textstyle{\green{\frac{5}{3}}\,
\frac{\Lambda^5\vert_0\,E^{10}\vert_0}{
\Lambda^3\vert_0\ \ \ \ \ \ \ }}}
\endaligned
\]
\[
\aligned
&
Q^{14}\big\vert_0
=
{\textstyle{\green{-\frac{25}{9}}\,
\frac{\Lambda^5\vert_0\,\Lambda^5\vert_0\,E^{10}\vert_0}{
\Lambda^3\vert_0\Lambda^3\vert_0\ \ \ \ \ \ \ \ }}}
\\
&
R^{15}\big\vert_0
=
{\textstyle{\green{-\frac{8}{3}}\,
\frac{\Lambda^5\vert_0\,D^6\vert_0\,E^{10}\vert_0}{
\Lambda^3\vert_0\,\Lambda^3\vert_0\ \ \ \ \ \ \ \, }}},
\endaligned
\]
with the 5 underlined bi-invariants being algebraically independent
and being considered as a transcendence basis, while the value of
$\Lambda^7 \big\vert_0$ comes from ``$\overset{ a}{ \equiv}$'' above;
the value of $D^8 \big \vert_0$ comes from ``$\overset{ b}{ \equiv}$''
above; the value of $N^{ 10} \big\vert_0$ comes from ``$\overset{ d}{
\equiv}$'' above; the value of $L^{ 12} \big\vert_0$ comes from
``$\overset{ r}{ \equiv}$'' below; the value of $Q^{ 14} \big \vert_0$
comes from ``$\overset{ q}{ \equiv}$'' below; and the value of $R^{
15} \big\vert_0$ comes from ``$\overset{ p}{ \equiv}$'' below. Then
one proceeds as in the proof of the lemma on
p.~\pageref{lemma-independence} to show mutual independence (details
will not be provided).

Importantly, the five new bi-invariants $M^8$, $E^{ 10}$, $L^{ 12}$,
$Q^{ 14}$ and $R^{ 15}$ again depend only upon the first three jet
components $\big( j^4 f_1, j^4 f_2, j^4 f_3\big)$, so that $W^{ 10}
\big\vert_0$ again will not intervene in the next ideal of
relations. In fact, all bi-invariants except $W^{ 10}$ 
live in dimension $n = 3$, and hence it
is enough to explore the structure of
${\sf UE}_4^3$.

\subsection*{Second loop of the algorithm}
Setting therefore $W^{ 10}$ apart, a Maple computation with the Degree
Reverse Lexicographic Ordering offers a reduced Gröbner basis for the
ideal of relations: 
\[
{\sf Ideal-Rel}\,
\bigg(
\aligned
&
\Lambda^3\big\vert_0,\ \
\Lambda^5\big\vert_0,\ \
\Lambda^7\big\vert_0,\ \
D^6\big\vert_0,\ \
D^8\big\vert_0,\ \
N^{10}\big\vert_0,\ \
\\
&\ \ \ \ 
M^8\big\vert_0,\ \
E^{10}\big\vert_0,\ \
L^{12}\big\vert_0,\ \
Q^{14}\big\vert_0,\ \
R^{15}\big\vert_0
\endaligned
\bigg)
\]
between our 11 bi-invariants restricted to $\{ f_1' = 0\}$, and this
basis consists of the 6 generators above together with the following
14 generators
\footnote{\, {\it See} {\sf dim-3-order-4-step-2-with-FGb.mw} at
\cite{ mer2008c}. Here again, the remainders behind a power
of $f_1'$ have all been computed and tested to know whether 
they belong to the algebra of the already known 11 bi-invariants. }:
\[
\aligned
0
&
\overset{g}{\equiv}
4\,D^8Q^{14}
-
5\,\Lambda^7R^{15}
\blue{-f_1'X^{21}},
\\
0
&
\overset{h}{\equiv}
24\,D^6Q^{14}
-
25\,\Lambda^5R^{15}
\green{+f_1'V^{19}},
\\
0
&
\overset{i}{\equiv}
L^{12}L^{12}
+
E^{10}Q^{14}
\green{-f_1'M^8R^{15}},
\\
0
&
\overset{j}{\equiv}
8\,N^{10}L^{12}
+
\Lambda^7R^{15}
\green{+f_1'X^{21}},
\\
0
&
\overset{k}{\equiv}
4\,D^8L^{12}
+
5\,\Lambda^5R^{15}
\green{-f_1'V^{19}},
\endaligned
\]
\[
\aligned
0
&
\overset{l}{\equiv}
8\,D^6L^{12}
+
5\,\Lambda^3R^{15}
\green{-{\textstyle{\frac{1}{3}}}\,f_1'U^{17}},
\\
0
&
\overset{m}{\equiv}
\Lambda^7L^{12}
+
\Lambda^5Q^{14}
\green{-2\,f_1'M^8N^{10}},
\\
0
&
\overset{n}{\equiv}
5\,\Lambda^5L^{12}
+
3\,\Lambda^3Q^{14}
\green{-f_1'D^8M^8},
\\
0
&
\overset{o}{\equiv}
8\,N^{10}E^{10}
+
\Lambda^5\,R^{15}
\blue{-f_1'V^{19}},
\\
0
&
\overset{p}{\equiv}
4\,D^8E^{10}
+
3\,\Lambda^3R^{15}
\blue{-f_1'U^{17}},
\endaligned
\]
\[
\aligned
0
&
\overset{q}{\equiv}
5\,\Lambda^7E^{10}
+
3\,\Lambda^3Q^{14}
\green{-6f_1'D^8M^8},
\\
0
&
\overset{r}{\equiv}
5\,\Lambda^5E^{10}
-
3\,\Lambda^3L^{12}
\green{-6\,f_1'D^6M^8},
\\
0
&
\overset{s}{\equiv}
8\,\Lambda^5N^{10}Q^{14}
-
\Lambda^7\Lambda^7R^{15}
\green{+f_1'Q^{14}Q^{14}+4\,f_1'N^{10}N^{10}M^8},
\\
0
&
\overset{t}{\equiv}
24\,\Lambda^3N^{10}Q^{14}
-
5\,\Lambda^5\Lambda^7R^{15}
\green{-5\,f_1'L^{12}Q^{14}+2\,f_1'M^8D^8N^{10}}.
\endaligned
\]
Here, three new bi-invariants appear: $U^{ 17}$, $V^{ 19}$ and $X^{
21}$, which are implicitly defined by the syzygies ``$\overset{
p}{\equiv}$'', ``$\overset{ o}{\equiv}$'', and ``$\overset{
g}{\equiv}$'', and we provide their explicit expression in terms of
$\Delta$ determinants\footnote{\, To be able do divide by $f_1'$, as
in~\cite{ mer2008a}, we sometimes need to replace $\Delta_{ 1, 2}^{ ',
\, '' }\, f_1'''$ by $-\Delta_{ 1, 2}^{ '', \, '''}\, f_1' + \Delta_{
1, 2}^{ ', \, '''}\, f_1''$, using the immediately checked syzygy: $0
\equiv \Delta_{ 1, 2}^{ '', \, '''}\, f_1' - \Delta_{ 1, 2}^{ ', \,
'''}\, f_1'' + \Delta_{ 1, 2}^{ ', \, '' }\, f_1'''$.
}:
\[
\aligned
U^{17}
&
=
\frac{4\,D^8E^{10}+3\,\Lambda^3R^{15}}{f_1'}
\\
&
=
15\,\Delta_{1,2,3}^{',\,'',\,''''}\,\Delta_{1,2,3}^{',\,'',\,''''}\,
\Delta_{1,2}^{',\,''}
-
36\,\Delta_{1,2,3}^{',\,''',\,''''}\,\Delta_{1,2,3}^{',\,'',\,'''}\,
\Delta_{1,2}^{',\,''}
-
\\
&\ \ \ \ \
-24\,\Delta_{1,2,3}^{',\,'',\,''''}\,
\Delta_{1,2,3}^{',\,'',\,'''}\,\Delta_{1,2}^{',\,'''}
+
144\,\Delta_{1,2,3}^{',\,'',\,'''}\,\Delta_{1,2,3}^{',\,'',\,'''}\,
\Delta_{1,2}^{'',\,'''},
\endaligned
\]
\[
\aligned
V^{19}
&
=
\frac{8\,N^{10}E^{10}+\Lambda^5R^{15}}{f_1'}
\\
&
=
24\,\Delta_{1,2,3}^{',\,''',\,''''}\,\Delta_{1,2,3}^{',\,'',\,''''}\,
\Delta_{1,2}^{',\,''}\,f_1'
-
60\,\Delta_{1,2,3}^{',\,''',\,''''}\,\Delta_{1,2,3}^{',\,'',\,'''}\,
\Delta_{1,2}^{',\,'''}\,f_1'
+
\\
&\ \ \ \ \
+\Delta_{1,2,3}^{',\,'',\,''''}\,\Delta_{1,2,3}^{',\,'',\,''''}\,
\Delta_{1,2}^{',\,'''}\,f_1'
-
75\,\Delta_{1,2,3}^{',\,'',\,''''}\,\Delta_{1,2,3}^{',\,'',\,''''}\,
\Delta_{1,2}^{',\,''}\,f_1''
+
\\
&\ \ \ \ \
+36\,\Delta_{1,2,3}^{',\,''',\,''''}\,\Delta_{1,2,3}^{',\,'',\,'''}\,
\Delta_{1,2}^{',\,''}\,f_1''
+
168\,\Delta_{1,2,3}^{',\,'',\,''''}\,\Delta_{1,2,3}^{',\,'',\,'''}\,
\Delta_{1,2}^{',\,'''}\,f_1''
-
\\
&\ \ \ \ \
-144\,\Delta_{1,2,3}^{',\,'',\,'''}\,\Delta_{1,2,3}^{',\,'',\,'''}\,
\Delta_{1,2}^{'',\,'''}\,f_1''
+
96\,\Delta_{1,2,3}^{',\,'',\,''''}\,\Delta_{1,2,3}^{',\,'',\,'''}\,
\Delta_{1,2}^{',\,''}\,f_1'''
-
\\
&\ \ \ \ \
-
240\,\Delta_{1,2,3}^{',\,'',\,'''}\,\Delta_{1,2,3}^{',\,'',\,'''}\,
\Delta_{1,2}^{',\,'''}\,f_1''',
\endaligned
\]
\[
\aligned
X^{21}
&
=
\frac{4\,D^8Q^{14}-5\,\Lambda^7R^{15}}{f_1'}
\\
&
=
-40\,\Delta_{1,2,3}^{',\,''',\,''''}\,\Delta_{1,2,3}^{',\,'',\,''''}\,
\Delta_{1,2}^{',\,'''}\,f_1'f_1'
-
4\,\Delta_{1,2,3}^{',\,'',\,''''}\,\Delta_{1,2,3}^{',\,'',\,''''}\,
\Delta_{1,2}^{',\,''''}\,f_1'f_1'
-
\\
&\ \ \ \ \
-4\,\Delta_{1,2,3}^{',\,'',\,''''}\,\Delta_{1,2,3}^{',\,'',\,''''}\,
\Delta_{1,2}^{'',\,'''}\,f_1'f_1'
+
60\,\Delta_{1,2,3}^{',\,''',\,''''}\,\Delta_{1,2,3}^{',\,'',\,'''}\,
\Delta_{1,2}^{',\,''''}\,f_1'f_1'
+
\\
&\ \ \ \ \
+240\,\Delta_{1,2,3}^{',\,''',\,''''}\,\Delta_{1,2,3}^{',\,'',\,'''}\,
\Delta_{1,2}^{'',\,'''}\,f_1'f_1'
+
130\,\Delta_{1,2,3}^{',\,'',\,''''}\,\Delta_{1,2,3}^{',\,'',\,''''}\,
\Delta_{1,2}^{',\,'''}\,f_1'f_1''
+
\\
&\ \ \ \ \
+120\,\Delta_{1,2,3}^{',\,''',\,''''}\,\Delta_{1,2,3}^{',\,'',\,''''}\,
\Delta_{1,2}^{',\,''}\,f_1'f_1''
-
168\,\Delta_{1,2,3}^{',\,'',\,''''}\,\Delta_{1,2,3}^{',\,'',\,'''}\,
\Delta_{1,2}^{',\,''''}\,f_1'f_1''
-
\endaligned
\]
\[
\aligned
&\ \ \ \ \
-668\,\Delta_{1,2,3}^{',\,'',\,''''}\,\Delta_{1,2,3}^{',\,'',\,'''}\,
\Delta_{1,2}^{'',\,'''}\,f_1'f_1''
-
360\,\Delta_{1,2,3}^{',\,''',\,''''}\,\Delta_{1,2,3}^{',\,'',\,'''}\,
\Delta_{1,2}^{',\,'''}\,f_1'f_1''
-
\\
&\ \ \ \ \
-160\,\Delta_{1,2,3}^{',\,'',\,''''}\,\Delta_{1,2,3}^{',\,'',\,'''}\,
\Delta_{1,2}^{',\,'''}\,f_1'f_1'''
+
240\,\Delta_{1,2,3}^{',\,'',\,'''}\,\Delta_{1,2,3}^{',\,'',\,'''}\,
\Delta_{1,2}^{',\,''''}\,f_1'f_1'''
+
\\
&\ \ \ \ \
+960\,\Delta_{1,2,3}^{',\,'',\,'''}\,\Delta_{1,2,3}^{',\,'',\,'''}\,
\Delta_{1,2}^{'',\,'''}\,f_1'f_1'''
-
375\,\Delta_{1,2,3}^{',\,'',\,''''}\,\Delta_{1,2,3}^{',\,'',\,''''}\,
\Delta_{1,2}^{',\,''}\,f_1''f_1''
+
\\
&\ \ \ \ \
+840\,\Delta_{1,2,3}^{',\,'',\,''''}\,\Delta_{1,2,3}^{',\,'',\,'''}\,
\Delta_{1,2}^{',\,'''}\,f_1''f_1''
+
180\,\Delta_{1,2,3}^{',\,''',\,''''}\,\Delta_{1,2,3}^{',\,'',\,'''}\,
\Delta_{1,2}^{',\,''}\,f_1''f_1''
+
\\
&\ \ \ \ \
+
144\,\Delta_{1,2,3}^{',\,'',\,'''}\,\Delta_{1,2,3}^{',\,'',\,'''}\,
\Delta_{1,2}^{',\,''''}\,f_1''f_1''
+
144\,\Delta_{1,2,3}^{',\,'',\,'''}\,\Delta_{1,2,3}^{',\,'',\,'''}\,
\Delta_{1,2}^{'',\,'''}\,f_1''f_1''
-
\\
&\ \ \ \ \
-1440\,\Delta_{1,2,3}^{',\,'',\,'''}\,\Delta_{1,2,3}^{',\,'',\,'''}\,
\Delta_{1,2}^{',\,'''}\,f_1''f_1'''
+
480\,\Delta_{1,2,3}^{',\,'',\,''''}\,\Delta_{1,2,3}^{',\,'',\,'''}\,
\Delta_{1,2}^{',\,''}\,f_1''f_1'''.
\endaligned
\]
Either a Maple computation or a glance at the syzygies ``$\overset{
7}{ \equiv}$'', ``$\overset{ 8}{ \equiv}$'', ``$\overset{ 9}{
\equiv}$'' below arriving in the third loop provides the values of
these two bi-invariants after setting $f_1' = 0$:
\[
\aligned
U^{17}\big\vert_0
&
=
{\textstyle{\green{12}}\,
\frac{D^6\vert_0\,D^6\vert_0\,M^8\vert_0}{\Lambda^3\vert_0
\ \ \ \ \ \ \ \ \ \ \ \ }
\green{+\frac{5}{3}}\,
\frac{E^{10}\vert_0\,E^{10}\vert_0}{\Lambda^3\vert_0
\ \ \ \ \ \ \ \ }},
\\
V^{19}\big\vert_0
&
=
{\textstyle{\green{\frac{25}{9}}\,
\frac{\Lambda^5\vert_0\,E^{10}\vert_0\,E^{10}\vert_0}{
\Lambda^3\vert_0\,\Lambda^3\vert_0\ \ \ \ \ \ \ \ }
\green{+4}\,
\frac{\Lambda^5\vert_0\,D^6\vert_0\,D^6\vert_0\,M^8\vert_0}{
\Lambda^3\vert_0\,\Lambda^3\vert_0\ \ \ \ \ \ \ \ \ \ \ \ \ \ }}},
\\
X^{21}\big\vert_0
&
=
{\textstyle{\green{-\frac{4}{3}}\,
\frac{\Lambda^5\vert_0\,\Lambda^5\vert_0\,
D^6\vert_0\,D^6\vert_0\,M^8\vert_0}{
\Lambda^3\vert_0\,\Lambda^3\vert_0\,\Lambda^3\vert_0
\ \ \ \ \ \ \ \ \ \ \ \ \ \ }
\green{-\frac{125}{27}}\,
\frac{\Lambda^5\vert_0\,\Lambda^5\vert_0\,E^{10}\vert_0\,E^{10}\vert_0}{
\Lambda^3\vert_0\,\Lambda^3\vert_0\,\Lambda^3\vert_0
\ \ \ \ \ \ \ \ \ }}}.
\endaligned
\]
Proceeding as in the lemma on p.~\pageref{lemma-independence}, one
checks patiently by hand that the 16 bi-invariants known so far:
\[
\boxed{
\aligned
&\ \ \ \
W^{10},\ \ \ \ \
f_1',\ \ \ \ \
\Lambda^3,\ \ \ \ \
\Lambda^5,\ \ \ \ \
\Lambda^7,\ \ \ \ \
D^6,\ \ \ \ \
D^8,\ \ \ \ \
N^{10},\ \ \ \ \ 
\\
&
M^8,\ \ \ \ \
E^{10},\ \ \ \ \
L^{12},\ \ \ \ \
Q^{14},\ \ \ \ \
R^{15},\ \ \ \ \
U^{17},\ \ \ \ \
V^{19},\ \ \ \ \
X^{21}
\endaligned}
\]
are mutually independent.

\subsection*{ Third loop of the algorithm} 
Again for the Degree Reverse Lexicographic ordering, setting $W^{ 10}$
apart, a Maple computation offers a reduced Gröbner basis for the
ideal of relations between the $14 = 15 - 1$ ($f_1'$ goes to zero)
restricted bi-invariants. The result consists of 50 generators
\footnote{\,{\it See} {\sf dim-3-order-4-step-3-with-FGb.mws} at
\cite{ mer2008c}. }. Taking the Lexicographic ordering instead:
\[
\small
\aligned
\Lambda^3\,>\,\Lambda^5\,>\,\Lambda^7\,>\,
D^6\,>\,D^8\,>\,
&
N^{10}\,>\,M^8\,>\,E^{10}\,>\,L^{12}\,>\,
\\
&
\,>\,Q^{14}\,>\,R^{15}\,>\,U^{17}\,>\,V^{19}\,>\,X^{21},
\endaligned
\]
one shows that the ideal of relations, in Gröbnerized form, contains
less equations\,\,---\,\,which is convenient\,\,---, namely the
following 41 equations\footnote{\, {\em See} {\sf
41-syzygies-dim-3-order-4.mw} at \cite{ mer2008c}.}, where we
underline the {\sf L}eading {\sf T}erm of each syzygy with the acronym
``${\sf LT}$'' appended\footnote{\, We recall that, in order to
appropriately read the ideal of relations between restricted
bi-invariants, one should set $f_1' = 0$, namely disregard the last
term(s) of each equation. We specify ``$\red{ +0}$'' when the
remainder behing a power of $f_1'$ vanishes identically.}:

\THEOREM

\[
\label{41-syzygies}
\aligned
0
&
\overset{1}{\equiv}
-5\,\Lambda^5\Lambda^5
+
3\,\underline{\Lambda^3\Lambda^7}_{{\scriptscriptstyle{\sf LT}}}
\green{-f_1'f_1'M^8},
\\
0
&
\overset{2}{\equiv}
-2\,\Lambda^5D^6+
\underline{\Lambda^3D^8}_{{\scriptscriptstyle{\sf LT}}}
\green{-{\textstyle{\frac{1}{3}}}f_1'\,E^{10}},
\\
0
&
\overset{3}{\equiv}
-\Lambda^7D^6
+
5\,\underline{\Lambda^3N^{10}}_{{\scriptscriptstyle{\sf LT}}}
\green{-f_1'L^{12}},
\\
0
&
\overset{4}{\equiv}
-5\,\Lambda^5E^{10}+3\,
\underline{\Lambda^3L^{12}}_{{\scriptscriptstyle{\sf LT}}}
\green{+6\,f_1'D^6M^8},
\\
0
&
\overset{5}{\equiv}
5\,\Lambda^7E^{10}
+
3\,\underline{\Lambda^3Q^{14}}_{{\scriptscriptstyle{\sf LT}}}
\green{-6\,f_1'D^8M^8},
\\
0
&
\overset{6}{\equiv}
4\,D^8E^{10}
+
3\,\underline{\Lambda^3R^{15}}_{{\scriptscriptstyle{\sf LT}}}
\green{-f_1'U^{17}},
\endaligned
\]
\[
\aligned
0
&
\overset{7}{\equiv}
-36\,D^6D^6M^8
-
5\,E^{10}E^{10}
+
3\,\underline{\Lambda^3U^{17}}_{{\scriptscriptstyle{\sf LT}}}
\red{+0},
\\
0
&
\overset{8}{\equiv}
-5\,E^{10}L^{12}
-
6\,D^6D^8M^8
+
3\,\underline{\Lambda^3V^{19}}_{{\scriptscriptstyle{\sf LT}}}
\red{+0},
\\
0
&
\overset{9}{\equiv}
5\,L^{12}L^{12}
+
3\,\underline{\Lambda^3X^{21}}_{{\scriptscriptstyle{\sf LT}}}
+
M^8D^8D^8
\red{+0},
\\
0
&
\overset{10}{\equiv}
-6\,\Lambda^7D^6
+
5\,\underline{\Lambda^5D^8}_{{\scriptscriptstyle{\sf LT}}}
\green{-f_1'L^{12}},
\\
0
&
\overset{11}{\equiv}
-\Lambda^7D^8
+
10\,\underline{\Lambda^5N^{10}}_{{\scriptscriptstyle{\sf LT}}}
\green{+f_1'Q^{14}},
\\
0
&
\overset{12}{\equiv}
\underline{\Lambda^5L^{12}}_{{\scriptscriptstyle{\sf LT}}}
-
\Lambda^7E^{10}
\green{+f_1'D^8M^8},
\endaligned
\]
\[
\aligned
0
&
\overset{13}{\equiv}
\Lambda^7L^{12}
+
\underline{\Lambda^5Q^{14}}_{{\scriptscriptstyle{\sf LT}}}
\green{-2\,f_1'M^8N^{10}},
\\
0
&
\overset{14}{\equiv}
8\,N^{10}E^{10}
+
\underline{\Lambda^5R^{15}}_{{\scriptscriptstyle{\sf LT}}}
\green{-f_1'V^{19}},
\\
0
&
\overset{15}{\equiv}
\underline{\Lambda^5U^{17}}_{{\scriptscriptstyle{\sf LT}}}
-
E^{10}L^{12}
-
6\,D^6D^8M^8
\red{+0},
\\
0
&
\overset{16}{\equiv}
\underline{\Lambda^5V^{19}}_{{\scriptscriptstyle{\sf LT}}}
-
M^8D^8D^8
-
L^{12}L^{12}
\green{+f_1'M^8R^{15}},
\\
0
&
\overset{17}{\equiv}
\underline{\Lambda^5X^{21}}_{{\scriptscriptstyle{\sf LT}}}
-
L^{12}Q^{14}
+
2\,D^8N^{10}M^8
\red{+0},
\\
0
&
\overset{18}{\equiv}
8\,N^{10}L^{12}
+
\underline{\Lambda^7R^{15}}_{{\scriptscriptstyle{\sf LT}}}
\green{+f_1'X^{21}},
\endaligned
\]
\[
\aligned
0
&
\overset{19}{\equiv}
-L^{12}L^{12}
+
\underline{\Lambda^7U^{17}}_{{\scriptscriptstyle{\sf LT}}}
\green{+f_1'},
-
5\,M^8D^8D^8
\red{+0},
\\
0
&
\overset{20}{\equiv}
L^{12}Q^{14}
+
\underline{\Lambda^7V^{19}}_{{\scriptscriptstyle{\sf LT}}}
-
10\,D^8M^8N^{10}
\red{+0},
\endaligned
\]
\[
\aligned
0
&
\overset{21}{\equiv}
20\,N^{10}N^{10}M^8
+
Q^{14}Q^{14}
+
\underline{\Lambda^7X^{21}}_{{\scriptscriptstyle{\sf LT}}}
\red{+0},
\\
0
&
\overset{22}{\equiv}
6\,\underline{D^6M^8R^{15}}_{{\scriptscriptstyle{\sf LT}}}
+
L^{12}U^{17}
-
E^{10}V^{19}
\red{+0},
\\
0
&
\overset{23}{\equiv}
5\,\underline{D^8M^8R^{15}}_{{\scriptscriptstyle{\sf LT}}}
-
Q^{14}U^{17}
-
L^{12}V^{19}
\red{+0},
\\
0
&
\overset{24}{\equiv}
10\,\underline{N^{10}M^8R^{15}}_{{\scriptscriptstyle{\sf LT}}}
-
Q^{14}V^{19}
+
L^{12}X^{21}
\red{+0},
\endaligned
\]
\[
\aligned
0
&
\overset{25}{\equiv}
5\,\underline{M^8R^{15}R^{15}}_{{\scriptscriptstyle{\sf LT}}}
+
V^{19}V^{19}
+
U^{17}X^{21}
\red{+0},
\\
0
&
\overset{26}{\equiv}
-D^8D^8
+
12\,\underline{D^6N^{10}}_{{\scriptscriptstyle{\sf LT}}}
\green{+f_1'R^{15}},
\\
0
&
\overset{27}{\equiv}
-5\,D^8E^{10}
+
6\,\underline{D^6L^{12}}_{{\scriptscriptstyle{\sf LT}}}
\green{+f_1'U^{17}},
\\
0
&
\overset{28}{\equiv}
3\,\underline{D^6Q^{14}}_{{\scriptscriptstyle{\sf LT}}}
+
25\,N^{10}E^{10}
\green{-3\,f_1'V^{19}},
\\
0
&
\overset{29}{\equiv}
5\,E^{10}R^{15}
-
D^8U^{17}
+
6\,\underline{D^6V^{19}}_{{\scriptscriptstyle{\sf LT}}}
\red{+0},
\\
0
&
\overset{30}{\equiv}
-3\,L^{12}R^{15}
+
N^{10}U^{17}
+
3\,\underline{D^6X^{21}}_{{\scriptscriptstyle{\sf LT}}}
\red{+0},
\endaligned
\]
\[
\aligned
0
&
\overset{31}{\equiv}
-10\,N^{10}E^{10}
+
\underline{D^8L^{12}}_{{\scriptscriptstyle{\sf LT}}}
\green{+f_1'V^{19}},
\\
0
&
\overset{32}{\equiv}
\underline{D^8Q^{14}}_{{\scriptscriptstyle{\sf LT}}}
+
10\,N^{10}L^{12}
\green{+f_1'X^{21}},
\\
0
&
\overset{33}{\equiv}
-2\,N^{10}U^{17}
+
\underline{D^8V^{19}}_{{\scriptscriptstyle{\sf LT}}}
+
L^{12}R^{15}
\red{+0},
\\
0
&
\overset{34}{\equiv}
Q^{14}R^{15}
+
2\,N^{10}V^{19}
+
\underline{D^8X^{21}}_{{\scriptscriptstyle{\sf LT}}}
\red{+0},
\\
0
&
\overset{35}{\equiv}
-2\,L^{12}N^{10}U^{17}
+
R^{15}L^{12}L^{12}
+
10\,\underline{V^{19}N^{10}E^{10}}_{{\scriptscriptstyle{\sf LT}}}
\green{-f_1'V^{19}V^{19}},
\\
0
&
\overset{36}{\equiv}
2\,N^{10}U^{17}Q^{14}
-
R^{15}L^{12}Q^{14}
+
10\,\underline{V^{19}N^{10}L^{12}}_{{\scriptscriptstyle{\sf LT}}}
\green{+f_1'V^{19}X^{21}},
\endaligned
\]
\[
\aligned
0
&
\overset{37}{\equiv}
10\,\underline{N^{10}L^{12}X^{21}}_{{\scriptscriptstyle{\sf LT}}}
-
R^{15}Q^{14}Q^{14}
-
2\,Q^{14}N^{10}V^{19}
\green{+f_1'X^{21}X^{21}},
\\
0
&
\overset{38}{\equiv}
2\,\underline{N^{10}U^{17}X^{21}}_{{\scriptscriptstyle{\sf LT}}}
-
X^{21}L^{12}R^{15}
+
V^{19}Q^{14}R^{15}
+
2\,N^{10}V^{19}V^{19}
\red{+0},
\\
0
&
\overset{39}{\equiv}
\underline{E^{10}Q^{14}}_{{\scriptscriptstyle{\sf LT}}}
+
L^{12}L^{12}
\green{-f_1'M^8R^{15}},
\\
0
&
\overset{40}{\equiv}
Q^{14}U^{17}
+
6\,L^{12}V^{19}
+
5\,\underline{E^{10}X^{21}}_{{\scriptscriptstyle{\sf LT}}}
\red{+0},
\\
0
&
\overset{41}{\equiv}
-6\,Q^{14}L^{12}V^{19}
-
Q^{14}Q^{14}U^{17}
+
5\,\underline{X^{21}L^{12}L^{12}}_{{\scriptscriptstyle{\sf LT}}}
\green{-5\,f_1'M^8R^{15}X^{21}}.
\endaligned
\]

\stopTHEOREM

Remarkably, each one of the 41 remainders behind a power of $f_1'$
belongs to the algebra of already known bi-invariants. No new
bi-invariant appears at this stage. In such a circumstance, 
according to the general theorem on
p.~\pageref{normal-syzygies}, we know
that our algorithm stops, so that we have gained the following
complete, quite nontrivial result.

\THEOREM

\smallskip\noindent\fbox{\bf THEOREM}\ \ 
\label{bi-invariant-4-4}
{\sf\em
In dimension $n = 4$ for jets of order $\kappa = 4$, the algebra ${\sf
UE}_4^4$ of jet polynomials ${\sf P} \big( j^4 f_1, j^4 f_2, j^4 f_3,
j^4f_4\big)$ invariant by reparametrization and invariant under the
unipotent action is generated by the 16 mutually independent
bi-invariants defined above:
\[
\boxed{
\aligned
&\ \ \ \
W^{10},\ \ \ \ \
f_1',\ \ \ \ \
\Lambda^3,\ \ \ \ \
\Lambda^5,\ \ \ \ \
\Lambda^7,\ \ \ \ \
D^6,\ \ \ \ \
D^8,\ \ \ \ \
N^{10},\ \ \ \ \ 
\\
&
M^8,\ \ \ \ \
E^{10},\ \ \ \ \
L^{12},\ \ \ \ \
Q^{14},\ \ \ \ \
R^{15},\ \ \ \ \
U^{17},\ \ \ \ \
V^{19},\ \ \ \ \
X^{21},
\endaligned}
\]
whose restriction to $\{ f_1' = 0 \}$ has a reduced gröbnerized ideal
of relations, for the Lexicographic ordering, which consists of the 41
syzygies written above.

Furthermore, any bi-invariant of weight $m$ writes uniquely in the
finite polynomial form:
\[
\aligned
{\sf P}\big(j^\kappa f\big)
=
&
\sum_{o,\,p}\,
(f_1')^o\,\big(W^{10}\big)^p\,
\sum_{(a,\dots,n)\in\N^{14}\backslash
(\square_1\cup\cdots\cup\square_{41})
\atop
3a+\cdots+21n=m-o-10p}\,
{\sf coeff}_{a,\dots,n,o,p}\,\cdot
\\
&\ \ \ \ \
\cdot
\big(\Lambda^3\big)^a\,
\big(\Lambda^5\big)^b\,
\big(\Lambda^7\big)^c\,
\big(D^6\big)^d\,
\big(D^8\big)^e\,
\big(N^{10}\big)^f
\big(M^8\big)^g\,
\big(E^{10}\big)^h\,
\\
&\ \ \ \ \ \ \
\big(L^{12}\big)^i\,
\big(Q^{14}\big)^j\,
\big(R^{15}\big)^k\,
\big(U^{17}\big)^l\,
\big(V^{19}\big)^m\,
\big(X^{21}\big)^n,
\endaligned
\]
with coefficients ${\sf coeff}_{a,\dots,n,o,p}$ subjected to no
restriction, where $\square_1$, \dots, $\square_{ 41}$ denote the
quadrants in $\N^{ 14}$ having vertex at the leading terms of the 41
syzygies in question.

Finally, in the preceding dimension $n = 3$ for jets of the same order
$\kappa = 4$, the algebra ${\sf UE}_4^3$ is generated by the same list
from which one removes only the four-dimensional Wronskian $W^{ 10}$,
the ideal of relations for the 15 restricted bi-invariants being
exactly the same, with an entirely similar normal form for
a general bi-invariant of weight $m$. 

}\medskip

\stopTHEOREM\medskip

As a consequence, by looking at the ${\sf GL}_4 ( \C)$-orbit of each
one of these 16 bi-invariants, we deduce a system of {\bf 2835}
generators for the algebra ${\sf E}_4^4$ of polynomials which are
invariant (only) by reparametrization.

\THEOREM

\smallskip\noindent\fbox{\bf THEOREM}\ \
{\sf\em
In dimension $n = 4$ for jets of order $\kappa = 4$, the algebra ${\sf
E}_4^4$ of jet polynomials ${\sf P} \big( j^4 f \big)$ invariant by
reparametrization is generated by the polarizations: 
\[
\boxed{
\aligned
&\ \ \ \ \ \ \ \
W^{10},\ \ \ \ \
f_i',\ \ \ \ \
\Lambda_{[i,j]}^3,\ \ \ \ \
\Lambda_{[i,j];\,\alpha}^5,\ \ \ \ \
\Lambda_{[i,j];\,\alpha,\beta}^7,\ \ \ \ \
D_{[i,j,k]}^6,
\\
&
D_{[i,j,k];\,\alpha}^8,\ \ \ \ \
N_{[i,j,k];\,\alpha,\beta}^{10},\ \ \ \ \ 
M_{[i,j],[k,l]}^8,\ \ \ \ \
E_{[i,j,k],[p,q]}^{10},\ \ \ \ \
L_{[i,j,k],[p,q];\,\alpha}^{12},
\\
&\ \ \ \ \ \ \ \ \ \ \
Q_{[i,j,k],[p,q];\,\alpha,\beta}^{14},\ \ \ \ \
R_{[i,j,k],[p,q,r];\,\alpha}^{15},\ \ \ \ \
U_{[i,j,k],[p,q,r],[s,t]}^{17},
\\
&\ \ \ \ \ \ \ \ \ \ \ \ \ \ \ \ \ \ \ 
V_{[i,j,k],[p,q,r],[s,t];\,\alpha}^{19},\ \ \ \ \
X_{[i,j,k],[p,q,r],[s,t];\,\alpha,\beta}^{21},
\endaligned}
\]
of the 16 bi-invariants $W^{ 10}$, $f_1'$, $\Lambda^3$, $\Lambda^5$,
$\Lambda^7$, $D^6$, $D^8$, $N^{ 10}$, $M^8$, $E^{ 10}$, $L^{ 12}$,
$Q^{ 14}$, $R^{ 15}$, $U^{ 17}$, $V^{ 19}$, $X^{ 21}$ generating 
the algebra ${\sf
UE}_4^4$ of bi-invariants; these polarized invariants are skew-symmetric
with respect to each collection of bracketed indices $[i,j,k]$,
$[p,q,r]$, $[s,t]$, and they are {\em explicitly} represented in terms
of $\Delta$-determinants by the following complete formulas:
\[
\aligned
W_{1,2,3,4}^{10}
&,
\\
f_i'
&,
\\
\Lambda_{[i,j]}^3
&
:=
\Delta_{i,j}^{',\,''},
\endaligned
\]
\[
\aligned
\\
\Lambda_{[i,j];\,\alpha}^5
&
:=
\Delta_{i,j}^{',\,'''}\,f_\alpha'
-
3\,\Delta_{i,j}^{',\,''}\,f_\alpha'',
\endaligned
\]
\[
\aligned
\Lambda_{[i,j];\,\alpha,\beta}^7
&
:=
\Delta_{i,j}^{',\,''''}\,f_\alpha' f_\beta'
+
4\,\Delta_{i,j}^{'',\,'''}\,f_\alpha'f_\beta'
-
5\Delta_{i,j}^{',\,'''}\,
\big(f_\alpha'f_\beta''+f_\alpha''f_\beta'\big)
+
\\
&\ \ \ \ \ \ \ \ \ \ \ \ \ \ \ \ \ \ \ \ \ \ \
\ \ \ \ \ \ \ \ \ \ \ \ \ \ \ \ \ \ \ \ \ \ \ \ \ \ \ \ 
\ \ \ \ \ \ \ \
+
15\,\Delta_{i,j}^{',\,''}\,f_\alpha''f_\beta'',
\\
D_{[i,j,k]}^6
&
:=
\Delta_{i,j,k}^{',\,'',\,'''},
\endaligned
\]
\[
\aligned
D_{[i,j,k];\,\alpha}^8
&
:=
\Delta_{i,j,k}^{',\,'',\,''''}\,f_\alpha'
-
6\,\Delta_{i,j,k}^{',\,'',\,'''}\,f_\alpha'',
\\
N_{[i,j,k];\,\alpha,\beta}^{10}
&
:=
\Delta_{i,j,k}^{',\,''',\,''''}\,f_\alpha'f_\beta'
-
\frac{3}{2}\,\Delta_{i,j,k}^{',\,'',\,''''}\,
\big(f_\alpha'f_\beta''+f_\alpha''f_\beta'\big)
+
\\
&\ \ \ \ \
+
2\,\Delta_{i,j,k}^{',\,'',\,'''}\,
\big(f_\alpha'f_\beta'''+f_\alpha'''f_\beta')
+
3\,\Delta_{i,j,k}^{',\,'',\,'''}\,f_\alpha''f_\beta'',
\endaligned
\]
\[
\aligned
M_{[i,j],[k,l]}^8
&
:=
3\,\Delta_{i,j}^{',\,''''}\,\Delta_{k,l}^{',\,''}
+
12\,\Delta_{i,j}^{'',\,'''}\,\Delta_{k,l}^{',\,''}
-
\\
&\ \ \ \ \
-
5\,\Delta_{i,j}^{',\,'''}\,\Delta_{k,l}^{',\,'''},
\\
E_{[i,j,k],[p,q]}^{10}
&
:=
3\,\Delta_{i,j,k}^{',\,'',\,''''}\,\Delta_{l,m}^{',\,''}
-
6\,\Delta_{i,j,k}^{',\,'',\,'''}\,\Delta_{l,m}^{',\,'''},
\\
L_{[i,j,k],[l,m];\,\alpha}^{12}
&
:=
5\,\Delta_{i,j,k}^{',\,'',\,''''}\,\Delta_{p,q}^{',\,'''}\,f_\alpha'
-
15\,\Delta_{i,j,k}^{',\,'',\,''''}\,\Delta_{p,q}^{',\,''}\,f_\alpha''
-
6\,\Delta_{i,j,k}^{',\,'',\,'''}\,\Delta_{p,q}^{',\,''''}\,f_\alpha'
-
\\
&\ \ \ \ \
-24\,\Delta_{i,j,k}^{',\,'',\,'''}\,\Delta_{p,q}^{'',\,'''}\,f_\alpha'
+
30\,\Delta_{i,j,k}^{',\,'',\,'''}\,\Delta_{p,q}^{',\,'''}\,f_\alpha',
\endaligned
\]
\[
\aligned
Q_{[i,j,k],[p,q];\,\alpha,\beta}^{14}
&
:=
-10\,\Delta_{i,j,k}^{',\,''',\,''''}\,\Delta_{p,q}^{',\,'''}\,
f_\alpha'f_\beta'
+
\Delta_{i,j,k}^{',\,'',\,''''}\,\Delta_{p,q}^{',\,''''}\,
f_\alpha'f_\beta'
+
\\
&\ \ \ \ \
+4\,\Delta_{i,j,k}^{',\,''',\,''''}\,\Delta_{p,q}^{'',\,'''}\,
f_\alpha'f_\beta'
+
+20\,\Delta_{i,j,k}^{',\,'',\,''''}\,\Delta_{p,q}^{',\,'''}\,
f_\alpha'f_\beta''
+
\\
&\ \ \ \ \
+30\,\Delta_{i,j,k}^{',\,''',\,''''}\,\Delta_{p,q}^{',\,''}\,
f_\alpha'f_\beta''
-
6\,\Delta_{i,j,k}^{',\,'',\,'''}\,\Delta_{p,q}^{',\,''''}\,
f_\alpha'f_\beta''
-
\\
&\ \ \ \ \
-24\,\Delta_{i,j,k}^{',\,'',\,'''}\,\Delta_{p,q}^{'',\,'''}\,
f_\alpha'f_\beta''
-
40\,\Delta_{i,j,k}^{',\,'',\,'''}\,\Delta_{p,q}^{',\,'''}\,
f_\alpha'f_\beta'''
-
\\
&\ \ \ \ \
-75\,\Delta_{i,j,k}^{',\,'',\,''''}\,\Delta_{p,q}^{',\,''}\,
f_\alpha''f_\beta''
+
30\,\Delta_{i,j,k}^{',\,'',\,'''}\,\Delta_{p,q}^{',\,'''}\,
f_\alpha''f_\beta''
+
\\
&\ \ \ \ \ 
+120\,\Delta_{i,j,k}^{',\,'',\,'''}\,\Delta_{p,q}^{',\,''}\,
f_\alpha''f_\beta''',
\\
R_{[i,j,k],[p,q,r];\,\alpha}^{15}
&
:=
\Delta_{i,j,k}^{',\,'',\,''''}\,\Delta_{p,q,r}^{',\,'',\,''''}\,
f_\alpha'
-
12\,\Delta_{i,j,k}^{',\,''',\,''''}\,\Delta_{p,q,r}^{',\,'',\,'''}\,
f_\alpha'
+
\\
&\ \ \ \ \
+
24\,\Delta_{i,j,k}^{',\,'',\,''''}\,\Delta_{p,q,r}^{',\,'',\,'''}\,
f_\alpha''
-
48\,\Delta_{i,j,k}^{',\,'',\,'''}\,
\Delta_{p,q,r}^{',\,'',\,'''}\,f_\alpha''',
\endaligned
\]
\[
\aligned
U_{[i,j,k],[p,q,r],[s,t]}^{17}
&
:=
15\,\Delta_{i,j,k}^{',\,'',\,''''}\,\Delta_{p,q,r}^{',\,'',\,''''}\,
\Delta_{s,t}^{',\,''}
-
36\,\Delta_{i,j,k}^{',\,''',\,''''}\,\Delta_{p,q,r}^{',\,'',\,'''}\,
\Delta_{s,t}^{',\,''}
-
\\
&\ \ \ \ \
-24\,\Delta_{i,j,k}^{',\,'',\,''''}\,
\Delta_{p,q,r}^{',\,'',\,'''}\,\Delta_{s,t}^{',\,'''}
+
144\,\Delta_{i,j,k}^{',\,'',\,'''}\,\Delta_{p,q,r}^{',\,'',\,'''}\,
\Delta_{s,t}^{'',\,'''},
\endaligned
\]
\[
\small
\aligned
V_{[i,j,k],[p,q,r],[s,t];\,\alpha}^{19}
&
:=
24\,\Delta_{i,j,k}^{',\,''',\,''''}\,\Delta_{p,q,r}^{',\,'',\,''''}\,
\Delta_{s,t}^{',\,''}\,f_\alpha'
-
60\,\Delta_{i,j,k}^{',\,''',\,''''}\,\Delta_{p,q,r}^{',\,'',\,'''}\,
\Delta_{s,t}^{',\,'''}\,f_\alpha'
+
\\
&\ \ \ \ \
+\Delta_{i,j,k}^{',\,'',\,''''}\,\Delta_{p,q,r}^{',\,'',\,''''}\,
\Delta_{s,t}^{',\,'''}\,f_\alpha'
-
75\,\Delta_{i,j,k}^{',\,'',\,''''}\,\Delta_{p,q,r}^{',\,'',\,''''}\,
\Delta_{s,t}^{',\,''}\,f_\alpha''
+
\\
&\ \ \ \ \
+36\,\Delta_{i,j,k}^{',\,''',\,''''}\,\Delta_{p,q,r}^{',\,'',\,'''}\,
\Delta_{s,t}^{',\,''}\,f_\alpha''
+
168\,\Delta_{i,j,k}^{',\,'',\,''''}\,\Delta_{p,q,r}^{',\,'',\,'''}\,
\Delta_{s,t}^{',\,'''}\,f_\alpha''
-
\endaligned
\]
\[
\small
\aligned
&\ \ \ \ \
-144\,\Delta_{i,j,k}^{',\,'',\,'''}\,\Delta_{p,q,r}^{',\,'',\,'''}\,
\Delta_{s,t}^{'',\,'''}\,f_\alpha''
+
96\,\Delta_{i,j,k}^{',\,'',\,''''}\,\Delta_{p,q,r}^{',\,'',\,'''}\,
\Delta_{s,t}^{',\,''}\,f_\alpha'''
-
\\
&\ \ \ \ \
-
240\,\Delta_{i,j,k}^{',\,'',\,'''}\,\Delta_{p,q,r}^{',\,'',\,'''}\,
\Delta_{s,t}^{',\,'''}\,f_\alpha''',
\endaligned
\]
\[
\small
\aligned
&
X_{[i,j,k],[p,q,r],[s,t];\,\alpha,\beta}^{21}
:=
\\
&
\ \ \ \ \
:=
-40\,\Delta_{i,j,k}^{',\,''',\,''''}\,\Delta_{p,q,r}^{',\,'',\,''''}\,
\Delta_{s,t}^{',\,'''}\,
f_\alpha'f_\beta'
-
4\,\Delta_{i,j,k}^{',\,'',\,''''}\,\Delta_{p,q,r}^{',\,'',\,''''}\,
\Delta_{s,t}^{',\,''''}\,
f_\alpha'f_\beta'
-
\\
&\ \ \ \ \
-4\,\Delta_{i,j,k}^{',\,'',\,''''}\,\Delta_{p,q,r}^{',\,'',\,''''}\,
\Delta_{s,t}^{'',\,'''}\,
f_\alpha'f_\beta'
+
60\,\Delta_{i,j,k}^{',\,''',\,''''}\,\Delta_{p,q,r}^{',\,'',\,'''}\,
\Delta_{s,t}^{',\,''''}\,
f_\alpha'f_\beta'
+
\\
&\ \ \ \ \
+240\,\Delta_{i,j,k}^{',\,''',\,''''}\,\Delta_{p,q,r}^{',\,'',\,'''}\,
\Delta_{s,t}^{'',\,'''}\,
f_\alpha'f_\beta'
+
130\,\Delta_{i,j,k}^{',\,'',\,''''}\,\Delta_{p,q,r}^{',\,'',\,''''}\,
\Delta_{s,t}^{',\,'''}\,
f_\alpha'f_\beta''
+
\\
&\ \ \ \ \
+120\,\Delta_{i,j,k}^{',\,''',\,''''}\,\Delta_{p,q,r}^{',\,'',\,''''}\,
\Delta_{s,t}^{',\,''}\,
f_\alpha'f_\beta''
-
168\,\Delta_{i,j,k}^{',\,'',\,''''}\,\Delta_{p,q,r}^{',\,'',\,'''}\,
\Delta_{s,t}^{',\,''''}\,
f_\alpha'f_\beta''
-
\\
&\ \ \ \ \
-668\,\Delta_{i,j,k}^{',\,'',\,''''}\,\Delta_{p,q,r}^{',\,'',\,'''}\,
\Delta_{s,t}^{'',\,'''}\,
f_\alpha'f_\beta''
-
360\,\Delta_{i,j,k}^{',\,''',\,''''}\,\Delta_{p,q,r}^{',\,'',\,'''}\,
\Delta_{s,t}^{',\,'''}\,
f_\alpha'f_\beta''
-
\endaligned
\]
\[
\small
\aligned
&\ \ \ \ \
-160\,\Delta_{i,j,k}^{',\,'',\,''''}\,\Delta_{p,q,r}^{',\,'',\,'''}\,
\Delta_{s,t}^{',\,'''}\,
f_\alpha'f_\beta'''
+
240\,\Delta_{i,j,k}^{',\,'',\,'''}\,\Delta_{p,q,r}^{',\,'',\,'''}\,
\Delta_{s,t}^{',\,''''}\,
f_\alpha'f_\beta'''
+
\\
&\ \ \ \ \
+960\,\Delta_{i,j,k}^{',\,'',\,'''}\,\Delta_{p,q,r}^{',\,'',\,'''}\,
\Delta_{s,t}^{'',\,'''}\,
f_\alpha'f_\beta'''
-
375\,\Delta_{i,j,k}^{',\,'',\,''''}\,\Delta_{p,q,r}^{',\,'',\,''''}\,
\Delta_{s,t}^{',\,''}\,
f_\alpha''f_\beta''
+
\\
&\ \ \ \ \
+840\,\Delta_{i,j,k}^{',\,'',\,''''}\,\Delta_{p,q,r}^{',\,'',\,'''}\,
\Delta_{s,t}^{',\,'''}\,
f_\alpha''f_\beta''
+
180\,\Delta_{i,j,k}^{',\,''',\,''''}\,\Delta_{p,q,r}^{',\,'',\,'''}\,
\Delta_{s,t}^{',\,''}\,
f_\alpha''f_\beta''
+
\\
&\ \ \ \ \
+
144\,\Delta_{i,j,k}^{',\,'',\,'''}\,\Delta_{p,q,r}^{',\,'',\,'''}\,
\Delta_{s,t}^{',\,''''}\,
f_\alpha''f_\beta''
+
144\,\Delta_{i,j,k}^{',\,'',\,'''}\,\Delta_{p,q,r}^{',\,'',\,'''}\,
\Delta_{s,t}^{'',\,'''}\,
f_\alpha''f_\beta''
-
\\
&\ \ \ \ \
-1440\,\Delta_{i,j,k}^{',\,'',\,'''}\,\Delta_{p,q,r}^{',\,'',\,'''}\,
\Delta_{s,t}^{',\,'''}\,
f_\alpha''f_\beta'''
+
480\,\Delta_{i,j,k}^{',\,'',\,''''}\,\Delta_{p,q,r}^{',\,'',\,'''}\,
\Delta_{s,t}^{',\,''}\,
f_\alpha''f_\beta''',
\endaligned
\]
where the roman indices satisfy $1 \leqslant i < j < k \leqslant 4$,
where $1 \leqslant p < q < r \leqslant 4$, where $1 \leqslant s < r
\leqslant 4$ and where the two greek indices $\alpha, \beta$ satisfy
$1\leqslant \alpha , \beta\leqslant 4$ without restriction and
finally the {\bf total number} of these invariants generating the
Demailly-Semple algebra ${\sf E}_4^4$ equals:
\[
\aligned
&
1
+
4+6+24+96+4+16+64
+
\\
&
+
36+24+96+384+64+96+384+1536
=
\fbox{\bf 2835}\,.
\endaligned
\]

Furthermore, in the preceding dimension $n = 3$ for jets of the same
order $\kappa = 4$, the Demailly-Semple algebra ${\sf E}_4^3$ is
generated by the analogous list from which one removes the
four-dimensional Wronskian $W_{1, 2, 3, 4}^{ 10}$ and in which the
triples of skew-symmetric indices $[i, j, k]$ and $[ p, q, r]$ are set
to $[1, 2, 3]$ while $[ p, q]$ satisfy $1 \leqslant p < q \leqslant 3$
and $\alpha, \beta$ satisfy $1 \leqslant \alpha, \beta \leqslant 3$
without restriction, whence the total number of generators of ${\sf
E}_4^3$ equals:
\[
\aligned
&
3+3+9+27+1+3+9
+
9+3+9+27+3+3+9+27
=
\fbox{\bf 145}\,.
\endaligned
\]
}

\stopTHEOREM

\section*{\S12.~Approximate Schur bundle decomposition 
of ${\sf E}_{4, m}^4 T_X^*$}
\label{Section-12}

\subsection*{ Finite generation}
Thus, we know from the preceding section that ${\sf UE}_4^4$ is
generated by the sixteen bi-invariant polynomials:
\[
\Lambda^3,\ \ \
\Lambda^5,\ \ \
\Lambda^7,\ \ \
D^6,\ \ \
D^8,\ \ \ 
N^{10},\ \ \
M^8,\ \ \
E^{10},\ \ \
L^{12},\ \ \
Q^{14},\ \ \
R^{15},\ \ \
U^{17},\ \ \
V^{19},\ \ \
X^{21},\ \ \
f_1',\ \ \ 
W^{10},
\]
whose weight appears as an exponent. A general polynomial in these 16
invariants writes:
\[
\small
\aligned
\sum\,\text{\sf coeff}\cdot\,
&
\big(\Lambda^3\big)^a\,
\big(\Lambda^5\big)^b\,
\big(\Lambda^7\big)^c\,
\big(D^6\big)^d\,
\big(D^8\big)^e\,
\big(N^{10}\big)^f\,
\big(M^8\big)^g\,
\big(E^{10}\big)^h\,
\big(L^{12}\big)^i\,
\\
&
\big(Q^{14}\big)^j\,
\big(R^{15}\big)^k\,
\big(U^{17}\big)^l\,
\big(V^{19}\big)^m\,
\big(X^{21}\big)^n\,
\big(f_1'\big)^o\,
\big(W^{10}\big)^p,
\endaligned
\] 
where $a, b, c, d, e, f, g, h, i, j, k, l, m, n, o$ and $p$ are
nonnegative integer exponents. We temporarily
use the letter $m$ which should not
make confusion with the weighting $m$ appearing in ${\sf UE}_{ \kappa,
m}^n$. When one requires that such a polynomial has weight $m$, the
sum should be restricted to exponents satisfying:
\[
m
=
3a+5b+7c+6d+8e+10f+8g+10h+12i+
14j+15k+17l+19m+21n+o+10p.
\]
When one furthermore restricts such a general polynomial to $\{ f_1' =
0\}$, one gets:
\[
\small
\aligned
\sum_{3a+5b+\cdots+21n+10p=m}\,
\text{\sf coeff}\cdot\,
&
\big(\Lambda^3\big\vert_0\big)^a\,
\big(\Lambda^5\big\vert_0\big)^b\,
\big(\Lambda^7\big\vert_0\big)^c\,
\big(D^6\big\vert_0\big)^d\,
\big(D^8\big\vert_0\big)^e\,
\big(N^{10}\big\vert_0\big)^f\,
\big(M^8\big\vert_0\big)^g\,
\big(E^{10}\big\vert_0\big)^h\,
\\
&
\big(L^{12}\big\vert_0\big)^i\,
\big(Q^{14}\big\vert_0\big)^j\,
\big(R^{15}\big\vert_0\big)^k\,
\big(U^{17}\big\vert_0\big)^l\,
\big(V^{19}\big\vert_0\big)^m\,
\big(X^{21}\big\vert_0\big)^n\,
\big(W^{10}\big\vert_0\big)^p.
\endaligned
\]
Next, let ${\sf Syz}_{ 41}$ denote the ideal of $\C \big[ \Lambda^3
\vert_0, \dots, X^{ 21} \vert_0 \big]$ generated by the 41
lexicographic syzygies written on p.~\pageref{41-syzygies} (in which
one sets $f_1' = 0$) holding between the ordered variables:
\[
\small
\aligned
&
\Lambda^3\big\vert_0\,>\,\Lambda^5\big\vert_0\,>\,
\Lambda^7\big\vert_0\,>\,
D^6\big\vert_0\,>\,D^8\big\vert_0\,>\,N^{10}\big\vert_0\,>\,
M^8\big\vert_0\,>\,E^{10}\big\vert_0\,>\,
\\
&
L^{12}\big\vert_0\,>\,
Q^{14}\big\vert_0\,>\,R^{15}\big\vert_0\,>\,U^{17}\big\vert_0\,>\,
V^{19}\big\vert_0\,>\,X^{21}\big\vert_0.
\endaligned
\]
We list in columns the 41 Leading Terms of these 41 syzygies:
\[
\aligned
\underline{\Lambda^3\vert_0\Lambda^7\vert_0}_{
{\scriptscriptstyle{\sf LT}}}:
&
\ \ \ \ \ \ a\geqslant 1,\ \ \ c\geqslant 1
\ \ \ \ \ \ \ \ \ \ \ \ \ \ \ \ \ \ \ \ \ \ \ \ \ \ \ \ \,
\underline{\Lambda^5\vert_0D^8\vert_0}_{
{\scriptscriptstyle{\sf LT}}}:
\ \ \ \ \ \ b\geqslant 1,\ \ \ e\geqslant 1
\\
\underline{\Lambda^3\vert_0D^8\vert_0}_{
{\scriptscriptstyle{\sf LT}}}:
&
\ \ \ \ \ \ a\geqslant 1,\ \ \ e\geqslant 1
\ \ \ \ \ \ \ \ \ \ \ \ \ \ \ \ \ \ \ \ \ \ \ \ \ \ \ \ \,
\underline{\Lambda^5\vert_0N^{10}\vert_0}_{
{\scriptscriptstyle{\sf LT}}}:
\ \ \ \ \ \ b\geqslant 1,\ \ \ f\geqslant 1
\\
\underline{\Lambda^3\vert_0N^{10}\vert_0}_{
{\scriptscriptstyle{\sf LT}}}:
&
\ \ \ \ \ \ a\geqslant 1,\ \ \ f\geqslant 1
\ \ \ \ \ \ \ \ \ \ \ \ \ \ \ \ \ \ \ \ \ \ \ \ \ \ \ \
\underline{\Lambda^5\vert_0L^{12}\vert_0}_{
{\scriptscriptstyle{\sf LT}}}:
\ \ \ \ \ \ b\geqslant 1,\ \ \ i\geqslant 1
\\
\underline{\Lambda^3\vert_0L^{12}\vert_0}_{
{\scriptscriptstyle{\sf LT}}}:
&
\ \ \ \ \ \ a\geqslant 1,\ \ \ i\geqslant 1
\ \ \ \ \ \ \ \ \ \ \ \ \ \ \ \ \ \ \ \ \ \ \ \ \ \ \ \ \
\underline{\Lambda^5\vert_0Q^{14}\vert_0}_{
{\scriptscriptstyle{\sf LT}}}:
\ \ \ \ \ \ b\geqslant 1,\ \ \ j\geqslant 1
\\
\underline{\Lambda^3\vert_0Q^{14}\vert_0}_{
{\scriptscriptstyle{\sf LT}}}:
&
\ \ \ \ \ \ a\geqslant 1,\ \ \ j\geqslant 1
\ \ \ \ \ \ \ \ \ \ \ \ \ \ \ \ \ \ \ \ \ \ \ \ \ \ \ \ \,
\underline{\Lambda^5\vert_0R^{15}\vert_0}_{
{\scriptscriptstyle{\sf LT}}}:
\ \ \ \ \ \ b\geqslant 1,\ \ \ k\geqslant 1
\\
\underline{\Lambda^3\vert_0R^{15}\vert_0}_{
{\scriptscriptstyle{\sf LT}}}:
&
\ \ \ \ \ \ a\geqslant 1,\ \ \ k\geqslant 1
\ \ \ \ \ \ \ \ \ \ \ \ \ \ \ \ \ \ \ \ \ \ \ \ \ \ \ \ \,
\underline{\Lambda^5\vert_0U^{17}\vert_0}_{
{\scriptscriptstyle{\sf LT}}}:
\ \ \ \ \ \ b\geqslant 1,\ \ \ l\geqslant 1
\\
\underline{\Lambda^3\vert_0U^{17}\vert_0}_{
{\scriptscriptstyle{\sf LT}}}:
&
\ \ \ \ \ \ a\geqslant 1,\ \ \ l\geqslant 1
\ \ \ \ \ \ \ \ \ \ \ \ \ \ \ \ \ \ \ \ \ \ \ \ \ \ \ \ \ \,
\underline{\Lambda^5\vert_0V^{19}\vert_0}_{
{\scriptscriptstyle{\sf LT}}}:
\ \ \ \ \ \ b\geqslant 1,\ \ \ m\geqslant 1
\\
\underline{\Lambda^3\vert_0V^{19}\vert_0}_{
{\scriptscriptstyle{\sf LT}}}:
&
\ \ \ \ \ \ a\geqslant 1,\ \ \ m\geqslant 1
\ \ \ \ \ \ \ \ \ \ \ \ \ \ \ \ \ \ \ \ \ \ \ \ \ \ \ \
\underline{\Lambda^5\vert_0X^{21}\vert_0}_{
{\scriptscriptstyle{\sf LT}}}:
\ \ \ \ \ \ b\geqslant 1,\ \ \ n\geqslant 1
\\
\underline{\Lambda^3\vert_0X^{21}\vert_0}_{
{\scriptscriptstyle{\sf LT}}}:
&
\ \ \ \ \ \ a\geqslant 1,\ \ \ n\geqslant 1
\endaligned
\]
\[
\aligned
\underline{\Lambda^7\vert_0R^{15}\vert_0}_{
{\scriptscriptstyle{\sf LT}}}:
&
\ \ \ \ \ \ c\geqslant 1,\ \ \ k\geqslant 1
\ \ \ \ \ \ \ \ \ \ \ \ \ \ \ \ \ \ \ \ \ \ \ \ \ \ \ \ \,
\underline{D^6\vert_0N^{10}\vert_0}_{
{\scriptscriptstyle{\sf LT}}}:
\ \ \ \ \ \ d\geqslant 1,\ \ \ f\geqslant 1
\\
\underline{\Lambda^7\vert_0U^{17}\vert_0}_{
{\scriptscriptstyle{\sf LT}}}:
&
\ \ \ \ \ \ c\geqslant 1,\ \ \ l\geqslant 1
\ \ \ \ \ \ \ \ \ \ \ \ \ \ \ \ \ \ \ \ \ \ \ \ \ \ \ \ \ \,
\underline{D^6\vert_0L^{12}\vert_0}_{
{\scriptscriptstyle{\sf LT}}}:
\ \ \ \ \ \ d\geqslant 1,\ \ \ i\geqslant 1
\\
\underline{\Lambda^7\vert_0V^{19}\vert_0}_{
{\scriptscriptstyle{\sf LT}}}:
&
\ \ \ \ \ \ c\geqslant 1,\ \ \ m\geqslant 1
\ \ \ \ \ \ \ \ \ \ \ \ \ \ \ \ \ \ \ \ \ \ \ \ \ \ \ \
\underline{D^6\vert_0Q^{14}\vert_0}_{
{\scriptscriptstyle{\sf LT}}}:
\ \ \ \ \ \ d\geqslant 1,\ \ \ j\geqslant 1
\\
\underline{\Lambda^7\vert_0X^{21}\vert_0}_{
{\scriptscriptstyle{\sf LT}}}:
&
\ \ \ \ \ \ c\geqslant 1,\ \ \ n\geqslant 1
\ \ \ \ \ \ \ \ \ \ \ \ \ \ \ \ \ \ \ \ \ \ \ \ \ \ \ \ \
\underline{D^6\vert_0V^{19}\vert_0}_{
{\scriptscriptstyle{\sf LT}}}:
\ \ \ \ \ \ d\geqslant 1,\ \ \ m\geqslant 1
\\
& 
\ \ \ \ \ \ \ \ \ \ \ \ \ \ \ \ \ \ \ \ \ \ \ \ \ \ \ \ \
\ \ \ \ \ \ \ \ \ \ \ \ \ \ \ \ \ \ \ \ \ \ \ \ \ \ \ \ \
\underline{D^6\vert_0X^{21}\vert_0}_{
{\scriptscriptstyle{\sf LT}}}:
\ \ \ \ \ \ d\geqslant 1,\ \ \ n\geqslant 1
\endaligned
\]
\[
\aligned
\underline{D^8\vert_0L^{12}\vert_0}_{
{\scriptscriptstyle{\sf LT}}}:
&
\ \ \ \ \ \ e\geqslant 1,\ \ \ i\geqslant 1
\ \ \ \ \ \ \ \ \ \ \ \ \ \ \ \ \ \ \ \ \ \ \ \ 
\underline{D^6\vert_0M^8\vert_0R^{15}\vert_0}_{
{\scriptscriptstyle{\sf LT}}}:
\ \ \ \ \ \ d\geqslant 1,\ \ \ g\geqslant 1,\ \ \ k\geqslant 1
\\
\underline{D^8\vert_0Q^{14}\vert_0}_{
{\scriptscriptstyle{\sf LT}}}:
&
\ \ \ \ \ \ e\geqslant 1,\ \ \ j\geqslant 1
\ \ \ \ \ \ \ \ \ \ \ \ \ \ \ \ \ \ \ \ \ \ \ 
\underline{D^8\vert_0M^8\vert_0R^{15}\vert_0}_{
{\scriptscriptstyle{\sf LT}}}:
\ \ \ \ \ \ e\geqslant 1,\ \ \ g\geqslant 1,\ \ \ k\geqslant 1
\\
\underline{D^8\vert_0V^{19}\vert_0}_{
{\scriptscriptstyle{\sf LT}}}:
&
\ \ \ \ \ \ e\geqslant 1,\ \ \ m\geqslant 1
\ \ \ \ \ \ \ \ \ \ \ \ \ \ \ \ \ \ \ \ \ \ 
\underline{N^{10}\vert_0M^8\vert_0R^{15}\vert_0}_{
{\scriptscriptstyle{\sf LT}}}:
\ \ \ \ \ \ f\geqslant 1,\ \ \ g\geqslant 1,\ \ \ k\geqslant 1
\\
\underline{D^8\vert_0X^{21}\vert_0}_{
{\scriptscriptstyle{\sf LT}}}:
&
\ \ \ \ \ \ e\geqslant 1,\ \ \ n\geqslant 1
\ \ \ \ \ \ \ \ \ \ \ \ \ \ \ \ \ \ \ \ \ \ \
\underline{M^8\vert_0R^{15}\vert_0R^{15}\vert_0}_{
{\scriptscriptstyle{\sf LT}}}:
\ \ \ \ \ \ g\geqslant 1,\ \ \ k\geqslant 2
\endaligned
\]

\[
\aligned
\underline{E^{10}\vert_0Q^{14}\vert_0}_{
{\scriptscriptstyle{\sf LT}}}:
&
\ \ \ \ \ \ h\geqslant 1,\ \ \ j\geqslant 1
\ \ \ \ \ \ \ \ \ \ \ \ \ \ \ \ \ \ \
\underline{N^{10}\vert_0E^{10}\vert_0V^{19}\vert_0}_{
{\scriptscriptstyle{\sf LT}}}:
\ \ \ \ \ \ f\geqslant 1,\ \ \ h\geqslant 1,\ \ \ m\geqslant 1
\\
\underline{E^{10}\vert_0X^{21}\vert_0}_{
{\scriptscriptstyle{\sf LT}}}:
&
\ \ \ \ \ \ h\geqslant 1,\ \ \ n\geqslant 1
\ \ \ \ \ \ \ \ \ \ \ \ \ \ \ \ \ \ \
\underline{N^{10}\vert_0L^{12}\vert_0V^{19}\vert_0}_{
{\scriptscriptstyle{\sf LT}}}:
\ \ \ \ \ \ f\geqslant 1,\ \ \ i\geqslant 1,\ \ \ m\geqslant 1
\\
&
\ \ \ \ \ \ \ \ \ \ \ \ \ \ \ \ \ \ \ \ \ \ \ \ \ \ \
\ \ \ \ \ \ \ \ \ \ \ \ \ \ \ \ \ \ \ \ \ \
\underline{N^{10}\vert_0L^{12}\vert_0X^{21}\vert_0}_{
{\scriptscriptstyle{\sf LT}}}:
\ \ \ \ \ \ f\geqslant 1,\ \ \ i\geqslant 1,\ \ \ n\geqslant 1
\\
\underline{L^{12}\vert_0L^{12}\vert_0X^{21}\vert_0}_{
{\scriptscriptstyle{\sf LT}}}:
&
\ \ \ \ \ \ i\geqslant 2,\ \ \ n\geqslant 1
\ \ \ \ \ \ \ \ \ \ \ \ \ \ \ \ \ \ \
\underline{N^{10}\vert_0U^{17}\vert_0X^{21}\vert_0}_{
{\scriptscriptstyle{\sf LT}}}:
\ \ \ \ \ \ f\geqslant 1,\ \ \ l\geqslant 1,\ \ \ n\geqslant 1
\endaligned
\]
If, by ${\sf LT} \big( {\sf Syz}_{ 41} \big)$, we denote the monomial
ideal of $\C \big[ \Lambda^3 \big\vert_0, \dots, X^{ 21} \big\vert_0
\big]$ generated by these 41 Leading Terms, a known elementary 
property of reduced Gröbner bases shows that:
\[
\C\big[
\Lambda^3\big\vert_0,\dots,X^{21}\big\vert_0\big]
\big/
{\sf Syz}_{41}
\simeq
\C\big[
\Lambda^3\big\vert_0,\dots,X^{21}\big\vert_0\big]
\big/
\text{\sf LT}\big({\sf Syz}_{41}\big).
\]
More suitably for our purposes, the theorem on 
p.~\pageref{bi-invariant-4-4} states that any
bi-invariant of weight $m$ writes uniquely under the form:
\[
\aligned
{\sf P}\big(j^\kappa f\big)
=
&
\sum_{o,\,p}\,
(f_1')^o\,\big(W^{10}\big)^p\,
\sum_{(a,b,\dots,n)\in\N^{14}\backslash
(\square_1\cup\cdots\cup\square_{41})
\atop
3a+\cdots+21n=m-o-10p}\,
{\sf coeff}_{a,\dots,n,o,p}\,\cdot
\\
&\ \ \ \ \
\cdot
\big(\Lambda^3\big)^a\,
\big(\Lambda^5\big)^b\,
\big(\Lambda^7\big)^c\,
\big(D^6\big)^d\,
\big(D^8\big)^e\,
\big(N^{10}\big)^f
\big(M^8\big)^g\,
\big(E^{10}\big)^h\,
\\
&\ \ \ \ \ \ \
\big(L^{12}\big)^i\,
\big(Q^{14}\big)^j\,
\big(R^{15}\big)^k\,
\big(U^{17}\big)^l\,
\big(V^{19}\big)^m\,
\big(X^{21}\big)^n,
\endaligned
\]
with coefficients ${\sf coeff}_{a,\dots,n,o,p}$ subjected to no
restriction, where $\square_1$, \dots, $\square_{ 41}$ denote the
quadrants in $\N^{ 14}$ having vertex at the leading terms of our 41
syzygies.

Our goal now is to compute an approximation of this general sum of
monomials which will suffice for our Euler-Poincaré characteristic
computations below.

A general monomial in $\C \big[ \Lambda^3 , \dots, X^{ 21}
\big]$ writes:
\[
\small
\aligned
\text{\sf Monomial}
&
=
\big(\Lambda^3\big)^a\,
\big(\Lambda^5\big)^b\,
\big(\Lambda^7\big)^c\,
\big(D^6\big)^d\,
\big(D^8\big)^e\,
\big(N^{10}\big)^f\,
\big(M^8\big)^g\,
\big(E^{10}\big)^h\,
\\
& \ \ \ \ \
\big(L^{12}\big)^i\,
\big(Q^{14}\big)^j\,
\big(R^{15}\big)^k\,
\big(U^{17}\big)^l\,
\big(V^{19}\big)^{m}\,
\big(X^{21}\big)^n.
\endaligned
\] 
Such a monomial {\em belongs} to the monomial ideal ${\sf LT} \big(
{\sf Syz}_{ 41} \big)$ {\em if and only if} it is a multiple of at
least one of the 41 Leading Terms. Equivalently, the 14-tuple of
integers $(a, \dots, n)$ belongs to at least one quadrant $\square_i$
with vertex the exponent of the leading term of the $i$-th syzygy.
For instance, being a multiple of $\Lambda^3 \Lambda^7$ occurs when
and only when $a \geqslant 1$ and $c \geqslant 1$. In fact, in our
complete list of the 41 leading terms above, just after each leading
Term, we have in advance written the condition that such a ${\sf
Monomial}$ be a multiple of it.

On the contrary, for ${\sf Monomial}$ {\em not to be a multiple} of
$\Lambda^3 \Lambda^7$, it is necessary and sufficient that $a = 0$ or
$c = 0$, and more generally, for it to belong to the relevant quotient
ideal:
\[
\C\big[
\Lambda^3,\dots,X^{21}\big]
\big/
\text{\sf LT}\big({\sf Syz}_{41}\big),
\]
it is necessary and sufficient that its 14-tuple exponent 
$\big( a, b, c, d, e, f, g, h, i, j, k, l, m, n \big) \in
\N^{ 14}$ belongs to the following {\em intersection} of 41 subsets of
$\N^{ 14}$:
\[
\big\{a=0\big\}\cup\big\{c=0\big\}
\bigcap
\big\{a=0\big\}\cup\big\{e=0\big\}
\bigcap
\cdots
\bigcap
\big\{f=0\}\cup\big\{l=0\big\}\cup\big\{n=0\big\}.
\]
To compute this intersection, we shall abbreviate for instance $\big\{
a = 0 \big\} \cup \big\{ c = 0 \big\}$ by $( a + c)$ with the symbol
``$+$'' denoting union, and with the intersection being denoted by an
unwritten multiplication symbol, so that we may develope for instance
the product of the first two terms as follows:
\[
\aligned
\big\{a=0\big\}\cup\big\{c=0\big\}
\bigcap
\big\{a=0\big\}\cup\big\{e=0\big\}
&
\equiv
(a+c)(a+e)
\\
&
=
aa+ae+ca+ce
\\
&
=
a+ce,
\endaligned
\]
and simplify it immediately, on understanding that the symbol $a$
represents $\big\{ a = 0 \}$, hence contains both $ae \equiv \big\{ a
= e = 0\}$ and $ca \equiv \big\{ c = a = 0\big\}$.

With such a convention, grouping by packages, we may compute the
intersections colum by column, starting with the first column
containing $\Lambda^3 \vert_0$:
\[
(a+c)(a+e)(a+f)(a+i)(a+j)(a+k)(a+l)(a+m)(a+n)
=
a+cefijklmn,
\]
and getting in sum nine ``words'' that we should further ``intersect'':
\[
\aligned
&
a+cefijklmn,
\\
&
b+efijklmn,
\\
&
c+klmn,
\\
&
d+fijmn,
\\
&
e+ijmn,
\\
&
(d+g+k)(e+g+k)(f+g+k)(g+k+k^1),
\\
&
h+jn,
\\
&
i+i^1+n,
\\
&
(f+h+m)(f+i+m)(f+i+n)(f+l+n).
\endaligned
\]
Here, the ``letter'' $k^1$ appearing at the end of the sixth line
means the subset $\big\{ k = 1 \big\}$ of $\N^{ 14}$, not to be
confused with $k \equiv \big\{ k = 0 \big\}$.
Let us develope step by step the sixth and
the ninth lines:
\[
\aligned
&
(d+g+k)(e+g+k)(f+g+k)(g+k+k^1)=
\\
&
(d+g+k)(e+g+k)(g+k+fk^1)=
\\
&
(d+g+k)(g+k+efk^1)=
\\
&
g+k+defk^1
\endaligned
\]
\[
\aligned
&
(f+h+m)(f+i+m)(f+i+n)(f+l+n)=
\\
&
(f+h+m)(f+i+m)(f+n+il)=
\\
&
(f+h+m)(f+il+in+mn)=
\\
&
f+hil+hin+mn+ilm.
\endaligned
\]
Now we compute the product of the 
lines 3, 4, 5, 7:
\[
\aligned
&
(c+klmn)(d+fijmn)(e+ijmn)(h+jn)=
\\
&
(c+klmn)(d+fijmn)(eh+ejn+ijmn)=
\\
&
(c+klmn)(deh+dejn+dijmn+fijmn)=
\\
&
cdeh+cdejn+cdijmn+cfijmn+dehklmn+dejklmn+dijklmn+fijklmn
\endaligned
\]
and the product of the lines 1 and 2:
\[
\aligned
ab+aefijklmn+cefijklmn,
\endaligned
\]
whence the product of the lines 1, 2, 3, 4, 5, 7 is:
\[
\aligned
&
abcdeh+abcdejn+abcdijmn+abcfijmn+abdehklmn+abdejklmn+
\\
&
+abdijklmn+abfijklmn+aefijklmn+cefijklmn.
\endaligned
\]
On the other hand, the product of the lines 9, 6, 8 is:
\[
\aligned
&
(f+hil+hin+mn+ilm)(g+k+defk^1)(i+i^1+n)=
\\
&
(f+hil+hin+mn+ilm)(gi+gi^1+gn+ik+i^1k+kn+defik^1+defi^1k^1+defk^1n)=
\endaligned
\]
When developing the latter product, sometimes words containing the
product $i i^1$ (or $k k^1$) might appear. But they denote the empty
set $\big\{ i = 0 \big\} \cap \big\{ i = 1 \big\}$, so they should be 
left out. The direct result of the product, before
any simplification, is:
\[
\aligned
&
=
fgi+fgi^1+fgn+fik+fi^1k+fkn+defik^1+defi^1k^1+defk^1n+
\\
&\ \ \ \ \
+
ghil+\emptyset+ghiln+hikl+\emptyset+hikln+defik^1l+\emptyset+defhik^1ln+
\\
&\ \ \ \ \
+
ghin+\emptyset+ghin+hikn+\emptyset+hikn+defhik^1n+\emptyset+defhik^1n+
\\
&\ \ \ \ \
+
gimn+gi^1mn+gmn+ikmn+i^1kmn+kmn+defik^1mn+defi^1k^1mn+defk^1mn+
\\
&\ \ \ \ \ 
+
gilm+\emptyset+gilmn+iklm+\emptyset+iklmn+defik^1lm+\emptyset+defik^1lmn,
\endaligned
\]
and after simplification:
\[
\aligned
&
=
fgi+fgi^1+fgn+fik+fi^1k+fkn+defik^1+defi^1k^1+defk^1n+
\\
&\ \ \ \ \
+ghil+hikl+ghin+hikn+gmn+kmn+gilm+iklm.
\endaligned
\]
The final multiplication shall be:
\[
\aligned
&
\Big(
abcdeh+abcdejn+abcdijmn+abcfijmn+abdehklmn+abdejklmn+
\\
&\ \ \ \ \
+abdijklmn+abfijklmnaefijklmn+cefijklmn\Big)\cdot
\\
&
\cdot\Big(fgi+fgi^1+fgn+fik+fi^1k+fkn+defik^1+defi^1k^1+defk^1n+
\\
&\ \ \ \ \
+ghil+hikl+ghin+hikn+gmn+kmn+gilm+iklm\Big),
\endaligned
\]
but we will not expand it completely.

\subsection*{ Twenty-four families of monomials}\
Instead, we will compute the product modulo words which contain more
than 9 letters. The reason why we do so will be appearent later. The
result then consists of 30 words of 9 letters:
\[
\boxed{
\aligned
{\sf A}:
&
\ \ \ \ \ \ \ \ \ \ \ \ \ \
abcdefghi
\ \ \ \ \ \ \ \ \ \ \ \ \ \ \ \ \ \ \ \ \ \
\ \ \ \ \ \ \ \ \ \ \ \ \ \ \ \ \ \ \ \ \ \ \
{\sf J}:
\ \ \ \ \ \ \ \ \ \ \ \ \ \
abcdegjmn
\\
{\sf A'}:
&
\ \ \ \ \ \ \ \ \ \ \ \ \ \
abcdefghi^1
\ \ \ \ \ \ \ \ \ \ \ \ \ \ \ \ \ \ \ \ \
\ \ \ \ \ \ \ \ \ \ \ \ \ \ \ \ \ \ \ \ \ \,
{\sf K}:
\ \ \ \ \ \ \ \ \ \ \ \ \ \
abcdehikl
\\
{\sf B}:
&
\ \ \ \ \ \ \ \ \ \ \ \ \ \
abcdefghn
\ \ \ \ \ \ \ \ \ \ \ \ \ \ \ \ \ \ \ \ \ \
\ \ \ \ \ \ \ \ \ \ \ \ \ \ \ \ \ \ \ \ \ \,
{\sf L}:
\ \ \ \ \ \ \ \ \ \ \ \ \ \
abcdehikn
\\
{\sf C}:
&
\ \ \ \ \ \ \ \ \ \ \ \ \ \
abcdefgjn
\ \ \ \ \ \ \ \ \ \ \ \ \ \ \ \ \ \ \ \ \
\ \ \ \ \ \ \ \ \ \ \ \ \ \ \ \ \ \ \ \ \ \,
{\sf M}:
\ \ \ \ \ \ \ \ \ \ \ \ \ \
abcdehkmn
\\
{\sf D}:
&
\ \ \ \ \ \ \ \ \ \ \ \ \ \
abcdefhik
\ \ \ \ \ \ \ \ \ \ \ \ \ \ \ \ \ \ \ \ \ \
\ \ \ \ \ \ \ \ \ \ \ \ \ \ \ \ \ \ \ \ \ \,
{\sf N}:
\ \ \ \ \ \ \ \ \ \ \ \ \ \
abcdejkmn
\\
{\sf D'}:
&
\ \ \ \ \ \ \ \ \ \ \ \ \ \
abcdefhi^1k
\ \ \ \ \ \ \ \ \ \ \ \ \ \ \ \ \ \ \ \ \
\ \ \ \ \ \ \ \ \ \ \ \ \ \ \ \ \ \ \ \ \
{\sf O}:
\ \ \ \ \ \ \ \ \ \ \ \ \ \
abcdgijmn
\\
{\sf D''}:
&
\ \ \ \ \ \ \ \ \ \ \ \ \ \
abcdefhik^1
\ \ \ \ \ \ \ \ \ \ \ \ \ \ \ \ \ \ \ \ \
\ \ \ \ \ \ \ \ \ \ \ \ \ \ \ \ \ \ \ \ \
{\sf P}:
\ \ \ \ \ \ \ \ \ \ \ \ \ \
abcdijkmn
\\
{\sf D'''}:
&
\ \ \ \ \ \ \ \ \ \ \ \ \ \
abcdefhi^1k^1
\ \ \ \ \ \ \ \ \ \ \ \ \ \ \ \ \ \ \ \
\ \ \ \ \ \ \ \ \ \ \ \ \ \ \ \ \ \ \ \
{\sf Q}:
\ \ \ \ \ \ \ \ \ \ \ \ \ \
abcfgijmn
\\
{\sf E}:
&
\ \ \ \ \ \ \ \ \ \ \ \ \ \
abcdefhkn
\ \ \ \ \ \ \ \ \ \ \ \ \ \ \ \ \ \ \ \ \
\ \ \ \ \ \ \ \ \ \ \ \ \ \ \ \ \ \ \ \ \ \,
{\sf R}:
\ \ \ \ \ \ \ \ \ \ \ \ \ \
abcfijkmn
\\
{\sf E'}:
&
\ \ \ \ \ \ \ \ \ \ \ \ \ \
abcdefhk^1n
\ \ \ \ \ \ \ \ \ \ \ \ \ \ \ \ \ \ \ \ \
\ \ \ \ \ \ \ \ \ \ \ \ \ \ \ \ \ \ \ \ \,
{\sf S}:
\ \ \ \ \ \ \ \ \ \ \ \ \ \
abdehklmn
\\
{\sf F}:
&
\ \ \ \ \ \ \ \ \ \ \ \ \ \
abcdefjkn
\ \ \ \ \ \ \ \ \ \ \ \ \ \ \ \ \ \ \ \ \
\ \ \ \ \ \ \ \ \ \ \ \ \ \ \ \ \ \ \ \ \ \
{\sf T}:
\ \ \ \ \ \ \ \ \ \ \ \ \ \
abdejklmn
\\
{\sf F'}:
&
\ \ \ \ \ \ \ \ \ \ \ \ \ \
abcdefjk^1n
\ \ \ \ \ \ \ \ \ \ \ \ \ \ \ \ \ \ \ \ \
\ \ \ \ \ \ \ \ \ \ \ \ \ \ \ \ \ \ \ \ \,
{\sf U}:
\ \ \ \ \ \ \ \ \ \ \ \ \ \
abdijklmn
\\
{\sf G}:
&
\ \ \ \ \ \ \ \ \ \ \ \ \ \
abcdeghil
\ \ \ \ \ \ \ \ \ \ \ \ \ \ \ \ \ \ \ \ \ \ \
\ \ \ \ \ \ \ \ \ \ \ \ \ \ \ \ \ \ \ \ \ \
{\sf V}:
\ \ \ \ \ \ \ \ \ \ \ \ \ \
abfijklmn
\\
{\sf H}:
&
\ \ \ \ \ \ \ \ \ \ \ \ \ \
abcdeghin
\ \ \ \ \ \ \ \ \ \ \ \ \ \ \ \ \ \ \ \ \ \ 
\ \ \ \ \ \ \ \ \ \ \ \ \ \ \ \ \ \ \ \ \
{\sf W}:
\ \ \ \ \ \ \ \ \ \ \ \ \ \
aefijklmn
\\
{\sf I}:
&
\ \ \ \ \ \ \ \ \ \ \ \ \ \
abcdeghimn
\ \ \ \ \ \ \ \ \ \ \ \ \ \ \ \ \ \ \ \
\ \ \ \ \ \ \ \ \ \ \ \ \ \ \ \ \ \ \ \ \,
{\sf X}:
\ \ \ \ \ \ \ \ \ \ \ \ \ \
cefijklmn
\endaligned}
\]
Recalling that the first 
word $abcdefghi$ for instance means the condition
$\big\{ a = b = c = d = e = f = g = h = i = 0 \big\}$ on the exponents
of a general monomial, we may therefore list in an extensive array the
24 families
${\sf A}$, ${\sf B}$, ${\sf C}$, 
${\sf D}$, ${\sf E}$, ${\sf F}$, 
${\sf G}$, ${\sf H}$, ${\sf I}$, 
${\sf J}$, ${\sf K}$, ${\sf L}$, 
${\sf M}$, ${\sf N}$, ${\sf O}$, 
${\sf P}$, ${\sf Q}$, ${\sf R}$, 
${\sf S}$, ${\sf T}$, ${\sf U}$,
${\sf V}$, ${\sf W}$, ${\sf X}$
of corresponding monomials, the subsidiary families ${\sf A}'$; ${\sf
D}'$, ${\sf D}''$, ${\sf D}'''$; ${\sf E}'$; ${\sf F}'$ being
considered as similar to ${\sf A}$; ${\sf D}$; ${\sf E}$; ${\sf F}$:

$\:$
\hspace{-5.75cm}
\begin{minipage}[t]{25cm}
\[
\tiny
\begin{array}{ccccccccccccccccc}
{\sf A}: 
& \bullet & \bullet & \bullet & \bullet & \bullet & \bullet & 
\bullet & \bullet & \bullet & \big(Q^{14}\big)^j & 
\big(R^{15}\big)^k & \big(U^{17}\big)^l & 
\big(V^{19}\big)^{m} & \big(X^{21}\big)^n &
(f_1')^o & \big(W^{10}\big)^p
\\
{\sf B}: 
& \bullet & \bullet & \bullet & \bullet & \bullet & \bullet & 
\bullet & \bullet & \big(L^{12}\big)^i & \big(Q^{14}\big)^j & 
\big(R^{15}\big)^k & \big(U^{17}\big)^l & 
\big(V^{19}\big)^{m} & \bullet &
(f_1')^o & \big(W^{10}\big)^p
\\
{\sf C}: 
& \bullet & \bullet & \bullet & \bullet & \bullet & \bullet & 
\bullet & \big(E^{10}\big)^h & \big(L^{12}\big)^i & \bullet & 
\big(R^{15}\big)^k & \big(U^{17}\big)^l & 
\big(V^{19}\big)^{m} & \bullet &
(f_1')^o & \big(W^{10}\big)^p
\\
{\sf D}: 
& \bullet & \bullet & \bullet & \bullet & \bullet & \bullet & 
\big(M^8\big)^g & \bullet & \bullet & \big(Q^{14}\big)^j & 
\bullet & \big(U^{17}\big)^l & 
\big(V^{19}\big)^{m} & \big(X^{21}\big)^n &
(f_1')^o & \big(W^{10}\big)^p
\\
{\sf E}: 
& \bullet & \bullet & \bullet & \bullet & \bullet & \bullet & 
\big(M^8\big)^g & \bullet & \big(L^{12}\big)^i & \big(Q^{14}\big)^j & 
\bullet & \big(U^{17}\big)^l & 
\big(V^{19}\big)^{m} & \bullet &
(f_1')^o & \big(W^{10}\big)^p
\\
{\sf F}: 
& \bullet & \bullet & \bullet & \bullet & \bullet & \bullet & 
\big(M^8\big)^g & \big(E^{10}\big)^h & 
\big(L^{12}\big)^i & \bullet & 
\bullet & \big(U^{17}\big)^l & 
\big(V^{19}\big)^{m} & \bullet &
(f_1')^o & \big(W^{10}\big)^p
\\
{\sf G}: 
& \bullet & \bullet & \bullet & \bullet & \bullet & \big(N^{10}\big)^f & 
\bullet & \bullet & 
\bullet & \big(Q^{14}\big)^j & 
\big(R^{15}\big)^k & \bullet & 
\big(V^{19}\big)^{m} & \big(X^{21}\big)^n &
(f_1')^o & \big(W^{10}\big)^p
\\
{\sf H}: 
& \bullet & \bullet & \bullet & \bullet & \bullet & \big(N^{10}\big)^f & 
\bullet & \bullet & 
\bullet & \big(Q^{14}\big)^j & 
\big(R^{15}\big)^k & \big(U^{17}\big)^l & 
\big(V^{19}\big)^{m} & \bullet &
(f_1')^o & \big(W^{10}\big)^p
\\
{\sf I}: 
& \bullet & \bullet & \bullet & \bullet & \bullet & \big(N^{10}\big)^f & 
\bullet & \bullet & 
\big(L^{12}\big)^i & \big(Q^{14}\big)^j & 
\big(R^{15}\big)^k & \big(U^{17}\big)^l & 
\bullet & \bullet &
(f_1')^o & \big(W^{10}\big)^p
\\
{\sf J}: 
& \bullet & \bullet & \bullet & \bullet & \bullet & \big(N^{10}\big)^f & 
\bullet & \big(E^{10}\big)^h & 
\big(L^{12}\big)^i & \bullet & 
\big(R^{15}\big)^k & \big(U^{17}\big)^l & 
\bullet & \bullet &
(f_1')^o & \big(W^{10}\big)^p
\\
{\sf K}: 
& \bullet & \bullet & \bullet & \bullet & \bullet & \big(N^{10}\big)^f & 
\big(M^8\big)^g & \bullet & 
\bullet & \big(Q^{14}\big)^j & 
\bullet & \bullet & 
\big(V^{19}\big)^{m} & \big(X^{21}\big)^n &
(f_1')^o & \big(W^{10}\big)^p
\\
{\sf L}: 
& \bullet & \bullet & \bullet & \bullet & \bullet & \big(N^{10}\big)^f & 
\big(M^8\big)^g & \bullet & 
\bullet & \big(Q^{14}\big)^j & 
\bullet & \big(U^{17}\big)^l & 
\big(V^{19}\big)^{m} & \bullet &
(f_1')^o & \big(W^{10}\big)^p
\\
{\sf M}: 
& \bullet & \bullet & \bullet & \bullet & \bullet & \big(N^{10}\big)^f & 
\big(M^8\big)^g & \bullet & 
\big(L^{12}\big)^i & \big(Q^{14}\big)^j & 
\bullet & \big(U^{17}\big)^l & 
\bullet & \bullet &
(f_1')^o & \big(W^{10}\big)^p
\\
{\sf N}: 
& \bullet & \bullet & \bullet & \bullet & \bullet & \big(N^{10}\big)^f & 
\big(M^8\big)^g & \big(E^{10}\big)^h & 
\big(L^{12}\big)^i & \bullet & 
\bullet & \big(U^{17}\big)^l & 
\bullet & \bullet &
(f_1')^o & \big(W^{10}\big)^p
\\
{\sf O}: 
& \bullet & \bullet & \bullet & \bullet & \big(D^8\big)^e & 
\big(N^{10}\big)^f & 
\bullet & \big(E^{10}\big)^h & 
\bullet & \bullet & 
\big(R^{15}\big)^k & \big(U^{17}\big)^l & 
\bullet & \bullet &
(f_1')^o & \big(W^{10}\big)^p
\\
{\sf P}: 
& \bullet & \bullet & \bullet & \bullet & \big(D^8\big)^e & 
\big(N^{10}\big)^f & 
\big(M^8\big)^g & \big(E^{10}\big)^h & 
\bullet & \bullet & 
\bullet & \big(U^{17}\big)^l & 
\bullet & \bullet &
(f_1')^o & \big(W^{10}\big)^p
\\
{\sf Q}: 
& \bullet & \bullet & \bullet & \big(D^6\big)^d & \big(D^8\big)^e & 
\bullet & 
\bullet & \big(E^{10}\big)^h & 
\bullet & \bullet & 
\big(R^{15}\big)^k & \big(U^{17}\big)^l & 
\bullet & \bullet &
(f_1')^o & \big(W^{10}\big)^p
\\
{\sf R}: 
& \bullet & \bullet & \bullet & \big(D^6\big)^d & \big(D^8\big)^e & 
\bullet & 
\big(M^8\big)^g & \big(E^{10}\big)^h & 
\bullet & \bullet & 
\bullet & \big(U^{17}\big)^l & 
\bullet & \bullet &
(f_1')^o & \big(W^{10}\big)^p
\\
{\sf S}: 
& \bullet & \bullet & \big(\Lambda^7\big)^c & \bullet & \bullet & 
\big(N^{10}\big)^f & 
\big(M^8\big)^g & \bullet & 
\big(L^{12}\big)^i & \big(Q^{14}\big)^j & 
\bullet & \bullet & 
\bullet & \bullet &
(f_1')^o & \big(W^{10}\big)^p
\\
{\sf T}: 
& \bullet & \bullet & \big(\Lambda^7\big)^c & \bullet & \bullet & 
\big(N^{10}\big)^f & 
\big(M^8\big)^g & \big(E^{10}\big)^h & 
\big(L^{12}\big)^i & \bullet & 
\bullet & \bullet & 
\bullet & \bullet &
(f_1')^o & \big(W^{10}\big)^p
\\
{\sf U}: 
& \bullet & \bullet & \big(\Lambda^7\big)^c & \bullet & 
\big(D^8\big)^e & 
\big(N^{10}\big)^f & 
\big(M^8\big)^g & \big(E^{10}\big)^h & 
\bullet & \bullet & 
\bullet & \bullet & 
\bullet & \bullet &
(f_1')^o & \big(W^{10}\big)^p
\\
{\sf V}: 
& \bullet & \bullet & \big(\Lambda^7\big)^c & \big(D^6\big)^d & 
\big(D^8\big)^e & 
\bullet & 
\big(M^8\big)^g & \big(E^{10}\big)^h & 
\bullet & \bullet & 
\bullet & \bullet & 
\bullet & \bullet &
(f_1')^o & \big(W^{10}\big)^p
\\
{\sf W}: 
& \bullet & \big(\Lambda^5\big)^b & 
\big(\Lambda^7\big)^c & \big(D^6\big)^d & 
\bullet & 
\bullet & 
\big(M^8\big)^g & \big(E^{10}\big)^h & 
\bullet & \bullet & 
\bullet & \bullet & 
\bullet & \bullet &
(f_1')^o & \big(W^{10}\big)^p
\\
{\sf X}: 
& \big(\Lambda^3\big)^a & \big(\Lambda^5\big)^b & \bullet &
\big(D^6\big)^d & \bullet & \bullet & \big(M^8\big)^g &
\big(E^{10}\big)^h & \bullet & \bullet & \bullet & \bullet & \bullet &
\bullet & \big(f_1'\big)^o & \big(W^{10}\big)^p
\end{array}
\]
\end{minipage}

\subsection*{ General Schur bundle decomposition of 
${\sf E}_{4, m}^4 T_X^*$} By general representation theory, the
polynomial action of ${\sf GL}_4 ( \C)$ decomposes in a certain direct
sum of irreducible Schur representations. What we call bi-invariants
correspond to vectors of highest weight for the ${\sf GL}_4 (
\C)$-representation. To each vector of highest weight corresponds one
and only one irreducible Schur representation. Such a vector of
highest weight is nothing else but a monomial:
\[
\small
\aligned
&
\big(\Lambda^3\big)^a\,
\big(\Lambda^5\big)^b\,
\big(\Lambda^7\big)^c\,
\big(D^6\big)^d\,
\big(D^8\big)^e\,
\big(N^{10}\big)^f\,
\big(M^8\big)^g\,
\big(E^{10}\big)^h\,
\\
&
\big(L^{12}\big)^i\,
\big(Q^{14}\big)^j\,
\big(R^{15}\big)^k\,
\big(U^{17}\big)^l\,
\big(V^{19}\big)^{m'}\,
\big(X^{21}\big)^n\,
(f_1')^o\,
\big(W^{10}\big),
\endaligned
\] 
with the usual condition on exponents: $3a + \cdots + 21 n + o + 10p =
m$ and $(a, \dots, n)$ belonging to the complement $\N^{ 14}
\big\backslash \big( \square_1 \cup \cdots \cup \square_{ 41} \big)$
of the 41 quadrants. From now on, we denote by 
$m'$ the exponent of $V^{ 19}$ to distinguish 
it from the weight $m$ of the bi-invariant.

To know what are the four integers $\ell_1, \ell_2, \ell_3,
\ell_4$ of the corresponding Schur representations $\Gamma^{
(\ell_1, \ell_2, \ell_3, \ell_4)} \C^4$, it suffices to
consider the diagonal matrices of ${\sf GL}_4 ( \C)$ of the form:
\[
{\sf x}
:=
\left(
\begin{array}{cccc}
x_1 & 0 & 0 & 0
\\ 
0 & x_2 & 0 & 0
\\
0 & 0 & x_3 & 0
\\
0 & 0 & 0 & x_4
\end{array}
\right),
\]
for which all vectors of highest weight are then just eigenvectors
having eigenvalue of the form $x_1^{ \ell_1} x_2^{ \ell_2} x_3^{
\ell_3} x_4^{ \ell_4}$. 

Here in our situation, coming back to the theorem which describes the
2835 generators of ${\sf E}_4^4$, we should at first write down our 16
bi-invariants under a form in which we emphasize the lower indices as
we did for the general invariants. This gives us the 
following more informative list:
\[
\boxed{
\aligned
&
\ell_{[1,2]}^3,\ \ \ \ \
\ell_{[1,2];\,1}^5,\ \ \ \ \
\ell_{[1,2];\,1,1}^7,\ \ \ \ \
D_{[1,2,3]}^6,\ \ \ \ \
D_{[1,2,3];\,1}^8,\ \ \ \ \
N_{[1,2,3];\,1,1}^{10},\ \ \ \ \
M_{[1,2],[1,2]}^8,\ \ \ \ \ 
\\
&
\ \ \ \ \ \ \ \ \ \ \ \ \
E_{[1,2,3],[1,2]}^{10},\ \ \ \ \
L_{[1,2,3],[1,2];\,1}^{12},\ \ \ \ \
Q_{[1,2,3],[1,2];\,1,1}^{14},\ \ \ \ \
R_{[1,2,3],[1,2,3];\,1}^{15},\ \ \ \ \
\\
&
U_{[1,2,3],[1,2,3],[1,2]}^{17},\ \ \ \ \
V_{[1,2,3],[1,2,3],[1,2];\,1}^{19},\ \ \ \ \
X_{[1,2,3],[1,2,3],[1,2];\,1,1}^{21},\ \ \ \ \
f_1',\ \ \ \ \
W_{[1,2,3,4]}^{10}
\endaligned
}\,.
\]
Then it is easy to realize that $\ell_1$, $\ell_2$, $\ell_3$,
$\ell_4$ just count the number of indices 1, 2, 3, 4 respectively
at the bottom of each invariant. Consequently, we have
the sixteen correspondences:
\[
\aligned
\big(\ell^3\big)^a:
&
\ \ \ \ \ \ \ \ \ \ \ \ \ \ \ \ \
\Gamma^{(a,a,0,0)}\C^4
\\
\big(\ell^5\big)^b:
&
\ \ \ \ \ \ \ \ \ \ \ \ \ \ \ \ \
\Gamma^{(2b,b,0,0)}\C^4
\\
\big(\ell^7\big)^c:
&
\ \ \ \ \ \ \ \ \ \ \ \ \ \ \ \ \
\Gamma^{(3c,c,0,0)}\C^4
\\
\big(D^6\big)^d:
&
\ \ \ \ \ \ \ \ \ \ \ \ \ \ \ \ \
\Gamma^{(d,d,d,0)}\C^4
\\
\big(D^8\big)^e:
&
\ \ \ \ \ \ \ \ \ \ \ \ \ \ \ \ \
\Gamma^{(2e,e,e,0)}\C^4
\endaligned
\]
\[
\aligned
\big(N^{10}\big)^f:
&
\ \ \ \ \ \ \ \ \ \ \ \ \ \ \ \ \
\Gamma^{(3f,f,f,0)}\C^4
\\
\big(M^8\big)^g:
&
\ \ \ \ \ \ \ \ \ \ \ \ \ \ \ \ \
\Gamma^{(2g,2g,0,0)}\C^4
\\
\big(E^{10}\big)^h:
&
\ \ \ \ \ \ \ \ \ \ \ \ \ \ \ \ \
\Gamma^{(2h,2h,h,0)}\C^4
\\
\big(L^{12}\big)^i:
&
\ \ \ \ \ \ \ \ \ \ \ \ \ \ \ \ \
\Gamma^{(3i,2i,i,0)}\C^4
\\
\big(Q^{14}\big)^j:
&
\ \ \ \ \ \ \ \ \ \ \ \ \ \ \ \ \
\Gamma^{(4j,2j,j,0)}\C^4
\endaligned
\]
\[
\aligned
\big(R^{15}\big)^k:
&
\ \ \ \ \ \ \ \ \ \ \ \ \ \ \ \ \
\Gamma^{(3k,2k,2k,0)}\C^4
\\
\big(U^{17}\big)^l:
&
\ \ \ \ \ \ \ \ \ \ \ \ \ \ \ \ \
\Gamma^{(3l,3l,2l,0)}\C^4
\\
\big(V^{19}\big)^{m'}:
&
\ \ \ \ \ \ \ \ \ \ \ \ \ \ \ \ \
\Gamma^{(4m',3m',2m',0)}\C^4
\\
\big(X^{21}\big)^n:
&
\ \ \ \ \ \ \ \ \ \ \ \ \ \ \ \ \
\Gamma^{(5n,3n,2n,0)}\C^4
\endaligned
\]
\[
\aligned
\big(f_1')^o:
&
\ \ \ \ \ \ \ \ \ \ \ \ \ \ \ \ \
\Gamma^{(o,0,0,0)}\C^4
\\
\big(W^{10}\big)^p:
&
\ \ \ \ \ \ \ \ \ \ \ \ \ \ \ \ \
\Gamma^{(p,p,p,p)}\C^4
\endaligned
\]
and it immediately follows that the Schur representation $\Gamma^{ (
\ell_1, \ell_2, \ell_3, \ell_4 )} \C^4$ which corresponds
to the general monomial written above has integers
$\ell_i$ given by:
\[
\left\{
\aligned
\ell_1
&
=
o+a+2b+3c+d+2e+3f+2g+2h+3i+4j+3k+3l+4m'+5n+p,
\\
\ell_2
&
=a+b+c+d+e+f+2g+2h+2i+2j+2k+3l+3m'+3n+p,
\\
\ell_3
&
=
d+e+f+h+i+j+2k+2l+2m'+2n+p,
\\
\ell_4
&
=
p.
\endaligned\right.
\]

By a direct application of the theorem on
p.~\pageref{bi-invariant-4-4} of \S11, we obtain an exact Schur bundle
decompositition of the graduate $m$-th part ${\sf E}_{ 4, m}^4 T_X^*$
of the Demailly-Semple bundle ${\sf E}_4^4 T_X^*$ on a complex
algebraic hypersurface $X \subset \P^5 ( \C)$.

\THEOREM

\smallskip\noindent\fbox{\bf THEOREM}\ \ 
\label{Schur-decomposition-4-4}
{\sf\em
In dimension $n = 4$ for jet order $\kappa = 4$, graduate $m$-th part
${\sf E}_{ 4, m}^4 T_X^*$ of the Demailly-Semple bundle ${\sf E}_4^4
T_X^* = \oplus_m \, {\sf E}_{ 4, m}^4 T_X^*$ on a complex algebraic
hypersurface $X \subset \P^5 ( \C)$ has the following decomposition in
direct sums of Schur bundles:
\[
\aligned
&
{\sf E}_{4,m}^4T_X^*
=
\bigoplus_{(a,b,\dots,n)\in\N^{14}\backslash
(\square_1\cup\cdots\cup\square_{41}) 
\atop
o+3a+\cdots+21n+10p=m}\,
\\
&
{\scriptsize
\Gamma
\left(
\aligned
o+a+2b+3c+d+2e+3f+2g+2h+3i+4j+3k+3l+4m'+5n+p&
\\
a+b+c+d+e+f+2g+2h+2i+2j+2k+3l+3m'+3n+p&
\\
d+e+f+h+i+j+2k+2l+2m'+2n+p&
\\
p&
\endaligned
\right)}\,T_X^*,
\endaligned
\]
where the 41 subsets $\square_i$ of $\N^{ 14}$ are precisely defined
by:
\[
\footnotesize
\aligned
&
\big\{a\geqslant 1,\,c\geqslant 1\big\},\ \ \ \ \
\big\{a\geqslant 1,\,e\geqslant 1\big\},\ \ \ \ \
\big\{a\geqslant 1,\,f\geqslant 1\big\},\ \ \ \ \
\big\{a\geqslant 1,\,i\geqslant 1\big\},
\\
&
\big\{a\geqslant 1,\,j\geqslant 1\big\},\ \ \ \ \
\big\{a\geqslant 1,\,k\geqslant 1\big\},\ \ \ \ \
\big\{a\geqslant 1,\,l\geqslant 1\big\},\ \ \ \ \
\big\{a\geqslant 1,\,m'\geqslant 1\big\},
\\
&
\big\{a\geqslant 1,\,n\geqslant 1\big\},\ \ \ \ \
\big\{b\geqslant 1,\,e\geqslant 1\big\},\ \ \ \ \
\big\{b\geqslant 1,\,f\geqslant 1\big\},\ \ \ \ \
\big\{b\geqslant 1,\,i\geqslant 1\big\},
\\
&
\big\{b\geqslant 1,\,j\geqslant 1\big\},\ \ \ \ \
\big\{b\geqslant 1,\,k\geqslant 1\big\},\ \ \ \ \
\big\{b\geqslant 1,\,l\geqslant 1\big\},\ \ \ \ \
\big\{b\geqslant 1,\,m'\geqslant 1\big\},
\\
&
\big\{b\geqslant 1,\,n\geqslant 1\big\},\ \ \ \ \
\big\{c\geqslant 1,\,k\geqslant 1\big\},\ \ \ \ \
\big\{c\geqslant 1,\,l\geqslant 1\big\},\ \ \ \ \
\big\{c\geqslant 1,\,m'\geqslant 1\big\},
\\
&
\big\{c\geqslant 1,\,n\geqslant 1\big\},\ \ \ \ \
\big\{d\geqslant 1,\,f\geqslant 1\big\},\ \ \ \ \
\big\{d\geqslant 1,\,i\geqslant 1\big\},\ \ \ \ \
\big\{d\geqslant 1,\,j\geqslant 1\big\},
\\
&
\big\{d\geqslant 1,\,m\geqslant 1\big\},\ \ \ \ \
\big\{d\geqslant 1,\,n\geqslant 1\big\},\ \ \ \ \
\big\{e\geqslant 1,\,i\geqslant 1\big\},\ \ \ \ \
\big\{e\geqslant 1,\,j\geqslant 1\big\},
\\
&
\big\{e\geqslant 1,\,m'\geqslant 1\big\},\ \ \ \ \
\big\{e\geqslant 1,\,n\geqslant 1\big\},\ \ \ \ \
\big\{d\geqslant 1,\,g\geqslant 1,\,k\geqslant 1\big\},\ \ \ \ \
\\
&
\big\{e\geqslant 1,\,g\geqslant 1,\,k\geqslant 1\big\},\ \ \ \ \
\big\{f\geqslant 1,\,g\geqslant 1,\,k\geqslant 1\big\},\ \ \ \ \
\big\{g\geqslant 1,\,k\geqslant 2\big\},\ \ \ \ \
\\
&
\big\{h\geqslant 1,\,j\geqslant 1\big\},\ \ \ \ \
\big\{h\geqslant 1,\,n\geqslant 1\big\},\ \ \ \ \
\big\{i\geqslant 2,\,n\geqslant 1\big\},\ \ \ \ \
\\
&
\big\{f\geqslant 1,\,h\geqslant 1,\,m'\geqslant 1\big\},\ \ \ \ \
\big\{f\geqslant 1,\,i\geqslant 1,\,m'\geqslant 1\big\},\ \ \ \ \
\big\{f\geqslant 1,\,i\geqslant 1,\,n\geqslant 1\big\},
\\
&
\big\{f\geqslant 1,\,l\geqslant 1,\,n\geqslant 1\big\}.
\endaligned
\]
In addition, in the preceding dimension $n = 3$ for jets of the same
order $\kappa = 4$, one has an entirely similar Schur bundle
decomposition of ${\sf E}_{ 4, m}^3 T_X^*$ for any $m$ in which one
removes $W^{ 10}$, one sets $p = 0$ and one removes the fourth
component $\ell_4$ of $\Gamma^{ ( \ell_1, \ell_2, \ell_3,
\ell_4)}$:
\[
\aligned
&
\ \ \ \ \ \ \ \ \ \ \ \ \ \
{\sf E}_{4,m}^3T_X^*
=
\bigoplus_{(a,b,\dots,n)\in\N^{14}\backslash
(\square_1\cup\cdots\cup\square_{41}) 
\atop
o+3a+\cdots+21n=m}\,
\\
&
{\footnotesize
\Gamma
\left(
\aligned
o+a+2b+3c+d+2e+3f+2g+2h+3i+4j+3k+3l+4m'+5n&
\\
a+b+c+d+e+f+2g+2h+2i+2j+2k+3l+3m'+3n&
\\
d+e+f+h+i+j+2k+2l+2m'+2n&
\endaligned
\right)}\,T_X^*.
\endaligned
\]
}

\stopTHEOREM

\subsection*{ Approximate Schur bundle decomposition}
We now come back to our 24 words of 9 letters and we
make three remarks which will simplify a bit the
further computations. 

\smallskip$\bullet$\
The full complement
$\N^{ 14} \big\backslash \big( \square_1 \cup \cdots \cup \square_{
41} \big)$ is slightly larger than the union of the 30 subsets of
$\N^{ 14}$ defined by ${\sf A}$, ${\sf A}'$, ${\sf B}$, \dots, ${\sf
W}$ ${\sf X}$, in the sense that it contains also a finite number of
subsets defined by equating to 0 (or to 1) more than 9 exponents.
These subsets will not contribute to the dominant term $m^{ 16}$ when
calculating the Euler-Poincaré characteristic of ${\sf E}_{ 4, m}^4
T_X^*$ and hence, they will at once be left out.

\smallskip$\bullet$\
The first family ${\sf A}$ corresponds to 
a general polynomial of the form:
\[
\small
\aligned
\sum_{o+14j+15k+17l+19m'+21n+10p=m}\,
{\sf A}_{j,k,l,m',n,o,p}\cdot
\big(Q^{14}\big)^j
&
\big(R^{15}\big)^k
\big(U^{17}\big)^l
\big(V^{19}\big)^{m'}
\\
&
\big(X^{21}\big)^n
\big(f_1'\big)^o
\big(W^{10}\big)^p.
\endaligned
\]
The second family ${\sf A}'$ corresponds to a 
general polynomial of the form:
\[
\small
\aligned
L^{12}\,
\sum_{o+14j+15k+17l+19m'+21n+10p=m-12}\,
{\sf A}_{j,k,l,m',n,o,p}'\cdot
\big(Q^{14}\big)^j
&
\big(R^{15}\big)^k
\big(U^{17}\big)^l
\big(V^{19}\big)^{m'}
\\
&
\big(X^{21}\big)^n
\big(f_1'\big)^o
\big(W^{10}\big)^p.
\endaligned
\]
It is entirely of the same type as ${\sf A}$, except that
the weight $m$ is
replaced by $m - 12$. We will see that its contribution to the
dominant $m^{ 16}$-term of the
Euler-Poincaré characteristic is exactly the same\footnote{\, The
argument will simply be that $(m - {\rm cst.})^{ 16} = m^{ 16} + {\rm O}
( m^{ 15})$ as $m \to \infty$. }, hence we will remove ${\sf A}'$ and
provide the family ${\sf A}$ with the multiplicity 2.
Similarly, ${\sf D}$, ${\sf E}$ and ${\sf F}$ will have multiplicity 4, 
2 and 2.

\smallskip$\bullet$\
The third (now second) family ${\sf B}$ corresponds to a 
general polynomial of the form:
\[
\small
\aligned
\sum_{o+12i+14j+15k+17l+19m'+10p=m}\,
{\sf B}_{i,j,k,l,m',o,p}\cdot
&
\big(L^{12}\big)^i
\big(Q^{14}\big)^j
\big(R^{15}\big)^k
\\
&
\big(U^{17}\big)^l
\big(V^{19}\big)^{m'}
\big(f_1'\big)^o
\big(W^{10}\big)^p,
\endaligned
\]
hence its intersection with the family ${\sf A}$ is nontrivial, 
consisting of polynomials of the form:
\[
\small
\aligned
\sum_{o+14j+15k+17l+19m'+10p=m}\,
\widetilde{\sf B}_{j,k,l,m',o,p}\cdot
\big(Q^{14}\big)^j
\big(R^{15}\big)^k
\big(U^{17}\big)^l
\big(V^{19}\big)^{m'}
\big(f_1'\big)^o
\big(W^{10}\big)^p.
\endaligned
\]
In principle, we should write the union of two overlapping families
${\sf A} \cup {\sf B}$ in the form of two non-intersecting families:
${\sf A} \cup \big( {\sf B} \backslash {\sf A} \big)$, but here again,
because the intersection ${\sf A} \cap {\sf B}$ is represented by the
word $abcdefghin$ which has $10 > 9$ letters, this intersection will
only contribute the Euler-Poincaré characteristic as an ${\rm O} ( m^{
15} )$, which will not perturb the dominant term $m^{ 16}$, as $m \to
\infty$. So we can consider the 24 remaining families (a bit of which
have multiplicities) without caring about
overlappings.

\smallskip

In summary, up to certain negligible sums of Schur bundles which will
not contribute to the dominant $m^{ 16}$-term while calculating the
Euler-Poincaré characteristic of ${\sf E}_{ 4, m}^4 T_X^*$, we have to
consider \red{\bf 24} direct sums of Schur bundles with
multiplicities, indexed from \green{${\sf A}$} up to \green{${\sf X}$}
in the roman alphabet:
\[
\footnotesize
\aligned
\blue{
\green{\underline{\sf A:}}\ \ \ \ \
\red{\bf 2\cdot}
\bigoplus_{m=o+14j+15k+17l+19m+21n+10p}\,\,\,
\Gamma
\left(
\aligned
o+4j+3k+3l+4m+5n+p&
\\
2j+2k+3l+3m+3n+p&
\\
j+2k+2l+2m+2n+p&
\\
p&
\endaligned
\right)\,T_X^*},
\endaligned
\]
\[
\footnotesize
\aligned
\blue{
\green{\underline{\sf B:}}\ \ \ \ \
\bigoplus_{m=o+12i+14j+15k+17l+19m+10p}\,\,\,
\Gamma
\left(
\aligned
o+3i+4j+3k+3l+4m+p
&
\\
2i+2j+2k+3l+3m+p
&
\\
i+j+2k+2l+2m+p
&
\\
p
&
\endaligned
\right)\,T_X^*},
\endaligned
\]
\[
\footnotesize
\aligned
\blue{
\green{\underline{\sf C:}}\ \ \ \ \
\bigoplus_{m=o+10h+12i+15k+17l+19m+10p}\,\,\,
\Gamma
\left(
\aligned
o+2h+3i+3k+3l+4m+p
&
\\
2h+2i+2k+3l+3m+p
&
\\
h+i+2k+2l+2m+p
&
\\
p
&
\endaligned
\right)\,T_X^*},
\endaligned
\]
\[
\footnotesize
\aligned
\blue{
\green{\underline{\sf D:}}\ \ \ \ \
\red{\bf 4\cdot}
\bigoplus_{m=o+8g+14j+17l+19m+21n+10p}\,\,\,
\Gamma
\left(
\aligned
o+2g+4j+3l+4m+5n+p
&
\\
2g+2j+3l+3m+3n+p
&
\\
j+2l+2m+2n+p
&
\\
p
&
\endaligned
\right)\,T_X^*},
\endaligned
\]
\[
\footnotesize
\aligned
\blue{
\green{\underline{\sf E:}}\ \ \ \ \
\red{\bf 2\cdot}
\bigoplus_{m=o+8g+12i+14j+17l+19m+10p}\,\,\,
\Gamma
\left(
\aligned
o+2g+3i+4j+3l+4m+p
&
\\
2g+2i+2j+3l+3m+p
&
\\
i+j+2l+2m+p
&
\\
p
&
\endaligned
\right)\,T_X^*},
\endaligned
\]
\[
\footnotesize
\aligned
\blue{
\green{\underline{\sf F:}}\ \ \ \ \
\red{\bf 2\cdot}
\bigoplus_{m=o+8g+10h+12i+17l+19m+10p}\,\,\,
\Gamma
\left(
\aligned
o+2g+2h+3i+3l+4m+p
&
\\
2g+2h+2i+3l+3m+p
&
\\
h+i+2l+2m+p
&
\\
p
&
\endaligned
\right)\,T_X^*},
\endaligned
\]
\[
\footnotesize
\aligned
\blue{
\green{\underline{\sf G:}}\ \ \ \ \
\bigoplus_{m=o+10f+14j+15k+19m+21n+10p}\,\,\,
\Gamma
\left(
\aligned
o+3f+4j+3k+4m+5n+p
&
\\
f+2j+2k+3m+3n+p
&
\\
f+j+2k+2m+2n+p
&
\\
p
&
\endaligned
\right)\,T_X^*},
\endaligned
\]
\[
\footnotesize
\aligned
\blue{
\green{\underline{\sf H:}}\ \ \ \ \
\bigoplus_{m=o+10f+14j+15k+17l+19m+10p}\,\,\,
\Gamma
\left(
\aligned
o+3f+4j+3k+3l+4m+p
&
\\
f+2j+2k+3l+3m+p
&
\\
f+j+2k+2l+2m+p
&
\\
p
&
\endaligned
\right)\,T_X^*},
\endaligned
\]
\[
\footnotesize
\aligned
\blue{
\green{\underline{\sf I:}}\ \ \ \ \
\bigoplus_{m=o+10f+12i+14j+15k+17l+10p}\,\,\,
\Gamma
\left(
\aligned
o+3f+3i+4j+3k+3l+p
&
\\
f+2i+2j+2k+3l+p
&
\\
f+i+j+2k+2l+p
&
\\
p
&
\endaligned
\right)\,T_X^*},
\endaligned
\]
\[
\footnotesize
\aligned
\blue{
\green{\underline{\sf J:}}\ \ \ \ \
\bigoplus_{m=o+10f+10h+12i+15k+17l+10p}\,\,\,
\Gamma
\left(
\aligned
o+3f+2h+3i+3k+3l+p
&
\\
f+2h+2i+2k+3l+p
&
\\
f+h+i+2k+2l+p
&
\\
p
&
\endaligned
\right)\,T_X^*},
\endaligned
\]
\[
\footnotesize
\aligned
\blue{
\green{\underline{\sf K:}}\ \ \ \ \
\bigoplus_{m=o+10f+8g+14j+19m+21n+10p}\,\,\,
\Gamma
\left(
\aligned
o+3f+2g+4j+4m+5n+p
&
\\
f+2g+2j+3m+3n+p
&
\\
f+j+2m+2n+p
&
\\
p
&
\endaligned
\right)\,T_X^*},
\endaligned
\]
\[
\footnotesize
\aligned
\blue{
\green{\underline{\sf L:}}\ \ \ \ \
\bigoplus_{m=o+10f+8g+14j+17l+19m+10p}\,\,\,
\Gamma
\left(
\aligned
o+3f+2g+4j+3l+4m+p
&
\\
f+2g+2j+3l+3m+p
&
\\
f+j+2l+2m+p
&
\\
p
&
\endaligned
\right)\,T_X^*},
\endaligned
\]
\[
\footnotesize
\aligned
\blue{
\green{\underline{\sf M:}}\ \ \ \ \
\bigoplus_{m=o+10f+8g+12i+14j+17l+10p}\,\,\,
\Gamma
\left(
\aligned
o+3f+2g+3i+4j+3l+p
&
\\
f+2g+2i+2j+3l+p
&
\\
f+i+j+2l+p
&
\\
p
&
\endaligned
\right)\,T_X^*},
\endaligned
\]
\[
\footnotesize
\aligned
\blue{
\green{\underline{\sf N:}}\ \ \ \ \
\bigoplus_{m=o+10f+8g+10h+12i+17l+10p}\,\,\,
\Gamma
\left(
\aligned
o+3f+2g+2h+3i+3l+p
&
\\
f+2g+2h+2i+3l+p
&
\\
f+h+i+2l+p
&
\\
p
&
\endaligned
\right)\,T_X^*},
\endaligned
\]
\[
\footnotesize
\aligned
\blue{
\green{\underline{\sf O:}}\ \ \ \ \
\bigoplus_{m=o+8e+10f+10h+15k+17l+10p}\,\,\,
\Gamma
\left(
\aligned
o+2e+3f+2h+3k+3l+p
&
\\
e+f+2h+2k+3l+p
&
\\
e+f+h+2k+2l+p
&
\\
p
&
\endaligned
\right)\,T_X^*},
\endaligned
\]
\[
\footnotesize
\aligned
\blue{
\green{\underline{\sf P:}}\ \ \ \ \
\bigoplus_{m=o+8e+10f+8g+10h+17l+10p}\,\,\,
\Gamma
\left(
\aligned
o+2e+3f+2g+2h+3l+p
&
\\
e+f+2g+2h+3l+p
&
\\
e+f+h+2l+p
&
\\
p
&
\endaligned
\right)\,T_X^*},
\endaligned
\]
\[
\footnotesize
\aligned
\blue{
\green{\underline{\sf Q:}}\ \ \ \ \
\bigoplus_{m=o+6d+8e+10h+15k+17l+10p}\,\,\,
\Gamma
\left(
\aligned
o+d+2e+2h+3k+3l+p
&
\\
d+e+2h+2k+3l+p
&
\\
d+e+h+2k+2l+p
&
\\
p
&
\endaligned
\right)\,T_X^*},
\endaligned
\]
\[
\footnotesize
\aligned
\blue{
\green{\underline{\sf R:}}\ \ \ \ \
\bigoplus_{m=o+6d+8e+8g+10h+17l+10p}\,\,\,
\Gamma
\left(
\aligned
o+d+2e+2g+2h+3l+p
&
\\
d+e+2g+2h+3l+p
&
\\
d+e+h+2l+p
&
\\
p
&
\endaligned
\right)\,T_X^*},
\endaligned
\]
\[
\footnotesize
\aligned
\blue{
\green{\underline{\sf S:}}\ \ \ \ \
\bigoplus_{m=o+7c+10f+8g+12i+14j+10p}\,\,\,
\Gamma
\left(
\aligned
o+3c+3f+2g+3i+4j+p
&
\\
c+f+2g+2i+2j+p
&
\\
f+i+j+p
&
\\
p
&
\endaligned
\right)\,T_X^*},
\endaligned
\]
\[
\footnotesize
\aligned
\blue{
\green{\underline{\sf T:}}\ \ \ \ \
\bigoplus_{m=o+7c+10f+8g+10h+12i+10p}\,\,\,
\Gamma
\left(
\aligned
o+3c+3f+2g+2h+3i+p
&
\\
c+f+2g+2h+2i+p
&
\\
f+h+i+p
&
\\
p
&
\endaligned
\right)\,T_X^*},
\endaligned
\]
\[
\footnotesize
\aligned
\blue{
\green{\underline{\sf U:}}\ \ \ \ \
\bigoplus_{m=o+7c+8e+10f+8g+10h+10p}\,\,\,
\Gamma
\left(
\aligned
o+3c+2e+3f+2g+2h+p
&
\\
c+e+f+2g+2h+p
&
\\
e+f+h+p
&
\\
p
&
\endaligned
\right)\,T_X^*},
\endaligned
\]
\[
\footnotesize
\aligned
\blue{
\green{\underline{\sf V:}}\ \ \ \ \
\bigoplus_{m=o+7c+6d+8e+8g+10h+10p}\,\,\,
\Gamma
\left(
\aligned
o+3c+d+2e+2g+2h+p
&
\\
c+d+e+2g+2h+p
&
\\
d+e+h+p
&
\\
p
&
\endaligned
\right)\,T_X^*},
\endaligned
\]
\[
\footnotesize
\aligned
\blue{
\green{\underline{\sf W:}}\ \ \ \ \
\bigoplus_{m=o+5b+7c+6d+8g+10h+10p}\,\,\,
\Gamma
\left(
\aligned
o+2b+3c+d+2g+2h+p
&
\\
b+c+d+2g+2h+p
&
\\
d+h+p
&
\\
p
&
\endaligned
\right)\,T_X^*},
\endaligned
\]
\[
\footnotesize
\aligned
\blue{
\green{\underline{\sf X:}}\ \ \ \ \
\bigoplus_{m=o+3a+5b+6d+8g+10h+10p}\,\,\,
\Gamma
\left(
\aligned
o+a+2b+d+2g+2h+p
&
\\
a+b+d+2g+2h+p
&
\\
d+h+p
&
\\
p
&
\endaligned
\right)\,T_X^*}.
\endaligned
\]
It is now time to speak of the asymptotic of the Euler characteristic
of a single Schur bundle.

\section*{ \S13.~Asymptotic expansion of the Euler characteristic
$\chi \big( \Gamma^{ (\ell_1, \ell_2, \dots, \ell_n)} T_X^* \big)$ }
\label{Section-13}

\subsection*{Euler-Poincaré characteristic of Schur bundles}
Let $X^n \subset \P^{ n+1 } ( \C)$ be a complex algebraic hypersurface
and denote by ${\sf c }_1, {\sf c }_2, \dots, {\sf c }_n$ be the Chern
classes ${\sf c }_k \big( T_X \big)$ of the tangent bundle $T_X$.
Each ${\sf c}_k$ may be represented by a smooth differential form of
bidegree $(k,k)$ on $X$. One thus assigns the weight $k$ to ${\sf
c}_k$. Because the total degrees of these forms are all even, 
the commutation relations ${\sf c}_{ k_1 } {\sf c}_{ k_2} = {\sf c}_{
k_2} {\sf c}_{k_1}$ hold for the cup product.

Every polynomial in the Chern classes:
\[
\sum_{k_1+\cdots+k_n=n}\,
{\sf coeff}\cdot
{\sf c}_{k_1}\,{\sf c}_{k_2}\cdots
{\sf c}_{k_n}
\]
which is homogeneous of degree $n = \dim X$ is represented by an $(n,
n)$-form on $X$, hence may be integrated. By a standard abuse of
language, such a polynomial is usually considered both as an $(n ,
n)$-form and as the purely numerical quantity:
\[
\int_X\,
\sum_{k_1+\cdots+k_n=n}\,
{\sf coeff}\,
{\sf c}_{k_1}\cdot{\sf c}_{k_2}\cdots
{\sf c}_{k_n}.
\]
For instance, if $d$ denotes the degree of $X$, one shows $\int_X \,
{\sf c}_1^n = d^{ n+1}$, a kind of relation often abbreviated ${\sf c}_1^n
= d^{ n+1}$. 

To speak in full generality (\cite{ dem1997, rou2006a, 
div2007}), the short exact sequence:
\[
0
\longrightarrow
T_X
\longrightarrow
T_{\P^{n+1}}\big\vert_X
\longrightarrow
\mathcal{O}_X(d)
\longrightarrow
0
\]
gives the relation ${\sc c}_\bullet \big( T_{\P^{n+1}}\big\vert_X
\big) = {\sc c}_\bullet \big( T_X\big) \cdot {\sc c}_\bullet \big(
\mathcal{ O}_X ( d) \big)$ between total Chern classes of the middle
term and of the two extreme ones, or more explicitly:
\[
(1+h)^{n+2}
=
\big[
1+{\sf c}_1+\cdots+{\sf c}_n
\big]
(1+dh),
\]
where $(1+h)^{ n+2}$ is the total Chern class of $\P^{ n+1}$ with $h =
{\sf c}_1 \big( \mathcal{ O}_{ \P^{ n+1}} ( 1) \big)$ being a $(1,
1)$-form. Consequently, by expanding both the left-hand and the
right-hand sides and by identifying terms of the same bidegree, we get
closed expressions for all the Chern classes.

\smallskip\noindent{\bf Lemma.}
{\em 
In terms of the hyperplane divisor $h = {\sf c}_1 \big( \mathcal{ O}_{
\P^{ n+1}} ( 1) \big)$ which satisfies $\int_X \, h^n = d = \deg X$,
the Chern classes ${\sf c}_k$ of $T_X$ are given by:
\[
{\sf c}_k
=
(-1)^k\,h^k\,
\Big(
d^k-
{\textstyle{\frac{(n+2)!}{1!\,\,(n+1)!}}}\,d^{k-1}
+\cdots+(-1)^k
{\textstyle{\frac{(n+2)!}{k!\,\,(n+2-k)!}}}
\Big).
\]
}\medskip

\proof
We indeed expand the two sides of the above relation between 
total Chern classes:
\[
1+
{\textstyle{\frac{(n+2)!}{1!\,(n+1)!}}}\,h
+\cdots+
{\textstyle{\frac{(n+2)!}{n!\,2!}}}\,h^n
=
1
+({\sc c}_1+dh)+({\sf c}_2+d{\sf c}_1h)
+\cdots+
({\sf c}_n+d{\sf c}_{n-1}h),
\]
on understanding that the forms $h^{ n+1}$, $h^{ n+2}$ and ${\sf c}_n
h$ of degree $> 2n$ vanish identically. Identifying forms of the same
bidegree yields the binomial-type recurrence relations: ${\sf c}_k =
\frac{ (n+2)!}{ k!\, (n+2-k)!}\, h^k - d {\sf c}_{ k-1} h$.
\endproof

It follows for instance as we said that ${\sf c}_1^n = (-1)^n d^{
n+1}$ and that ${\sf c}_1^{ n-2} {\sf c}_2 = (-1)^{ n-2} d \big( d-
\frac{ (n+2) !}{ (n+1)! \, 1 !} \big)^{ n-2} \big( d^2 - \frac{ (n+2)
!}{ (n+1)! \, 1 !} \, d + \frac{ (n+2) !}{ n! \, 2 !} \big)$ are
numerical quantities.

\smallskip

Following~\cite{ hirz1966}, one introduces the formal factorization:
\[
1
+
{\sf c}_1\,x
+
{\sf c}_2\,x^2
+\cdots+
{\sf c}_n\,x^n
=
\prod_{0\leqslant i\leqslant n}\,
\big(1+{\sf a}_i\,x\big),
\]
using new formal symbols ${\sf a}_i$ whose elementary
symmetric functions regive the Chern classes ${\sf c}_k$:
\[
{\sf c}_k
=
\sum_{1\leqslant i_1<i_2<\cdots<i_k\leqslant n}\,
{\sf a}_{i_1}\,{\sf a}_{i_2}\cdots{\sf a}_{i_k},
\]
so that any polynomial ${\sf P} \big( {\sf a}_1, \dots, {\sf a}_n
\big)$ in the ${\sf a}_i$ which is invariant under all permutations of
its arguments may in fact be expressed in terms of the ${\sf c}_k$.
Every such a symmetric ${\sf P} \big( {\sf a}_1, \dots, {\sf a}_n
\big)$ which is homogeneous of degree $n$ may thus be considered as a
numerical quantity, after integration.

\smallskip\noindent{\bf Proposition.}
(\cite{ hirz1966, rou2004})
{\em 
The Euler-Poincaré characteristic:
\[
\chi\Big(X,\,\Gamma^{(\ell_1,\dots,\ell_n)}\,T_X\Big)
=
\sum_{i=0}^n\,(-1)^i\,\dim H^i\big(X,\,
\Gamma^{(\ell_1,\dots,\ell_n)}\,T_X\big)
\] 
of an arbitrary Schur bundle $\Gamma^{ (\ell_1, \ell_2, \dots,
\ell_n)} \, T_X$ with $\ell_1 \geqslant \ell_2 \geqslant \cdots
\geqslant \ell_n$ is given as {\rm (}the integral over $X$ of{\rm )}
the rewriting by means of the ${\sf c}_k$ of all the terms which are
homogeneous of degree $n$ with respect to ${\sf a}_1, \dots, {\sf
a}_n$ in the expansion of the {\rm (}symmetric{\rm )} quotient:
\[
\left\vert
\begin{array}{ccc}
e^{{\sf a}_1\ell_1'} & \cdots & e^{{\sf a}_1\ell_n'}
\\
\vdots & \ddots & \vdots
\\
e^{{\sf a}_n\ell_1'} & \cdots & e^{{\sf a}_n\ell_n'}
\end{array}
\right\vert\,
\bigg/\,
{\scriptsize{
\left\vert
\begin{array}{ccc}
e^{(n-1){\sf a}_1} & \cdots & 1
\\
\vdots & \ddots & \vdots
\\
e^{(n-1){\sf a}_n} & \cdots & 1
\end{array}
\right\vert}},
\]
in which one has abbreviated for notational condensation:
\[
\ell_1'
:=
\ell_1+n-1,\ \ \ \ \ \ 
\ell_2'
:=
\ell_2+n-2,\,\dots\dots,\ \ \
\ell_n'
:=
\ell_n.
\]
}\medskip

We shall admit this result.
In fact, the well known Van der Monde determinant
yields an approximate expression of the denominator:
\[
\aligned
{\scriptsize{
\left\vert
\begin{array}{ccc}
e^{(n-1){\sf a}_1} & \cdots & 1
\\
\vdots & \ddots & \vdots
\\
e^{(n-1){\sf a}_n} & \cdots & 1
\end{array}
\right\vert}}
&
=
\prod_{1\leqslant i<j\leqslant n}\,
\big(e^{{\sf a}_i}-e^{{\sf a}_j}\big)
\\
&
=
\prod_{1\leqslant i<j\leqslant n}\,
\big({\sf a}_i-{\sf a}_j\big)
\cdot
\big[1+R({\sf a}_1,\dots,{\sf a}_n)\big],
\endaligned
\]
where the remainder $R ( {\sf a})$ denotes a local holomorphic
function which vanishes at the origin. Because the determinant at the
numerator also visibly vanishes whenever one ${\sf a}_{ i_1}$ is equal
to another ${\sf a}_{ i_2}$, for some two distinct indices $i_1$ and
$i_2$, this numerator also is a multiple, as a holomorphic function,
of the same product $\prod_{ 1 \leqslant i < j \leqslant n}\, \big(
{\sf a}_i - {\sf a}_j \big)$. Consequently, when one expands
simultaneously the numerator and the denominator, the two products
should cancel out:
\[
\frac{\prod_{i<j}\,({\sf a}_i-{\sf a}_j)\,
\big[S({\sf a},\,\ell')\big]}{
\prod_{i<j}\,({\sf a}_i-{\sf a}_j)\big[1+R({\sf a})\big]}
=
S({\sf a},\,\ell')
\Big[
1
-
R({\sf a})
+
R({\sf a})^2
-
R({\sf a})^3
+\cdots
\Big]
\]
and one should obtain a power series in which only the homogeneous
terms of degree $n$ in the ${\sf a}_i$ are relevant. Getting a partial
{\em explicit} expression of the result is our next goal.

\subsection*{ Asymptotic expansion of 
the Euler-Poincaré characteristic of $\Gamma^{ ( \ell_1, \ell_2,
\dots, \ell_n ) }\, T_X$} A {\sl partition} of $n$ is any sequence:
\[
\lambda
=
\big(\lambda_1,\,\lambda_2,\,\dots,\,\lambda_n\big)
\]
of non-negative integers listed in decreasing order:
\[
\lambda_1\geqslant\lambda_2\geqslant\cdots\geqslant\lambda_n,
\]
whose total sum equals $n$:
\[
\lambda_1+\lambda_2
+\cdots+
\lambda_n
=
n.
\]
The {\sl diagram} of a partition $\lambda = ( \lambda_1, \lambda_2,
\dots, \lambda_n)$ in the real plane consists of $\lambda_1$ squares
of length one placed above $\lambda_2$ squares of length one, {\it
etc.}, placed above $\lambda_n$ squares of length one, all horizontal
series of squares being justified to the left along a fixed vertical
line; some figures appear below. The {\sl conjugate} of a partition
$\lambda$ is the partition $\lambda^c = ( \lambda_1^c, \lambda_2^c,
\dots, \lambda_n^c)$ whose diagram is obtained from the diagram of
$\lambda$ by reflecting it across its main diagonal. Hence
$\lambda_i^c$ is the number of squares in the $i$-th column of
$\lambda$, or equivalently $\lambda_i^c = \text{\rm Card}\, \big\{ j :
\, \lambda_j \geqslant i \big\}$.

\THEOREM

\smallskip\noindent\fbox{\bf THEOREM}
{\sf\em
The terms of highest order with respect to $\vert \ell \vert = \max_{
1 \leqslant i \leqslant n } \, \ell_i$ in the Euler-Poincaré
characteristic of the Schur bundle $\Gamma^{ ( \ell_1, \ell_2, \dots,
\ell_n )} \, T_X$ are homogeneous of order ${\rm O} \big( \vert \ell
\vert^{ \frac{ n ( n+1)}{ 2}} \big)$ and they are given by a sum of
$\ell_i'$-determinants indexed by all the partitions $(\lambda_1,
\dots, \lambda_n)$ of $n$:
\[
\small
\aligned
&
\chi\Big(X,\,\,
\Gamma^{(\ell_1,\ell_2,\dots,\ell_n)}\,T_X\Big)
=
\\
&
=
\sum_{\lambda\,\text{\rm partition of}\,\,n}\,
\frac{
{\sf C}_{\lambda^c}}{(\lambda_1+n-1)!\,\cdots\,\lambda_n!}\,
\left\vert
\begin{array}{cccc}
{\ell_1'}^{\lambda_1+n-1} & {\ell_2'}^{\lambda_1+n-1} &
\cdots & {\ell_n'}^{\lambda_1+n-1}
\\
{\ell_1'}^{\lambda_2+n-2} & {\ell_2'}^{\lambda_2+n-2} &
\cdots & {\ell_n'}^{\lambda_2+n-2}
\\
\vdots & \vdots & \ddots & \vdots
\\
{\ell_1'}^{\lambda_n} & {\ell_2'}^{\lambda_n} & 
\cdots & {\ell_n'}^{\lambda_n}
\end{array}
\right\vert
+
\\
&\ \ \ \ \ \ \ \ \ \
+
{\rm O}\Big(
\vert\ell\vert^{\frac{n(n+1)}{2}-1}
\Big),
\endaligned
\]
where $\ell_i ' := \ell_i + n - i$ for notational brevity, with
coefficients ${\sf C}_{\lambda^c}$ being expressed in terms of the
Chern classes ${\sf c}_k \big( T_X \big) = {\sf c}_k$ of $T_X$
by means of {\em Giambelli's determinantal expression} depending
upon the {\em conjugate} partition $\lambda^c$:
\[
{\sf C}_{\lambda^c}
=
{\sf C}_{(\lambda_1^c,\dots,\lambda_n^c)}
=
\left\vert
\begin{array}{cccccc}
{\sf c}_{\lambda_1^c} & {\sf c}_{\lambda_1^c+1} &
{\sf c}_{\lambda_1^c+2} & \cdots & {\sf c}_{\lambda_1^c+n-1}
\\
{\sf c}_{\lambda_2^c-1} & {\sf c}_{\lambda_2^c} &
{\sf c}_{\lambda_2^c+1} & \cdots & {\sf c}_{\lambda_2^c+n-2}
\\
\vdots & \vdots & \vdots & \ddots & \vdots
\\
{\sf c}_{\lambda_n^c-n+1} & {\sf c}_{\lambda_n^c-n+2} &
{\sf c}_{\lambda_n^c-n+3} & \cdots & 
{\sf c}_{\lambda_n^c}
\end{array}
\right\vert,
\]
on understanding by convention that ${\sf c}_k := 0$ 
for $k< 0$ or $k > n$, and that ${\sf c}_0 := 1$. }

\stopTHEOREM

In fact, replacing the $\ell_i'$ by the $\ell_i$ everywhere in the
framed formula would be harmless, because the difference between any
two corresponding determinants is easily seen to be an ${\rm O} \big(
\vert\ell\vert^{\frac{n(n+1)}{2}-1} \big)$, neglected in the
remainder.

\smallskip

We give two expanded instances of this general formula. Firstly, in
dimension $n = 3$, there are only three partitions of $3$, namely
$3+0+0$, $2+1+0$ and $1+1+1$, along which we draw the diagram of the
conjugate partitions $1+1+1$, $2 + 1 + 1$ and $3+0+0$ together with
the corresponding Giambelli determinants:

\begin{center}
\begin{picture}(0,0)%
\includegraphics{3-diagram.pstex}%
\end{picture}%
\setlength{\unitlength}{4144sp}%
\begingroup\makeatletter\ifx\SetFigFont\undefined%
\gdef\SetFigFont#1#2#3#4#5{%
  \reset@font\fontsize{#1}{#2pt}%
  \fontfamily{#3}\fontseries{#4}\fontshape{#5}%
  \selectfont}%
\fi\endgroup%
\begin{picture}(4589,837)(829,-561)
\put(1211,-136){\makebox(0,0)[lb]{\smash{{\SetFigFont{12}{14.4}{\familydefault}{\mddefault}{\updefault}{\color[rgb]{0,0,0}$\left\vert\begin{array}{ccc}{\sf c}_1&{\sf c}_2&{\sf c}_3\\1&{\sf c}_1&{\sf c}_2\\0&1&{\sf c}_1\end{array}\right\vert$}%
}}}}
\put(4686,-466){\makebox(0,0)[lb]{\smash{{\SetFigFont{12}{14.4}{\familydefault}{\mddefault}{\updefault}{\color[rgb]{0,0,0}$\left\vert\begin{array}{ccc}{\sf c}_3&0&0\\0&1&0\\0&0&1\end{array}\right\vert$}%
}}}}
\put(3109,-503){\makebox(0,0)[lb]{\smash{{\SetFigFont{12}{14.4}{\familydefault}{\mddefault}{\updefault}{\color[rgb]{0,0,0}$\left\vert\begin{array}{ccc}{\sf c}_2&{\sf c}_3&0\\1&{\sf c}_1&{\sf c}_2\\0&0&1\end{array}\right\vert$}%
}}}}
\end{picture}%

\end{center}

\noindent
so that we can write down in great details the leading terms, for
$\vert \ell \vert \to \infty$, of the Euler-Poincaré characteristic:
\[
\aligned
&
\chi\big(
X,\,
\Gamma^{(\ell_1,\ell_2,\ell_3)}\,T_X
\big)
=
\\
&
=
\frac{{\sf c}_1^3-2\,{\sf c}_1{\sf c}_2+{\sf c}_3}{0!\,\,1!\,\,5!}\
\left\vert
\begin{array}{ccc}
1\, & 1\, & 1\,
\\
\ell_1\, & \ell_2\, & \ell_3\,
\\
\ell_1^5\, & \ell_2^5\, & \ell_3^5\,
\end{array}
\right\vert
+
\frac{{\sf c}_1{\sf c}_2-{\sf c}_3}{0!\,\,2!\,\,4!}\
\left\vert
\begin{array}{ccc}
1\, & 1\, & 1\,
\\
\ell_1^2\, & \ell_2^2\, & \ell_3^2\,
\\
\ell_1^4\, & \ell_2^4\, & \ell_3^4\,
\end{array}
\right\vert
+
\\
&
\ \ \ \ \ \ \ \
+
\frac{{\sf c}_3}{1!\,\,2!\,\,3!}\
\left\vert
\begin{array}{ccc}
\ell_1\, & \ell_2\, & \ell_3\,
\\
\ell_1^2\, & \ell_2^2\, & \ell_3^2\,
\\
\ell_1^3\, & \ell_2^3\, & \ell_3^3\,
\end{array}
\right\vert
+
{\rm O}\big(\vert\ell\vert^5\big).
\endaligned
\]
Secondly, in dimension $n = 4$, there are five partitions of $4$,
namely $4$, $3+1$, $2+2$, $2+1+1$ and $1+1+1+1$ along which we again
draw the diagram of the conjugate partition together with the
corresponding Giambelli determinants:

\smallskip

\begin{picture}(0,0)%
\includegraphics{4-diagram.pstex}%
\end{picture}%
\setlength{\unitlength}{4144sp}%
\begingroup\makeatletter\ifx\SetFigFont\undefined%
\gdef\SetFigFont#1#2#3#4#5{%
  \reset@font\fontsize{#1}{#2pt}%
  \fontfamily{#3}\fontseries{#4}\fontshape{#5}%
  \selectfont}%
\fi\endgroup%
\begin{picture}(4804,2330)(514,-2014)
\put(837,-215){\makebox(0,0)[lb]{\smash{{\SetFigFont{12}{14.4}{\familydefault}{\mddefault}{\updefault}{\color[rgb]{0,0,0}$\left\vert\begin{array}{cccc}{\sf c}_1&{\sf c}_2&{\sf c}_3&{\sf c}_4\\1&{\sf c}_1&{\sf c}_2&{\sf c}_3\\0&1&{\sf c}_1&{\sf c}_2\\0&0&1&{\sf c}_1\end{array}\right\vert$}%
}}}}
\put(2782,-548){\makebox(0,0)[lb]{\smash{{\SetFigFont{12}{14.4}{\familydefault}{\mddefault}{\updefault}{\color[rgb]{0,0,0}$\left\vert\begin{array}{cccc}{\sf c}_2&{\sf c}_3&{\sf c}_4&0\\1&{\sf c}_1&{\sf c}_2&{\sf c}_3\\0&1&{\sf c}_1&{\sf c}_2\\0&0&0&1\end{array}\right\vert$}%
}}}}
\put(5318,-255){\makebox(0,0)[lb]{\smash{{\SetFigFont{12}{14.4}{\familydefault}{\mddefault}{\updefault}{\color[rgb]{0,0,0}$\left\vert\begin{array}{cccc}{\sf c}_2&{\sf c}_3&{\sf c}_4&0\\{\sf c}_1&{\sf c}_2&{\sf c}_3&{\sf c}_4\\0&0&1&{\sf c}_1\\0&0&0&1\end{array}\right\vert$}%
}}}}
\put(3862,-1956){\makebox(0,0)[lb]{\smash{{\SetFigFont{12}{14.4}{\familydefault}{\mddefault}{\updefault}{\color[rgb]{0,0,0}$\left\vert\begin{array}{cccc}{\sf c}_4&0&0&0\\0&1&{\sf c}_1&{\sf c}_2\\0&0&1&{\sf c}_1\\0&0&0&1\end{array}\right\vert$}%
}}}}
\put(1638,-1899){\makebox(0,0)[lb]{\smash{{\SetFigFont{12}{14.4}{\familydefault}{\mddefault}{\updefault}{\color[rgb]{0,0,0}$\left\vert\begin{array}{cccc}{\sf c}_3&{\sf c}_4&0&0\\1&{\sf c}_1&{\sf c}_2&{\sf c}_3\\0&0&1&{\sf c}_1\\0&0&0&1\end{array}\right\vert$}%
}}}}
\end{picture}%

\bigskip\noindent
so that we can write down in length the asymptotic of the
Euler-Poincaré characteristic also in this case, of major interest to 
us:
\[
\label{chi-4-4}
\boxed{
\aligned
&
\chi\big(
X,\,
\Gamma^{(\ell_1,\ell_2,\ell_3,\ell_4)}\,T_X
\big)
=
\\
&
=
\frac{{\sf c}_1^4-3\,{\sf c}_1^2{\sf c}_2+{\sf c}_2^2
+2\,{\sf c}_1{\sf c}_3-{\sf c}_4}{0!\,\,1!\,\,2!\,\,7!}\
\left\vert
\begin{array}{cccc}
1\, & 1\, & 1\, & 1\,
\\
\ell_1\, & \ell_2\, & \ell_3\, & \ell_4\,
\\
\ell_1^2\, & \ell_2^2\, & \ell_3^2\, & \ell_4^2\,
\\
\ell_1^7\, & \ell_2^7\, & \ell_3^7\, & \ell_4^7\,
\end{array}
\right\vert
+
\\
&
+
\frac{{\sf c}_1^2{\sf c}_2-{\sf c}_2^2-{\sf c}_1{\sf c}_3+{\sf c}_4}
{0!\,\,1!\,\,3!\,\,6!}\
\left\vert
\begin{array}{cccc}
1\, & 1\, & 1\, & 1\,
\\
\ell_1\, & \ell_2\, & \ell_3\, & \ell_4\,
\\
\ell_1^3\, & \ell_2^3\, & \ell_3^3\, & \ell_4^3\,
\\
\ell_1^6\, & \ell_2^6\, & \ell_3^6\, & \ell_4^6\,
\end{array}
\right\vert
+
\frac{-{\sf c}_1{\sf c}_3+{\sf c}_2^2}
{0!\,\,1!\,\,4!\,\,5!}\
\left\vert
\begin{array}{cccc}
1\, & 1\, & 1\, & 1\,
\\
\ell_1\, & \ell_2\, & \ell_3\, & \ell_4\,
\\
\ell_1^4\, & \ell_2^4\, & \ell_3^4\, & \ell_4^4\,
\\
\ell_1^5\, & \ell_2^5\, & \ell_3^5\, & \ell_4^5\,
\end{array}
\right\vert
+
\\
&
+
\frac{{\sf c}_1{\sf c}_3-{\sf c}_4}{0!\,\,2!\,\,3!\,\,5!}\
\left\vert
\begin{array}{cccc}
1\, & 1\, & 1\, & 1\,
\\
\ell_1^2\, & \ell_2^2\, & \ell_3^2\, & \ell_4^2\,
\\
\ell_1^3\, & \ell_2^3\, & \ell_3^3\, & \ell_4^3\,
\\
\ell_1^5\, & \ell_2^5\, & \ell_3^5\, & \ell_4^5\,
\end{array}
\right\vert
+
\frac{{\sf c}_4}{1!\,\,2!\,\,3!\,\,4!}\
\left\vert
\begin{array}{cccc}
\ell_1\, & \ell_2\, & \ell_3\, & \ell_4\,
\\
\ell_1^2\, & \ell_2^2\, & \ell_3^2\, & \ell_4^2\,
\\
\ell_1^3\, & \ell_2^3\, & \ell_3^3\, & \ell_4^3\,
\\
\ell_1^4\, & \ell_2^4\, & \ell_3^4\, & \ell_4^4\,
\end{array}
\right\vert
+
{\rm O}\big(\vert\ell\vert^9\big).
\endaligned}
\]

\proof[Proof of the general theorem] Taking the proposition for
granted, we start by expanding plainly in Taylor series the
exponentials of the numerator determinant:
\[
\aligned
\left\vert
\begin{array}{ccc}
e^{{\sf a}_1\ell_1'} & \cdots & e^{{\sf a}_1\ell_n'}
\\
\vdots & \ddots & \vdots
\\
e^{{\sf a}_n\ell_1'} & \cdots & e^{{\sf a}_n\ell_n'}
\end{array}
\right\vert
&
=
\left\vert
\begin{array}{ccc}
\sum_{\mu\geqslant 0}\,\frac{(\ell_1')^\mu}{\mu!}\,{\sf a}_1^\mu
& \cdots &
\sum_{\mu\geqslant 0}\,\frac{(\ell_n')^\mu}{\mu!}\,{\sf a}_1^\mu
\\
\vdots & \ddots & \vdots
\\
\sum_{\mu\geqslant 0}\,\frac{(\ell_1')^\mu}{\mu!}\,{\sf a}_n^\mu
& \cdots &
\sum_{\mu\geqslant 0}\,\frac{(\ell_n')^\mu}{\mu!}\,{\sf a}_n^\mu
\end{array}
\right\vert
\\
&
=
\sum_{\mu_1,\mu_2,\dots,\mu_n\geqslant 0}\,
\frac{(\ell_1')^{\mu_1}}{\mu_1!}\,
\frac{(\ell_2')^{\mu_2}}{\mu_2!}\,\cdots\,
\frac{(\ell_n')^{\mu_n}}{\mu_n!}\,
\left\vert
\begin{array}{ccc}
{\sf a}_1^{\mu_1} & \cdots & {\sf a}_1^{\mu_n}
\\
\vdots & \ddots & \vdots
\\
{\sf a}_n^{\mu_1} & \cdots & {\sf a}_n^{\mu_n}
\end{array}
\right\vert
\endaligned
\]
and we then develope the result by multilinearity. According to what
has already been noticed after the proposition, dividing this last sum
by the determinant at the denominator amounts to multiplying it by
$\big[ 1\big/ \prod_{ i < j }\, \big( {\sf a}_i - {\sf a}_j \big)
\big] \cdot\big[ 1+ \sum_{ k\geqslant 1 }\, (-1)^k \, R( {\sf a})^k
\big]$, so we obtain:
\[
\aligned
&
\chi
\Big(
X,\,\Gamma^{(\ell_1,\dots,\ell_n)}\,T_X
\Big)
=
\sum_{\mu_1,\mu_2,\dots,\mu_n\geqslant 0}\,
\frac{(\ell_1')^{\mu_1}}{\mu_1!}\,
\frac{(\ell_2')^{\mu_2}}{\mu_2!}\,\cdots\,
\frac{(\ell_n')^{\mu_n}}{\mu_n!}\cdot
\\
&
\cdot
\text{\sf homogeneous n-th part of}\Bigg(
\frac{1}{\prod_{i<j}\big({\sf a}_i-{\sf a}_j\big)}\,
\left\vert
\begin{array}{ccc}
{\sf a}_1^{\mu_1} & \cdots & {\sf a}_1^{\mu_n}
\\
\vdots & \ddots & \vdots
\\
{\sf a}_n^{\mu_1} & \cdots & {\sf a}_n^{\mu_n}
\end{array}
\right\vert
\Big[
1+{\rm O}_1({\sf a})
\Big]\Bigg),
\endaligned
\]
where we have gathered all terms $-R ( {\sf a} ) + R( {\sf a} )^2 -
\cdots$ simply as a remainder ${\rm O}_1({\sf a})$ vanishing at ${\sf
a} = 0$. The order at ${\sf a} = 0$ of the Van der Monde denominator
$\prod_{ i < j }\, \big( {\sf a}_i - {\sf a}_j\big)$ is equal to
$\frac{ n ( n-1)}{ 2}$, while the order of the determinant $\big\vert
{\sf a}_i^{\mu_j} \big\vert$ equals $\mu_1 + \cdots +
\mu_n$. Consequently, when selecting in the sum $\sum_{ \mu_1, \dots,
\mu_n \geqslant 0}$ only homogeneous terms of order $n$ with respect
to ${\sf a}$, one must consider:

\begin{itemize}

\smallskip\item[$\bullet$]
all terms with $\mu_1 + \cdots + \mu_n = n + \frac{ n ( n-1)}{ 2} =
\frac{ n ( n+1)}{ 2}$ if the determinant is multiplied by the term 1
inside the last brackets; with respect to the $\ell_i'$, this then gives
terms which are homogeneous of degree $\frac{ n ( n+1)}{ 2}$;

\smallskip\item[$\bullet$]
some appropriate terms with $\mu_1 + \cdots + \mu_n < \frac{ n (
n+1)}{ 2}$ if the determinant is multiplied by some nonzero monomial
belonging to the remainder ${\rm O}_1 ( {\sf a})$; with respect to the
$\ell_i'$, this then gives terms in ${\rm O} \big( \vert \ell' \vert
\big)^{ \frac{ n ( n+1)}{ 2} - 1}$, and we announced in the theorem
that we should neglect them.

\end{itemize}\smallskip

\noindent
As a result, we may therefore represent as follows the principal terms
of the Euler-Poincaré characteristic, considered 
asymptotically for $\vert \ell \vert \to \infty$:
\[
\aligned
\chi
\Big(
X,\,\Gamma^{(\ell_1,\dots,\ell_n)}\,T_X
\Big)
&
=
\sum_{\mu_1+\cdots+\mu_n=\frac{n(n+1)}{2}\atop
\mu_1,\dots,\mu_n\geqslant 0}\,\,
\frac{(\ell_1')^{\mu_1}}{\mu_1!}\,
\frac{(\ell_2')^{\mu_2}}{\mu_2!}\,\cdots\,
\frac{(\ell_n')^{\mu_n}}{\mu_n!}\cdot
\\
&\ \ \ \ \
\cdot
\frac{1}{\prod_{i<j}\big({\sf a}_i-{\sf a}_j\big)}\,
\left\vert
\begin{array}{ccc}
{\sf a}_1^{\mu_1} & \cdots & {\sf a}_1^{\mu_n}
\\
\vdots & \ddots & \vdots
\\
{\sf a}_n^{\mu_1} & \cdots & {\sf a}_n^{\mu_n}
\end{array}
\right\vert
+
{\rm O}\Big(\vert\ell'\vert^{\frac{n(n+1)}{2}-1}\Big).
\endaligned
\]
Whenever there exist two equal exponents $\mu_{ i_1 } = \mu_{ i_2 }$
for two distinct indices $i_1 \neq i_2$, the determinant obviously
vanishes. So in the sum, one may assume the $\mu_i$ to be pairwise
distinct. Furthermore, for any $n$-tuple $(\mu_1, \dots, \mu_n)$ of
pairwise distinct $\mu_i$, there exists a unique permutation $\sigma
\in \mathfrak{ S}_n$ rearranging them in decreasing order: $\mu_{
\sigma ( 1)} > \mu_{ \sigma ( 2)} > \cdots > \mu_{ \sigma ( n)}$.
Consequently, we can split as follows the sum to be considered:
\[
\aligned
\chi
\Big(
X,\,\Gamma^{(\ell_1,\dots,\ell_n)}\,T_X
\Big)
&
=
\sum_{\sigma\in\mathfrak{S}_n}\,
\sum_{\mu_1+\cdots+\mu_n=\frac{n(n+1)}{2}\atop
\mu_1>\cdots>\mu_n\geqslant 0}\,
\frac{(\ell_1')^{\mu_{\sigma(1)}}}{\mu_{\sigma(1)}!}\cdots\,
\frac{(\ell_n')^{\mu_{\sigma(n)}}}{\mu_{\sigma(n)}!}\cdot
\\
&\ \ \ \ \ 
\cdot
\frac{1}{\prod_{i<j}\big({\sf a}_i-{\sf a}_j\big)}\,
\left\vert
\begin{array}{ccc}
{\sf a}_1^{\mu_{\sigma(1)}} & \cdots & {\sf a}_1^{\mu_{\sigma(n)}}
\\
\vdots & \ddots & \vdots
\\
{\sf a}_n^{\mu_{\sigma(1)}} & \cdots & {\sf a}_n^{\mu_{\sigma(n)}}
\end{array}
\right\vert
+
{\rm O}\Big(\vert\ell'\vert^{\frac{n(n+1)}{2}-1}\Big).
\endaligned
\]
Finally, one easily convinces oneself that there is a one-to-one
correspondence between the $n$-tuples $\mu = (\mu_1, \dots, \mu_n)$ as
above with $\mu_1 > \cdots > \mu_n \geqslant 0$ and $\mu_1 + \cdots +
\mu_n = \frac{ n ( n+1)}{ 2}$ on the one hand, and on the other hand,
the partitions $\lambda = ( \lambda_1, \lambda_2, \dots, \lambda_n)$
of $n$, namely with $\lambda_1 \geqslant \cdots \geqslant \lambda_n
\geqslant 0$ and $\lambda_1 + \cdots + \lambda_n = n$, a
correspondence which is simply given by:
\[
\mu_i\longmapsto\lambda_i:=\mu_i-n+i
\ \ \ \ \ \
\text{\rm and has obvious inverse}
\ \ \ \ \ \ 
\lambda_i\longmapsto\mu_i:=\lambda_i+n-i.
\]
Taking account of the skew-symmetry $\big\vert {\sf a}_i^{ \mu_{
\sigma(j) }} \big\vert = {\rm sgn} (\sigma)\, \big\vert {\sf
a}_i^{\mu_j} \big \vert$, we thus obtain an almost final asymptotic
representation of the Euler-Poincaré characteristic:
\[
\aligned
&
\chi
\Big(
X,\,\Gamma^{(\ell_1,\dots,\ell_n)}\,T_X
\Big)
=
\sum_{\sigma\in\mathfrak{S}_n}\,
\sum_{\lambda_1+\cdots+\lambda_n=n\atop
\lambda_1\geqslant\cdots\geqslant\lambda_n\geqslant 0}\,
\,\,
\frac{(\ell_{\sigma^{-1}(1)}')^{
\lambda_1+n-1}}{(\lambda_1+n-1)!}\,\cdots\,
\frac{(\ell_{\sigma^{-1}(n)}')^{\lambda_n}}{\lambda_n!}\cdot
{\rm sgn}(\sigma)\cdot
\\
&
\cdot
\frac{1}{\prod_{i<j}\big({\sf a}_i-{\sf a}_j\big)}\,
\left\vert
\begin{array}{ccc}
{\sf a}_1^{\lambda_1+n-1} & \cdots & {\sf a}_1^{\lambda_n}
\\
\vdots & \ddots & \vdots
\\
{\sf a}_n^{\lambda_1+n-1} & \cdots & {\sf a}_n^{\lambda_n}
\end{array}
\right\vert
+
{\rm O}\Big(\vert\ell'\vert^{\frac{n(n+1)}{2}-1}\Big).
\endaligned
\]
To conclude the proof of the theorem, using ${\rm sgn} ( \sigma^{ -1})
= {\rm sgn} ( \sigma)$, it now suffices only to observe the compulsory
reconstitution of $\ell'$-determinants:
\[
\small
\aligned
&
\sum_{\sigma\in\mathfrak{S}_n}\,
{\rm sgn}\,({\sigma^{-1}})\cdot
\frac{(\ell_{\sigma^{-1}(1)}')^{
\lambda_1+n-1}}{(\lambda_1+n-1)!}\,\cdots\,
\frac{(\ell_{\sigma^{-1}(n)}')^{\lambda_n}}{\lambda_n!}
=
\\
&
\ \ \ \ \ \ \ \ \ \ \ \ \ \ \ 
=
\frac{1}{(\lambda_1+n-1)!\cdots\lambda_n!}\cdot
\left\vert
\begin{array}{ccc}
{\ell_1'}^{\lambda_1+n-1} & \cdots & {\ell_n'}^{\lambda_1+n-1}
\\
\vdots & \ddots & \vdots
\\
{\ell_1'}^{\lambda_n} & \cdots & {\ell_n'}^{\lambda_n},
\end{array}
\right\vert,
\endaligned
\]
and also to recognize the Schur polynomials:
\[
S_\lambda({\sf a})
=
S_{(\lambda_1,\dots,\lambda_n)}({\sf a})
=
\frac{1}{\prod_{i<j}\,
\big({\sf a}_i-{\sf a}_j\big)}\,
\left\vert
\begin{array}{cccc}
{\sf a}_1^{\lambda_1+n-1} & {\sf a}_1^{\lambda_2+n-2}
& \cdots & {\sf a}_1^{\lambda_n}
\\
{\sf a}_2^{\lambda_1+n-1} & {\sf a}_2^{\lambda_2+n-2}
& \cdots & {\sf a}_2^{\lambda_n}
\\
\vdots & \vdots & \ddots & \vdots
\\
{\sf a}_n^{\lambda_1+n-1} & {\sf a}_n^{\lambda_2+n-2}
& \cdots & {\sf a}_n^{\lambda_n}
\end{array}
\right\vert
\] 
indexed by the partitions of $n$, which according to Giambelli's
formulas (Appendix~A of \cite{ fuha1991}), are expressed in terms of
the elementary symmetric functions ${\sf c}_k = \sum_{ 1 \leqslant
i_1 < \cdots < i_k \leqslant n}\, {\sf a}_{ i_1} \cdots {\sf a}_{
i_k}$ of the ${\sf a}_i$ by means of the specific determinants written
and exemplified above. Thus, the proof is achieved.
\endproof

\subsection*{ Computation of the Euler-Poincaré characteristic of
${\sf E}_{ 4, m}^4 T_X^*$} As is known, duality shows that the
cotangent bundle $T_X^*$ has Chern classes ${\sf c}_k \big( T_X^*
\big)$ related to those of $T_X$ by the relations:
\[
{\sf c}_k^*
:=
{\sf c}_k\big(T_X^*\big)
=
(-1)^k\,{\sf c}_k\big(T_X)
=
(-1)^k\,{\sf c}_k.
\]
Consequently, the dual Giambelli determinants satisfy ${\sf C}_{
\lambda^c}^* = (-1)^n {\sf C}_{\lambda^c}$, because all monomials
${\sf c}_{ \mu_1}^* \cdots {\sf c}_{ \mu_n}^*$ have total weight
$\mu_1 + \cdots + \mu_n = n$ and we therefore deduce:
\[
\chi\big(
X,\,\Gamma^{(\ell_1,\dots,\ell_n)}T_X^*\big)
=
(-1)^n\,
\chi\big(
X,\,\Gamma^{(\ell_1,\dots,\ell_n)}T_X\big).
\]
When considering Demailly-Semple and Schur bundles, everything shall be
expressed in terms of Chern classes of $T_X$ (not of $T_X^*$).

\section*{\S14.~Euler characteristic calculations}
\label{Section-14}

\subsection*{ Explaining the final calculations on an example}
We may now come back to our 24 sums of Schur bundles (with
multiplicities). Consider for instance the family ${\sf A}$. 
In it, we have:
\[
\left\{
\aligned
\ell_1
&
=
o+4j+3k+3l+4m'+5n+p,
\\
\ell_2
&
=
2j+2k+3l+3m'+3n+p,
\\
\ell_3
&
=
j+2k+2l+2m'+2n+p,
\\
\ell_4
&
=
p.
\endaligned\right.
\]
But since sums of weight should be equal to $m$:
\[
o+14j+15k+17l+19m'+21n+10p
=
m,
\]
we may eliminate $o$ and this provides $\ell_1$ with the value:
\[
\ell_1
=
m-10j-12k-14l-15m'-16n-9p,
\]
while $\ell_2$, $\ell_3$ and $\ell_4$ where at the beginning
independent of $o$. The Euler-Poincaré characteristic being additive,
we have:
\[
\chi\Big(X,\,\oplus_{\sf A}\,
\Gamma^{(\ell_1,\ell_2,\ell_3,\ell_4)}T_X^*
\Big)
=
\sum_{o+14j+15k+17l+19m'+21n+10p=m}\,
\chi\big(X,\,\Gamma^{(\ell_1,\ell_2,\ell_3,\ell_4)}T_X^*\big).
\]
Furthermore, according to the formula written on p.~\pageref{chi-4-4},
the dominant term of the Euler-Poincaré characteristic, as $\vert \ell
\vert \to \infty$, of a single Schur bundle in such a sum is given, in
terms of the Chern classes ${\sf c}_k$ of $T_X$, by:
\[
\aligned
\chi
\big(X,\,
\Gamma^{(\ell_1,\ell_2,\ell_3,\ell_4)}T_X^*
\big)
&
=
\frac{{\sf c}_1^4-3\,{\sf c}_1^2{\sf c}_2+{\sf c}_2^2+
2\,{\sf c}_1{\sf c}_3-{\sf c}_4}{0!\,1!\,2!\,7!}\,
\Delta_{0127}
+
\\
&\ \ \ \ \
+
\frac{{\sf c}_1^2{\sf c}_2-{\sf c}_2^2-{\sf c}_1{\sf c}_3+
{\sf c}_4}{0!\,1!\,3!\,6!}\,
\Delta_{0136}
+
\frac{-{\sf c}_1{\sf c}_3+{\sf c}_2^2}{0!\,1!\,4!\,5!}\,
\Delta_{0145}
+
\\
&\ \ \ \ \
+
\frac{{\sf c}_1{\sf c}_3-{\sf c}_4}{0!\,2!\,3!\,5!}\,
\Delta_{0235}
+
\frac{{\sf c}_4}{1!\,2!\,3!\,4!}\,
\Delta_{1234}
\\
&\ \ \ \ \ \
+
{\rm O}\big(\vert\ell\vert^9\big),
\endaligned
\]
on understanding that, in the five determinants:
\[
\aligned
\Delta_{0137}
&
:=
\left\vert
\begin{array}{cccc}
1\, & 1\, & 1\, & 1\,
\\
\ell_1\, & \ell_2\, & \ell_3\, & \ell_4\,
\\
\ell_1^2\, & \ell_2^2\, & \ell_3^2\, & \ell_4^2\,
\\
\ell_1^7\, & \ell_2^7\, & \ell_3^7\, & \ell_4^7\,
\end{array}
\right\vert,
\ \ \ \ \ \ \ \ \ \ \ \ \ \ \ \ \
\Delta_{0136}
:=
\left\vert
\begin{array}{cccc}
1\, & 1\, & 1\, & 1\,
\\
\ell_1\, & \ell_2\, & \ell_3\, & \ell_4\,
\\
\ell_1^3\, & \ell_2^3\, & \ell_3^3\, & \ell_4^3\,
\\
\ell_1^6\, & \ell_2^6\, & \ell_3^6\, & \ell_4^6\,
\end{array}
\right\vert,
\\
\Delta_{0145}
&
:=
\left\vert
\begin{array}{cccc}
1\, & 1\, & 1\, & 1\,
\\
\ell_1\, & \ell_2\, & \ell_3\, & \ell_4\,
\\
\ell_1^4\, & \ell_2^4\, & \ell_3^4\, & \ell_4^4\,
\\
\ell_1^5\, & \ell_2^5\, & \ell_3^5\, & \ell_4^5\,
\end{array}
\right\vert,
\ \ \ \ \ \ \ \ \ \ \ \ \ \ \ \ \
\Delta_{0235}
:=
\left\vert
\begin{array}{cccc}
1\, & 1\, & 1\, & 1\,
\\
\ell_1^2\, & \ell_2^2\, & \ell_3^2\, & \ell_4^2\,
\\
\ell_1^3\, & \ell_2^3\, & \ell_3^3\, & \ell_4^3\,
\\
\ell_1^5\, & \ell_2^5\, & \ell_3^5\, & \ell_4^5\,
\end{array}
\right\vert,
\\
\Delta_{1234}
&
:=
\left\vert
\begin{array}{cccc}
\ell_1\, & \ell_2\, & \ell_3\, & \ell_4\,
\\
\ell_1^2\, & \ell_2^2\, & \ell_3^2\, & \ell_4^2\,
\\
\ell_1^3\, & \ell_2^3\, & \ell_3^3\, & \ell_4^3\,
\\
\ell_1^4\, & \ell_2^4\, & \ell_3^4\, & \ell_4^4\,
\end{array}
\right\vert,
\endaligned
\]
one should substitute the above values for $\ell_1$, $\ell_2$,
$\ell_3$ and $\ell_4$ in terms of $j$, $k$, $l$, $m'$, $n$ and $p$.

On the other hand, it is well known that the dominant term of a
multiple sum is given by an integral, so that we have to
compute\footnote{\, The $\Delta$ determinants being of degree 10 in
the $\ell_i$, the presence of six integrals entails that the result is
$m^{ 16}$ times a fractional constant plus an ${\rm O} \big( m^{ 15}
\big)$. If there would be 5 or less integrals, this would leave us
with an ${\rm O} \big( m^{ 15} \big)$, negligible in comparison with
$m^{ 16}$ as $m \to \infty$. By this remark we therefore justify why
we considered only the approximate Schur bundle decomposition of ${\sf
E}_{ 4, m}^4 T_X^*$ in the \S12. }:
\[
\aligned
&
\int_0^{\frac{(m-15k-17l-19m'-21n-10p)}{14}}\,dj\,
\int_0^{\frac{(m-17l-19m'-21n-10p)}{15}}\,dk\,
\int_0^{\frac{(m-19m'-21n-10p)}{17}}\,dl\,
\\
&
\int_0^{\frac{(m-21n-10p)}{19}}\,dm'\,
\int_0^{\frac{(m-10p)}{21}}\,dn\,
\int_0^{\frac{m}{10}}\,dp\,\,\,
\left\{
\aligned
&
\Delta_{0127}
\\
&
\Delta_{0136}
\\
&
\Delta_{0145}
\\
&
\Delta_{0235}
\\
&
\Delta_{1234}
\endaligned\right.
\endaligned
\]
It happens that all the five integrals are equal to $m^{ 16}$ times a
fractional number. A computation with the help of Maple yields the
values of these five fractional numbers, which, we guess, would be
quite uneasy to get by hand:
\[
\aligned
{\sf A}_{0127}
&
=
\frac{157423754766863651482110063939631617713614267}{
7470130440549849070995762660822781685545412418720000000000000},
\\
{\sf A}_{0136}
&
=
\frac{285224611253902544589491011638457808537315047}{
34860608722565962331313559083839647865878591287360000000000000}
\\
{\sf A}_{0145}
&
=
\frac{10306128852122999807705628256770676631371801}{
5229091308384894349697033862575947179881788693104000000000000}
\\
{\sf A}_{0235}
&
=
\frac{2097522233626513305099611552292506537139247}{
2376859685629497431680469937534521445400813042320000000000000}
\\
{\sf A}_{1234}
&
=
\frac{20051359515371820286197508247902844353}{
2102485347748339169995992868230447983547822240000000000000}
\\
\endaligned
\]

\subsection*{ End of the computation}
Similarly, for the other 23 families, we compute these 5-tuples of
rational numbers and at the end, we make the 
addition\footnote{\, {\em See} {\sf new-riemann-roch-4-4.mws} at~\cite{
mer2008c}. }:
\[
\aligned
{\sf Coeff}_{0127} 
&
=
2\,{\sf A}_{0127}+{\sf B}_{0127}+{\sf C}_{0127}+4\,{\sf D}_{0127}+2\,{\sf E}_{0127}+2\,{\sf F}_{0127}+{\sf G}_{0127}+{\sf H}_{0127}+
\\
&\ \ \ \ \ 
+{\sf I}_{0127}+{\sf J}_{0127}+{\sf K}_{0127}+{\sf L}_{0127}+{\sf M}_{0127}+{\sf N}_{0127}+{\sf O}_{0127}+{\sf P}_{0127}+
\\
&\ \ \ \ \
+{\sf Q}_{0127}+{\sf R}_{0127}+{\sf S}_{0127}+{\sf T}_{0127}+{\sf U}_{0127}+{\sf V}_{0127}+{\sf W}_{0127}+{\sf X}_{0127}
\\
&
=
\frac{2127566277536547206644157}{
65144733745232853829877760000000000000},
\endaligned
\]
\[
\aligned
{\sf Coeff}_{_{0136}} 
&
=
2\,{\sf A}_{0136}+{\sf B}_{0136}+{\sf C}_{0136}+4\,{\sf D}_{0136}+2\,{\sf E}_{0136}+2\,{\sf F}_{0136}+{\sf G}_{0136}+{\sf H}_{0136}+
\\
&\ \ \ \ \
+{\sf I}_{0136}+{\sf J}_{0136}+{\sf K}_{0136}+{\sf L}_{0136}+{\sf M}_{0136}+{\sf N}_{0136}+{\sf O}_{0136}+{\sf P}_{0136}+
\\
&\ \ \ \ \
+{\sf Q}_{0136}+{\sf R}_{0136}+{\sf S}_{0136}+{\sf T}_{0136}+{\sf U}_{0136}+{\sf V}_{0136}+{\sf W}_{0136}+{\sf X}_{0136}
\\
&
=
\frac{52676407087143116547997}{4053450099703377571636838400000000000},\endaligned
\]
\[
\aligned
{\sf Coeff}_{_{0145}} 
&
=
2\,{\sf A}_{0145}+{\sf B}_{0145}+{\sf C}_{0145}+4\,{\sf D}_{0145}+2\,{\sf E}_{0145}+2\,{\sf F}_{0145}+{\sf G}_{0145}+{\sf H}_{0145}+
\\
&\ \ \ \ \
+{\sf I}_{0145}+{\sf J}_{0145}+{\sf K}_{0145}+{\sf L}_{0145}+{\sf M}_{0145}+{\sf N}_{0145}+{\sf O}_{0145}+{\sf P}_{0145}+
\\
&\ \ \ \ \ 
+{\sf Q}_{0145}+{\sf R}_{0145}+{\sf S}_{0145}+{\sf T}_{0145}+{\sf U}_{0145}+{\sf V}_{0145}+{\sf W}_{0145}+{\sf X}_{0145}
\\
&
=
\frac{164685282124542664946051}{50668126246292219645460480000000000000},
\endaligned
\]
\[
\aligned
{\sf Coeff}_{_{0235}}
&
=
2\,{\sf A}_{0235}+{\sf B}_{0235}+{\sf C}_{0235}+4\,{\sf D}_{0235}+2\,{\sf E}_{0235}+2\,{\sf F}_{0235}+{\sf G}_{0235}+{\sf H}_{0235}+
\\
&\ \ \ \ \ 
+{\sf I}_{0235}+{\sf J}_{0235}+{\sf K}_{0235}+{\sf L}_{0235}+{\sf M}_{0235}+{\sf N}_{0235}+{\sf O}_{0235}+{\sf P}_{0235}+{\sf Q}_{0235}+
\\
&\ \ \ \ \ 
+{\sf R}_{0235}+{\sf S}_{0235}+{\sf T}_{0235}+{\sf U}_{0235}+{\sf V}_{0235}+{\sf W}_{0235}+{\sf X}_{0235}
\\
&
= 
\frac{122298240743566105217737}{114003284054157494202286080000000000000},
\endaligned
\]
\[
\aligned
{\sf Coeff}_{_{1234}} 
&
=
2\,{\sf A}_{1234}+{\sf B}_{1234}+{\sf C}_{1234}+4\,{\sf D}_{1234}+2\,{\sf E}_{1234}+2\,{\sf F}_{1234}+{\sf G}_{1234}+{\sf H}_{1234}+
\\
&\ \ \ \ \
+{\sf I}_{1234}+{\sf J}_{1234}+{\sf K}_{1234}+{\sf L}_{1234}+{\sf M}_{1234}+{\sf N}_{1234}+{\sf O}_{1234}+{\sf P}_{1234}+{\sf Q}_{1234}+
\\
&\ \ \ \ \
+{\sf R}_{1234}+{\sf S}_{1234}+{\sf T}_{1234}+{\sf U}_{1234}+{\sf V}_{1234}+{\sf W}_{1234}+{\sf X}_{1234}
\\
&
=
\frac{1429957461022772407321}{130289467490465707659755520000000000000}.
\endaligned
\]
Coming back to the Euler-Poincaré characteristic we therefore get:
\[
\aligned
\chi\big(X,\,{\sf E}_{4,m}^4T_X^*\big)
&
=
\frac{{\sf c}_1^4-3\,{\sf c}_1^2{\sf c}_2+{\sf c}_2^2+
2\,{\sf c}_1{\sf c}_3-{\sf c}_4}{0!\,1!\,2!\,7!}\,
{\sf Coeff}_{0127}
+
\\
&\ \ \ \ \
+
\frac{{\sf c}_1^2{\sf c}_2-{\sf c}_2^2-{\sf c}_1{\sf c}_3+
{\sf c}_4}{0!\,1!\,3!\,6!}\,
{\sf Coeff}_{0136}
+
\frac{-{\sf c}_1{\sf c}_3+{\sf c}_2^2}{0!\,1!\,4!\,5!}\,
{\sf Coeff}_{0145}
+
\\
&\ \ \ \ \
+
\frac{{\sf c}_1{\sf c}_3-{\sf c}_4}{0!\,2!\,3!\,5!}\,
{\sf Coeff}_{0235}
+
\frac{{\sf c}_4}{1!\,2!\,3!\,4!}\,
{\sf Coeff}_{1234}
+
{\rm O}\big(m^{15}\big)
\\
&
=
m^{16}\bigg(
\frac{2127566277536547206644157}{
656658916151947166605167820800000000000000}\,
{\sf c}_1^4-
\\
&\ \ \ \ \ \ \ \ \ \ \ \ \ \ \
-
\frac{139915351328310309504209}{
20846314798474513225560883200000000000000}\,
{\sf c}_1^2{\sf c}_2+
\\
&\ \ \ \ \ \ \ \ \ \ \ \ \ \ \
+
\frac{18230301659778006701051}{
13401202370447901359289139200000000000000}\,
{\sf c}_2^2+
\\
&\ \ \ \ \ \ \ \ \ \ \ \ \ \ \
+
\frac{405575296543809994270429}{
131331783230389433321033564160000000000000}\,
{\sf c}_1{\sf c}_3-
\\
&\ \ \ \ \ \ \ \ \ \ \ \ \ \ \
-
\frac{6163697191750462398371}{
6566589161519471666051678208000000000000}\,
{\sf c}_4\bigg)
+
\\
&
+
{\rm O}\big(m^{15}\big).
\endaligned
\]
In terms of the degree:
\[
\boxed{
\aligned
\chi\big(X,\,{\sf E}_{4,m}^4T_X^*\big)
&
=
\frac{m^{16}}{1313317832303894333210335641600000000000000}\,
\cdot\,d\,\cdot
\\
&\ \ \ \ \
\cdot
\big(
50048511135797034256235\,d^4
-
\\
&\ \ \ \ \ \ \ \ \ \
-
6170606622505955255988786\,d^3
-
\\
&\ \ \ \ \ \ \ \ \ \
-
928886901354141153880624704\,d
+
\\
&\ \ \ \ \ \ \ \ \ \
+
141170475250247662147363941\,d^2
+
\\
&\ \ \ \ \ \ \ \ \ \
+
1624908955061039283976041114
\big)
\\
&\ \ \ \ \ \
+
{\rm O}\big(m^{15}\big).
\endaligned}
\]
The four roots of the 4-th degree numerator in parentheses are:
\[
2.794353346\cdots,\ \ \ \ \ \ \ 
6.784939538\cdots,\ \ \ \ \ \ \
17.86618823\cdots,\ \ \ \ \ \ \
95.84703014\cdots,
\]
hence in conclusion, the characteristic is positive
for all degrees $d \geqslant 96$.

\subsection*{ Jets of order $\kappa = 4$ in dimension $n = 3$}
For a hypersurface $X^3 \subset \P_4 ( \C)$ of degree $d$, thanks to a
similar but quicker Maple computation\footnote{\, {\em See} {\sf
new-riemann-roch-3-4.mws} at~\cite{ mer2008c}. }, one obtains the
asymptotic:
\[
\aligned
\chi\big(X,\,{\sf E}_{4,m}^3T_X^*\big)
&
=
m^{11}\bigg(
-\frac{78181453985171}{2013023350054886400000000}\,{\sf c}_1^3
+
\\
&\ \ \ \ \ \ \ \ \ \ \ \ \ \ \ 
+
\frac{3780346214152789}{343555985076033945600000000}\,{\sf c}_3
-
\\
&\ \ \ \ \ \ \ \ \ \ \ \ \ \ \ 
-
\frac{46223512567695359}{1030667955228101836800000000}\,
{\sf c}_1{\sf c}_2
\bigg)
+
\\
&\ \ \ \ \
+
{\rm O}\big(m^{10}\big),
\endaligned
\]
and then in terms of the degree $d$ of the hypersurface $X$:
\[
\boxed{
\aligned
\chi\big(X,\,{\sf E}_{4,m}^3T_X^*\big)
&
=
\frac{m^{11}}{206133591045620367360000000}\,
\cdot\,d\,\cdot\,
\\
&\ \ \ \ \ \ \ \ \ \
\cdot\,\big(
1029286103034112\,d^3
-
38980726828290305\,d^2
+
\\
&\ \ \ \ \ \ \ \ \ \ \ \ \ \ \
+
299551055917162501\,d
-
561169562618151944
\big)
\endaligned}\,.
\]
The three roots of the third degree numerator in parentheses are:
\[
2.852373090\cdots,\ \ \ \ \ \ \ 
6.765004304\cdots,\ \ \ \ \ \ \ 
28.25423742,
\]
hence in conclusion, the Euler-Poincaré characteristic of ${\sf E}_{
4, m}^3 T_X^*$ is positive in all degrees $d \geqslant 29$ as $m \to
\infty$. This condition improves the condition $d \geqslant 43$
obtained in~\cite{ rou2006a} for the positivity of $\chi\big( X, \,
{\sf E}_{ 3, m}^3 T_X^*\big)$ as $m \to \infty$.

\subsection*{ Existence of sections}
Finally, in order to get positivity of the dimension $h^0$ of the
vector space of sections of ${\sf E}_{ 4, m}^3 T_X^*$, it would
suffice, in the trivial minoration:
\[
\aligned
h^0\big(X,\,{\sf E}_{4,m}^3T_X^*\big)
\geqslant
\chi\big(X,\,{\sf E}_{4,m}^3T_X^*\big)
-
h^2\big(X,\,{\sf E}_{4,m}^3T_X^*\big),
\endaligned
\] 
stemming from the definition $\chi = h^0 - h^1 + h^2 - h^3$, to
possess a good majoration of $h^2$. This main task is achieved in \cite{
rou2006b, rou2007b}: for each Schur bundle, one has:
\[
h^2
\big(
X,\,\Gamma^{(\ell_1,\ell_2,\ell_3)}T_X^*
\big)
\leqslant
d(d+13)\,
\frac{3(\ell_1+\ell_2+\ell_3)^3}{2}\,
(\ell_1-\ell_2)(\ell_1-\ell_3)(\ell_2-\ell_3)
+
{\rm O}\big(\vert\ell\vert^5\big).
\]
When summing up our 24 sums of Schur bundles (with multiplicities), 
a Maple computation provides: 
\[
h^2
\big(
X,\,\Gamma^{(\ell_1,\ell_2,\ell_3)}T_X^*
\big)
\leqslant
d(d+13)\,
\frac{342988705758851}{29822568148961280000000}\,m^{11}
+
{\rm O}\big(m^{10}\big).
\]
Finally, one sees that $\chi$ minus this upper bound for $h^2$ is
positive, for $m \to \infty$, in all degrees $d \geqslant 72$. This
last condition on the degre insuring the existence of invariant jet
differentials improves the condition $d \geqslant 97$ obtained
in~\cite{ rou2006b} and appears to be slightly better than the
condition $d \geqslant 74$ obtained recently in~\cite{ div2007}.

\vfill\end{document}